\documentclass[a4paper,12pt,reqno]{amsart}

\usepackage{amsmath,amssymb,amsthm}
\usepackage{mathrsfs}
\usepackage{amsmath}
\usepackage{amssymb,amsmath}
\usepackage{graphicx}
\usepackage{latexsym}
\usepackage{color}
\usepackage{graphicx}
\usepackage{mathrsfs}
\usepackage{enumerate}
\usepackage[abbrev]{amsrefs}
\usepackage[T1]{fontenc}
\usepackage{graphicx}
\usepackage{fancyhdr}
\usepackage{CJK}
\usepackage{cases}
\usepackage{hyperref}
\hypersetup{
colorlinks=true,
linkcolor=blue,
anchorcolor=blue,
citecolor=blue}
\allowdisplaybreaks[4]

\theoremstyle{plain}
\newtheorem{thm}{Theorem}[section]
\newtheorem{lemm}[thm]{Lemma}
\newtheorem{prop}[thm]{Proposition}
\newtheorem{cor}[thm]{Corollary}

\theoremstyle{definition}

\newtheorem{rem}[thm]{Remark}

\makeatletter

\@addtoreset{equation}{section}
\makeatother

\newcommand{\f}{\frac}
\newcommand{\beq}{\begin{equation}}
\newcommand{\eeq}{\end{equation}}

\newcommand{\Grad}{\nabla}
\newcommand{\p}{\partial}

\newcommand{\xh}{x_{\rm h}}
\newcommand{\yh}{y_{\rm h}}
\newcommand{\uh}{u_{\rm h}}
\newcommand{\vh}{v_{\rm h}}
\newcommand{\Bh}{B_{\rm h}}
\newcommand{\alphah}{\alpha_{\rm h}}
\newcommand{\nablah}{\nabla_{\rm h}}
\newcommand{\Gh}{G_{\rm h}}
\newcommand{\Gv}{G_{\rm v}}
\newcommand{\Deltah}{\Delta_{\rm h}}
\newcommand{\R}{\mathbb{R}}
\newcommand{\Nhu}{\mathcal{U}^{\rm h}}
\newcommand{\Nvu}{\mathcal{U}^{\rm v}}
\newcommand{\Nhb}{\mathcal{B}^{\rm h}}
\newcommand{\Nvb}{\mathcal{B}^{\rm v}}
\newcommand{\Nhd}{\mathcal{D}^{\rm h}}
\newcommand{\Nvd}{\mathcal{D}^{\rm v}}
\newcommand{\Nhe}{\mathcal{E}^{\rm h}}
\newcommand{\Nve}{\mathcal{E}^{\rm v}}

\begin{document}
\title[Anisotropic MHD system]
{Large time behavior of the solutions to 3D incompressible MHD system with horizontal dissipation or horizontal magnetic diffusion}
\author{ Yang Li}

\address[]{(Yang Li)    School of Mathematical Sciences and Center of Pure Mathematics, Anhui University, 230601, People's Republic of China}
\email{lynjum@163.com}
\date{\today}

\begin{abstract}
In this paper, we consider the asymptotic behavior of global solutions to 3D anisotropic incompressible MHD systems. For the 3D MHD system with horizontal dissipation and full magnetic diffusion, it is shown that $\uh(t)$ decays at the rate of $O(t^{-(1-\f{1}{p})})$, $u_3(t)$ decays at the rate of $O(t^{-\f{3}{2}(1-\f{1}{p})})$ and $B(t)$ decays at the rate of $O(t^{-\f{3}{2}(1-\f{1}{p})-\f{1}{2}})$. Furthermore, we give the asymptotic expansion of solutions. We prove that the leading term of $\uh(t)$ is a combination of linear solution and two integrals from nonlinear coupling effects, while for $u_3(t)$ the leading term is given by only the linear solution without the influence of magnetic field. Though the dissipation of velocity is weak, we show that the full magnetic diffusion is robust enough to keep the asymptotic expansion of magnetic field basically expected. However, the magnetic field turns out to affect the higher order asymptotic expansions of $u_3(t)$.
For the 3D MHD system with full dissipation and horizontal magnetic diffusion, it is shown that $\Bh(t)$ decays at the rate of $O(t^{-(1-\f{1}{p})})$, $B_3(t)$ decays at the rate of $O(t^{-\f{3}{2}(1-\f{1}{p})})$, while $u(t)$ decays at the rate of $O(  t^{-\f{9}{8}(1-\f{1}{p})-\f{1}{2}}  \log(2+t)        )   $. Moreover, we give the asymptotic expansion of solutions. We conclude that the leading term of $\Bh(t)$ is a combination of linear solution and two integrals from nonlinear coupling effects between velocity field and magnetic field, while for $B_3(t)$ the leading term is given by only the linear solution without the influence of velocity field. For higher order asymptotic expansion of solutions, we show that the velocity field does affect the asymptotic expansion of $B_3(t)$. Moreover, the nonlinear integrals of magnetic field turn out to be the second leading terms of higher order asymptotic expansion of $u(t)$. The obtained results reveal the crucial role played by the magnetic field in the asymptotic analysis of anisotropic incompressible MHD equations and are expected to enhance the understanding of incompressible MHD equations from the viewpoint of mathematics.
\end{abstract}
\maketitle

\tableofcontents

\section{Introduction and the main results}
\subsection{Background and motivation}

The mutual interactions between the electrically conducting fluids and magnetic field are described by magnetohydrodynamic equations. The theory of MHD finds its wide applications in astrophysics, plasmas, controlled nuclear fusion and engineering, among many others.

Assuming that the fluids are incompressible and the density of fluids is normalized to be one unit, one gets the following incompressible MHD system in three dimensions (see \cite{ST83}):
\begin{align}\label{MHD-1}
    \begin{cases}
    \partial_t u - \nu \Delta u + ( u \cdot \nabla ) u - ( B \cdot \nabla ) B + \nabla p = 0, & t>0 , x \in \mathbb{R}^3,\\
    \partial_t B -  \mu \Delta B + ( u \cdot \nabla ) B - ( B \cdot \nabla ) u = 0, & t>0 , x \in \mathbb{R}^3,\\
    \nabla \cdot u = \nabla \cdot B = 0 , & t \geqslant 0, x \in \mathbb{R}^3,\\
    u(0,x) = u_0(x), B(0,x) = B_0(x),
    & x \in \mathbb{R}^3.
    \end{cases}
\end{align}
Here the unknowns $u=(u_1,u_2,u_3)$ and $B=(B_1,B_2,B_3)$ represent the velocity field and magnetic field respectively. $t>0$ and $x\in \R^3$ signify the temporal and spatial variables. $p$ is the scalar pressure. $\nu \geq 0$ is the kinematic viscosity coefficient and $\mu \geq 0$ is the resistivity coefficient acting as the magnetic diffusion.

Below we first recall some known results on the stability and large time behavior of system (\ref{MHD-1}) according to the sign of $\nu,\mu$.
\begin{enumerate}[(i)]
\item{$\nu>0,\mu>0$. Mathematical results towards this type of system is nowadays classical. For large initial data, it admits global regular solutions in 2D and global existence of weak solutions in 3D. It also admits global small smooth solutions in 3D. We refer to \cite{ST83} for instance. The decay estimates in $L^2$-norm goes back to Kozono \cite{Ko89} and Mohgooner, Sarayker \cite{MS1}. Schonbek et al. \cite{SSS} proved the decay estimates in $L^2$-norm for both upper and lower bounds. We refer to \cites{GZ1,HH1} for more decay estimates in the whole space or half space. Notice that the velocity field and magnetic field always admit the same decay rate in these work.
}

\item{$\nu>0,\mu=0$. This type of equations has received lots of attention in recent years due to its physical importance and mathematical challenges. Suppose there exists a strong bounded solution, Agapito and Schonbek \cite{AS1} showed that the diffusion in the velocity is sufficient to prevent compensate oscillations between the two energies. Fefferman et al. \cites{FMDRR1,FMDRR2} proved the local well-posedness in Sobolev spaces. Lin et al. \cite{LXZ1} proved the global existence and uniqueness of small smooth solutions in $\R^2$; Xu and Zhang \cite{XZ2} proved the similar result of global small smooth solutions in $\R^3$. Both papers made use of the Lagrangian coordinates and imposed the so-called admissibility conditions on the initial magnetic field. Later, Abidi and Zhang \cite{AZ1} removed the admissibility condition and the large time decay rate of the solution is also obtained in $H^2(\R^3)$ if the initial magnetic field is a constant vector. In Eulerian coordinates, Pan et al. \cite{PZZ1} proved the global small smooth solutions on the periodic boxes $\mathbb{T}^3$, by imposing the odevity conditions of initial data. Tan and Wang \cite{TW1} obtained the global well-posedness and almost exponential decay rates in an infinite slab $\R^2 \times (0,1)$ via the two-tier energy method. In $\R^2$, Ren et al. \cite{RWXZ1} proved the global existence and decay estimates of smooth solution in $L^2$-norm by the anisotropic Littlewood–Paley analysis.

}

\item{$\nu=0,\mu>0$. Again by imposing the odevity conditions of initial data, Zhou and Zhu \cite{ZZ3} proved the global classical solutions on the periodic domain $\mathbb{T}^2$. Wei and Zhang \cite{WZ} obtained the global well-posedness provided that the initial data is sufficiently small in $H^4(\mathbb{T}^2)$ and the integral of initial magnetic field over $\mathbb{T}^2$ vanishes. Lei and Zhou \cite{LZ1} established the BKM's blowup criterion in $\R^2,\R^3$ and the global existence of weak solutions in $\R^2$.

}

\item{$\nu=0,\mu=0$. In the pioneering work, Bardos et al. \cite{BSS1} proved the global small smooth solutions in the presence of a strong background magnetic field. We refer to Cai, Lei \cite{CL1} and He et al. \cite{HXY1} for more recent global dynamics towards this system. Xu \cite{Xu1} recently considered the well-posedness and vanishing
domain thickness limit in three-dimensional thin domains $\R^2 \times (-\delta,\delta)$, $\delta>0$.

}

\end{enumerate}

Next we recall some known results on the incompressible MHD system with partial dissipation and magnetic diffusion. For clarity, we for the moment rewrite $\nu \Delta u$ in (\ref{MHD-1}) by $\nu_1 \p_1^2 u +\nu_2 \p_2^2 u+ \nu_3 \p_3^2 u$ and $\mu \Delta B$ by $\mu_1 \p_1^2 B +\mu_2 \p_2^2 B +\mu_3 \p_3^2 B$ with $\nu_i, \mu_i$ being non-negative constants. In the pioneering work, Cao and Wu \cite{CW1} proved the unique global classical solution to 2D incompressible MHD system with either $\nu_1>0,\nu_2=0,\mu_1=0,\mu_2>0$ or $\nu_1=0,\nu_2>0,\mu_1>0,\mu_2=0$. See also \cites{LWZ1,LJWY1} for more recent progress on the global stability and decay estimates of 2D anisotropic MHD system in the presence of strong background magnetic field. Wang and Wang \cite{WW1} proved the global small smooth solutions in the standard Sobolev space to 3D incompressible MHD system with either $\nu_1=\nu_2>0,\nu_3=0,\mu_1=0,\mu_2=\mu_3>0$ or $\nu_1=\nu_2>0,\nu_3=0,\mu_1=\mu_2>0,\mu_3=0$; see also Yue and Zhong \cite{YZ1} for the global well-posedness in anisotropic Sobolev space in case of $\nu_1=\nu_2>0,\nu_3=0,\mu_1=\mu_2>0,\mu_3=0$ and Shang and Zhai \cite{SZ22} for the decay estimates in $L^2$-norm for the same system. Very recently, Wu and Zhu \cite{WZ} proved the global stability for 3D incompressible MHD system with $\nu_1=\nu_2>0, \nu_3=0,\mu_1=\mu_2=0,\mu_3>0$ around a strong background magnetic field.

In view of the available literature, we find that there are few results on the large time behavior of solutions to the 3D incompressible MHD system with partial dissipation and magnetic diffusion, and in particular on the asymptotic profile of solutions. Specifically, it is not clear up to now how the mutual interactions between the velocity field and magnetic field affect the decay rates of global solutions and furthermore the asymptotic expansion of solutions to the 3D anisotropic incompressible MHD equations. This motivates the present work.

As is well-known in geophysical fluid dynamics, a suitable model for the approximation of turbulent diffusion involves mainly on the horizontal diffusion, while the vertical diffusion is negligible in some setting, see for instance Chapter 4 in the monograph \cite{PED} for more physical discussions. We also refer to \cites{CM1,CL2,CR1} for more physical backgrounds, simulations and numerical results on the incompressible MHD equations with anisotropic viscosity or anisotropic resistivity. As a consequence, in this paper, we consider two classes of anisotropic incompressible MHD equations in 3D. Without loss of generality, from now on the viscosity coefficient and resistivity coefficient are assumed to be one unit.

The first one that we consider in this paper is the incompressible MHD system with horizontal dissipation and full magnetic diffusion. The governing equations read as
\begin{align}\label{MHD-2}
    \begin{cases}
    \partial_t u -  \Deltah u + ( u \cdot \nabla ) u - ( B \cdot \nabla ) B + \nabla p = 0, & t>0 , x \in \mathbb{R}^3,\\
    \partial_t B -  \Delta B + ( u \cdot \nabla ) B - ( B \cdot \nabla ) u = 0, & t>0 , x \in \mathbb{R}^3,\\
    \nabla \cdot u = \nabla \cdot B = 0 , & t \geqslant 0, x \in \mathbb{R}^3,\\
    u(0,x) = u_0(x), B(0,x) = B_0(x),
    & x \in \mathbb{R}^3.
    \end{cases}
\end{align}
Observe that when $B=0$ system (\ref{MHD-2}) reduces to the well-known 3D incompressible Navier-Stokes system with horizontal dissipation. System of this type appears often in geophysical fluids. We next recall some previous results about this system. Global well-posedness in anisotropic Sobolev space goes back to the work of Chemin et al. \cite{CDGG1} and Iftimie \cite{If1}. We refer to \cites{CZ1,PZ1,ZT1,ZF1} for more results on the global well-posedness in anisotropic Besov space. See also \cites{GY22,LPZ1} for more recent progress. It should be emphasized that the large time behavior of this model appeared very recently. Ji et al. \cite{JWY} obtained the decay estimates in $L^2$-norms by imposing the smallness of initial velocity in the intersection of Sobolev space and negative ones; the regularity condition on initial data was soon relaxed by Xu and Zhang \cite{XZ}. Fujii \cite{F} later considered the decay estimates in the general $L^p$-norms and further obtained the asymptotic expansion of solutions. However, the evolution of magnetic field and its influence upon the fluid velocity make the asymptotic analysis for (\ref{MHD-2}) more involved. Moreover, as pointed out by Ji et al. \cite{JWY} and Lin et al. \cite{LJWY1}, the classical methods like the Fourier splitting technique does not seem to be effective for partially dissipated equations when it comes to the decay estimates of solutions.

The second one that we consider in this paper is the incompressible MHD system with full dissipation and horizontal magnetic diffusion. The governing equations read as
\begin{align}\label{MHD-3}
    \begin{cases}
    \partial_t u -  \Delta u + ( u \cdot \nabla ) u - ( B \cdot \nabla ) B + \nabla p = 0, & t>0 , x \in \mathbb{R}^3,\\
    \partial_t B -  \Deltah B + ( u \cdot \nabla ) B - ( B \cdot \nabla ) u = 0, & t>0 , x \in \mathbb{R}^3,\\
    \nabla \cdot u = \nabla \cdot B = 0 , & t \geqslant 0, x \in \mathbb{R}^3,\\
    u(0,x) = u_0(x), B(0,x) = B_0(x),
    & x \in \mathbb{R}^3.
    \end{cases}
\end{align}
When $B=0$ system (\ref{MHD-3}) reduces to the classical 3D incompressible Navier-Stokes equations. We shall not review the known results on global well-posedness. As to the decay estimates of solutions, we refer to Miyakawa \cites{MY1,MY2} for a series of work using Hardy spaces. Furthermore, Fujigaki and Miyakawa \cite{FM} obtained the asymptotic profiles of solutions. Obviously, the mutual interactions between the velocity field and magnetic field make the system (\ref{MHD-3}) more involved and one needs to explore the inner structures. In particular, the equation of magnetic field loses the dissipation in $x_3$-direction, bringing extra difficulty in analyzing the large time behavior of solutions.


\subsection{The main results}\label{ma}
The main results of this paper are summarized as the following four theorems.
Our first result is concerned with the decay estimates of global strong solutions in $L^p$-norm and the asymptotic profiles of solutions to the incompressible MHD system with horizontal dissipation and full magnetic diffusion.
\begin{thm}\label{thm1}
Let $s \geq 5$ be an integer. Then there exists $\epsilon_1=\epsilon_1(s)>0$ such that
\begin{enumerate}[(1)]
\item{
For any $u_0,B_0\in X^s(\R^3)$ satisfying $\nabla \cdot u_0=\nabla \cdot B_0=0, |x|B_0(x)\in L^1(\R^3)$, $\| (u_0,B_0) \|_{X^s} +\| |x|B_0(x) \|_{L^1_x} \leq \epsilon_1$, there exists a unique global strong solution $(u,B)\in C([0,\infty);X^s(\R^3))$ to \eqref{MHD-2} obeying the following decay estimates:
\beq\label{ma-1}
\| \nabla ^{\alpha}   \uh (t) \|_{L^p} \leq C t^{ -(1-\f{1}{p}) -\f{|\alphah|}{2}        }  \| (u_0,B_0) \|_{   \overline{X^{s} }  },
\eeq
\beq\label{ma-2}
\| \nablah ^{\alphah}   u_3 (t) \|_{L^p} \leq C t^{ -\f{3}{2}(1-\f{1}{p}) -\f{|\alphah|}{2}        }  \|(u_0,B_0)\|_{   \overline{X^{s} }  },
\eeq
\beq\label{ma-3}
\| \nabla ^{\alpha}   B (t) \|_{L^p} \leq C t^{ -\f{3}{2}(1-\f{1}{p}) -\f{1}{2}-\f{|\alpha|}{2}        }  \| (u_0,B_0) \|_{   \overline{X^{s} }  },
\eeq
for any $1\leq p \leq \infty,t>0$ and $\alpha=(\alphah,\alpha_3)\in (\mathbb{N} \cup \{ 0\})^2 \times (\mathbb{N} \cup \{ 0\})$ with $|\alpha|\leq 1$. Here we set
\[
\| (u_0,B_0) \|_{   \overline{X^{s} }  }:=
\| (u_0,B_0) \|_{ X^{s}  }+ \| |x|B_0(x) \|_{L^1_x} .
\]
}
\item{Furthermore, there hold the asymptotic profiles of solutions. For any $1\leq p \leq \infty$,
\begin{equation}\label{ma-4}
     \begin{split}
         \lim_{t\rightarrow \infty} t^{1-\f{1}{p}}  \Big\|  \uh (t,x)  & \ -G_{\rm h}(t,\xh) \int_{\R^2} u_{0,{\rm h}}  (y_{\rm h},x_3)  d y_{\rm h} \\
    & \  + G_{\rm h}(t,\xh)  \int_0^{\infty} \int_{\R^2}  \p_3 (u_3 \uh) (\tau, y_{\rm h},x_3)
    d y_{\rm h} d \tau\\
    & \  - G_{\rm h}(t,\xh)  \int_0^{\infty} \int_{\R^2}  \p_3 (B_3 \Bh) (\tau, y_{\rm h},x_3)
    d y_{\rm h} d \tau  \Big\|_{L^p_x}=0.      \\
     \end{split}
\end{equation}
For any $1\leq p < \infty$,
\begin{equation}\label{ma-5}
\lim_{t \rightarrow \infty} t^{ \f{3}{2} (1-\f{1}{p})   }
\Big \|
u_3(t,x)- G_{\rm h}(t,\xh) \int_{\R^2} u_{0,3} (y_{\rm h},x_3)  d y_{\rm h}
\Big\|_{L^p_x}=0.
\end{equation}
}
\end{enumerate}
\end{thm}

\begin{rem}
\begin{enumerate}[(1)]
\item{
The asymptotic limit (\ref{ma-4}) shows that the leading term of $\uh(t)$ is the combination of linear solution and two integrals resulting from nonlinear effects. In particular, the magnetic field does affect the leading term.
}
\item{
It follows from the asymptotic limit (\ref{ma-5}) that the leading term of $u_3(t)$ is given only by the linear solution without the influence of magnetic field. Notice also that this limit fails if $p=\infty$.
}
\item{
From the decay estimate of magnetic field (\ref{ma-3}) we see that for any $1 \leq p \leq \infty$,
\begin{equation}\nonumber
\lim_{t \rightarrow \infty} t^{ \f{3}{2} (1-\f{1}{p})   }
 \|
B (t,x)\|_{L^p_x}=0;
\end{equation}
whence, upon utilizing the basic fact that $\int_{\R^3} B_0(y) dy=0$ (see the proof of Lemma \ref{lm3-2}),
\begin{equation}\nonumber
\lim_{t \rightarrow \infty} t^{ \f{3}{2} (1-\f{1}{p})   }
 \left \|
B (t,x)-G(t,x) \int_{\R^3} B_0(y) dy  \right\|_{L^p_x}=0.
\end{equation}
This implies that the leading term of $B(t)$ is given by the trivial linear solution.
}

\end{enumerate}
\end{rem}

Our second result treats the decay estimates of global strong solutions in $L^p$-norm and the asymptotic profiles of solutions to the incompressible MHD system with full dissipation and horizontal magnetic diffusion.
\begin{thm}\label{thm2}
Let $s \geq 5$ be an integer. Then there exists $\delta_1=\delta_1(s)>0$ such that
\begin{enumerate}[(1)]
\item{
For any $u_0,B_0\in X^s(\R^3)$ satisfying $\nabla \cdot u_0=\nabla \cdot B_0=0, |x| u_0(x) \in L^1(\R^3),\nablah u_0 \in L^1(\R^3)$, $\| (u_0,B_0) \|_{X^s} + \| |x|u_0(x) \|_{L^1_x} +\| \nablah u_0 \|_{L^1} \leq \delta_1$, there exists a unique global strong solution $(u,B)$ to \eqref{MHD-3} obeying the following decay estimates:
\beq\label{ma-8}
\| \nabla ^{\alpha}   \Bh (t) \|_{L^p} \leq C t^{ -(1-\f{1}{p}) -\f{|\alphah|}{2}        }  \| (u_0,B_0) \|_{X^{s}_{\#}  },
\eeq
\beq\label{ma-9}
\| \nablah ^{\alphah}   B_3 (t) \|_{L^p} \leq C t^{ -\f{3}{2}(1-\f{1}{p}) -\f{|\alphah|}{2}        }  \|(u_0,B_0)\|_{X^{s}_{\#}  },
\eeq
\beq\label{ma-10}
\| \nabla ^{\alpha}   u (t) \|_{L^p} \leq C (1+t)^{ -\f{9}{8}(1-\f{1}{p}) -\f{1+|\alpha|}{2}        } \log(2+t) \| (u_0,B_0) \|_{X^{s}_{\#}  },
\eeq
for any $1\leq p \leq \infty,t>0$ and $\alpha=(\alphah,\alpha_3)\in (\mathbb{N} \cup \{ 0\})^2 \times (\mathbb{N} \cup \{ 0\})$ with $|\alpha|\leq 1$. Here we denoted by
\[
\| (u_0,B_0) \|_{X^{s}_{\#}  }:=
\| (u_0,B_0) \|_{ X^{s}  }+ \| |x|u_0(x) \|_{L^1_x}
+\| \nablah u_0 \|_{L^1}.
\]
}
\item{Furthermore, there hold the asymptotic profiles of solutions. For any $1\leq p \leq \infty$,
\begin{equation}\label{ma-10-1}
     \begin{split}
         \lim_{t\rightarrow \infty} t^{1-\f{1}{p}}  \Big\|  \Bh (t,x)  & \ -G_{\rm h}(t,\xh) \int_{\R^2} B_{0,{\rm h}}  (y_{\rm h},x_3)  d y_{\rm h} \\
    & \  + G_{\rm h}(t,\xh)  \int_0^{\infty} \int_{\R^2}  \p_3 (u_3 \Bh) (\tau, y_{\rm h},x_3)
    d y_{\rm h} d \tau\\
    & \  - G_{\rm h}(t,\xh)  \int_0^{\infty} \int_{\R^2}  \p_3 (B_3 \uh) (\tau, y_{\rm h},x_3)
    d y_{\rm h} d \tau  \Big\|_{L^p_x}=0.      \\
     \end{split}
\end{equation}
For any $1\leq p < \infty$,
\begin{equation}\label{ma-10-2}
\lim_{t \rightarrow \infty} t^{ \f{3}{2} (1-\f{1}{p})   }
\Big \|
B_3(t,x)- G_{\rm h}(t,\xh) \int_{\R^2} B_{0,3} (y_{\rm h},x_3)  d y_{\rm h}
\Big\|_{L^p_x}=0.
\end{equation}
}
\end{enumerate}

\end{thm}

\begin{rem}
\begin{enumerate}[(1)]
\item{
We know from the asymptotic limit (\ref{ma-10-1}) that the leading term of $\Bh(t)$ is the combination of linear solution and two integrals resulting from nonlinear coupling effects. In particular, the velocity field does affect the leading term.
}
\item{
It follows from the asymptotic limit (\ref{ma-10-2}) that the leading term of $B_3(t)$ is given by only the linear solution without the influence of velocity field. However, this limit fails if $p=\infty$.
}
\item{
In light of the decay estimate of velocity field (\ref{ma-10}), we see that for any $1 \leq p \leq \infty$,
\begin{equation}\nonumber
\lim_{t \rightarrow \infty} t^{ \f{9}{8} (1-\f{1}{p})   }
 \|
u (t,x)\|_{L^p_x}=0;
\end{equation}
whence, upon realizing the basic fact that $\int_{\R^3} u_0(y) dy=0$ (see the proof of Lemma \ref{lm3-2}),
\begin{equation}\nonumber
\lim_{t \rightarrow \infty} t^{ \f{9}{8} (1-\f{1}{p})   }
 \left \|
u (t,x)-G(t,x) \int_{\R^3} u_0(y) dy  \right\|_{L^p_x}=0.
\end{equation}
This implies that the leading term of $u(t)$ is given by the trivial linear solution. On the one hand, we see that
\[
 \left \|
u (t,x)-G(t,x) \int_{\R^3} u_0(y) dy  \right\|_{L^p_x}
\]
decays to zero faster than $t^{ -\f{9}{8} (1-\f{1}{p})   } $. On the other hand, for system (\ref{MHD-2}) we obtain that
\[
 \left \|
B (t,x)-G(t,x) \int_{\R^3} B_0(y) dy  \right\|_{L^p_x}
\]
decays to zero faster than $t^{ -\f{3}{2} (1-\f{1}{p})   } $. These different decay rates reveal the distinct structures between the incompressible MHD systems (\ref{MHD-2}) and (\ref{MHD-3}).
}

\end{enumerate}
\end{rem}

Given additionally the spatial decay assumption for the initial velocity field, we are able to strengthen the asymptotic limits and obtain higher order asymptotic expansion of the vertical component of velocity field. For the incompressible MHD system with horizontal dissipation and full magnetic diffusion, we have
\begin{thm}\label{thm3}
Let $s \geq 9$ be an integer. Then there exists $\epsilon_2=\epsilon_2(s)>0$ such that for any $u_0,B_0\in X^s(\R^3)$ satisfying $\nabla \cdot u_0=\nabla \cdot B_0=0, |x| B_0(x)\in L^1(\R^3)$, $\| (u_0,B_0) \|_{ \overline{X^s} } \leq \epsilon_2$ and $|\xh| u_0(x) \in L^1(\R^2_{\xh};(L^1 \cap L^{\infty})(\R_{x_3})  )$, the unique global solution $(u,B)$ to \eqref{MHD-2} obeys
\begin{enumerate}[(1)]
\item{ For any $2<p \leq \infty$, it holds that
\begin{equation}\label{ma11}
     \begin{split}
          \Big\|  \uh (t,x)  & \ -G_{\rm h}(t,\xh) \int_{\R^2} u_{0,{\rm h}}  (y_{\rm h},x_3)  d y_{\rm h} \\
    & \  + G_{\rm h}(t,\xh)  \int_0^{\infty} \int_{\R^2}  \p_3 (u_3 \uh) (\tau, y_{\rm h},x_3)
    d y_{\rm h} d \tau\\
    & \  - G_{\rm h}(t,\xh)  \int_0^{\infty} \int_{\R^2}  \p_3 (B_3 \Bh) (\tau, y_{\rm h},x_3)
    d y_{\rm h} d \tau  \Big\|_{L^p_x}      \\
   \leq  &\  C \| (u_0,B_0) \|_{\widetilde{X^s}      } t^{ -(1-\f{1}{p}) -\f{1}{2}   }
   \log t   ,    \\
     \end{split}
\end{equation}
for any $t \geq 2$. Here we set
\[
\| (u_0,B_0) \|_{\widetilde{X^s}      }:=
\| (u_0,B_0) \|_{X^s     }  + \| |\xh| u_0 (x)     \|_{ L^1(   \R^2_{\xh};(L^1\cap L^{\infty})  (\R_{x_3})  )          }
+\| |x|B_0(x) \|_{L^1_x}.
\]
}
\item{For any $1\leq p \leq \infty$, it holds that
\begin{equation}\label{ma12}
     \begin{split}
       \Big \|
u_3(t,x)-  &\   G_{\rm h}(t,\xh) \int_{\R^2} u_{0,3} (y_{\rm h},x_3)  d y_{\rm h}
\Big\|_{L^p_x}  \\
    \leq  & \ \left\{\begin{aligned}
& C  \| (u_0,B_0) \|_{\widetilde{X^s}      }   t^{ -\f{3}{2} (1-\f{1}{p}) -\f{1}{2p}    }      ,\,   1<p \leq \infty       \\
& C  \| (u_0,B_0) \|_{\widetilde{X^s}      } t^{-\f{1}{2}} \log t  ,\,\,\,\, p=1     \\
\end{aligned}\right.
     \end{split}
\end{equation}
for any $t \geq 2$.
}
\item{ For any $1<p\leq \infty$, it holds that
\begin{equation}\label{ma13}
     \begin{split}
         \lim_{t\rightarrow \infty} t^  {  \f{3}{2}(1-\f{1}{p}) +\f{1}{2p} }  \Big\|  u_3 (t,x)  & \ -G_{\rm h}(t,\xh) \int_{\R^2} u_{0,3}  (y_{\rm h},x_3)  d y_{\rm h} \\
    & \  + \nablah G_{\rm h}(t,\xh)  \cdot  \int_{\R^2}  \yh u_{0,3} ( y_{\rm h},x_3)
    d y_{\rm h} \\
    & \  -  \nablah  G_{\rm h}(t,\xh)  \cdot \int_0^{\infty} \int_{\R^2}   (u_3 \uh) (\tau, y_{\rm h},x_3)  d y_{\rm h} d \tau          \\
    &\  +  \nablah  G_{\rm h}(t,\xh)  \cdot \int_0^{\infty} \int_{\R^2}   (B_3 \Bh) (\tau, y_{\rm h},x_3)  d y_{\rm h} d \tau     \Big\|_{L^p_x}=0.
     \end{split}
\end{equation}

}
\item{ For any $1\leq p \leq \infty$, it holds that
\begin{equation}\label{ma14}
     \begin{split}
         \lim_{t\rightarrow \infty} t^  {  \f{3}{2}(1-\f{1}{p}) +\f{1}{2} }  \Big\|  B_i (t,x)  & \ +  \nabla G(t,x)  \cdot  \int_{\R^3}   y B_{0,i} ( y)
    d y \\
    & \  +  \sum_{j=1}^3 \p_j  G(t,x) \int_0^{\infty} \int_{\R^3}   (B_i u_j ) (\tau, y )  d y d \tau          \\
     & \  - \sum_{j=1}^3 \p_j  G(t,x) \int_0^{\infty} \int_{\R^3}   (u_i B_j ) (\tau, y )  d y d \tau     \Big\|_{L^p_x}=0       \\
     \end{split}
\end{equation}
for any $i=1,2,3$.
}
\end{enumerate}
\end{thm}

In the last theorem, we show the higher order asymptotic expansions of solutions for the incompressible MHD system with full dissipation and horizontal magnetic diffusion. 
\begin{thm}\label{thm4}
Let $s \geq 9$ be an integer. Then there exists $\delta_2=\delta_2(s)>0$ such that for any $u_0,B_0\in X^s(\R^3)$ satisfying $\nabla \cdot u_0=\nabla \cdot B_0=0, |x| u_0(x)\in L^1(\R^3)$, $\| (u_0,B_0) \|_{X^s _{\#} } \leq \delta_2$ and $|\xh| B_0(x) \in L^1(\R^2_{\xh};(L^1 \cap L^{\infty})(\R_{x_3})  )$, the unique global solution $(u,B)$ to \eqref{MHD-3} obeys
\begin{enumerate}[(1)]
\item{For any $1\leq p \leq \infty$, it holds that
\begin{equation}\label{ma16}
     \begin{split}
       \Big \|
B_3(t,x)-  &\   G_{\rm h}(t,\xh) \int_{\R^2} B_{0,3} (y_{\rm h},x_3)  d y_{\rm h}
\Big\|_{L^p_x}  \\
    \leq  & \ C  \| (u_0,B_0) \|_{X^s_{\# \#}    }  t^{ -\f{3}{2} (1-\f{1}{p}) -\f{1}{2p}    }
     \end{split}
\end{equation}
for any $t \geq 2$. Here we set
\[
\| (u_0,B_0) \|_{X^s_{\# \#}    }:=
\| (u_0,B_0) \|_{X^s     }  + \| |\xh| B_0 (x)     \|_{ L^1(   \R^2_{\xh};(L^1\cap L^{\infty})  (\R_{x_3})  )          }
+\| |x| u_0(x) \|_{L^1_x}.
\]
}
\item{ For any $1\leq p\leq \infty$, it holds that
\begin{equation}\label{ma17}
     \begin{split}
         \lim_{t\rightarrow \infty} t^  {  \f{3}{2}(1-\f{1}{p}) +\f{1}{2p} }  \Big\|  B_3 (t,x)  & \ -G_{\rm h}(t,\xh) \int_{\R^2} B_{0,3}  (y_{\rm h},x_3)  d y_{\rm h} \\
    & \  + \nablah G_{\rm h}(t,\xh)  \cdot  \int_{\R^2}  \yh B_{0,3} ( y_{\rm h},x_3)
    d y_{\rm h} \\
    & \  +  \nablah  G_{\rm h}(t,\xh)  \cdot \int_0^{\infty} \int_{\R^2}   (B_3 \uh) (\tau, y_{\rm h},x_3)  d y_{\rm h} d \tau          \\
    &\  -  \nablah  G_{\rm h}(t,\xh)  \cdot \int_0^{\infty} \int_{\R^2}   (u_3 \Bh) (\tau, y_{\rm h},x_3)  d y_{\rm h} d \tau     \Big\|_{L^p_x}=0.
     \end{split}
\end{equation}

}
\item{ For any $1\leq p \leq \infty$ and any $0<\sigma< \f{1}{2}$, there hold the higher order asymptotic expansion of velocity field:
\begin{equation}\label{ma18}
     \begin{split}
         \lim_{t\rightarrow \infty} t^  {  \f{3}{2}(1-\f{1}{p}) +\f{1}{2} -\sigma }  \Big\|  u_1 (t,x)  & \ +  \nabla G(t,x)  \cdot  \int_{\R^3}    y u_{0,1} ( y)
    d y \\
    & \  -   \int_0^{t} e^{(t-\tau)\Delta } \nablah  \cdot (B_1 \Bh ) (\tau )   d \tau          \\
     & \  - \sum_{k,l=1}^2   \int_0^{t}    \p_1 \p_k \p_l N(t-\tau) \ast  (B_k B_l ) (\tau)  d \tau     \Big\|_{L^p_x}=0 ,      \\
     \end{split}
\end{equation}
\begin{equation}\label{ma19}
     \begin{split}
       \lim_{t\rightarrow \infty} t^  {  \f{3}{2}(1-\f{1}{p}) +\f{1}{2} -\sigma }   \Big\|  u_2 (t,x)  & \ +  \nabla G(t,x)  \cdot  \int_{\R^3}    y u_{0,2} ( y)
    d y \\
    & \  -  \int_0^{t} e^{(t-\tau)\Delta } \nablah  \cdot (B_2 \Bh ) (\tau )   d \tau          \\
     & \  - \sum_{k,l=1}^2   \int_0^{t}    \p_2 \p_k \p_l N(t-\tau) \ast  (B_k B_l ) (\tau)  d \tau     \Big\|_{L^p_x}=0 ,      \\
     \end{split}
\end{equation}
\begin{equation}\label{ma20}
     \begin{split}
       \lim_{t\rightarrow \infty} t^  {  \f{3}{2}(1-\f{1}{p}) +\f{1}{2} -\sigma }   \Big\|  u_3 (t,x)  & \ +  \nabla G(t,x)  \cdot  \int_{\R^3}    y u_{0,3} ( y)
    d y \\
     & \  - \sum_{k,l=1}^2   \int_0^{t}    \p_3 \p_k \p_l N(t-\tau) \ast  (B_k B_l ) (\tau)  d \tau     \Big\|_{L^p_x}=0 .       \\
     \end{split}
\end{equation}
}
Here, $N(t,x)$ is the function with the following representation:
\[
		N(t,x)=\int_0^{\infty}     (4 \pi  (t+s) ) ^{-\f{3}{2}}    e^{-\frac{|x|^2}{4 (t+s)}}    ds
		=\int_0^{\infty}  G( t+s,x)    ds.
\]
\end{enumerate}

\end{thm}
\begin{rem}
\begin{enumerate}[(1)]

\item{
In contrast with (\ref{ma16}) for $B_3(t)$, we are unable to show the corresponding asymptotic limit for $\Bh(t)$. We include some calculations in Section \ref{app01} that are responsible for this. In particular, 
the estimate (\ref{hm23}) gives no useful information and prevents us from getting the desired estimate.
}

\item{
Due to the slowly decay estimates of magnetic field and the structure of equations (\ref{MHD-3}), we see from (\ref{ma18})-(\ref{ma20}) that the leading terms of higher order asymptotic expansion of velocity field consist mainly of the nonlinear effects from magnetic field, apart from the linear solution.
}

\end{enumerate}
\end{rem}


The procedure of proof is roughly split into the following four steps.
\begin{itemize}
\item{
We reformulate the equations in suitable integral forms. This is crucial for obtaining the decay estimates of solutions in the general $L^p$-norms, for any $1 \leq p \leq \infty$.
}

\item{
We analyze in details for all the Duhamel terms and derive the asymptotic expansions of solutions for some slowly decaying Duhamel terms.
}

\item{
Based on the above estimates and the estimates of linearized equations, we show the decay estimates of nonlinear solutions by the bootstrapping argument, together with the asymptotic profiles of nonlinear solutions.
}

\item{
Given more spatial decay assumptions of initial data, we are able to provide more elaborate asymptotic limits and moreover on higher order asymptotic expansions of solutions.
}

\end{itemize}

Several new ingredients of this paper are summarized below.
\begin{itemize}
\item {For the incompressible MHD system with horizontal dissipation and full magnetic diffusion (\ref{MHD-2}), together with the incompressible MHD system with full dissipation and horizontal magnetic diffusion (\ref{MHD-3}), we establish the decay estimates of global solutions in the general $L^p$-norms, for any $1\leq p \leq \infty$. Moreover, these decay estimates are optimal in the sense that they coincide with the corresponding linearized equations, with the only exception of $u(t)$ to the system (\ref{MHD-3}), cf. (\ref{ma-10}).
}

\item{
In order to understand the interplay between velocity field and magnetic field for the incompressible MHD system with horizontal dissipation and full magnetic diffusion (\ref{MHD-2}), we focus attention on the asymptotic expansions of solutions. We show that the leading term of $\uh(t)$ consists of the linear solution combined with the nonlinear effects from both velocity field and magnetic field, while the leading term of $u_3(t)$ is given by only the linear solution. Furthermore, the decay rates for these limits are shown under suitable assumptions on the initial data. However, for the higher order asymptotic expansions, we shall prove that magnetic field does affect the asymptotic expansion of $u_3(t)$. Moreover, the higher order asymptotic expansion of magnetic field will also be established.
}

\item{
For the sake of comprehending the mutual interactions between velocity field and magnetic field for the incompressible MHD system with full dissipation and horizontal magnetic diffusion (\ref{MHD-3}), we again concentrate on the asymptotic expansions of solutions. We shall show that the asymptotic profile of $\Bh(t)$ is given by the combination of linear solution and two integrals from the nonlinear coupling between velocity field and magnetic field, while the leading term of $B_3(t)$ consists only of the linear solution. Nevertheless, for higher order asymptotic expansions, it turns out that the velocity field does affect the asymptotic expansion of $B_3(t)$. Moreover, the nonlinear integrals of magnetic field become the second leading terms of asymptotic expansion of $u(t)$, mainly due to the special structure of equations.
}
\end{itemize}

The rest of this paper is arranged as follows. In Section \ref{li}, we recall some classical decay estimates related to the heat kernels in 2D and 3D; furthermore, we also give the decay estimates for the linearized equations. Section \ref{hd} is devoted to the anisotropic incompressible MHD system (\ref{MHD-2}). Precisely speaking, we first reformulate the equations in an integral form, then study the decay estimates for the Duhamel terms, next consider the asymptotic expansions for the Duhamel terms and finally finish the decay estimates of nonlinear solutions mainly via the bootstrapping argument. The asymptotic profiles of solutions follow largely from the decay estimates and asymptotic expansion of Duhemel terms, together with the asymptotic analysis for linear solutions. Following the same spirit, we handle the anisotropic incompressible MHD system (\ref{MHD-3}) in Section \ref{hm}. Under suitable assumptions of initial data, we show some higher order asymptotic expansions of solutions in Section \ref{hie}. The paper ends up with the conclusion and some technical remarks.

In the end of this section, we list some notations that will be frequently used throughout this paper.
{\flushleft{\bf{Notations.} }} For any vector $v\in \R^3$, we denote by $v_{\rm h}$ the horizontal components of $v$. The operator $\Deltah$ means the Laplacian operator acting only on horizontal variables, namely $\Deltah:= \p_1^2+\p_2^2$. Similarly, $\nablah:=(\p_1,\p_2)$ signifies the gradient operator acting only on horizontal variables.
We denote by $G(t,x)$ the 3D Gaussian:
\[
G(t,x):= (4\pi t)^{-\f{3}{2}} e^{ -\f{|x|^2}{4t}        },\,\,
(t,x)\in (0,\infty)\times \R^3,
\]
and $\Gh(t,\xh)$ the 2D Gaussian:
\[
\Gh(t,\xh):= (4\pi t)^{-1}  e^{ -\f{|\xh|^2}{4t}        },\,\,
(t,\xh)\in (0,\infty)\times \R^2.
\]
Let $1\leq p,q \leq \infty$. $L^p(\R^3)$ is the usual Lebesgue space with norm $\| \cdot \|_{L^p_x}$ or simply $\|\cdot \|_{L^p}$. $L^{p}_{\rm h} L^{q}_{\rm v} (\R^3)$ signifies the anisotropic Lebesgue space $L^p(\R^2_{\xh};L^q(\R_{x_3}))$ with norm $\| f \|_{L^{p}_{\rm h} L^{q}_{\rm v}  }:=  \|     \| f(\xh,x_3)  \|_{L^q_{x_3}}        \|_{L^p_{\xh}}        $. For brevity, we will adopt the following abbreviation of norms:
\begin{equation}\nonumber
     \begin{split}
         \| f \| _{  L^1_{\xh}(W^{1,1}\cap W^{1,\infty})_{x_3}         }  & \  :=
\| f \| _{  L^1(  {\R^2_{\xh} } ; (W^{1,1}\cap W^{1,\infty})(\R_{x_3})       )  } , \\
   \| f \|_{ L^1_{\xh} W^{2,1}_{x_3}     }  & \  :=
\| f \|_{ L^1 (\R^2_{\xh};  W^{2,1}(\R_{x_3})   ) }. \\
     \end{split}
\end{equation}
For any $s\in \mathbb{N}$, we set
\[
X^s(\R^3):= H^s(\R^3)\cap   L^1(  {\R^2_{\xh} } ; (W^{1,1}\cap W^{1,\infty})(\R_{x_3})       ) .
\]
The same letter $C$ denotes various generic positive constants changing from line to line. Sometimes we use $C_1,C_2$ to distinguish them.

\section{Linear estimates}\label{li}

In this section, we shall recall several well-known decay estimates related to heat kernels. As an immediate application, we show the decay estimates for the linearized equations under consideration in $L^p$-norms.

\subsection{Preliminary lemmas}

\begin{lemm}\label{lm1}
\begin{enumerate}[(1)]
\item{
For any $\alpha=(\alphah, \alpha_3) \in (\mathbb{N} \cup \{ 0\})^2 \times (\mathbb{N} \cup \{ 0\})$ and $m\in \mathbb{N} \cup \{ 0\} $, there exists $C>0$ depending only on $\alpha$ and $m$ such that in $3D$
\[
\| |x|^m \nabla^{\alpha} G(t,x) \|_{ L^p_x }
\leq C  t^{ -\f{3}{2}  (1-\f{1}{p}) -\f{ |\alpha|  }{2} +\f{m}{2}         }
\]
and in $2D$
\[
\| |\xh |^m \nablah^{\alphah} \Gh(t,\xh) \|_{ L^p_{\rm h } }
\leq C  t^{ -  (1-\f{1}{p}) -\f{ |\alphah|  }{2} +\f{m}{2}         }
\]
hold for any $1\leq p \leq \infty$ and $t>0$. In particular,
\[
\|\nabla^{\alpha}  e ^{t \Deltah } f \|_{  L^{p_2}_{\rm h} L^{q}_{\rm v}      }
\leq C t^{  -(\f{1}{p_1}-\f{1}{p_2})   - \f{ |\alphah|  }{2}          }
\| \p_3^{\alpha_3} f \|_{  L^{p_1}_{\rm h} L^{q}_{\rm v}      }
\]
for any $1\leq p_1\leq p_2 \leq \infty,1\leq q \leq \infty$ and $t>0$ and all $f$ obeying $\p_3^{\alpha_3} f \in  L^{p_1}_{\rm h} L^{q}_{\rm v} (\R^3)  $;
\begin{equation}\nonumber
     \begin{split}
         \||\xh| \nabla^{\alpha}  e ^{t \Deltah } f \|_{  L^{p}_{\rm h} L^{q}_{\rm v}      }
\leq   & \  C t^{  -(\f{1}{p_1}-\f{1}{p})   - \f{ |\alphah|  }{2} +\f{1}{2}         }
\| \p_3^{\alpha_3} f \|_{  L^{p_1}_{\rm h} L^{q}_{\rm v}      } \\
     & \ +C        t^{  -(\f{1}{p_2}-\f{1}{p})   - \f{ |\alphah|  }{2}  }
\| |\xh|  \p_3^{\alpha_3} f \|_{  L^{p_2}_{\rm h} L^{q}_{\rm v}      }        \\
     \end{split}
\end{equation}
for any $1\leq p_1, p_2 \leq p \leq \infty,1\leq q \leq \infty$ and $t>0$ and all $f$ obeying $\p_3^{\alpha_3} f \in L^{p_1}_{\rm h} L^{q}_{\rm v} (\R^3) ,  |\xh|\p_3^{\alpha_3} f(x) \in  L^{p_2}_{\rm h} L^{q}_{\rm v} (\R^3)$.
}
\item{
Suppose that $m=0,1$ and $1\leq p,q \leq \infty$. For any $f$ obeying $|x|^{m} f(x)\in L^1(\R^3)$, it holds the asymptotic expansion in 3D:
\[
\lim_{t\rightarrow \infty}
t^{ \f{3}{2}(1-\f{1}{p})+\f{m}{2}  }
\left\|   e ^{t \Delta } f (x) - \sum_{|\alpha|\leq m}  \nabla^{\alpha} G (t,x)
\int_{\R^3}      (-y )^{\alpha} f(y) d y \right \|_{  L^{p}_{x}      }=0;
\]
for any $f$ obeying $|\xh|^m f(x) \in L^1( \R^2_{\xh};(L^1\cap L^{\infty})(\R_{x_3})   )
$, it holds the asymptotic expansion in 2D:
\[
\lim_{t\rightarrow \infty}
t^{ (1-\f{1}{p})+\f{m}{2}  }
\left\|   e ^{t \Deltah } f (x) - \sum_{|\alphah|\leq m}  \nablah^{\alphah} \Gh (t,\xh)
\int_{\R^2}      (-\yh )^{\alphah} f(\yh,x_3) d \yh \right \|_{  L^{p}_{\rm h} L^{q}_{\rm v}      }=0.
\]
}
\item{
Let $1\leq p,q \leq \infty$. For any $f$ obeying $|x|f(x) \in L^1(\R^3)$, it holds in 3D:
\[
\left\|   e ^{t \Delta } f (x) -  G (t,x)
\int_{\R^3}    f(y) d y \right \|_{  L^{p}_{x}      }
\leq C  t^{ -\f{3}{2}(1-\f{1}{p})-\f{1}{2}  } \| |x|f(x) \|_{L^1_x}
;
\]
for any $f$ obeying $|\xh|f(x) \in L^{1}_{\rm h} L^{q}_{\rm v} (\R^3) $, it holds in 2D:
\[
\left\|   e ^{t \Deltah } f (x) -  \Gh (t,\xh)
\int_{\R^2}    f(\yh,x_3) d \yh \right \|_{  L^{p}_{\rm h} L^{q}_{\rm v}      }
\leq C  t^{ -(1-\f{1}{p})-\f{1}{2}  } \| |\xh|f(x) \|_{L^1_{\rm h}  L^q_{\rm v}}.
\]
}
\end{enumerate}
\end{lemm}

The proof of item $(1)$ in Lemma \ref{lm1} is standard. For the proof of items $(2),(3)$ in Lemma \ref{lm1} we may refer to \cite{FM} for the details.

To proceed, we recall the enhanced dissipation mechanism for $e^{t \Deltah} v_{0,3}$ with $v_0$ integrable and divergence free. For anisotropic incompressible Navier-Stokes system, this celebrated mechanism was initiated in the work of Xu and Zhang \cite{XZ} in $L^2$ framework and further studied by Fujii
\cite{F} in $L^p$ framework as well as in the asymptotic expansions of solutions.

\begin{lemm}\label{lm2}
\begin{enumerate}[(1)]
\item{
Let $1\leq p,q \leq \infty$. Assume that $v_0=(v_{0,1},v_{0,2},v_{0,3}) \in L^1(\R^3)$ obeying $\Grad \cdot v_0=0$. For any $t>0$ it holds that
\[
\|  e^{t \Deltah} v_{0,3} \|_{L^p_{\rm h}  L^q_{\rm v}}
\leq C t^{ -(1-\f{1}{p})    -\f{1}{2} (1-\f{1}{q})       }  \| v_0 \|_{L^1_x}.
\]
Assume furthermore that $|\xh| v_0\in L^1(\R^3)$, we have
\[
\|  |\xh|     e^{t \Deltah} v_{0,3} \|_{L^1_{\rm h}  L^{\infty}_{\rm v}}
\leq C   \| (1+|\xh|)v_0(x) \|_{L^1_x}.
\]
}
\item{
Let $1\leq p \leq \infty$ and $m=0,1$. Suppose that $\Grad \cdot v_0=0$ and $|\xh|^m v_0(x) \in L^1( \R^2_{\xh};(L^1\cap L^{\infty})(\R_{x_3})   ) $. Then there exists a non-negative remainder $\mathcal{R}_{p,m}(t)$ obeying $ \mathcal{R}_{p,m}(t) \rightarrow 0$ as time tends to infinity and the enhanced asymptotic expansion:
\begin{equation}\nonumber
     \begin{split}
         \Big\|   e^{t \Deltah} v_{0,3}(x) -    & \ \sum_{|\alphah|\leq m}  \nablah^{\alphah} \Gh (t,\xh)
\int_{\R^2}      (-\yh )^{\alphah} v_{0,3}(\yh,x_3) d \yh \Big\|_{L^p_x}  \\
    \leq  & \  C  t ^{  -\f{3}{2} (1-\f{1}{p}) -\f{m}{2}       }
        \mathcal{R}_{p,m}(t)^{\f{1}{p}}   \|  |\xh|^m v_0(x)  \|_{L^1_x}^{1-\f{1}{p}} .            \\
     \end{split}
\end{equation}
}
\end{enumerate}
\end{lemm}

We refer to \cite{F} for a detailed proof.

When reformulating the equations in integral form, we will encounter $K(t,x)$ defined as
\[
	K(t,x):=\int_0^{\infty}\frac{e^{-\frac{|\xh|^2}{4 (t+s)}}}{4 \pi  (t+s)}\frac{e^{-\frac{x_3^2}{4 s}}}{(4 \pi  s)^{\frac{1}{2}}}ds
		=\int_0^{\infty}\Gh( t+s,\xh)\Gv( s,x_3)ds.
\]
The following lemma then gives the decay estimates related to $K(t,x)$.
\begin{lemm}\label{lm3}
Let $1\leq p,q \leq \infty$, $(\beta,\gamma)\in (\mathbb{N}\cup \{ 0\})^2\times (\mathbb{N}\cup \{ 0\})$ and $m \in \mathbb{N}\cup \{ 0\}$.
\begin{enumerate}[(1)]
\item{
Suppose that $|\beta|+\gamma> \f{2}{p}+\f{1}{q}-1+m$. It holds for any $t>0$ that
\[
\Big\|  |\xh|^{m} \nablah^{\beta} (-\Deltah)^{  \f{\gamma}{2}   } K(t,x)  \Big\|_{L^p_{\rm h}  L^q_{\rm v}}
\leq C t ^{  -(1-\f{1}{p}) -\f{1}{2} (1-\f{1}{q}) -\f{ |\beta|+\gamma-2       }{2}                                +\f{m}{2}        }  .
\]
}
\item{
Suppose that $1\leq p_1,p_2\leq p, 1\leq q_1,q_2\leq q$ obeying
\[
\left\{
\begin{aligned}
\left( \f{1}{p_1} -\f{1}{p} \right) +\f{1}{2} \left( \f{1}{q_1} -\f{1}{q} \right) +\f{|\beta|+\gamma-3}{2} & >0, \\
\left( \f{1}{p_2} -\f{1}{p} \right) +\f{1}{2} \left( \f{1}{q_2} -\f{1}{q} \right) +\f{|\beta|+\gamma-2}{2} & >0. \\
\end{aligned}
\right.
\]
It holds for any $t>0$ that
\begin{equation}\nonumber
     \begin{split}
      \Big\|  |\xh| \nablah^{\beta} (-\Deltah)^{  \f{\gamma}{2}   } K\ast f(x)  \Big\|_{L^p_{\rm h}  L^q_{\rm v}}   \leq  & \ C t ^{  -(\f{1}{p_1}-\f{1}{p}) -\f{1}{2} (\f{1}{q_1}-\f{1}{q}) -\f{ |\beta|+\gamma-3       }{2}                                     } \| f \|_{L^{p_1}_{\rm h}  L^{q_1}_{\rm v}}   \\
 & \  +C   t ^{  -(\f{1}{p_2}-\f{1}{p}) -\f{1}{2} (\f{1}{q_2}-\f{1}{q}) -\f{ |\beta|+\gamma-2       }{2}                                     } \| |\xh|f (x)\|_{L^{p_2}_{\rm h}  L^{q_2}_{\rm v}}. \\
     \end{split}
\end{equation}
}
\end{enumerate}
\end{lemm}

The proof of Lemma \ref{lm3} follows from Lemma \ref{lm1} and we refer to \cite{F} for a detailed proof. Roughly speaking, the lemma states that the operator $\nablah^{\beta} (-\Deltah)^{  \f{\gamma}{2}   } K \ast$ decays like the 3D heat kernel with horizontal derivative of order $|\beta|+\gamma-2$; this is crucial for the asymptotic analysis of the incompressible MHD system with horizontal dissipation and full magnetic diffusion \ref{MHD-2}.

When handling the incompressible MHD system with full dissipation and horizontal magnetic diffusion, we shall encounter the following function
\begin{equation}\nonumber
     \begin{split}
         N(t,x)= & \  \int_0^{\infty}     (4 \pi  (t+s) ) ^{-\f{3}{2}}    e^{-\frac{|x|^2}{4 (t+s)}}    ds
		  \\
    = & \  \int_0^{\infty}  G( t+s,x)    ds.   \\
     \end{split}
\end{equation}
The next lemma gives the decay estimates of $N(t,x)$ that is crucial in the nonlinear estimates. Remark that these estimates have been used essentially in the context of anisotropic Navier-Stokes system, see \cite{FM}. For convenience, we include a proof.
\begin{lemm}\label{lm3-1}
Let $1\leq p \leq \infty, \alpha \in (\mathbb{N}\cup \{ 0\})^3,m\in \mathbb{N}\cup \{ 0\} $ obeying
\[
\f{3}{2}\left(  1-\f{1}{p} \right)  + \f{|\alpha|}{2} -\f{m }{2} >1.
\]
Then there exists a generic $C>0$ such that for any $t>0$
\begin{equation}\label{li-13-1}
\|  |x|^m \Grad^{\alpha} N(t,x)    \|_{L^p} \leq C t^{  -   \f{3}{2}(  1-\f{1}{p} ) - \f{|\alpha|}{2} +  \f{m }{2} +1         }.
\end{equation}
In particular, in case of $m=1$ we have
\begin{equation}\label{li-13-2}
        \|  |x| \Grad^{\alpha} N(t) \ast f(x)    \|_{L^p}  \leq
        C t ^{  -   \f{3}{2}(  1-\f{1}{p} ) - \f{|\alpha|}{2} +  \f{3 }{2}         }
        \| f \|_{L^1}
       +     C t ^{  -   \f{3}{2}(  1-\f{1}{p} ) - \f{|\alpha|}{2} + 1        }
       \| |x|f(x) \|_{L^1}
\end{equation}
for any $t>0$ and any $f$ satisfying $f\in L^1(\R^3),|x|f(x) \in L^1(\R^3)$.
\end{lemm}
{\it Proof.} By Lemma \ref{lm1},
\begin{equation}\nonumber
     \begin{split}
         \|  |x|^m \Grad^{\alpha} N(t,x)    \|_{L^p}  \leq  & \
         \int_0^{\infty}  \|  |x|^m \Grad^{\alpha} G(t+s,x)    \|_{L^p}  ds
         \\
      \leq & \   \int_0^{\infty} (t+s)   ^{  -   \f{3}{2}(  1-\f{1}{p} ) - \f{|\alpha|}{2} +  \f{m }{2}         }  ds          \\
        =  & \  C t      ^{  -   \f{3}{2}(  1-\f{1}{p} ) - \f{|\alpha|}{2} +  \f{m }{2}        +1 } ,     \\
     \end{split}
\end{equation}
which proves (\ref{li-13-1}). Next, it is easily seen that
\begin{equation}\nonumber
     \begin{split}
        |x| |\Grad^{\alpha} N(t) \ast f (x)|  \leq & \  \int_{\R^3} |x-y|  |\Grad^{\alpha} N(t,x-y) | |f(y)|  dy                     \\
   & \  +    \int_{\R^3}  |\Grad^{\alpha} N(t,x-y) |   |y|  |f(y)|  dy .       \\
     \end{split}
\end{equation}
A direct application of Young inequality gives (\ref{li-13-2}).  $\Box$

The next lemma aims to reveal the enhanced decay rate for a class of weighted and divergence-free data, compared with the standard heat kernel. The result is classical and appeared sparsely in \cites{FM,MY1,MY2}. Due to its importance in our paper, we would like to provide a proof for convenience of the reader.
\begin{lemm}\label{lm3-2}
Let $1\leq p \leq \infty$. Assume that $v_0=(v_{0,1},v_{0,2},v_{0,3})$ satisfies $\nabla \cdot v_0=0$ and
\[
\int_{\R^3} (1+|x|) |v_0 (x)| dx <\infty.
\]
Then there exists a generic $C>0$ such that
\[
\|  e^{t \Delta} v_0 (x)\|_{L^p_x}
\leq C t^{  -\f{3}{2} (1-\f{1}{p})-\f{1}{2} }
\| |x|v_0(x) \|_{L^1_x}  .
\]

\end{lemm}
{\it Proof.} We first claim that
\[
\int_{\R^3} v_0 (x) dx=0.
\]
Indeed, for any $j=1,2,3$, it holds that
\[
\int_{\R^3} v_{0,j}(x) dx =
\int_{\R^3} v_0  \cdot \Grad x_j dx=
-\int_{\R^3} x_j \Grad \cdot v_0 dx=0.
\]
Using this claim, we see
\begin{equation}\nonumber
     \begin{split}
       e^{t \Delta} v_0 (x) = & \  \int_{\R^3} G(t,x-y) v_0 (y) dy                     \\
  = & \ \int_{\R^3} \left( G(t,x-y) -G(t,x) \right) v_0 (y) dy      \\
 = &\  \int_{\R^3} \int_0^1 -y (\Grad G) (t,x-\theta y) v_0 (y) d\theta dy;                \\
     \end{split}
\end{equation}
whence
\begin{equation}\nonumber
     \begin{split}
     \|  e^{t \Delta} v_0 (x)\|_{L^p_x}  \leq & \  C \int_{\R^3}  \int_0^1  \| (\Grad G) (t,x-\theta y)  \|_{L^p_x}  |y| |v_0(y)| d \theta dy       \\
  \leq  & \ C   t^{  -\f{3}{2} (1-\f{1}{p})-\f{1}{2} }
\| |y|v_0(y) \|_{L^1_y} .              \\
     \end{split}
\end{equation}
This completes the proof.    $\Box$
\begin{rem}\label{rk-2}
Let the assumptions of Lemma \ref{lm3-2} be in force. For any $\alpha \in (\mathbb{N} \cup \{ 0\})^3 $, the same argument gives
\[
\|  \Grad^{\alpha} e^{t \Delta} v_0 (x)\|_{L^p_x}
\leq C t^{  -\f{3}{2} (1-\f{1}{p})-\f{|\alpha|}{2}-\f{1}{2} }
\| |x|v_0(x) \|_{L^1_x}  .
\]
\end{rem}


We recall the well-known Gagliardo-Nirenberg interpolation inequality that will be used in the sequel.
\begin{lemm}\label{lm5}
Let $1\leq q \leq \infty$ and $j,k \in \mathbb{N},j<k$, and
\begin{equation*}
either \,\,\,\,
\left\{\begin{aligned}
& r=1  \\
& \f{j}{k} \leq \theta \leq 1           \\
\end{aligned}\right.
\,\,\,\, or \,\,\,\,
\left\{\begin{aligned}
& 1<r<\infty  \\
& \f{j}{k} \leq \theta < 1   .        \\
\end{aligned}\right.
\end{equation*}
Set
\[
\f{1}{p}  =\f{j}{n}  +\theta \left(\f{1}{r}  -\f{k}{n}          \right) +\f{1-\theta}{q}.
\]
There exists an absolute positive constant $C$ such that
\begin{equation}\label{li-13-5}
\| \Grad^j f \|_{L^p(\R^n)} \leq C  \| \Grad^k f \|_{L^r (\R^n)}^{\theta} \| f \|_{L^q(\R^n)}^{1-\theta}
\end{equation}
for any $f\in L^{q}(\R^n) \cap \dot{W}^{k,r}(\R^n)  $.

\end{lemm}

\subsection{Decay estimates for the linearized equations}\label{se-l}
As a direct consequence of the above lemmas, we present the decay estimates for the linearized equations. For the incompressible MHD system (\ref{MHD-2}) with horizontal dissipation and full magnetic diffusion, we obtain by dropping the nonlinear terms that
\begin{align}\label{li-14}
    \begin{cases}
    \partial_t u -  \Deltah u + \nabla p = 0, & t>0 , x \in \mathbb{R}^3,\\
    \partial_t B -  \Delta B  = 0, & t>0 , x \in \mathbb{R}^3,\\
    \nabla \cdot u = \nabla \cdot B = 0 , & t \geqslant 0, x \in \mathbb{R}^3,\\
    u(0,x) = u_0(x), B(0,x) = B_0(x),
    & x \in \mathbb{R}^3.
    \end{cases}
\end{align}
It follows easily that
\[
u(t)=e^{t\Deltah} u_0,\,\,\,\, B(t)=e^{t \Delta} B_0;
\]
whence for any $1\leq p \leq \infty$ and $\alpha=(\alphah,\alpha_3)\in (\mathbb{N} \cup \{ 0\})^2 \times (\mathbb{N} \cup \{ 0\})$ with $|\alpha|\leq 1$, with the help of Lemmas \ref{lm1}-\ref{lm2},
\begin{equation}\label{li-15}
     \begin{split}
         \| e^{t \Deltah} u_{0,\rm h} \|_{L^p}   & \  =O  ( t^{-(1-\f{1}{p})}  ),\,\,   \| \Grad^{\alpha} e^{t \Deltah} u_{0,\rm h} \|_{L^p} =O  ( t^{-(1-\f{1}{p})-\f{|\alphah|}{2}  }  ), \\
   \| e^{t \Deltah} u_{0,3} \|_{L^p} & \   =O ( t^{-\f{3}{2}(1-\f{1}{p})} ), \,\,   \| \nablah^{\alphah} e^{t \Deltah} u_{0,3} \|_{L^p} =O  ( t^{-\f{3}{2}(1-\f{1}{p})-\f{|\alphah|}{2}  }  ),  \\
       \| e^{t\Delta} B_0 \|_{L^p}  & \ =O (    t^{-\f{3}{2}(1-\f{1}{p}) -\f{1}{2}}      ) ,\,\,
         \| \Grad^{\alpha}e^{t\Delta} B_0 \|_{L^p}  =O (    t^{-\f{3}{2}(1-\f{1}{p}) -\f{1}{2} -\f{|\alpha|}{2}}      ).          \\
     \end{split}
\end{equation}
It should be remarked that the decay estimate of $\| e^{t\Delta} B_0 \|_{L^p}$ is $t^{-\f{1}{2}}$-faster than the 3D heat kernel. This is the consequence of our assumptions of initial magnetic field, see Lemma \ref{lm3-2} and Remark \ref{rk-2}. In Section \ref{hd}, we shall prove that, through a careful analysis, the nonlinear solution obeys the same decay estimates and furthermore derive the asymptotic profile of solution.

In a similar manner, for the incompressible MHD system (\ref{MHD-3}) with full dissipation and horizontal magnetic diffusion, we obtain again by dropping the nonlinear terms that
\begin{align}\label{li-16}
    \begin{cases}
    \partial_t u -  \Delta u + \nabla p = 0, & t>0 , x \in \mathbb{R}^3,\\
    \partial_t B -  \Deltah B  = 0, & t>0 , x \in \mathbb{R}^3,\\
    \nabla \cdot u = \nabla \cdot B = 0 , & t \geqslant 0, x \in \mathbb{R}^3,\\
    u(0,x) = u_0(x), B(0,x) = B_0(x),
    & x \in \mathbb{R}^3.
    \end{cases}
\end{align}
It follows that
\[
u(t)=e^{t\Delta} u_0,\,\,\,\, B(t)=e^{t \Deltah} B_0;
\]
whence for any $1\leq p \leq \infty$, with the help of Lemmas \ref{lm1}-\ref{lm2},
\begin{equation}\label{li-17}
     \begin{split}
         \| e^{t \Deltah} B_{0,\rm h} \|_{L^p}   & \  =O  ( t^{-(1-\f{1}{p})}  ),\,\,
         \| \Grad^{\alpha} e^{t \Deltah} B_{0,\rm h} \|_{L^p}  =O  ( t^{-(1-\f{1}{p})-\f{|\alphah|}{2}}  ), \\
   \| e^{t \Deltah} B_{0,3} \|_{L^p} & \   =O ( t^{-\f{3}{2}(1-\f{1}{p})} ), \,\,
          \| \nablah^{\alphah}e^{t \Deltah} B_{0,3} \|_{L^p}  =O ( t^{-\f{3}{2}(1-\f{1}{p})-\f{|\alphah|}{2}} ),            \\
       \| e^{t\Delta} u_0 \|_{L^p}  & \ =O (    t^{-\f{3}{2}(1-\f{1}{p}) -\f{1}{2}}      )
       ,\,\,    \| \Grad^{\alpha} e^{t\Delta} u_0 \|_{L^p}    =O (    t^{-\f{3}{2}(1-\f{1}{p}) -\f{1}{2}-\f{|\alpha|}{2}}      )           .          \\
     \end{split}
\end{equation}
As the previous case, the decay estimate of $\| e^{t\Delta} u_0 \|_{L^p}$ is $t^{-\f{1}{2}}$-faster than the 3D heat kernel due to the assumptions of initial velocity field, see Lemma \ref{lm3-2} and Remark \ref{rk-2} for the details. Furthermore, we deduce from our assumption of initial velocity that\footnote{Notice that $\| \nabla^{\alpha} e^{t\Delta} u_0 \|_{L^p}= \|\nabla^{\alpha} G(t) \ast u_0  \|_{L^p} \leq \| G(t) \|_{L^1} \| \nabla^{\alpha}  u_0 \|_{L^p} \leq \| \nabla^{\alpha}  u_0 \|_{L^p}$. Our assumptions of initial velocity in Theorem \ref{thm2} ensure that $\nabla^{\alpha}  u_0 \in L^1(\R^3) \cap L^{\infty}(\R^3)$ for any $\alpha \in (\mathbb{N} \cup \{ 0\})^3, |\alpha| \leq 1$.
}
\begin{equation}\label{li-17-1}
\| \nabla^{\alpha} e^{t\Delta} u_0 \|_{L^p} \leq \infty, \,\,\text{ for any  }  \alpha \in (\mathbb{N} \cup \{ 0\})^3, |\alpha| \leq 1.
\end{equation}
Therefore, we actually have
\[
\| \Grad^{\alpha} e^{t\Delta} u_0 \|_{L^p}  \leq C_{\rhd} (1+t)^{-\f{3}{2}(1-\f{1}{p}) -\f{1}{2}-\f{|\alpha|}{2}}         .
\]
Here $C_{\rhd}$ consists of the norms in (\ref{li-17-1}). Our main goal in Section \ref{hm} is to confirm that the nonlinear solution satisfies the same decay estimate for $\Bh(t),B_3(t)$, while $u(t)$ admits a slower decay estimate than the linear estimate. Moreover, we shall derive the asymptotic profile of solution.

\section{Large time behavior of the MHD system with horizontal dissipation and full magnetic diffusion}\label{hd}
\subsection{Reformulation of the equations}
In this section, we consider the incompressible MHD system with horizontal dissipation and full magnetic diffusion:
\begin{align}\label{MHD-2-1}
    \begin{cases}
    \partial_t u -  \Deltah u + ( u \cdot \nabla ) u - ( B \cdot \nabla ) B + \nabla p = 0, & t>0 , x \in \mathbb{R}^3,\\
    \partial_t B -   \Delta B + ( u \cdot \nabla ) B - ( B \cdot \nabla ) u = 0, & t>0 , x \in \mathbb{R}^3,\\
    \nabla \cdot u = \nabla \cdot B = 0 , & t \geqslant 0, x \in \mathbb{R}^3,\\
    u(0,x) = u_0(x), B(0,x) = B_0(x),
    & x \in \mathbb{R}^3,
    \end{cases}
\end{align}
The corresponding integral equation of \eqref{MHD-2-1} is given by
\begin{align}
    &u(t) = e^{ t \Deltah} u_0
    - \int_0^t e^{ (t - \tau) \Deltah} \mathbb{P} \nabla \cdot ( u \otimes  u - B \otimes B ) (\tau) d\tau,\label{Int-2-1-u}\\
    &B(t) = e^{ t \Delta} B_0
    - \int_0^t e^{ (t - \tau) \Delta} \nabla \cdot ( B \otimes  u - u \otimes B ) (\tau) d\tau,\label{Int-2-1-B}
\end{align}
where $\mathbb{P}$ is the Helmholtz projection on $\mathbb{R}^3$.
Then, by the similar argument as in Section 2 of \cite{F},
we have the following decomposition of the nonlinear terms of \eqref{Int-2-1-u}-\eqref{Int-2-1-B}:
\begin{prop}\label{prop1}
    Let $(u,B)$ be the solution to the integral equation \eqref{Int-2-1-u}-\eqref{Int-2-1-B}.
    Then, the integral equation for $u$ is decomposed as
\beq\label{hd4}
        \uh(t) = e^{ t \Deltah} u_{0,{\rm h}} + \sum_{m=1}^5 \Big( \Nhu_m[u](t) -  \Nhu_m[B](t) \Big)
\eeq
\beq\label{hd5}
        u_3(t) = e^{ t \Deltah} u_{0,3} + \sum_{m=1}^3 \Big( \Nvu_m[u](t) -  \Nvu_m[B](t) \Big),
\eeq
where
	\begin{align*}
		\Nhu_1[v](t)&:=-\int_0^te^{ (t-\tau)\Deltah}\partial_{3}(v_3\vh)(\tau)d\tau,\\
		\Nhu_2[v](t)&:=-\int_0^te^{ (t-\tau)\Deltah}\nablah \cdot (\vh\otimes\vh)(\tau)d\tau,\\
		\Nhu_3[v](t)&:=\int_0^t\nablah e^{ (t-\tau)\Deltah}(v_3(\tau)^2)d\tau,\\
		\Nhu_4[v](t)&:=-\sum_{k,l=1}^2\int_0^t\nablah\partial_{k}\partial_{l}K(t-\tau)*(v_k v_l)(\tau)d\tau,\\
		\Nhu_5[v](t)&:=2\sum_{k=1}^2\int_0^t\nablah \partial_k(-\Deltah)^{\frac{1}{2}}\widetilde{K}(t-\tau)*(v_3 v_k)(\tau)d\tau\\
		&\qquad +\int_0^t\nablah \Deltah K(t-\tau)*(v_3(\tau)^2)d\tau
	\end{align*}
	and
	\begin{align*}
		\Nvu_1[v](t)&:=\int_0^te^{(t-\tau)\Deltah}\nablah\cdot(v_3\vh)(\tau)d\tau,\\
		\Nvu_2[v](t)&:=\sum_{k,l=1}^2\int_0^t(-\Deltah)^{\frac{1}{2}}\partial_{k}\partial_{l}\widetilde{K}(t-\tau)*(v_k v_l)(\tau)d\tau,\\
		\Nvu_3[v](t)&:=2\sum_{k=1}^2\int_0^t\partial_k\Deltah K(t-\tau)*(v_3v_k)(\tau)d\tau\\
		&\qquad+\int_0^t(-\Deltah)^{\frac{3}{2}}\widetilde{K}(t-\tau)*(v_3(\tau)^2)d\tau
	\end{align*}
	for $v = u, B$.
	Here, $K(t,x)$ and $\widetilde{K}(t,x)$ are the functions with the following representations:
	\begin{align*}
		&K(t,x)=\int_0^{\infty}\frac{e^{-\frac{|\xh|^2}{4 (t+s)}}}{4 \pi  (t+s)}\frac{e^{-\frac{x_3^2}{4 s}}}{(4 \pi  s)^{\frac{1}{2}}}ds
		=\int_0^{\infty}\Gh( t+s,\xh)\Gv( s,x_3)ds,\\
		&\widetilde{K}(t,x)={\rm sgn}(x_3)K(t,x).
	\end{align*}
    For the integral equation of $B$, it holds
\beq\label{hd6}
        \Bh(t) = e^{ t \Delta} B_{0,{\rm h}} + \sum_{m=1}^3 \Nhb_m[u,B](t),
\eeq
\beq\label{hd7}
        B_3(t) = e^{ t \Delta} B_{0,3} + \sum_{m=1}^3 \Nvb_m[u,B](t),
\eeq
where
    \begin{equation}\nonumber
     \begin{split}
         \Nhb_1[u,B](t) & \  :=-\int_0^t e^{(t-\tau)\Delta} \nabla \cdot ( \Bh \otimes u )(\tau)d\tau, \\
     \Nhb_2[u,B](t) & \  :=\int_0^t e^{(t-\tau)\Delta} \nablah \cdot( \uh \otimes \Bh )(\tau)d\tau, \\
     \Nhb_3[u,B](t)    & \ :=\int_0^t e^{(t-\tau)\Delta} \partial_3 ( B_3 \uh )(\tau)d\tau,\\
     \end{split}
\end{equation}
and
    \begin{align*}
        \Nvb_1[u,B](t)&:=-\int_0^t e^{(t-\tau)\Delta} \nabla \cdot ( B_3  u )(\tau)d\tau,\\
        \Nvb_2[u,B](t)&:=\int_0^t e^{(t-\tau)\Delta} \nablah \cdot (u_3  \Bh )(\tau)d\tau,\\
        \Nvb_3[u,B](t)&:=\int_0^t e^{(t-\tau)\Delta} \partial_3 ( u_3 B_3)(\tau)d\tau.
    \end{align*}
\end{prop}

\subsection{Decay estimates for the Duhamel terms}
In this subsection, we show the decay estimates for the Duhamel terms in (\ref{hd4})-(\ref{hd7}). To this end, we assume that $s\in \mathbb{N},0<T\leq \infty, A, A_{\ast} \geq 0$ and make the following ansatzes:
\begin{enumerate}[(i)]
\item {
$u,B\in C([0,\infty);X^s(\R^3)),\Grad \cdot u=\Grad \cdot B=0$ and
\beq\label{hd8}
\| (u,B)(t) \|_{H^s} \leq A,\,\,     \| (u,B)(t) \| _{  L^1_{\xh}(W^{1,1}\cap W^{1,\infty})_{x_3}         }
\leq A(1+t)
\eeq
for any $t>0$.
}
\item {
For any $1\leq p \leq \infty$,
\begin{equation}\label{hd9}
     \begin{split}
      \| \nabla ^{\alpha}   \uh (t) \|_{L^p} & \ \leq C A t^{ -(1-\f{1}{p}) -\f{|\alphah|}{2}        }  \\
    \| \nablah ^{\alphah}   u_3 (t) \|_{L^p} & \   \leq C A t^{ -\f{3}{2}(1-\f{1}{p}) -\f{|\alphah|}{2}        }  , \\
      \| \nabla ^{\alpha}   B (t) \|_{L^p}    & \  \leq C A t^{ -\f{3}{2}(1-\f{1}{p}) -\f{1}{2}-\f{|\alpha|}{2}        }, \\
     \end{split}
\end{equation}
for any $0<t<T$ and any $\alpha=(\alphah,\alpha_3)\in (\mathbb{N} \cup \{ 0\})^2 \times (\mathbb{N} \cup \{ 0\})$ with $|\alpha|\leq 1$.
}
\item{
For any $0<t <T$ and $\alpha=(\alphah,\alpha_3)\in (\mathbb{N} \cup \{ 0\})^2 \times (\mathbb{N} \cup \{ 0\})$ with $|\alpha|\leq 1$,
\begin{equation}\label{hd9-1}
     \begin{split}
        \| \Grad^{\alpha} u (t) \|_{L^{\infty}_{\rm h}  L^{1}_{\rm v}} & \  \leq A t^{-1-\f{|\alphah|}{2}}, \\
     \| \Grad^{\alpha} B (t) \|_{L^{\infty}_{\rm h}  L^{1}_{\rm v}} & \  \leq A t^{-1-\f{1+|\alpha|}{2}}.  \\
     \end{split}
\end{equation}
}
\item{
For any $t>0$,
\begin{equation}\label{hd9-2}
     \begin{split}
        \| |\xh|\uh (t,x) \|_{L^{1}_{\rm h}  L^{\infty}_{\rm v}} & \  \leq A_{\ast} (1+t)^{  \f{1}{2}  },  \\
     \| |\xh| u_3 (t,x) \|_{L^{1}_{\rm h}  L^{\infty}_{\rm v}} & \  \leq A_{\ast},  \\
         \| |x| \Bh (t,x) \|_{L^{1}_{\rm h}  L^{\infty}_{\rm v}} & \  \leq A_{\ast} (1+t)^{  \f{1}{2}  },  \\
           \| |x| B_3 (t,x) \|_{L^{1}_{\rm h}  L^{\infty}_{\rm v}} & \  \leq A_{\ast} (1+t)^{  \f{1}{2}  }.  \\
     \end{split}
\end{equation}
}
\end{enumerate}

Based on the ansatzes above, one gets the following
\begin{lemm}\label{lem6}
Let $(u,B)$ be subject to the ansatzes $(i),(ii)$ with $s\geq 3$. There exists a generic constant $C>0$ such that for any $0<t<T$ and any $1\leq p \leq \infty$ we have the decay estimates for velocity:
\begin{equation}\label{hd10}
     \begin{split}
         \| \uh(t) \|_{L^p},\| \p_3 \uh (t)\|_{L^p} & \  \leq  CA (1+t)^{-(1-\f{1}{p})}   ,          \\
     \|  \nablah \uh (t)\|_{L^p}  & \  \leq    \left\{\begin{aligned}
& CA (1+t)^{  -(1-\f{1}{p})-\f{1}{2}           }\,\,\, \text{if  } 2\leq p\leq \infty         \\
& CA t^{  -(1-\f{1}{p})-\f{1}{2}           }\,\,\, \text{if  }   1\leq p <2,         \\
\end{aligned}\right.                           \\
       \| u_3 ( t) \|_{L^p}  & \  \leq      CA (1+t)^{  -\f{3}{2}   (1-\f{1}{p})         }      ,           \\
        \| \p_3 u_3 ( t) \|_{L^p}  & \  \leq      CA (1+t)^{  -   (1-\f{1}{p})  -\f{1}{2}       }      ,           \\
          \|  \nablah u_3 (t)\|_{L^p}  & \   \leq    \left\{\begin{aligned}
& CA (1+t)^{  -\f{3}{2} (1-\f{1}{p})-\f{1}{2}           }\,\,\, \text{if  } 2\leq p\leq \infty         \\
& CA t^{  -  \f{3}{2}  (1-\f{1}{p})-\f{1}{2}           }\,\,\, \text{if  }   1\leq p <2,         \\
\end{aligned}\right.                           \\
     \| \uh (t) \|_{L^1_{\rm h}  L^{\infty}_{\rm v}}    & \  \leq  CA  ,  \\
     \| u_3 (t) \|_{L^1_{\rm h}  L^{\infty}_{\rm v}}    & \  \leq  CA (1+t)^{-\f{1}{2} },  \\
     \end{split}
\end{equation}
and the decay estimates for magnetic field:
\begin{equation}\label{hd11}
     \begin{split}
         \| B(t) \|_{L^p} & \   \leq CA (1+t )^{  -\f{3}{2} (1-\f{1}{p})-\f{1}{2}           },  \\
        \|  \nablah B (t)\|_{L^p}  & \   \leq    \left\{\begin{aligned}
& CA (1+t)^{  -\f{3}{2} (1-\f{1}{p})-1          }\,\,\, \text{if  } 2\leq p\leq \infty         \\
& CA t^{  -  \f{3}{2}  (1-\f{1}{p})-1           }\,\,\, \text{if  }   1\leq p <2,         \\
\end{aligned}\right.                           \\
        \| \p_3 B(t) \|_{L^p} & \   \leq CA (1+t )^{  -\f{3}{2} (1-\f{1}{p})-1           },  \\
       \| B(t) \|_{L^1_{\rm h}  L^{\infty}_{\rm v}}   & \ \leq CA (1+t)^{-1}.
     \end{split}
\end{equation}
\end{lemm}

Observe that the decay estimates in Lemma \ref{lem6} follow directly from the ansatzes (\ref{hd8})-(\ref{hd9}). The details are thus omitted. We are in a position to give the decay estimates for the Duhamel terms.
\begin{lemm}\label{lem7}
Let $(u,B)$ be subject to the ansatzes $(i),(ii)$ with $s\geq 5$. There exists a generic constant $C>0$ such that for any $0<t<T$ and any $1\leq p \leq \infty$ we have the decay estimates for the Duhamel terms associated with velocity equation:
\begin{equation}\label{hd12}
     \begin{split}
         \| \Grad^{\alpha} \Nhu_1[u](t)   \|_{L^p}     & \ \leq CA^2 (1+t )^{ -(1-\f{1}{p}) -\f{|\alphah|}{2}   } ,  \\
         \| \Grad^{\alpha} \Nhu_2[u](t)   \|_{L^p}     & \ \leq CA^2 (1+t )^{ -(1-\f{1}{p}) -\f{1+|\alphah|}{2}   } \log(2+t) ,        \\
         \| \Grad^{\alpha} \Nhu_3[u](t)   \|_{L^p}     & \ \leq CA^2 (1+t )^{ -(1-\f{1}{p}) -\f{1+|\alphah|}{2}   } ,        \\
         \| \Grad^{\alpha} \Nhu_4[u](t)   \|_{L^p}     & \ \leq CA^2 (1+t )^{ -\f{9}{8}(1-\f{1}{p}) -\f{1+|\alphah|}{2}   }\log(2+t) ,        \\
         \| \Grad^{\alpha} \Nhu_5[u](t)   \|_{L^p}     & \ \leq CA^2 (1+t )^{ -\f{9}{8}(1-\f{1}{p}) -\f{1+|\alphah|}{2}   } ,        \\
         \| \Grad^{\alpha} \Nhu_1[B](t)   \|_{L^p}     & \ \leq CA^2 (1+t ) ^{  -(1-\f{1}{p})    -\f{ |\alphah|   }{2}             } ,        \\
            \| \Grad^{\alpha} \Nhu_2[B](t)   \|_{L^p}     & \ \leq CA^2 (1+t )^{   -(1-\f{1}{p})    -\f{ 1+|\alphah|   }{2}       } ,        \\
               \| \Grad^{\alpha} \Nhu_3[B](t)   \|_{L^p}     & \ \leq CA^2 (1+t ) ^{   -(1-\f{1}{p})    -\f{ 1+|\alphah|   }{2}       },        \\
                  \| \Grad^{\alpha} \Nhu_4[B](t)   \|_{L^p}     & \ \leq CA^2 (1+t ) ^{   -\f{3}{2}(1-\f{1}{p})    -\f{ 1+|\alphah|   }{2}       },        \\
                     \| \Grad^{\alpha} \Nhu_5[B](t)   \|_{L^p}     & \ \leq CA^2 (1+t ) ^{   -\f{3}{2}(1-\f{1}{p})    -\f{ 1+|\alphah|   }{2}       },        \\
     \end{split}
\end{equation}

\begin{equation}\label{hd13}
     \begin{split}
        \| \Grad^{\alpha} \Nvu_1[u](t)   \|_{L^p} & \  \leq CA^2 (1+t )^{ -(1-\f{1}{p}) -\f{1+|\alphah|}{2}   }  , \\
     \| \Grad^{\alpha} \Nvu_2[u](t)   \|_{L^p} & \  \leq CA^2 (1+t )^{ -\f{9}{8}(1-\f{1}{p}) -\f{1+|\alphah|}{2}   } \log(2+t)  , \\
        \| \Grad^{\alpha} \Nvu_3[u](t)   \|_{L^p} & \  \leq CA^2 (1+t )^{ -\f{9}{8}(1-\f{1}{p}) -\f{1+|\alphah|}{2}   }   , \\
         \| \Grad^{\alpha} \Nvu_1[B](t)   \|_{L^p} & \  \leq CA^2 (1+t ) ^{   -(1-\f{1}{p})    -\f{ 1+|\alphah|   }{2}       }  , \\
         \| \Grad^{\alpha} \Nvu_2[B](t)   \|_{L^p} & \  \leq CA^2 (1+t )^{   -\f{3}{2}(1-\f{1}{p})    -\f{ 1+|\alphah|   }{2}       }   , \\
         \| \Grad^{\alpha} \Nvu_3[B](t)   \|_{L^p} & \  \leq CA^2 (1+t )^{   -\f{3}{2}(1-\f{1}{p})    -\f{ 1+|\alphah|   }{2}       }   , \\
     \end{split}
\end{equation}
for any $\alpha=(\alphah,\alpha_3)\in (\mathbb{N} \cup \{ 0\})^2 \times (\mathbb{N} \cup \{ 0\})$ with $|\alpha|\leq 1$ and the decay estimates for the Duhamel terms associated with magnetic equation:
\begin{equation}\label{hd14}
     \begin{split}
         \| \Grad^{\alpha} \Nhb_1[u,B](t)   \|_{L^p} & \           \leq CA^2 (1+t ) ^{  -\f{3}{2}(1-\f{1}{p})    -\f{ 1+|\alpha|   }{2}             }  , \\
    \| \Grad^{\alpha} \Nhb_2[u,B](t)   \|_{L^p} & \           \leq CA^2 (1+t ) ^{  -\f{3}{2}(1-\f{1}{p})    -\f{ 1+|\alpha|   }{2}             }  , \\
          \| \Grad^{\alpha} \Nhb_3[u,B](t)   \|_{L^p} & \           \leq CA^2 (1+t ) ^{  -\f{3}{2}(1-\f{1}{p})    -\f{ 1+|\alpha|   }{2}             }  , \\
     \end{split}
\end{equation}
\begin{equation}\label{hd15}
     \begin{split}
         \| \Grad^{\alpha} \Nvb_1[u,B](t)   \|_{L^p} & \           \leq CA^2 (1+t )  ^{  -\f{3}{2}(1-\f{1}{p})    -\f{ 1+|\alpha|   }{2}             } , \\
    \| \Grad^{\alpha} \Nvb_2[u,B](t)   \|_{L^p} & \           \leq CA^2 (1+t ) ^{  -\f{3}{2}(1-\f{1}{p})    -\f{ 1+|\alpha|   }{2}             }  , \\
          \| \Grad^{\alpha} \Nvb_3[u,B](t)   \|_{L^p} & \           \leq CA^2 (1+t )^{  -\f{3}{2}(1-\f{1}{p})    -\f{ 1+|\alpha|   }{2}             }   , \\
     \end{split}
\end{equation}
for any $\alpha=(\alphah,\alpha_3)\in (\mathbb{N} \cup \{ 0\})^2 \times (\mathbb{N} \cup \{ 0\})$ with $|\alpha|\leq 1$.
\end{lemm}
{\it Proof.} Notice that the first five estimates in (\ref{hd12}) and the first three estimates in (\ref{hd13}) have been proved in \cite{F}. In the first part of the proof, we show the decay estimates in (\ref{hd12})-(\ref{hd13}) involved with magnetic field. Upon recalling the definition of these Duhamel terms, it suffices to show
\begin{equation}\label{hd16}
\left\|   \Grad^{\alpha} \int_0^t e^{ (t-\tau)\Deltah  }\p_3 (B_3\Bh ) (\tau  )     d \tau \right\|_{L^p}
\leq CA^2 (1+t)^{  -(1-\f{1}{p})    -\f{ |\alphah|   }{2}             },
\end{equation}
\begin{equation}\label{hd17}
\left\|   \Grad^{\alpha} \int_0^t e^{ (t-\tau)\Deltah  }\nablah (B_k B_l ) (\tau  )     d \tau \right\|_{L^p}
\leq CA^2 (1+t)^{   -(1-\f{1}{p})    -\f{ 1+|\alphah|   }{2}       },
\end{equation}
\begin{equation}\label{hd18}
\left\|   \Grad^{\alpha} \int_0^t   K_{\beta,\gamma}^{ (m)} (t-\tau) \ast (B_k B_l ) (\tau  )     d \tau \right\|_{L^p}
\leq CA^2 (1+t)^{   -\f{3}{2}(1-\f{1}{p})    -\f{ 1+|\alphah|   }{2}       },
\end{equation}
for any $k,l\in \{1,2,3\}$ and $m\in \{1,2\}$. $(\beta,\gamma)  \in (\mathbb{N} \cup \{ 0\})^2 \times (\mathbb{N} \cup \{ 0\})$ obeying $|\beta|+\gamma=3$. Here and in the sequel we adopt the following notations:
\[
K_{\beta,\gamma}^{ (1)} (t,x):= \nablah^{\beta}  (-\Deltah )^{ \f{\gamma}{2}  } K(t,x),\,\,\,\,
K_{\beta,\gamma}^{ (2)} (t,x):=\text{sgn}(x_3) \nablah^{\beta}  (-\Deltah )^{ \f{\gamma}{2}  } K(t,x).
\]
We start with the proof of (\ref{hd16}) for the case $p=1$. Notice that
\begin{equation}\nonumber
     \begin{split}
    \| \p_3 (B_3 \Bh) (\tau)\|_{L^1}        \leq  & \    \|\p_3 B_3 (\tau)  \|_{L^{\infty}} \|  \Bh (\tau) \|  _{L^1}   +
    \| B_3 (\tau) \|_{L^{\infty}} \| \p_3 \Bh   (\tau) \|  _{L^1}     \\
    \leq  & \    CA^2 (1+\tau)^{-3} ,        \\
     \end{split}
\end{equation}
\begin{equation}\nonumber
     \begin{split}
    \| \p_3^2 (B_3 \Bh) (\tau) \|_{L^1}        \leq  & \    \|\p_3^2 B_3 (\tau)  \|_{L^{2}} \|  \Bh (\tau) \|  _{L^2}   +
   2 \| \p_3 B_3 (\tau) \|_{L^{2}} \| \p_3 \Bh (\tau) \|  _{L^2}     \\
   &\ +  \| B_3 (\tau)  \|_{L^{2}} \|   \p_3^2 \Bh (\tau) \|  _{L^2}  \\
 \leq   &\   CA^2 (1+\tau) ^{  -\f{5}{4}          };
     \end{split}
\end{equation}
whence
\begin{equation}\nonumber
     \begin{split}
     \left\|   \Grad^{\alpha} \int_0^t e^{ (t-\tau)\Deltah  }\p_3 (B_3\Bh ) (\tau  )     d \tau \right\|_{L^1}   \leq  & \  C  \int_0^t  \| \nablah^{\alphah} \Gh (t-\tau) \|_{L^1(\R^2)} \| \p_3^{\alpha_3+1} (B_3 \Bh) (\tau) \|     _{L^1}         d \tau                     \\
  \leq  & \  CA^2          \int_0^  {  \f{t}{2} } (t-\tau)^{  -\f{|\alphah|}{2}         } (1+\tau) ^{  -\f{5}{4} } d \tau                 \\
         & \ +     CA^2          \int_{  \f{t}{2} } ^t (t-\tau)^{  -\f{|\alphah|}{2}         } (1+\tau) ^{  -\f{5}{4}          }    d\tau                 \\
     \leq     & \   CA^2 t^{    -\f{|\alphah|}{2}           }
     \end{split}
\end{equation}
for $t \geq 1$ and
\begin{equation}\nonumber
     \begin{split}
     \left\|   \Grad^{\alpha} \int_0^t e^{ (t-\tau)\Deltah  }\p_3 (B_3\Bh ) (\tau  )     d \tau \right\|_{L^1}   \leq  & \  C  \int_0^t  \| \nablah^{\alphah} \Gh (t-\tau) \|_{L^1(\R^2)} \| \p_3^{\alpha_3+1} (B_3 \Bh) (\tau) \|     _{L^1}         d \tau                     \\
  \leq  & \  CA^2          \int_0^  {  t } (t-\tau)^{  -\f{|\alphah|}{2}         } (1+\tau) ^{  -\f{5}{4} } d \tau                  \\
     \leq     & \   CA^2
     \end{split}
\end{equation}
for $0<t \leq 1$. This gives (\ref{hd16}) with $p=1$. For the case $p=\infty$, we divide it into two classes: either $\alpha_3=0$ or $\alpha_3=1$. If $\alpha_3=0$,
\begin{equation}\nonumber
     \begin{split}
        \| \p_3(B_3 \Bh) (\tau) \|_{L^1_{\rm h}  L^{\infty}_{\rm v}}  \leq  & \   \|  \p_3 B_3 (\tau) \| _{L^{\infty}}
        \| \Bh ( \tau ) \|_{L^1_{\rm h}  L^{\infty}_{\rm v}}   +
        \| B_3 ( \tau ) \|_{L^1_{\rm h}  L^{\infty}_{\rm v}}    \|  \p_3  \Bh  (\tau) \| _{L^{\infty}}              \\
    \leq  & \ CA^2 (1+\tau) ^{ -\f{7}{2}     } ,      \\
     \end{split}
\end{equation}
\begin{equation}\nonumber
     \begin{split}
        \| \p_3(B_3 \Bh) (\tau) \|_{L^{\infty}}  \leq  & \   \|  \p_3 B_3 (\tau) \| _{L^{\infty}}
        \| \Bh ( \tau ) \|_{L^{\infty}}  +
        \| B_3 ( \tau ) \|_{L^{\infty}}   \|  \p_3  \Bh  (\tau) \| _{L^{\infty}}              \\
    \leq  & \ CA^2 (1+\tau) ^{ -\f{9}{2}     } ;      \\
     \end{split}
\end{equation}
whence
\begin{equation}\nonumber
     \begin{split}
     \left\|   \nablah^{\alphah} \int_0^t e^{ (t-\tau)\Deltah  }\p_3 (B_3\Bh ) (\tau  )     d \tau \right\|_{L^{\infty}}   \leq  & \  C   \int_0^  {  \f{t}{2} }  \| \nablah^{\alphah} \Gh (t-\tau) \|_{L^{\infty}(\R^2)} \| \p_3  (B_3 \Bh) (\tau) \|     _{L^1_{\rm h}  L^{\infty}_{\rm v}}         d \tau                     \\
         & \ +     C          \int_{  \f{t}{2} } ^t        \| \nablah^{\alphah} \Gh (t-\tau) \|_{L^{1}(\R^2)} \| \p_3  (B_3 \Bh) (\tau) \|    _{L^{\infty}}         d \tau         \\
     \leq     & \   CA^2 \int_0^  {  \f{t}{2} } (t-\tau)^{ -1-\f{|\alphah|}{2}   }    (1+\tau) ^{ -\f{7}{2}     } d\tau \\
     &\ + CA^2  \int_{  \f{t}{2} } ^t    (t-\tau)^{ -\f{|\alphah|}{2}   }    (1+\tau) ^{ -\f{9}{2}     } d\tau \\
\leq      &\  CA^2 t ^{ -1-\f{|\alphah|}{2}   }
     \end{split}
\end{equation}
for $t \geq 1$ and
\begin{equation}\nonumber
     \begin{split}
     \left\|  \nablah^{\alphah}  \int_0^t e^{ (t-\tau)\Deltah  }\p_3 (B_3\Bh ) (\tau  )     d \tau \right\|_{L^{\infty}}      \leq  & \  C  \int_0^t  \| \nablah^{\alphah} \Gh (t-\tau) \|_{L^1(\R^2)} \| \p_3 (B_3 \Bh) (\tau) \|     _{L^{\infty}}        d \tau                     \\
     \leq     & \   CA^2    \int_0^  {  t }   (t-\tau)^{  -\f{|\alphah|}{2}         } (1+\tau) ^{  -\f{9}{2} } d \tau   \\
     \leq &\ CA^2
     \end{split}
\end{equation}
for $0<t \leq 1$. This gives (\ref{hd16}) with $p=\infty,\alpha_3=0$. If $\alpha_3=1$,
\begin{equation}\nonumber
     \begin{split}
        \| \p_3^2(B_3 \Bh) (\tau)  \|_{L^1_{\rm h}  L^{\infty}_{\rm v}}  \leq     & \  \| \p_3^2 B_3 (\tau) \|_{L^2_{\rm h}  L^{\infty}_{\rm v}}   \| \Bh (\tau)  \|_{L^2_{\rm h}  L^{\infty}_{\rm v}}
        + 2 \| \p_3 B_3 (\tau) \|_{L^2_{\rm h}  L^{\infty}_{\rm v}}    \| \p_3  \Bh (\tau)  \|_{L^2_{\rm h}  L^{\infty}_{\rm v}}             \\
         & \  +      \|  B_3 (\tau) \|_{L^2_{\rm h}  L^{\infty}_{\rm v}}   \|  \p_3^2 \Bh (\tau)  \|_{L^2_{\rm h}  L^{\infty}_{\rm v}}                             \\
      \leq    & \ C  \Big(   \| \p_3^2 B_3 (\tau) \|_{L^2} ^{ \f{1}{2}  } \| \p_3^3 B_3 (\tau) \|_{L^2} ^{ \f{1}{2}  }               \|  \Bh(\tau) \| _{L^2}  ^{ \f{1}{2}  }   \| \p_3  \Bh(\tau) \| _{L^2}  ^{ \f{1}{2}  }       \\
         & \  +   \| \p_3 B_3 (\tau) \|_{L^2} ^{ \f{1}{2}  } \| \p_3^2 B_3 (\tau) \|_{L^2} ^{ \f{1}{2}  }                   \| \p_3 \Bh (\tau) \|_{L^2} ^{ \f{1}{2}  } \| \p_3^2 \Bh (\tau) \|_{L^2} ^{ \f{1}{2}  }       \\
         &\ +        \|  B_3 (\tau) \|_{L^2} ^{ \f{1}{2}  } \| \p_3 B_3 (\tau) \|_{L^2} ^{ \f{1}{2}  }                   \| \p_3^2 \Bh (\tau) \|_{L^2} ^{ \f{1}{2}  } \| \p_3^3 \Bh (\tau) \|_{L^2} ^{ \f{1}{2}  }   \Big)                          \\
       \leq   &\   CA^2  (1+\tau)^{- \f{7}{4}  }  ,      \\
     \end{split}
\end{equation}
\begin{equation}\nonumber
     \begin{split}
        \| \p_3^2(B_3 \Bh) (\tau)  \|_{L^2_{\rm h}  L^{\infty}_{\rm v}}  \leq     & \  \| \p_3^2 B_3 (\tau) \|_{L^2_{\rm h}  L^{\infty}_{\rm v}}   \| \Bh (\tau)  \|_{L^{\infty}}
        + 2 \| \p_3 B_3 (\tau) \|_{L^{\infty}}     \| \p_3  \Bh (\tau)  \|_{L^2_{\rm h}  L^{\infty}_{\rm v}}             \\
         & \  +      \|  B_3 (\tau) \|_{L^{\infty}}    \|  \p_3^2 \Bh (\tau)  \|_{L^2_{\rm h}  L^{\infty}_{\rm v}}                             \\
      \leq    & \ C  \Big(   \| \p_3^2 B_3 (\tau) \|_{L^2} ^{ \f{1}{2}  } \| \p_3^3 B_3 (\tau) \|_{L^2} ^{ \f{1}{2}  }              \| \Bh (\tau)  \|_{L^{\infty}}         \\
         & \  +  \| \p_3 B_3 (\tau) \|_{L^{\infty}}                 \| \p_3 \Bh (\tau) \|_{L^2} ^{ \f{1}{2}  } \| \p_3^2 \Bh (\tau) \|_{L^2} ^{ \f{1}{2}  }       \\
         &\ +        \|  B_3 (\tau) \|_{L^{\infty}}                   \| \p_3^2 \Bh (\tau) \|_{L^2} ^{ \f{1}{2}  } \| \p_3^3 \Bh (\tau) \|_{L^2} ^{ \f{1}{2}  }   \Big)                          \\
       \leq   &\   CA^2  (1+\tau)^{- 2 }  ;      \\
     \end{split}
\end{equation}
whence \footnote{Due to the nice decay estimates of magnetic field, the Gagliardo-Nirenberg inequality is not necessary when estimating $\| \p_3^2(B_3 \Bh) (\tau)  \|_{L^1_{\rm h}  L^{\infty}_{\rm v}}, \| \p_3^2(B_3 \Bh) (\tau)  \|_{L^2_{\rm h}  L^{\infty}_{\rm v}}$. In other words, the regularity assumption for initial magnetic field could be relaxed to $X^s(\R^3)$ with $s \geq 3$. }
\begin{equation}\nonumber
     \begin{split}
     \left\|   \p_3 \int_0^t e^{ (t-\tau)\Deltah  }\p_3 (B_3\Bh ) (\tau  )     d \tau \right\|_{L^{\infty}}   \leq  & \  C   \int_0^  {  \f{t}{2} }  \|  \Gh (t-\tau) \|_{L^{\infty}(\R^2)} \| \p_3^2  (B_3 \Bh) (\tau) \|     _{L^1_{\rm h}  L^{\infty}_{\rm v}}         d \tau                     \\
         & \ +     C          \int_{  \f{t}{2} } ^t        \|  \Gh (t-\tau) \|_{L^{2}(\R^2)} \| \p_3^2  (B_3 \Bh) (\tau) \|   _{L^2_{\rm h}  L^{\infty}_{\rm v}}         d \tau         \\
     \leq     & \   CA^2 \int_0^  {  \f{t}{2} } (t-\tau)^{ -1  }    (1+\tau) ^{ -\f{7}{4}     } d\tau \\
     &\ + CA^2  \int_{  \f{t}{2} } ^t    (t-\tau)^{ -\f{1}{2}   }    (1+\tau) ^{ -2    } d\tau \\
\leq      &\  CA^2 t ^{ -1  }
     \end{split}
\end{equation}
for $t \geq 1$ and
\begin{equation}\nonumber
     \begin{split}
     \left\|  \p_3  \int_0^t e^{ (t-\tau)\Deltah  }\p_3 (B_3\Bh ) (\tau  )     d \tau \right\|_{L^{\infty}}      \leq  & \  C  \int_0^t  \| \Gh (t-\tau) \|_{L^2(\R^2)} \| \p_3^2 (B_3 \Bh) (\tau) \|   _{L^2_{\rm h}  L^{\infty}_{\rm v}}        d \tau                     \\
     \leq     & \   CA^2    \int_0^  {  t }   (t-\tau)^{  -\f{ 1  }{2}         } (1+\tau) ^{  -2 } d \tau   \\
     \leq &\ CA^2
     \end{split}
\end{equation}
for $0<t \leq 1$. This gives (\ref{hd16}) with $p=\infty,\alpha_3=1$. Combining the estimates above and the interpolation inequality, we have verified (\ref{hd16}) completely.

We proceed with the proof of (\ref{hd17}) for the case $p=1$. Obviously,
\begin{equation}\nonumber
     \begin{split}
        \|  \p_3^{\alpha_3} (B_k B_l) (\tau)  \|_{L^1}  & \  \leq   2\| \p_3^{\alpha_3}  B(\tau) \|_{L^1}  \| B(\tau) \|_{L^{\infty}}     \\
       &\   \leq    \left\{\begin{aligned}
&   \| B(\tau) \|_{L^1}  \| B(\tau) \|_{L^{\infty}}             \,\,\, \text{if  } \alpha_3=0        \\
&  \| \p_3 B(\tau) \|_{L^1}  \| B(\tau) \|_{L^{\infty}}  \,\,\, \text{if  }   \alpha_3=1        \\
\end{aligned}\right.                           \\
         &\   \leq    \left\{\begin{aligned}
&   CA^2 (1+\tau)^{-\f{5}{2}}            \,\,\, \text{if  } \alpha_3=0        \\
&   CA^2 (1+\tau)^{-3}   \,\,\, \text{if  }   \alpha_3=1        \\
\end{aligned}\right.                           \\
&\ \leq CA^2 (1+\tau)^{-\f{5}{2}}  ,
     \end{split}
\end{equation}
\[
  \|  \Grad^{\alpha}   (B_k B_l) (\tau)  \|_{L^1}    \leq   2\|  \Grad^{\alpha}   B(\tau) \|_{L^{\infty}}   \| B(\tau) \|_{L^1}
  \leq
  CA^2 (1+\tau)^{-\f{5}{2}} ;
\]
whence
\begin{equation}\nonumber
     \begin{split}
       \left\|   \Grad^{\alpha} \int_0^t e^{ (t-\tau)\Deltah  }\nablah (B_k B_l ) (\tau  )     d \tau \right\|_{L^1}           \leq  & \  C   \int_0^  {  \f{t}{2} }  \| \nablah ^{\alphah} \nablah \Gh (t-\tau) \|_{L^{1}(\R^2)} \| \p_3^{\alpha_3}  (B_k B_l ) (\tau) \| _{L^1}        d \tau                     \\
         & \ +     C          \int_{  \f{t}{2} } ^t        \| \nablah \Gh (t-\tau) \|_{L^{1}(\R^2)} \| \Grad^{\alpha}  (B_k   B_l) (\tau) \|   _{L^1}          d \tau         \\
     \leq     & \   CA^2 \int_0^  {  \f{t}{2} } (t-\tau)^{ -\f{  1+|\alphah|  }{  2 }  }    (1+\tau) ^{ -\f{5}{2}     } d\tau \\
     &\ + CA^2  \int_{  \f{t}{2} } ^t    (t-\tau)^{ -\f{1}{2}   }    (1+\tau) ^{ -\f{5}{2}     } d\tau \\
\leq      &\  CA^2 t^{ -\f{  1+|\alphah|  }{  2 }  }
 \end{split}
\end{equation}
for $t \geq 1$ and
\begin{equation}\nonumber
     \begin{split}
     \left\|   \Grad^{\alpha}  \int_0^t e^{ (t-\tau)\Deltah  }\nablah (B_k  B_l  ) (\tau  )     d \tau \right\|_{L^{1}}      \leq  & \  C  \int_0^t   \| \nablah \Gh (t-\tau) \|_{L^{1}(\R^2)} \| \Grad^{\alpha}  (B_k   B_l) (\tau) \|   _{L^1}          d \tau                     \\
     \leq     & \   CA^2    \int_0^  {  t }    (t-\tau)^{ -\f{1}{2}   }    (1+\tau) ^{ -\f{5}{2}     } d\tau   \\
     \leq &\ CA^2
     \end{split}
\end{equation}
for $0<t \leq 1$. This gives (\ref{hd17}) with $p=1$. For the case $p=\infty$, one easily finds that
\begin{equation}\nonumber
     \begin{split}
         \| \p_3^{\alpha_3} (B_k B_l)(\tau) \|_{L^1_{\rm h}  L^{\infty}_{\rm v}}   &\   \leq C \|  \p_3^{\alpha_3} B(\tau) \| _{L^{\infty}}  \| B(\tau) \|  _{L^1_{\rm h}  L^{\infty}_{\rm v}}                    \\
     & \  \leq   CA^2 (1+\tau) ^{-3},             \\
     \end{split}
\end{equation}
\begin{equation}\nonumber
     \begin{split}
         \| \Grad^{\alpha} (B_k B_l)(\tau) \|_{L^{\infty}}    &\   \leq C \|  \Grad^{\alpha}  B(\tau) \| _{L^{\infty}}  \| B(\tau) \|  _{L^{\infty}}                  \\
     & \  \leq   CA^2 (1+\tau) ^{-4};            \\
     \end{split}
\end{equation}
whence
\begin{equation}\nonumber
     \begin{split}
       \left\|   \Grad^{\alpha} \int_0^t e^{ (t-\tau)\Deltah  }\nablah (B_k B_l ) (\tau  )     d \tau \right\| _{L^{\infty}}              \leq  & \  C   \int_0^  {  \f{t}{2} }  \| \nablah ^{\alphah} \nablah \Gh (t-\tau) \|_{L^{\infty}(\R^2)} \| \p_3^{\alpha_3}  (B_k B_l ) (\tau) \| _{L^1_{\rm h}  L^{\infty}_{\rm v}}        d \tau                     \\
         & \ +     C          \int_{  \f{t}{2} } ^t        \| \nablah \Gh (t-\tau) \|_{L^{1}(\R^2)} \| \Grad^{\alpha}  (B_k   B_l) (\tau) \|   _{L^{\infty}}          d \tau         \\
     \leq     & \   CA^2 \int_0^  {  \f{t}{2} } (t-\tau)^{ -1-\f{  1+|\alphah|  }{  2 }  }    (1+\tau) ^{ -3     } d\tau \\
     &\ + CA^2  \int_{  \f{t}{2} } ^t    (t-\tau)^{ -\f{1}{2}   }    (1+\tau) ^{ -4     } d\tau \\
\leq      &\  CA^2 t^{-1 -\f{  1+|\alphah|  }{  2 }  }
 \end{split}
\end{equation}
for $t \geq 1$ and
\begin{equation}\nonumber
     \begin{split}
     \left\|   \Grad^{\alpha}  \int_0^t e^{ (t-\tau)\Deltah  }\nablah (B_k  B_l  ) (\tau  )     d \tau \right\|_{L^{\infty}}       \leq  & \  C  \int_0^t   \| \nablah \Gh (t-\tau) \|_{L^{1}(\R^2)} \| \Grad^{\alpha}  (B_k   B_l) (\tau) \|   _{L^{\infty}}          d \tau                     \\
     \leq     & \   CA^2    \int_0^  {  t }    (t-\tau)^{ -\f{1}{2}   }    (1+\tau) ^{ -4   } d\tau   \\
     \leq &\ CA^2
     \end{split}
\end{equation}
for $0<t \leq 1$. This gives (\ref{hd17}) with $p=\infty$.  Combining the estimates above and the interpolation inequality, we have verified (\ref{hd17}) for any $1\leq p \leq \infty$. In the same spirit, we next handle (\ref{hd18}) with $p=1$. It follows that
\begin{equation}\nonumber
     \begin{split}
       \left\|   \Grad^{\alpha} \int_0^t K_{\beta,\gamma}^{(m)}(t-\tau) \ast (B_k B_l ) (\tau  )     d \tau \right\|_{L^1}           \leq  & \  C   \int_0^  {  \f{t}{2} }  \| \nablah ^{\alphah} K_{\beta,\gamma}^{(m)}(t-\tau)  \|_{L^{1}} \| \p_3^{\alpha_3}  (B_k B_l ) (\tau) \| _{L^1}        d \tau                     \\
         & \ +     C          \int_{  \f{t}{2} } ^t        \| K_{\beta,\gamma}^{(m)}(t-\tau)  \|_{L^{1}} \| \Grad^{\alpha}  (B_k   B_l) (\tau) \|   _{L^1}          d \tau         \\
     \leq     & \   CA^2 \int_0^  {  \f{t}{2} } (t-\tau)^{ -\f{  1+|\alphah|  }{  2 }  }    (1+\tau) ^{ -\f{5}{2}     } d\tau \\
     &\ + CA^2  \int_{  \f{t}{2} } ^t    (t-\tau)^{ -\f{1}{2}   }    (1+\tau) ^{ -\f{5}{2}     } d\tau \\
\leq      &\  CA^2 t^{ -\f{  1+|\alphah|  }{  2 }  }
 \end{split}
\end{equation}
for $t \geq 1$ and
\begin{equation}\nonumber
     \begin{split}
     \left\|   \Grad^{\alpha}  \int_0^t  K_{\beta,\gamma}^{(m)}(t-\tau) \ast (B_k B_l ) (\tau  )     d \tau \right\|_{L^{1}}      \leq  & \  C  \int_0^t    \|  K_{\beta,\gamma}^{(m)}(t-\tau)  \|_{L^{1}} \| \Grad^{\alpha}  (B_k B_l ) (\tau) \| _{L^1}        d \tau                     \\
     \leq     & \   CA^2    \int_0^  {  t }    (t-\tau)^{ -\f{1}{2}   }    (1+\tau) ^{ -\f{5}{2}     } d\tau   \\
     \leq &\ CA^2
     \end{split}
\end{equation}
for $0<t \leq 1$. This gives (\ref{hd18}) with $p=1$. In case of $p=\infty$, we deduce from Lemma \ref{lm3} that
\begin{equation}\nonumber
     \begin{split}
       \left\|   \Grad^{\alpha} \int_0^t K_{\beta,\gamma}^{(m)}(t-\tau) \ast (B_k B_l ) (\tau  )     d \tau \right\|_{L^{\infty}}           \leq  & \  C   \int_0^  {  \f{t}{2} }  \| \nablah ^{\alphah} K_{\beta,\gamma}^{(m)}(t-\tau)  \|_{L^{\infty}}    \| \p_3^{\alpha_3}  (B_k B_l ) (\tau) \| _{L^1}        d \tau                     \\
         & \ +     C          \int_{  \f{t}{2} } ^t        \| K_{\beta,\gamma}^{(m)}(t-\tau)  \|_{L^{1}} \| \Grad^{\alpha}  (B_k   B_l) (\tau) \|  _{L^{\infty}}             d \tau         \\
     \leq     & \   CA^2 \int_0^  {  \f{t}{2} } (t-\tau)^{ -\f{3} {2}-\f{  1+|\alphah|  }{  2 }  }    (1+\tau) ^{ -\f{5}{2}     } d\tau \\
     &\ + CA^2  \int_{  \f{t}{2} } ^t    (t-\tau)^{ -\f{1}{2}   }    (1+\tau) ^{ -4    } d\tau \\
\leq      &\  CA^2 t^{ -\f{3}{2}-\f{  1+|\alphah|  }{  2 }  }
 \end{split}
\end{equation}
for $t \geq 1$ and
\begin{equation}\nonumber
     \begin{split}
     \left\|   \Grad^{\alpha}  \int_0^t  K_{\beta,\gamma}^{(m)}(t-\tau) \ast (B_k B_l ) (\tau  )     d \tau \right\|_{L^{\infty}}      \leq  & \  C  \int_0^t    \|  K_{\beta,\gamma}^{(m)}(t-\tau)  \|_{L^{1}} \| \Grad^{\alpha}  (B_k B_l ) (\tau) \| _{L^{\infty}}       d \tau                     \\
     \leq     & \   CA^2    \int_0^  {  t }    (t-\tau)^{ -\f{1}{2}   }    (1+\tau) ^{ -4    } d\tau   \\
     \leq &\ CA^2
     \end{split}
\end{equation}
for $0<t \leq 1$. This gives (\ref{hd18}) with $p=\infty$. In view of the estimates above and the interpolation inequality, we have shown (\ref{hd18}) for any $1\leq p \leq \infty$.

The second part of the proof is about the decay estimates for the Duhamel terms associated with magnetic equation, namely (\ref{hd14})-(\ref{hd15}). Realizing the expressions of these Duhamel terms\footnote{Roughly speaking, $B=(\Bh,B_3),u=(\uh,u_3)$. Thus $B \otimes u$ contains four terms: $\Bh \cdot \uh, \Bh u_3,B_3 \uh,B_3 u_3$; the same is true for $u\otimes B$. Essentially, there are only two terms: $B_l u_k,B_l u_3$ for any $k\in \{ 1,2 \},l \in \{1,2,3\}$. Finally, since $\Grad =(\nablah,\p_3)$, one sees that $\Grad \cdot (B\otimes u -u \otimes B)$ contains only four terms: $\nablah (B_l u_k),\p_3 (B_l u_k), \nablah (B_l u_3), \p_3 (B_l u_3)$ for any $k\in \{ 1,2 \},l \in \{1,2,3\}$.}, it suffices to prove
\begin{equation}\label{hd19}
\left\|   \Grad^{\alpha} \int_0^t e^{ (t-\tau)\Delta  }\p_3 (B_3  u_k ) (\tau  )     d \tau \right\|_{L^p}
\leq CA^2 (1+t)^{  -\f{3}{2}(1-\f{1}{p})    -\f{ 1+|\alpha|   }{2}             },
\end{equation}
\begin{equation}\label{hd20}
\left\|   \Grad^{\alpha} \int_0^t e^{ (t-\tau)\Delta  }\p_3 (B_l  u_3 ) (\tau  )     d \tau \right\|_{L^p}
\leq CA^2 (1+t)^{  -\f{3}{2}(1-\f{1}{p})    -\f{ 1+|\alpha|   }{2}             }    ,
\end{equation}
\begin{equation}\label{hd21}
\left\|   \Grad^{\alpha} \int_0^t e^{ (t-\tau)\Delta  }\nablah (u_3  B_k ) (\tau  )     d \tau \right\|_{L^p}
\leq CA^2 (1+t)^{  -\f{3}{2}(1-\f{1}{p})    -\f{ 1+|\alpha|   }{2}             },
\end{equation}
\begin{equation}\label{hd22}
\left\|   \Grad^{\alpha} \int_0^t e^{ (t-\tau)\Delta  }\nablah (u_k  B_l ) (\tau  )     d \tau \right\|_{L^p}
\leq CA^2 (1+t)^{  -\f{3}{2}(1-\f{1}{p})    -\f{ 1+|\alpha|   }{2}             },
\end{equation}
for any $k\in \{ 1,2 \},l \in \{1,2,3\}$. We begin with the proof of (\ref{hd19}). Observe that for any $k\in \{ 1,2 \}$
\begin{equation}\nonumber
     \begin{split}
         \| (B_3 u_k) (\tau) \|_{L^1}  & \  \leq  \|B_3(\tau)  \| _{L^{\infty}}  \| u_k(\tau) \|  _{L^1}          \\
    & \  \leq CA^2  (1+\tau)^{ -2 }  ,        \\
     \end{split}
\end{equation}
\begin{equation}\nonumber
     \begin{split}
        \| \p_3(B_3 u_k) (\tau) \|     _{L^1}   & \  \leq
        \| \p_3 B_3 (\tau) \|_{L^{\infty}} \| u_k (\tau) \|  _{L^1}        +
        \|B_3(\tau)  \| _{L^{\infty}}  \|\p_3 u_k(\tau) \|  _{L^1}    \\
         & \     \leq   CA^2   (1+\tau)^{ -2 }  ;          \\
     \end{split}
\end{equation}
whence
\begin{equation}\nonumber
     \begin{split}
     \left\|   \Grad^{\alpha} \int_0^t e^{ (t-\tau)\Delta  }\p_3 (B_3 u_k ) (\tau  )     d \tau \right\|_{L^1}   \leq  & \  C  \int_0^{\f{t}{2}}  \| \p_3 \nabla^{\alpha} G (t-\tau) \|_{L^1}     \| (B_3 u_k) (\tau) \|     _{L^1}         d \tau                     \\
     &\ +  \int_{  \f{t}{2} } ^t \|  \nabla^{\alpha} G (t-\tau) \|_{L^1}     \| \p_3(B_3 u_k) (\tau) \|     _{L^1}         d \tau    \\
  \leq  & \  CA^2          \int_0^  {  \f{t}{2} } (t-\tau)^{  -\f{1+|\alpha|}{2}         } (1+\tau) ^{  -2 } d \tau                 \\
         & \ +     CA^2          \int_{  \f{t}{2} } ^t (t-\tau)^{  -\f{|\alpha|}{2}         } (1+\tau) ^{   -2       }    d\tau                 \\
     \leq     & \   CA^2 t^{    -\f{1+|\alpha|}{2}           }
     \end{split}
\end{equation}
for $t \geq 1$ and
\begin{equation}\nonumber
     \begin{split}
     \left\|   \Grad^{\alpha} \int_0^t e^{ (t-\tau)\Delta  }\p_3 (B_3 u_k ) (\tau  )     d \tau \right\|_{L^1}   \leq  & \  C  \int_0^t  \|  \nabla^{\alpha} G (t-\tau) \|_{L^1}     \| \p_3(B_3 u_k) (\tau) \|     _{L^1}         d \tau                 \\
  \leq  & \  CA^2          \int_0^  {  t } (t-\tau)^{  -\f{|\alpha|}{2}         } (1+\tau) ^{ -2 } d \tau                  \\
     \leq     & \   CA^2
     \end{split}
\end{equation}
for $0<t \leq 1$. This gives (\ref{hd19}) with $p=1$. Observe next that
\begin{equation}\nonumber
     \begin{split}
        \| \p_3(B_3 u_k) (\tau) \|     _{L^{\infty}}    & \  \leq
        \| \p_3 B_3 (\tau) \|_{L^{\infty}} \| u_k (\tau) \|  _{L^{\infty}}        +
        \|B_3(\tau)  \| _{L^{\infty}}  \|\p_3 u_k(\tau) \|   _{L^{\infty}}     \\
         & \     \leq   CA^2   (1+\tau)^{ -3 }  ;          \\
     \end{split}
\end{equation}
whence
\begin{equation}\nonumber
     \begin{split}
       \left\|   \Grad^{\alpha} \int_0^t e^{ (t-\tau)\Delta  }\p_3 (B_3 u_k ) (\tau  )     d \tau \right\| _{L^{\infty}}              \leq  & \  C   \int_0^  {  \f{t}{2} }  \| \p_3 \nabla ^{\alpha}  G (t-\tau) \|_{L^{\infty}} \|   (B_3 u_k ) (\tau) \| _{L^1}        d \tau                     \\
         & \ +     C          \int_{  \f{t}{2} } ^t        \| \Grad^{\alpha} G (t-\tau) \|_{L^{1}} \| \p_3  (B_3   u_k) (\tau) \|   _{L^{\infty}}          d \tau         \\
     \leq     & \   CA^2 \int_0^  {  \f{t}{2} } (t-\tau)^{ -\f{3}{2}-\f{  1+|\alpha|  }{  2 }  }    (1+\tau) ^{ -2     } d\tau \\
     &\ + CA^2  \int_{  \f{t}{2} } ^t    (t-\tau)^{ -\f{|\alpha|}{2}   }    (1+\tau) ^{ -3     } d\tau \\
\leq      &\  CA^2 t^{-\f{3}{2} -\f{  1+|\alpha|  }{  2 }  }
 \end{split}
\end{equation}
for $t \geq 1$ and
\begin{equation}\nonumber
     \begin{split}
     \left\|   \Grad^{\alpha} \int_0^t e^{ (t-\tau)\Delta  }\p_3 (B_3 u_k ) (\tau  )     d \tau \right\|_{L^{\infty}}   \leq  & \  C  \int_0^t  \|  \nabla^{\alpha} G (t-\tau) \|_{L^1}     \| \p_3(B_3 u_k) (\tau) \|     _{L^{\infty}}         d \tau                 \\
  \leq  & \  CA^2          \int_0^  {  t } (t-\tau)^{  -\f{|\alpha|}{2}         } (1+\tau) ^{ -3 } d \tau                  \\
     \leq     & \   CA^2
     \end{split}
\end{equation}
for $0<t \leq 1$. This gives (\ref{hd19}) with $p=\infty$. We thus obtain (\ref{hd19}) via the interpolation inequality.

Next we turn to (\ref{hd20}). Observe that for any $l\in \{ 1,2 ,3\}$
\begin{equation}\nonumber
     \begin{split}
         \| (B_l u_3) (\tau) \|_{L^1}  & \  \leq  \|B_l(\tau)  \| _{L^{\infty}}  \| u_3(\tau) \|  _{L^1}          \\
    & \  \leq CA^2  (1+\tau)^{ -2 }  ,        \\
     \end{split}
\end{equation}
\begin{equation}\nonumber
     \begin{split}
        \| \p_3(B_l u_3) (\tau) \|     _{L^1}   & \  \leq
        \| \p_3 B_l (\tau) \|_{L^{\infty}} \| u_3 (\tau) \|  _{L^1}        +
        \|B_l(\tau)  \| _{L^{\infty}}  \|\p_3 u_3(\tau) \|  _{L^1}    \\
         & \     \leq   CA^2   (1+\tau)^{ -\f{5}{2} }  ;          \\
     \end{split}
\end{equation}
whence
\begin{equation}\nonumber
     \begin{split}
     \left\|   \Grad^{\alpha} \int_0^t e^{ (t-\tau)\Delta  }\p_3 (B_l u_3 ) (\tau  )     d \tau \right\|_{L^1}   \leq  & \  C  \int_0^{\f{t}{2}}  \| \p_3 \nabla^{\alpha} G (t-\tau) \|_{L^1}     \| (B_l u_3) (\tau) \|     _{L^1}         d \tau                     \\
     &\ +  \int_{  \f{t}{2} } ^t \|  \nabla^{\alpha} G (t-\tau) \|_{L^1}     \| \p_3(B_l u_3) (\tau) \|     _{L^1}         d \tau    \\
  \leq  & \  CA^2          \int_0^  {  \f{t}{2} } (t-\tau)^{  -\f{1+|\alpha|}{2}         } (1+\tau) ^{  -2 } d \tau                 \\
         & \ +     CA^2          \int_{  \f{t}{2} } ^t (t-\tau)^{  -\f{|\alpha|}{2}         } (1+\tau) ^{   -\f{5}{2}       }    d\tau                 \\
     \leq     & \   CA^2 t^{    -\f{1+|\alpha|}{2}           }
     \end{split}
\end{equation}
for $t \geq 1$ and
\begin{equation}\nonumber
     \begin{split}
     \left\|   \Grad^{\alpha} \int_0^t e^{ (t-\tau)\Delta  }\p_3 (B_l u_3 ) (\tau  )     d \tau \right\|_{L^1}   \leq  & \  C  \int_0^t  \|  \nabla^{\alpha} G (t-\tau) \|_{L^1}     \| \p_3(B_l u_3) (\tau) \|     _{L^1}         d \tau                 \\
  \leq  & \  CA^2          \int_0^  {  t } (t-\tau)^{  -\f{|\alpha|}{2}         } (1+\tau) ^{ -\f{5}{2} } d \tau                  \\
     \leq     & \   CA^2
     \end{split}
\end{equation}
for $0<t \leq 1$. This gives (\ref{hd20}) with $p=1$. Notice also that
\begin{equation}\nonumber
     \begin{split}
        \| \p_3(B_l u_3) (\tau) \|     _{L^{\infty}}    & \  \leq
        \| \p_3 B_l (\tau) \|_{L^{\infty}} \| u_3 (\tau) \|  _{L^{\infty}}        +
        \|B_l(\tau)  \| _{L^{\infty}}  \|\p_3 u_3(\tau) \|   _{L^{\infty}}     \\
         & \     \leq   CA^2   (1+\tau)^{ -\f{7}{2} }  ;          \\
     \end{split}
\end{equation}
whence
\begin{equation}\nonumber
     \begin{split}
       \left\|   \Grad^{\alpha} \int_0^t e^{ (t-\tau)\Delta  }\p_3 (B_l u_3 ) (\tau  )     d \tau \right\| _{L^{\infty}}              \leq  & \  C   \int_0^  {  \f{t}{2} }  \| \p_3 \nabla ^{\alpha}  G (t-\tau) \|_{L^{\infty}} \|   (B_l u_3 ) (\tau) \| _{L^1}        d \tau                     \\
         & \ +     C          \int_{  \f{t}{2} } ^t        \| \Grad^{\alpha} G (t-\tau) \|_{L^{1}} \| \p_3  (B_l   u_3) (\tau) \|   _{L^{\infty}}          d \tau         \\
     \leq     & \   CA^2 \int_0^  {  \f{t}{2} } (t-\tau)^{ -\f{3}{2}-\f{  1+|\alpha|  }{  2 }  }    (1+\tau) ^{ -2     } d\tau \\
     &\ + CA^2  \int_{  \f{t}{2} } ^t    (t-\tau)^{ -\f{|\alpha|}{2}   }    (1+\tau) ^{ -\f{7}{2}     } d\tau \\
\leq      &\  CA^2 t^{-\f{3}{2} -\f{  1+|\alpha|  }{  2 }  }
 \end{split}
\end{equation}
for $t \geq 1$ and
\begin{equation}\nonumber
     \begin{split}
     \left\|   \Grad^{\alpha} \int_0^t e^{ (t-\tau)\Delta  }\p_3 (B_3 u_k ) (\tau  )     d \tau \right\|_{L^{\infty}}   \leq  & \  C  \int_0^t  \|  \nabla^{\alpha} G (t-\tau) \|_{L^1}     \| \p_3(B_3 u_k) (\tau) \|     _{L^{\infty}}         d \tau                 \\
  \leq  & \  CA^2          \int_0^  {  t } (t-\tau)^{  -\f{|\alpha|}{2}         } (1+\tau) ^{ -\f{7}{2} } d \tau                  \\
     \leq     & \   CA^2
     \end{split}
\end{equation}
for $0<t \leq 1$. This gives (\ref{hd20}) with $p=\infty$. (\ref{hd20}) follows immediately via the interpolation inequality.

We proceed with (\ref{hd21}) using similar techniques. For any $k \in \{ 1,2 \}$
\begin{equation}\nonumber
     \begin{split}
         \| ( u_3  B_k ) (\tau) \|_{L^1}  & \  \leq   \| u_3(\tau) \|  _{L^1}     \|B_k(\tau)  \| _{L^{\infty}}       \\
    & \  \leq CA^2  (1+\tau)^{ -2 }  ,        \\
     \end{split}
\end{equation}
\begin{equation}\nonumber
     \begin{split}
        \| \nablah ( u_3 B_k) (\tau) \|     _{L^1}   & \  \leq
        \| \nablah u_3 (\tau) \|_{L^{\infty}} \| B_k (\tau) \|  _{L^1}        +
        \| u_3(\tau) \|  _{L^1}  \|\nablah B_k(\tau)  \| _{L^{\infty}}     \\
         & \     \leq   CA^2   (1+\tau)^{ -\f{5}{2} }  ;          \\
     \end{split}
\end{equation}
whence
\begin{equation}\nonumber
     \begin{split}
     \left\|   \Grad^{\alpha} \int_0^t e^{ (t-\tau)\Delta  }\nablah ( u_3 B_k ) (\tau  )     d \tau \right\|_{L^1}   \leq  & \  C  \int_0^{\f{t}{2}}  \| \nablah \nabla^{\alpha} G (t-\tau) \|_{L^1}     \| ( u_3 B_k) (\tau) \|  _{L^1}         d \tau                     \\
     &\ +  \int_{  \f{t}{2} } ^t \|  \nabla^{\alpha} G (t-\tau) \|_{L^1}     \| \nablah ( u_3  B_k ) (\tau) \|     _{L^1}         d \tau    \\
  \leq  & \  CA^2          \int_0^  {  \f{t}{2} } (t-\tau)^{  -\f{1+|\alpha|}{2}         } (1+\tau) ^{  -2 } d \tau                 \\
         & \ +     CA^2          \int_{  \f{t}{2} } ^t (t-\tau)^{  -\f{|\alpha|}{2}         } (1+\tau) ^{   -\f{5}{2}       }    d\tau                 \\
     \leq     & \   CA^2 t^{    -\f{1+|\alpha|}{2}           }
     \end{split}
\end{equation}
for $t \geq 1$ and
\begin{equation}\nonumber
     \begin{split}
     \left\|   \Grad^{\alpha} \int_0^t e^{ (t-\tau)\Delta  }\nablah ( u_3 B_k ) (\tau  )     d \tau \right\|_{L^1}   \leq  & \  C  \int_0^t  \|  \nabla^{\alpha} G (t-\tau) \|_{L^1}     \| \nablah( u_3 B_k) (\tau) \|     _{L^1}         d \tau                 \\
  \leq  & \  CA^2          \int_0^  {  t } (t-\tau)^{  -\f{|\alpha|}{2}         } (1+\tau) ^{ -\f{5}{2} } d \tau                  \\
     \leq     & \   CA^2
     \end{split}
\end{equation}
for $0<t \leq 1$. This gives (\ref{hd21}) with $p=1$. Seeing that
\begin{equation}\nonumber
     \begin{split}
        \| \nablah ( u_3 B_k ) (\tau) \|     _{L^{\infty}}    & \  \leq
        \| \nablah B_k (\tau) \|_{L^{\infty}} \| u_3 (\tau) \|  _{L^{\infty}}        +
        \|B_k(\tau)  \| _{L^{\infty}}  \|\nablah u_3(\tau) \|   _{L^{\infty}}     \\
         & \     \leq   CA^2   (1+\tau)^{ -4 }  ;          \\
     \end{split}
\end{equation}
whence
\begin{equation}\nonumber
     \begin{split}
       \left\|   \Grad^{\alpha} \int_0^t e^{ (t-\tau)\Delta  }\nablah ( u_3 B_k ) (\tau  )     d \tau \right\| _{L^{\infty}}              \leq  & \  C   \int_0^  {  \f{t}{2} }  \| \nablah \nabla ^{\alpha}  G (t-\tau) \|_{L^{\infty}} \|   ( u_3 B_k) (\tau) \| _{L^1}        d \tau                     \\
         & \ +     C          \int_{  \f{t}{2} } ^t        \| \Grad^{\alpha} G (t-\tau) \|_{L^{1}} \| \nablah   ( u_3 B_k)(\tau) \|   _{L^{\infty}}          d \tau         \\
     \leq     & \   CA^2 \int_0^  {  \f{t}{2} } (t-\tau)^{ -\f{3}{2}-\f{  1+|\alpha|  }{  2 }  }    (1+\tau) ^{ -2     } d\tau \\
     &\ + CA^2  \int_{  \f{t}{2} } ^t    (t-\tau)^{ -\f{|\alpha|}{2}   }    (1+\tau) ^{ -4    } d\tau \\
\leq      &\  CA^2 t^{-\f{3}{2} -\f{  1+|\alpha|  }{  2 }  }
 \end{split}
\end{equation}
for $t \geq 1$ and
\begin{equation}\nonumber
     \begin{split}
     \left\|   \Grad^{\alpha} \int_0^t e^{ (t-\tau)\Delta  }\nablah ( u_3 B_k ) (\tau  )     d \tau \right\|_{L^{\infty}}   \leq  & \  C  \int_0^t  \|  \nabla^{\alpha} G (t-\tau) \|_{L^1}     \| \nablah ( u_3 B_k ) (\tau) \|     _{L^{\infty}}         d \tau                 \\
  \leq  & \  CA^2          \int_0^  {  t } (t-\tau)^{  -\f{|\alpha|}{2}         } (1+\tau) ^{ - 4 } d \tau                  \\
     \leq     & \   CA^2
     \end{split}
\end{equation}
for $0<t \leq 1$. This gives (\ref{hd21}) with $p=\infty$. (\ref{hd21}) follows readily via the interpolation inequality.

Finally, we show (\ref{hd22}). For any $k \in \{ 1,2 \},l \in \{ 1,2,3 \} $
\begin{equation}\nonumber
     \begin{split}
         \| ( u_k  B_l ) (\tau) \|_{L^1}  & \  \leq   \| u_k(\tau) \|  _{L^1}     \|B_l(\tau)  \| _{L^{\infty}}       \\
    & \  \leq CA^2  (1+\tau)^{ -2 }  ,        \\
     \end{split}
\end{equation}
\begin{equation}\nonumber
     \begin{split}
        \| \nablah ( u_k B_l) (\tau) \|     _{L^1}   & \  \leq
        \| \nablah B_l (\tau) \|_{L^{\infty}} \| u_k (\tau) \|  _{L^1}        +
        \| B_l (\tau) \|  _{L^{\infty}}    \|\nablah u_k (\tau)  \| _{L^1}       \\
         & \     \leq   CA^2   (1+\tau)^{ -\f{5}{2} }  ;          \\
     \end{split}
\end{equation}
whence
\begin{equation}\nonumber
     \begin{split}
     \left\|   \Grad^{\alpha} \int_0^t e^{ (t-\tau)\Delta  }\nablah ( u_k B_l ) (\tau  )     d \tau \right\|_{L^1}   \leq  & \  C  \int_0^{\f{t}{2}}  \| \nablah \nabla^{\alpha} G (t-\tau) \|_{L^1}     \| ( u_k B_l) (\tau) \|  _{L^1}         d \tau                     \\
     &\ +  \int_{  \f{t}{2} } ^t \|  \nabla^{\alpha} G (t-\tau) \|_{L^1}     \| \nablah ( u_k  B_l ) (\tau) \|     _{L^1}         d \tau    \\
  \leq  & \  CA^2          \int_0^  {  \f{t}{2} } (t-\tau)^{  -\f{1+|\alpha|}{2}         } (1+\tau) ^{  -2 } d \tau                 \\
         & \ +     CA^2          \int_{  \f{t}{2} } ^t (t-\tau)^{  -\f{|\alpha|}{2}         } (1+\tau) ^{   -\f{5}{2}       }    d\tau                 \\
     \leq     & \   CA^2 t^{    -\f{1+|\alpha|}{2}           }
     \end{split}
\end{equation}
for $t \geq 1$ and
\begin{equation}\nonumber
     \begin{split}
     \left\|   \Grad^{\alpha} \int_0^t e^{ (t-\tau)\Delta  }\nablah ( u_3 B_k ) (\tau  )     d \tau \right\|_{L^1}   \leq  & \  C  \int_0^t  \|  \nabla^{\alpha} G (t-\tau) \|_{L^1}     \| \nablah( u_3 B_k) (\tau) \|     _{L^1}         d \tau                 \\
  \leq  & \  CA^2          \int_0^  {  t } (t-\tau)^{  -\f{|\alpha|}{2}         } (1+\tau) ^{ -\f{5}{2} } d \tau                  \\
     \leq     & \   CA^2
     \end{split}
\end{equation}
for $0<t \leq 1$. This gives (\ref{hd22}) with $p=1$. In light of
\begin{equation}\nonumber
     \begin{split}
        \| \nablah ( u_k B_l ) (\tau) \|     _{L^{\infty}}    & \  \leq
        \| \nablah B_l (\tau) \|_{L^{\infty}} \| u_k (\tau) \|  _{L^{\infty}}        +
        \|B_l(\tau)  \| _{L^{\infty}}  \|\nablah u_k(\tau) \|   _{L^{\infty}}     \\
         & \     \leq   CA^2   (1+\tau)^{ -\f{7}{2} }  ;          \\
     \end{split}
\end{equation}
whence
\begin{equation}\nonumber
     \begin{split}
       \left\|   \Grad^{\alpha} \int_0^t e^{ (t-\tau)\Delta  }\nablah ( u_k B_l ) (\tau  )     d \tau \right\| _{L^{\infty}}              \leq  & \  C   \int_0^  {  \f{t}{2} }  \| \nablah \nabla ^{\alpha}  G (t-\tau) \|_{L^{\infty}} \|   ( u_k B_l) (\tau) \| _{L^1}        d \tau                     \\
         & \ +     C          \int_{  \f{t}{2} } ^t        \| \Grad^{\alpha} G (t-\tau) \|_{L^{1}} \| \nablah   ( u_k B_l)(\tau) \|   _{L^{\infty}}          d \tau         \\
     \leq     & \   CA^2 \int_0^  {  \f{t}{2} } (t-\tau)^{ -\f{3}{2}-\f{  1+|\alpha|  }{  2 }  }    (1+\tau) ^{ -2     } d\tau \\
     &\ + CA^2  \int_{  \f{t}{2} } ^t    (t-\tau)^{ -\f{|\alpha|}{2}   }    (1+\tau) ^{ -\f{7}{2}    } d\tau \\
\leq      &\  CA^2 t^{-\f{3}{2} -\f{  1+|\alpha|  }{  2 }  }
 \end{split}
\end{equation}
for $t \geq 1$ and
\begin{equation}\nonumber
     \begin{split}
     \left\|   \Grad^{\alpha} \int_0^t e^{ (t-\tau)\Delta  }\nablah ( u_k B_l ) (\tau  )     d \tau \right\|_{L^{\infty}}   \leq  & \  C  \int_0^t  \|  \nabla^{\alpha} G (t-\tau) \|_{L^1}     \| \nablah ( u_k B_l ) (\tau) \|     _{L^{\infty}}         d \tau                 \\
  \leq  & \  CA^2          \int_0^  {  t } (t-\tau)^{  -\f{|\alpha|}{2}         } (1+\tau) ^{ - \f{7}{2} } d \tau                  \\
     \leq     & \   CA^2
     \end{split}
\end{equation}
for $0<t \leq 1$. This gives (\ref{hd22}) with $p=\infty$. (\ref{hd21}) follows as before via the interpolation inequality. The proof of Lemma \ref{lem7} is finished completely.         $\Box$

Given more regular initial data, we are able to show the decay estimates for the Duhamel terms in $L^{\infty}_{\rm h}  L^{1}_{\rm v}$. These estimates will play crucial role in higher order asymptotic expansions of solutions.
\begin{lemm}\label{lem8}
Let $(u,B)$ be subject to the ansatzes $(i),(ii),(iii)$ with $s\geq 9$. There exists a generic constant $C>0$ such that for any $0<t<T$ and any $\alpha=(\alphah,\alpha_3)\in (\mathbb{N} \cup \{ 0\})^2 \times (\mathbb{N} \cup \{ 0\})$ with $|\alpha|\leq 1$ we have the decay estimates for the Duhamel terms associated with velocity equation:
\begin{equation}\label{hd23}
     \begin{split}
         \|\Grad^{\alpha}  \Nhu_m[u](t)  \|_{L^{\infty}_{\rm h}  L^{1}_{\rm v}} & \  \leq  CA^2 t^{ -1-\f{|\alphah|}{2}  },        \\
   \|\Grad^{\alpha}  \Nvu_n[u](t)  \|_{L^{\infty}_{\rm h}  L^{1}_{\rm v}} & \  \leq  CA^2 t^{ -1-\f{|\alphah|}{2}  } ,       \\
    \|\Grad^{\alpha}  \Nhu_m[B](t)  \|_{L^{\infty}_{\rm h}  L^{1}_{\rm v}} & \  \leq  CA^2 t^{ -1-\f{|\alphah|}{2}  },        \\
   \|\Grad^{\alpha}  \Nvu_n[B](t)  \|_{L^{\infty}_{\rm h}  L^{1}_{\rm v}} & \  \leq  CA^2 t^{ -1-\f{|\alphah|}{2}  } ,       \\
     \end{split}
\end{equation}
for any $m\in \{ 1,2,3,4,5 \},n\in \{ 1,2,3 \}$ and the decay estimates for the Duhamel terms associated with magnetic equation:
\begin{equation}\label{hd24}
     \begin{split}
 \|\Grad^{\alpha}  \Nhb_n[u,B](t)  \|_{L^{\infty}_{\rm h}  L^{1}_{\rm v}} & \  \leq  CA^2 t^{ -1-\f{1+|\alpha|}{2}  } ,       \\
 \|\Grad^{\alpha}  \Nvb_n[u,B](t)  \|_{L^{\infty}_{\rm h}  L^{1}_{\rm v}} & \  \leq  CA^2 t^{ -1-\f{1+|\alpha|}{2}  } ,       \\
     \end{split}
\end{equation}
for any $n\in \{ 1,2,3 \}$.
\end{lemm}
{\it Proof.} Observe that the first two estimates in (\ref{hd23}) have been proved in \cite{F}, see Lemma 4.3 therein. It should be remarked that the assumption $s \geq 9$ comes mainly from these two estimates, due to the weak dissipation mechanism in the velocity equation. To show the last two estimates in (\ref{hd23}) involved with magnetic field, it suffices to verify
\begin{equation}\label{hd25}
\left\|   \Grad^{\alpha} \int_0^t e^{ (t-\tau)\Deltah  }\p_3 (B_3\Bh ) (\tau  )     d \tau \right\|_{L^{\infty}_{\rm h}  L^{1}_{\rm v}}
\leq  CA^2 t^{ -1-\f{|\alphah|}{2}  }            ,
\end{equation}
\begin{equation}\label{hd26}
\left\|   \Grad^{\alpha} \int_0^t e^{ (t-\tau)\Deltah  }\nablah (B_k B_l ) (\tau  )     d \tau \right\|_{L^{\infty}_{\rm h}  L^{1}_{\rm v}}
\leq  CA^2 t^{ -1-\f{|\alphah|}{2}  }      ,
\end{equation}
\begin{equation}\label{hd27}
\left\|   \Grad^{\alpha} \int_0^t   K_{\beta,\gamma}^{ (m)} (t-\tau) \ast (B_k B_l ) (\tau  )     d \tau \right\|_{L^{\infty}_{\rm h}  L^{1}_{\rm v}}
\leq    CA^2 t^{ -1-\f{|\alphah|}{2}  }       ,
\end{equation}
for any $k,l\in \{1,2,3\}$ and $m\in \{1,2\}$. $(\beta,\gamma)  \in (\mathbb{N} \cup \{ 0\})^2 \times (\mathbb{N} \cup \{ 0\})$ obeying $|\beta|+\gamma=3$. We begin with (\ref{hd25}). For the case $\alpha_3=0$,
\begin{equation}\nonumber
     \begin{split}
        \| \p_3(B_3 \Bh) \|_{L^{\infty}_{\rm h}  L^{1}_{\rm v}}  & \  \leq  \| \p_3 B_3 \|_{L^{\infty}}    \| \Bh \|_{L^{\infty}_{\rm h}  L^{1}_{\rm v}}   +    \|  B_3 \|_{L^{\infty}}    \| \p_3  \Bh \|_{L^{\infty}_{\rm h}  L^{1}_{\rm v}}        \\
    & \  \leq  CA^2 \left(  (1+\tau)^{-\f{5}{2}} \tau^{-\f{3}{2}} +  (1+\tau)^{-2} \tau^{-2 }      \right) ;        \\
     \end{split}
\end{equation}
whence
\begin{equation}\nonumber
     \begin{split}
     \Big\|   \nablah^{\alphah}  &\  \int_0^{ t } e^{ (t-\tau)\Deltah  }\p_3 (B_3\Bh ) (\tau  )     d \tau \Big\|_{L^{\infty}_{\rm h}  L^{1}_{\rm v}}    \\
     \leq  & \    \int_0^{\f{t}{2}}  \| \nablah^{\alphah} \Gh (t-\tau) \|_{L^{\infty}(\R^2)} \| \p_3 (B_3 \Bh) (\tau) \|     _{L^1}         d \tau                     \\
     &\ + \int_{  \f{t}{2} } ^t  \| \nablah^{\alphah} \Gh (t-\tau) \|_{L^{1}(\R^2)}  \| \p_3 (B_3 \Bh) (\tau) \|  _{L^{\infty}_{\rm h}  L^{1}_{\rm v}}  d \tau         \\
  \leq  & \  CA^2          \int_0^  {  \f{t}{2} } (t-\tau)^{  -1-\f{|\alphah|}{2}         } (1+\tau) ^{  -3 } d \tau                 \\
         & \ +     CA^2          \int_{  \f{t}{2} } ^t (t-\tau)^{  -\f{|\alphah|}{2}         } \left(  (1+\tau)^{-\f{5}{2}} \tau^{-\f{3}{2}} +  (1+\tau)^{-2} \tau^{-2 }      \right)   d\tau                 \\
     \leq     & \   CA^2 t^{  -1  -\f{|\alphah|}{2}           }. \\
     \end{split}
\end{equation}
For the case $\alpha_3=1$,
\begin{equation}\nonumber
     \begin{split}
        \| \p_3^2(B_3 \Bh) \|_{L^{\infty}_{\rm h}  L^{1}_{\rm v}} \leq  & \    \| \p_3^2 B_3 \|_{L^{\infty}} \|\Bh \| _{L^{\infty}_{\rm h}  L^{1}_{\rm v}}   +2 \| \p_3 B_3 \|_{L^{\infty}}   \| \p_3 \Bh \|_{L^{\infty}_{\rm h}  L^{1}_{\rm v}}          \\
        &\ +  \| \p_3^2 \Bh \|_{L^{\infty}}   \| B_3 \| _{L^{\infty}_{\rm h}  L^{1}_{\rm v}}  \\
       \leq  &\  C \Big( \| B_3 \|_{H^7} ^{\f{2}{7}}    \|\p_3 B_3 \|_{L^{\infty}} ^{\f{5}{7}}
       \|\Bh \| _{L^{\infty}_{\rm h}  L^{1}_{\rm v}}  + \| \p_3 B_3 \|_{L^{\infty}}   \| \p_3 \Bh \|_{L^{\infty}_{\rm h}  L^{1}_{\rm v}}  \\
       &\ +   \| \Bh \|_{H^7} ^{\f{2}{7}}    \|\p_3 \Bh \|_{L^{\infty}} ^{\f{5}{7}}
       \|B_3 \| _{L^{\infty}_{\rm h}  L^{1}_{\rm v}}   \Big)      \\
  \leq  & \    CA^2 \left(  (1+\tau)^{ - \f{25}{14}   } \tau^{-\f{3}{2}} +  (1+\tau)^{-\f{5}{2}} \tau^{-2 }      \right) ;        \\
     \end{split}
\end{equation}
whence
\begin{equation}\nonumber
     \begin{split}
     \left\|  \p_3 \int_0^{ t } e^{ (t-\tau)\Deltah  }\p_3 (B_3\Bh ) (\tau  )     d \tau \right\|_{L^{\infty}_{\rm h}  L^{1}_{\rm v}}   \leq  & \    \int_0^{\f{t}{2}}  \|  \Gh (t-\tau) \|_{L^{\infty}(\R^2)} \| \p_3^2 (B_3 \Bh) (\tau) \|     _{L^1}         d \tau                     \\
     &\ + \int_{  \f{t}{2} } ^t  \| \Gh (t-\tau) \|_{L^{1}(\R^2)}  \| \p_3^2 (B_3 \Bh) (\tau) \|  _{L^{\infty}_{\rm h}  L^{1}_{\rm v}}  d \tau         \\
  \leq  & \  CA^2          \int_0^  {  \f{t}{2} } (t-\tau)^{  -1     } (1+\tau) ^{  -\f{5}{4} } d \tau                 \\
         & \ +     CA^2          \int_{  \f{t}{2} } ^t
        \left(  (1+\tau)^{ - \f{25}{14}   } \tau^{-\f{3}{2}}+  (1+\tau)^{-\f{5}{2}} \tau^{-2 }      \right)   d\tau                 \\
     \leq     & \   CA^2 t^{  -1           }. \\
     \end{split}
\end{equation}
The proof of (\ref{hd25}) is thus finished. We next come to (\ref{hd26}). In view of
\begin{equation}\nonumber
     \begin{split}
        \| \Grad^{\alpha} (B_k B_l)(\tau) \|_{L^{\infty}_{\rm h}  L^{1}_{\rm v}}  \leq  & \  C  \|  B(\tau)  \| _{L^{\infty}_{\rm h}  L^{1}_{\rm v}}  \|  \Grad^{\alpha}  B(\tau) \|_{L^{\infty}}      \\
    \leq & \   CA^2 \tau^{-\f{3}{2}} (1+\tau)^{-1-\f{|\alphah|}{2}} ,         \\
     \end{split}
\end{equation}
\begin{equation}\nonumber
     \begin{split}
        \| \p_3^{\alpha_3} (B_k B_l)(\tau) \|_{L^{2}_{\rm h}  L^{1}_{\rm v}}  \leq  & \    \|  \p_3^{\alpha_3} (B_k B_l)(\tau) \| _{ L^{1}} ^{ \f{1}{2}  }    \|  \p_3^{\alpha_3} (B_k B_l)(\tau) \|_{L^{\infty}_{\rm h}  L^{1}_{\rm v}} ^{ \f{1}{2}  }      \\
    \leq & \   CA^2 \tau^{-\f{3}{4}} (1+\tau)^{-\f{9}{4}} ;        \\
     \end{split}
\end{equation}
whence
\begin{equation}\nonumber
     \begin{split}
     \left\|   \nabla^{\alpha} \int_0^{ t } e^{ (t-\tau)\Deltah  } \nablah (B_k B_l ) (\tau  )     d \tau \right\|_{L^{\infty}_{\rm h}  L^{1}_{\rm v}}   \leq  & \    \int_0^{\f{t}{2}}  \| \nablah^{\alphah} \nablah \Gh (t-\tau) \|_{L^{2}(\R^2)} \| \p_3 ^{\alpha_3}(B_k  B_l) (\tau) \|   _{L^{2}_{\rm h}  L^{1}_{\rm v}}        d \tau                     \\
     &\ + \int_{  \f{t}{2} } ^t  \| \nablah^{\alphah} \Gh (t-\tau) \|_{L^{1}(\R^2)}  \| \nabla^{\alpha} (B_k  B_l) (\tau) \|  _{L^{\infty}_{\rm h}  L^{1}_{\rm v}}  d \tau         \\
  \leq  & \  CA^2          \int_0^  {  \f{t}{2} } (t-\tau)^{  -1-\f{|\alphah|}{2}         } (1+\tau) ^{  -\f{9}{4} } \tau^{-\f{3}{4}} d \tau                 \\
         & \ +     CA^2          \int_{  \f{t}{2} } ^t (t-\tau)^{  -\f{1}{2}         }   (1+\tau)^{-1-\f{|\alphah|}{2}} \tau^{-\f{3}{2}}    d\tau                 \\
     \leq     & \   CA^2 t^{  -1  -\f{|\alphah|}{2}           }. \\
     \end{split}
\end{equation}
This verifies (\ref{hd26}). Exactly in the same way, we treat (\ref{hd27}) as follows:
\begin{equation}\nonumber
     \begin{split}
     \left\|   \nabla^{\alpha} \int_0^{ t }  K_{\beta,\gamma}^{ (m)} \ast  (B_k B_l ) (\tau  )     d \tau \right\|_{L^{\infty}_{\rm h}  L^{1}_{\rm v}}   \leq  & \    \int_0^{\f{t}{2}}  \| \nablah^{\alphah}  K_{\beta,\gamma}^{ (m)}(t-\tau)\| _{L^{2}_{\rm h}  L^{1}_{\rm v}}    \| \p_3 ^{\alpha_3}(B_k  B_l) (\tau) \|   _{L^{2}_{\rm h}  L^{1}_{\rm v}}        d \tau                     \\
     &\ + \int_{  \f{t}{2} } ^t  \| K_{\beta,\gamma}^{ (m)}(t-\tau)  \|_{L^{1}}  \| \nabla^{\alpha} (B_k  B_l) (\tau) \|  _{L^{\infty}_{\rm h}  L^{1}_{\rm v}}  d \tau         \\
  \leq  & \  CA^2          \int_0^  {  \f{t}{2} } (t-\tau)^{  -1-\f{|\alphah|}{2}         } (1+\tau) ^{  -\f{9}{4} } \tau^{-\f{3}{4}} d \tau                 \\
         & \ +     CA^2          \int_{  \f{t}{2} } ^t (t-\tau)^{  -\f{1}{2}         }   (1+\tau)^{-1  -\f{|\alphah|}{2} } \tau^{-\f{3}{2}}    d\tau                 \\
     \leq     & \   CA^2 t^{  -1  -\f{|\alphah|}{2}           }. \\
     \end{split}
\end{equation}

The second part of the proof is about the decay estimates for the Duhamel terms associated with magnetic equation, namely (\ref{hd24}). As in Lemma \ref{lem7}, it suffices to prove
\begin{equation}\label{hd28}
\left\|   \Grad^{\alpha} \int_0^t e^{ (t-\tau)\Delta  }\p_3 (B_3  u_k ) (\tau  )     d \tau \right\|_{L^{\infty}_{\rm h}  L^{1}_{\rm v}}
\leq CA^2 t^{  -1   -\f{ 1+|\alpha|   }{2}             },
\end{equation}
\begin{equation}\label{hd29}
\left\|   \Grad^{\alpha} \int_0^t e^{ (t-\tau)\Delta  }\p_3 (B_l  u_3 ) (\tau  )     d \tau \right\|_{L^{\infty}_{\rm h}  L^{1}_{\rm v}}
\leq CA^2 t^{  -1    -\f{ 1+|\alpha|   }{2}             }    ,
\end{equation}
\begin{equation}\label{hd30}
\left\|   \Grad^{\alpha} \int_0^t e^{ (t-\tau)\Delta  }\nablah (u_3  B_k ) (\tau  )     d \tau \right\|_{L^{\infty}_{\rm h}  L^{1}_{\rm v}}
\leq CA^2 t^{  -1   -\f{ 1+|\alpha|   }{2}             },
\end{equation}
\begin{equation}\label{hd31}
\left\|   \Grad^{\alpha} \int_0^t e^{ (t-\tau)\Delta  }\nablah (u_k  B_l ) (\tau  )     d \tau \right\|_{L^{\infty}_{\rm h}  L^{1}_{\rm v}}
\leq CA^2 t^{  -1   -\f{ 1+|\alpha|   }{2}             },
\end{equation}
for any $k\in \{ 1,2 \},l \in \{1,2,3\}$. To show (\ref{hd28}), we notice that for any $k\in \{ 1,2 \}$
\begin{equation}\nonumber
     \begin{split}
        \| \p_3(B_3 u_k) \|_{L^{\infty}_{\rm h}  L^{1}_{\rm v}}  & \  \leq  \| \p_3 B_3 \|_{L^{\infty}_{\rm h}  L^{1}_{\rm v}}      \| u_k \|_{L^{\infty}}    +    \|  B_3 \|_{L^{\infty}_{\rm h}  L^{1}_{\rm v}}       \| \p_3  u_k \|   _{L^{\infty}}     \\
    & \  \leq  CA^2 \left(  (1+\tau)^{-1} \tau^{-2} +  (1+\tau)^{-1} \tau^{-\f{3}{2} }      \right) ;        \\
     \end{split}
\end{equation}
yielding
\begin{equation}\nonumber
     \begin{split}
     \Big \|   \nabla^{\alpha}   &\        \int_0^{ t } e^{ (t-\tau)\Delta  }\p_3 (B_3 u_k ) (\tau  )     d \tau \Big\|_{L^{\infty}_{\rm h}  L^{1}_{\rm v}}   \\
     \leq  & \    \int_0^{\f{t}{2}}  \| \p_3 \nabla^{\alpha} G (t-\tau) \|_{L^{\infty}_{\rm h}  L^{1}_{\rm v}} \|  (B_3 u_k) (\tau) \|     _{L^1}         d \tau                     \\
     &\ + \int_{  \f{t}{2} } ^t  \| \nabla^{\alpha} G (t-\tau) \|_{L^{1}}  \| \p_3 (B_3 u_k) (\tau) \|  _{L^{\infty}_{\rm h}  L^{1}_{\rm v}}  d \tau          \\
  \leq  & \  CA^2          \int_0^  {  \f{t}{2} } (t-\tau)^{  -1-\f{1+ |\alpha|}{2}         } (1+\tau) ^{  -2 } d \tau                 \\
         & \ +     CA^2          \int_{  \f{t}{2} } ^t (t-\tau)^{  -\f{|\alpha|}{2}         } \left[  (1+\tau)^{-1} \tau^{-2} +  (1+\tau)^{-1} \tau^{-\f{3}{2} }      \right]   d\tau                 \\
     \leq     & \   CA^2 t^{  -1  -\f{1+|\alpha|}{2}           }. \\
     \end{split}
\end{equation}
Analogously we prove (\ref{hd29}). For any $l\in \{ 1,2,3 \}$, it holds
\begin{equation}\nonumber
     \begin{split}
        \| \p_3(B_l u_3) \|_{L^{\infty}_{\rm h}  L^{1}_{\rm v}}  & \  \leq  \| \p_3 B_l \|_{L^{\infty}_{\rm h}  L^{1}_{\rm v}}      \| u_3 \|_{L^{\infty}}    +    \|  B_l \|_{L^{\infty}_{\rm h}  L^{1}_{\rm v}}       \| \p_3  u_3 \|   _{L^{\infty}}     \\
    & \  \leq  CA^2 \left(  (1+\tau)^{-\f{3}{2}} \tau^{-2} +  (1+\tau)^{-\f{3}{2}} \tau^{-\f{3}{2} }      \right) ;        \\
     \end{split}
\end{equation}
yielding
\begin{equation}\nonumber
     \begin{split}
     \Big\|   \nabla^{\alpha}  &\          \int_0^{ t } e^{ (t-\tau)\Delta  }\p_3 (B_l u_3 ) (\tau  )     d \tau \Big\|_{L^{\infty}_{\rm h}  L^{1}_{\rm v}}  \\
     \leq  & \    \int_0^{\f{t}{2}}  \| \p_3 \nabla^{\alpha} G (t-\tau) \|_{L^{\infty}_{\rm h}  L^{1}_{\rm v}} \|  (B_l u_3) (\tau) \|     _{L^1}         d \tau                     \\
     &\ + \int_{  \f{t}{2} } ^t  \| \nabla^{\alpha} G (t-\tau) \|_{L^{1}}  \| \p_3 (B_l u_3) (\tau) \|  _{L^{\infty}_{\rm h}  L^{1}_{\rm v}}  d \tau          \\
  \leq  & \  CA^2          \int_0^  {  \f{t}{2} } (t-\tau)^{  -1-\f{1+ |\alpha|}{2}         } (1+\tau) ^{  -2 } d \tau                 \\
         & \ +     CA^2          \int_{  \f{t}{2} } ^t (t-\tau)^{  -\f{|\alpha|}{2}         } \left[  (1+\tau)^{-\f{3}{2}} \tau^{-2} +  (1+\tau)^{-\f{3}{2}} \tau^{-\f{3}{2} }     \right]   d\tau                 \\
     \leq     & \   CA^2 t^{  -1  -\f{1+|\alpha|}{2}           }. \\
     \end{split}
\end{equation}
Following the same spirit, we show (\ref{hd30}). For any $k \in \{ 1,2 \}$, it holds
\begin{equation}\nonumber
     \begin{split}
        \| \nablah ( u_3  B_k) \|_{L^{\infty}_{\rm h}  L^{1}_{\rm v}}  & \  \leq  \| B_k \|_{L^{\infty}_{\rm h}  L^{1}_{\rm v}}      \|\nablah  u_3 \|_{L^{\infty}}    +    \| \nablah  B_k \|_{L^{\infty}_{\rm h}  L^{1}_{\rm v}}       \|   u_3 \|   _{L^{\infty}}     \\
    & \  \leq  CA^2 \left(  (1+\tau)^{-2} \tau^{-\f{3}{2}} +  (1+\tau)^{-\f{3}{2}} \tau^{-2 }      \right) ;        \\
     \end{split}
\end{equation}
yielding
\begin{equation}\nonumber
     \begin{split}
     \Big\|   \nabla^{\alpha}   &\  \int_0^{ t } e^{ (t-\tau)\Delta  }\nablah ( u_3 B_k ) (\tau  )     d \tau \Big\|_{L^{\infty}_{\rm h}  L^{1}_{\rm v}}  \\
     \leq  & \    \int_0^{\f{t}{2}}  \| \p_3 \nabla^{\alpha} G (t-\tau) \|_{L^{\infty}_{\rm h}  L^{1}_{\rm v}} \|  ( u_3 B_k ) (\tau) \|     _{L^1}         d \tau                     \\
     &\ + \int_{  \f{t}{2} } ^t  \| \nabla^{\alpha} G (t-\tau) \|_{L^{1}}  \| \nablah ( u_3 B_k) (\tau) \|  _{L^{\infty}_{\rm h}  L^{1}_{\rm v}}  d \tau          \\
  \leq  & \  CA^2          \int_0^  {  \f{t}{2} } (t-\tau)^{  -1-\f{1+ |\alpha|}{2}         } (1+\tau) ^{  -2 } d \tau                 \\
         & \ +     CA^2          \int_{  \f{t}{2} } ^t (t-\tau)^{  -\f{|\alpha|}{2}         } \left[   (1+\tau)^{-2} \tau^{-\f{3}{2}} +  (1+\tau)^{-\f{3}{2}} \tau^{-2 }     \right]   d\tau                 \\
     \leq     & \   CA^2 t^{  -1  -\f{1+|\alpha|}{2}           }. \\
     \end{split}
\end{equation}
Finally, we verify (\ref{hd31}). For any $k \in \{ 1,2 \},l\in \{ 1,2,3 \}$, it holds
\begin{equation}\nonumber
     \begin{split}
        \| \nablah ( u_k  B_l) \|_{L^{\infty}_{\rm h}  L^{1}_{\rm v}}  & \  \leq  \| \nablah B_l \|_{L^{\infty}_{\rm h}  L^{1}_{\rm v}}      \| u_k \|_{L^{\infty}}    +    \|   B_l \|_{L^{\infty}_{\rm h}  L^{1}_{\rm v}}       \|  \nablah  u_k \|   _{L^{\infty}}     \\
    & \  \leq  CA^2 \left(  (1+\tau)^{-1} \tau^{-2} +  (1+\tau)^{-\f{3}{2}} \tau^{-2 }      \right) ;        \\
     \end{split}
\end{equation}
yielding
\begin{equation}\nonumber
     \begin{split}
     \Big\|   \nabla^{\alpha}  &\       \int_0^{ t } e^{ (t-\tau)\Delta  }\nablah ( u_k B_l ) (\tau  )     d \tau \Big\|_{L^{\infty}_{\rm h}  L^{1}_{\rm v}}  \\
     \leq  & \    \int_0^{\f{t}{2}}  \| \p_3 \nabla^{\alpha} G (t-\tau) \|_{L^{\infty}_{\rm h}  L^{1}_{\rm v}} \|  ( u_k B_l ) (\tau) \|     _{L^1}         d \tau                     \\
     &\ + \int_{  \f{t}{2} } ^t  \| \nabla^{\alpha} G (t-\tau) \|_{L^{1}}  \| \nablah ( u_k B_l ) (\tau) \|  _{L^{\infty}_{\rm h}  L^{1}_{\rm v}}  d \tau          \\
  \leq  & \  CA^2          \int_0^  {  \f{t}{2} } (t-\tau)^{  -1-\f{1+ |\alpha|}{2}         } (1+\tau) ^{  -2 } d \tau                 \\
         & \ +     CA^2          \int_{  \f{t}{2} } ^t (t-\tau)^{  -\f{|\alpha|}{2}         } \left[   (1+\tau)^{-1} \tau^{-2} +  (1+\tau)^{-\f{3}{2}} \tau^{-2 }    \right]   d\tau                 \\
     \leq     & \   CA^2 t^{  -1  -\f{1+|\alpha|}{2}           }. \\
     \end{split}
\end{equation}
This completes the proof of Lemma \ref{lem8}.  $\Box$

\subsection{Asymptotic expansions for the Duhamel terms}
In this subsection, we present the asymptotic expansions for the Duhamel terms $\Nhu_1[u],\Nhu_1[B],\Nvu_1[u],\Nvu_1[B]$.
\begin{lemm}\label{lem9}
Let $(u,B)$ be subject to the ansatzes $(i),(ii)$ with $s\geq 3,A>0,T=\infty$. Then there hold the asymptotic limits for any $1\leq p \leq \infty$:
\beq\label{hd32}
\lim_{t\rightarrow \infty} t^{1-\f{1}{p}}
\left\|
\Nhu_1[u] (t,x)+ G_{\rm h}(t,\xh)  \int_0^{\infty} \int_{\R^2}  \p_3 (u_3 \uh) (\tau, y_{\rm h},x_3)
    d y_{\rm h}   d \tau
\right\|_{L^p_x} =0,
\eeq
\beq\label{hd33}
\lim_{t\rightarrow \infty} t^{1-\f{1}{p}}
\left\|
\Nhu_1[B] (t,x)+ G_{\rm h}(t,\xh)  \int_0^{\infty} \int_{\R^2}  \p_3 (B_3 \Bh) (\tau, y_{\rm h},x_3)
    d y_{\rm h}   d \tau
\right\|_{L^p_x} =0.
\eeq
Suppose further that $(u,B)$ be subject to the ansatzes $(iii),(iv)$. Then for any $1<p\leq \infty$ there exists $C>0$ depending only on $p$ such that
\begin{equation}\label{hd34}
     \begin{split}
\Big\|
\Nhu_1[u] (t,x)+
&\ G_{\rm h}(t,\xh)  \int_0^{\infty} \int_{\R^2}  \p_3 (u_3 \uh) (\tau, y_{\rm h},x_3)
    d y_{\rm h}  d \tau
\Big\|_{L^p_x} \\
    \leq & \  CA ( A+ A_{\ast} )t^{-(1-\f{1}{p})-\f{1}{2}} \log (2+t),        \\
     \end{split}
\end{equation}
and for any $2<p \leq \infty$  there exists $C>0$ depending only on $p$ such that
\begin{equation}\label{hd35}
     \begin{split}
\Big\|
\Nhu_1[B] (t,x)+
&\ G_{\rm h}(t,\xh)  \int_0^{\infty} \int_{\R^2}  \p_3 (B_3 \Bh) (\tau, y_{\rm h},x_3)
    d y_{\rm h}  d \tau
\Big\|_{L^p_x} \\
    \leq & \  CA ( A+ A_{\ast} )t^{-(1-\f{1}{p})-\f{1}{2}} (1+  t^{ -\f{1}{2} }  ).        \\
     \end{split}
\end{equation}
\end{lemm}
{\it Proof.} Noticing that (\ref{hd32}) and (\ref{hd34}) have been proved in \cite{F}, we are left to show (\ref{hd33}) and (\ref{hd35}). Inspired by \cites{FM,F}, we make the following decomposition
\[
\Nhu_1[B] (t,x)+ G_{\rm h}(t,\xh)  \int_0^{\infty} \int_{\R^2}  \p_3 (B_3 \Bh) (\tau, y_{\rm h},x_3)
    d y_{\rm h}   d \tau=-
\sum_{m=1}^4 \mathcal{I}_m (t,x),
\]
with
\begin{equation}\nonumber
     \begin{split}
        \mathcal{I}_1 (t,x) =& \  \int_0^{\f{t}{2}}   \int_{\R^2}
   \Big( \Gh(t-\tau,\xh-\yh)-\Gh(t,\xh-\yh)    \Big) \p_3(B_3 \Bh)(\tau,\yh,x_3)  d y_{\rm h}   d \tau,   \\
   \mathcal{I}_2 (t,x) = & \    \int_0^{\f{t}{2}}   \int_{\R^2}
   \Big( \Gh(t,\xh-\yh)-\Gh(t,\xh)    \Big) \p_3(B_3 \Bh)(\tau,\yh,x_3)  d y_{\rm h}   d \tau,
   \\
       \mathcal{I}_3 (t,x) = & \  \int_{\f{t}{2}}^{t} e^{(t-\tau)\Deltah} \p_3(B_3 \Bh)(\tau,x) d\tau,       \\
       \mathcal{I}_4 (t,x) = & \  -\Gh(t,\xh) \int_{\f{t}{2}}^{\infty} \int_{\R^2} \p_3(B_3 \Bh)(\tau,\yh,x_3)  d y_{\rm h}   d \tau.   \\
     \end{split}
\end{equation}
In view of Lemma \ref{lem6}\footnote{In what follows, for any $1\leq p\leq \infty$, we denote by $p'$ the H\"{o}lder conjugate exponent of $p$.},
\begin{equation}\nonumber
     \begin{split}
        \| \p_3(B_3 \Bh)(\tau)  \|_{L^{1}_{\rm h}  L^{p}_{\rm v}} \leq  & \
        \| \p_3 B_3 (\tau) \|_{L^p}   \| \Bh(\tau)   \|_{L^{p'}_{\rm h}  L^{\infty}_{\rm v}}
        +\| \p_3 \Bh (\tau) \|_{L^p}   \| B_3(\tau)   \|_{L^{p'}_{\rm h}  L^{\infty}_{\rm v}}
        \\
   \leq  & \       \| \p_3 B_3 (\tau) \|_{L^p}   \| \Bh(\tau)   \|_{L^{1}_{\rm h}  L^{\infty}_{\rm v}}^{1-\f{1}{p}}  \| \Bh(\tau)   \|_{  L^{\infty} }  ^{\f{1}{p}}                     \\
         & \ +   \| \p_3 \Bh (\tau) \|_{L^p}   \| B_3(\tau)   \|_{L^{1}_{\rm h}  L^{\infty}_{\rm v}}^{1-\f{1}{p}}  \| B_3 (\tau)   \|_{  L^{\infty} }  ^{\f{1}{p}}                     \\
        \leq  & \  CA^2 (1+\tau)^{     -3    }. \\
     \end{split}
\end{equation}
It remains to control $\mathcal{I}_m (m=1,2,3,4)$ suitably. For $\mathcal{I}_1$, we may rewrite
\[
 \mathcal{I}_1 (t,x) =  - \int_0^{\f{t}{2}}   \int_{\R^2} \int_0^1
 \tau (\p_t \Gh) (t-\theta\tau,\xh-\yh) \p_3(B_3 \Bh)(\tau,\yh,x_3)   d \theta d y_{\rm h}   d \tau;
\]
whence
\begin{equation}\nonumber
     \begin{split}
        \|  \mathcal{I}_1(t) \|_{L^p} \leq  & \ \int_0^{\f{t}{2}}
        \int_0^1  \tau \| (\p_t \Gh)(t-\theta \tau)  \|_{L^p(\R^2)}
        \|  \p_3 (B_3 \Bh)(\tau)  \|_{L^{1}_{\rm h}  L^{p}_{\rm v}} d \theta d \tau
        \\
    \leq  & \  CA^2     \int_0^{\f{t}{2}}
        \int_0^1  \tau   (t-\theta \tau)^{ -(1-\f{1}{p})    -1    } (1+\tau)^{-3} d \theta d \tau         \\
        \leq  & \ CA^2 t^{    -(1-\f{1}{p})    -1    }   \int_0^{\f{t}{2}}  (1+\tau)^{-2} d\tau           \\
        \leq  & \   CA^2 t^{    -(1-\f{1}{p})    -1    } .   \\
     \end{split}
\end{equation}
Next we consider $\mathcal{I}_2$. Through a change of variable $\xh \mapsto t^{\f{1}{2}}  \xi_{\rm h} $,
\[
        \|  \mathcal{I}_2(t) \|_{L^p}
        \leq
t^{-(1-\f{1}{p})}   \int_0^{\f{t}{2}}   \int_{\R^2}
\| \Gh(1,\cdot - t^{-\f{1}{2}} \yh ) -\Gh(1,\cdot) \| _{L^p(\R^2)}
\| \p_3(B_3 \Bh) (\tau, \yh, \cdot) \|_{L^p(\R)}   d y_{\rm h}   d \tau .
\]
Seeing that $\p_3(B_3 \Bh)  \in L^1(0,\infty;  L^{1}_{\rm h}  L^{p}_{\rm v}(\R^3)    )      $ and invoking Lebesgue's dominated convergence theorem, it follows that $\lim_{t\rightarrow \infty }t^{(1-\f{1}{p})} \|  \mathcal{I}_2(t) \|_{L^p}=0  $. Finally, we handle $\mathcal{I}_3$ and $\mathcal{I}_4$ together.
\begin{equation}\nonumber
     \begin{split}
         \|  \mathcal{I}_3(t) \|_{L^p}+ \|  \mathcal{I}_4(t) \|_{L^p} \leq
         & \       \int_{\f{t}{2}}^{t}   \|  \Gh(t-\tau) \|_{L^p(\R^2)}        \|  \p_3(B_3 \Bh)(\tau ) \|_{     L^{1}_{\rm h}  L^{p}_{\rm v}   }    d \tau              \\
   & \  +    \|  \Gh(t) \|_{L^p(\R^2)}   \int_{\f{t}{2}}^{\infty}
            \|  \p_3(B_3 \Bh)(\tau ) \|_{     L^{1}_{\rm h}  L^{p}_{\rm v}   }  d \tau                              \\
     \leq     & \ CA^2  \int_{\f{t}{2}}^{t}   (t-\tau)^{ -(1-\f{1}{p}) } (1+\tau)^{  -3     } d\tau + CA^2     t^{ -(1-\f{1}{p}) } \int_{\f{t}{2}}^{\infty}  (1+\tau)^{  -3     } d \tau      \\
        \leq  & \   CA^2   t^{    -(1-\f{1}{p})    -1    }.               \\
     \end{split}
\end{equation}
Combining the estimates above, we obtain (\ref{hd33}).

The rest of the proof is devoted to (\ref{hd35}). It is easy to see that
\begin{equation}\nonumber
     \begin{split}
        \|  |\yh| \p_3(B_3 \Bh)(\tau ) \|_{L^1}  \leq & \   \| \p_3 B_3(\tau) \| _{     L^{\infty}_{\rm h}  L^{1}_{\rm v}   }
        \| |\yh|  \Bh (\tau)  \| _{     L^{1}_{\rm h}  L^{\infty}_{\rm v}   }       \\
    & \ +   \| \p_3 \Bh (\tau) \| _{     L^{\infty}_{\rm h}  L^{1}_{\rm v}   }
        \| |\yh|  B_3 (\tau)  \| _{     L^{1}_{\rm h}  L^{\infty}_{\rm v}   }      \\
       \leq   & \  CA A_{\ast}  \tau^{-2} (1+\tau)^{     \f{1}{2}  } , \\
     \end{split}
\end{equation}
\begin{equation}\nonumber
     \begin{split}
        \|  |\yh| \p_3(B_3 \Bh)(\tau ) \|_{     L^{1}_{\rm h}  L^{\infty}_{\rm v}   }  \leq & \   \| \p_3 B_3(\tau) \| _{   L^{\infty}  }
        \| |\yh|  \Bh (\tau)  \| _{     L^{1}_{\rm h}  L^{\infty}_{\rm v}   }       \\
    & \ +   \| \p_3 \Bh (\tau) \| _{     L^{\infty}   }
        \| |\yh|  B_3 (\tau)  \| _{     L^{1}_{\rm h}  L^{\infty}_{\rm v}   }      \\
       \leq   & \  CA A_{\ast}   (1+\tau)^{-\f{5}{2} } , \\
     \end{split}
\end{equation}
which immediately yields
\begin{equation}\nonumber
     \begin{split}
         \|  |\yh| \p_3(B_3 \Bh)(\tau ) \|_{     L^{1}_{\rm h}  L^{p}_{\rm v}   }  \leq       & \      \|  |\yh| \p_3(B_3 \Bh)(\tau ) \|_{L^1} ^{ \f{1}{p} }     \|  |\yh| \p_3(B_3 \Bh)(\tau ) \|_{     L^{1}_{\rm h}  L^{\infty}_{\rm v}   }  ^{ 1-\f{1}{p} }             \\
    \leq  & \  CA A_{\ast}   \tau^{-\f{2}{p}   } (1+\tau)^{  -2+ \f{5}{2}  \f{1}{p}       } .    \\
     \end{split}
\end{equation}
As before, we reformulate $\mathcal{I}_2$ as
\[
\mathcal{I}_2(t,x)=
- \int_0^{\f{t}{2}  } \int_{\R^2} \int_0^1
(\nablah \Gh) (t,\xh-\theta \yh) \cdot \yh \p_3(B_3 \Bh)(\tau,\yh,x_3)
d \theta d \yh d \tau;
\]
whence
\begin{equation}\nonumber
     \begin{split}
         \|  \mathcal{I}_2(t) \|_{L^p} \leq  & \
         \int_0^{\f{t}{2}  } \int_{\R^2} \int_0^1
         \|   (\nablah \Gh) (t,\cdot-\theta \yh)      \|_{L^p(\R^2)}
         \|    |\yh| \p_3(B_3 \Bh)(\tau,\yh,\cdot)   \|_{L^p(\R)}d \theta d \yh d \tau
         \\
  \leq  & \  C t^{  -(1-\f{1}{p})-\f{1}{2}     }  \int_0^{\f{t}{2}  }
   \|    |\yh| \p_3(B_3 \Bh)(\tau)   \|_{     L^{1}_{\rm h}  L^{p}_{\rm v}   }  d \tau
  \\
        \leq  & \  CA A_{\ast}t^{  -(1-\f{1}{p})-\f{1}{2}     }
        \int_0^{\f{t}{2}  }   \tau^{-\f{2}{p}   } (1+\tau)^{   -2+ \f{5}{2}  \f{1}{p}      }    d \tau
        \\
     \leq     & \   CA A_{\ast}t^{  -(1-\f{1}{p})-\f{1}{2}     }.    \\
     \end{split}
\end{equation}
It should be emphasized that the hypotheses $2<p\leq \infty$ was essentially used in the last step. Revisiting $\mathcal{I}_m$ with $m=1,3,4$ in the first part, we arrive at (\ref{hd35}). The proof of Lemma \ref{lem9} is completely finished.         $\Box$

Employing similar arguments as in the previous lemma, we next consider the asymptotic expansions for $\Nvu_1[u]$ and $\Nvu_1[B]$. Notice that the asymptotic expansion of $\Nvu_1[u]$ was used in \cite{F} without a proof, here we would like to provide a detailed proof for completeness.
\begin{lemm}\label{lem10}
Let $(u,B)$ be subject to the ansatzes $(i),(ii)$ with $s\geq 3,A>0,T=\infty$. Then there hold the asymptotic limits for any $1\leq p \leq \infty$:
\beq\label{hd36}
\lim_{t\rightarrow \infty} t^{(1-\f{1}{p})+\f{1}{2}}
\left\|
\Nvu_1[u] (t,x)-\nablah G_{\rm h}(t,\xh) \cdot  \int_0^{\infty} \int_{\R^2}   (u_3 \uh) (\tau, y_{\rm h},x_3)
    d y_{\rm h}   d \tau
\right\|_{L^p_x} =0,
\eeq
\beq\label{hd37}
\lim_{t\rightarrow \infty} t^{(1-\f{1}{p}) + \f{1}{2}}
\left\|
\Nvu_1[B] (t,x)-  \nablah G_{\rm h}(t,\xh) \cdot \int_0^{\infty} \int_{\R^2}   (B_3 \Bh) (\tau, y_{\rm h},x_3)
    d y_{\rm h}   d \tau
\right\|_{L^p_x} =0.
\eeq

\end{lemm}
{\it Proof.} We begin with the proof of (\ref{hd36}). In analogy with \cites{FM,F}, we make the following decomposition
\[
\Nvu_1[u] (t,x)-\nablah G_{\rm h}(t,\xh)  \int_0^{\infty} \int_{\R^2}   (u_3 \uh) (\tau, y_{\rm h},x_3)
    d y_{\rm h}   d \tau=
\sum_{m=1}^4 \mathcal{I}_m (t,x),
\]
with
\begin{equation}\nonumber
     \begin{split}
        \mathcal{I}_1 (t,x) =& \  \int_0^{\f{t}{2}}   \int_{\R^2}  \nablah
   \Big( \Gh(t-\tau,\xh-\yh)-\Gh(t,\xh-\yh)    \Big) (u_3 \uh)(\tau,\yh,x_3)  d y_{\rm h}   d \tau,   \\
   \mathcal{I}_2 (t,x) = & \    \int_0^{\f{t}{2}}   \int_{\R^2} \nablah
   \Big( \Gh(t,\xh-\yh)-\Gh(t,\xh)    \Big) (u_3 \uh)(\tau,\yh,x_3)  d y_{\rm h}   d \tau,
   \\
       \mathcal{I}_3 (t,x) = & \  \int_{\f{t}{2}}^{t} e^{(t-\tau)\Deltah} \nablah (u_3 \uh)(\tau,x) d\tau,       \\
       \mathcal{I}_4 (t,x) = & \  -\nablah \Gh(t,\xh) \int_{\f{t}{2}}^{\infty} \int_{\R^2} (u_3 \uh)(\tau,\yh,x_3)  d y_{\rm h}   d \tau.   \\
     \end{split}
\end{equation}
Observe that
\begin{equation}\nonumber
     \begin{split}
        \| (u_3 \uh)(\tau)  \|_{L^{1}_{\rm h}  L^{p}_{\rm v}} \leq  & \
\|  \uh (\tau) \|_{L^p}   \| u_3(\tau)   \|_{L^{p'}_{\rm h}  L^{\infty}_{\rm v}}
        \\
   \leq  & \    \|  \uh (\tau) \|_{L^p}   \| u_3(\tau)   \|_{L^{1}_{\rm h}  L^{\infty}_{\rm v}}^{1-\f{1}{p}}  \| u_3 (\tau)   \|_{  L^{\infty} }  ^{\f{1}{p}}                     \\
        \leq  & \  CA^2 (1+\tau)^{     -\f{3}{2}   }. \\
     \end{split}
\end{equation}
Rewriting $\mathcal{I}_1$ via
\[
 \mathcal{I}_1 (t,x) =  - \int_0^{\f{t}{2}}   \int_{\R^2} \int_0^1
 \tau (\p_t \nablah \Gh) (t-\theta\tau,\xh-\yh)   (u_3 \uh)(\tau,\yh,x_3)   d \theta d y_{\rm h}   d \tau;
\]
whence
\begin{equation}\nonumber
     \begin{split}
        \|  \mathcal{I}_1(t) \|_{L^p} \leq  & \ \int_0^{\f{t}{2}}
        \int_0^1  \tau \| (\p_t \nablah \Gh)(t-\theta \tau)  \|_{L^p(\R^2)}
        \|   (u_3 \uh)(\tau)  \|_{L^{1}_{\rm h}  L^{p}_{\rm v}} d \theta d \tau
        \\
    \leq  & \  CA^2     \int_0^{\f{t}{2}}
        \int_0^1  \tau   (t-\theta \tau)^{ -(1-\f{1}{p})    -\f{3}{2}    } (1+\tau)^{-\f{3}{2}} d \theta d \tau         \\
        \leq  & \ CA^2 t^{    -(1-\f{1}{p})    -\f{3}{2}     }   \int_0^{\f{t}{2}}  (1+\tau)^{-\f{1}{2}} d\tau           \\
        \leq  & \   CA^2 t^{    -(1-\f{1}{p})    -1    } .   \\
     \end{split}
\end{equation}
Performing a change of variable $\xh \mapsto t^{\f{1}{2}}  \xi_{\rm h} $,
\[
        \|  \mathcal{I}_2(t) \|_{L^p}
        \leq
t^{-(1-\f{1}{p})-\f{1}{2}}   \int_0^{\f{t}{2}}   \int_{\R^2}
\|\nablah \Gh(1,\cdot - t^{-\f{1}{2}} \yh ) -\nablah \Gh(1,\cdot) \| _{L^p(\R^2)}
\| (u_3 \uh) (\tau, \yh, \cdot) \|_{L^p(\R)}   d y_{\rm h}   d \tau .
\]
Seeing that $(u_3 \uh)  \in L^1(0,\infty;  L^{1}_{\rm h}  L^{p}_{\rm v}(\R^3)    )      $ and invoking Lebesgue's dominated convergence theorem, it follows that $\lim_{t\rightarrow \infty }t^{(1-\f{1}{p}) +\f{1}{2} } \|  \mathcal{I}_2(t) \|_{L^p}=0  $. To proceed, we handle $\mathcal{I}_3$ and $\mathcal{I}_4$ together.
\begin{equation}\nonumber
     \begin{split}
         \|  \mathcal{I}_3(t) \|_{L^p}+ \|  \mathcal{I}_4(t) \|_{L^p} \leq
         & \       \int_{\f{t}{2}}^{t}   \| \nablah \Gh(t-\tau) \|_{L^p(\R^2)}        \|  (u_3 \uh)(\tau ) \|_{     L^{1}_{\rm h}  L^{p}_{\rm v}   }    d \tau              \\
   & \  +    \|  \nablah \Gh(t) \|_{L^p(\R^2)}   \int_{\f{t}{2}}^{\infty}
            \|  (u_3 \uh)(\tau ) \|_{     L^{1}_{\rm h}  L^{p}_{\rm v}   }  d \tau                              \\
     \leq     & \ CA^2  \int_{\f{t}{2}}^{t}   (t-\tau)^{ -(1-\f{1}{p})-\f{1}{2} } (1+\tau)^{  -\f{3}{2}     } d\tau      \\
     &\ + CA^2     t^{ -(1-\f{1}{p}) -\f{1}{2}} \int_{\f{t}{2}}^{\infty}  \tau^{  -\f{3}{2}     } d \tau  \\
        \leq  & \   CA^2   t^{    -(1-\f{1}{p})    -1    }.               \\
     \end{split}
\end{equation}
The obtained estimates above allow us to verify (\ref{hd36}) readily. Following exactly the same way, we prove (\ref{hd37}). The details are thus omitted. The proof of Lemma \ref{lem10} is completely finished.                           $\Box$

To handle the asymptotic expansion of magnetic field, we have
\begin{lemm}\label{lem11}
Let $(u,B)$ be subject to the ansatzes $(i),(ii)$ with $s\geq 3,A>0,T=\infty$. Then there hold the asymptotic limits for any $1\leq p \leq \infty$ and any $i\in \{ 1,2,3 \}$:
\beq\label{hd38}
\lim_{t\rightarrow \infty} t^{\f{3}{2}(1-\f{1}{p})+\f{1}{2}}
\left\|
\sum_{j=1}^3 \int_0^t e^{(t-\tau)\Delta } \p_j (B_i u_j)(\tau) d \tau
-  \sum_{j=1}^3 \p_j G(t,x)  \int_0^{\infty} \int_{\R^3}   (B_i u_j) (\tau, y) d y   d \tau
\right\|_{L^p_x} =0,
\eeq
\beq\label{hd39}
\lim_{t\rightarrow \infty} t^{\f{3}{2}(1-\f{1}{p})+\f{1}{2}}
\left\|
\sum_{j=1}^3 \int_0^t e^{(t-\tau)\Delta } \p_j (u_i B_j )(\tau) d \tau
-  \sum_{j=1}^3 \p_j G(t,x)  \int_0^{\infty} \int_{\R^3}   (u_i B_j) (\tau, y) d y   d \tau
\right\|_{L^p_x} =0.
\eeq
\end{lemm}
{\it Proof.} We shall only present the proof of (\ref{hd38}), since the other is carried out in the same way. Following \cites{FM,F}, it holds the decomposition
\[
\sum_{j=1}^3 \int_0^t e^{(t-\tau)\Delta } \p_j (B_i u_j)(\tau) d \tau
- \sum_{j=1}^3 \p_j G(t,x)  \int_0^{\infty} \int_{\R^3}   (B_i u_j) (\tau, y) d y   d \tau
=\sum_{m=1}^4 \mathcal{I}_m (t,x),
\]
with
\begin{equation}\nonumber
     \begin{split}
        \mathcal{I}_1 (t,x) =& \  \sum_{j=1}^3 \int_0^{\f{t}{2} } \int_{\R^3} \p_j \Big(
        G(t-\tau,x-y)-G(t,x-y)
        \Big)  (B_i u_j)(\tau, y) d y d\tau  ,   \\
   \mathcal{I}_2 (t,x) = & \ \sum_{j=1}^3    \int_0^{\f{t}{2}}   \int_{\R^3}
   \p_j \Big(  G(t,x-y)-G(t,x)    \Big) (B_i u_j)(\tau, y) d y d\tau   ,
   \\
       \mathcal{I}_3 (t,x) = & \  \sum_{j=1}^3   \int_{\f{t}{2}}^{t} e^{(t-\tau)\Delta}  \p_3(B_i u_j)(\tau,x) d\tau,       \\
       \mathcal{I}_4 (t,x) = & \  -\sum_{j=1}^3 \p_j G(t,x) \int_{\f{t}{2}}^{\infty} \int_{\R^3} (B_i u_j)(\tau, y) d y d\tau  .   \\
     \end{split}
\end{equation}
It is easily noticed that
\[
\| (B_i u_j)(\tau)  \|_{L^1} \leq \| B_i (\tau) \|_{L^{\infty}}
\| u_j(\tau) \|_{L^1}
\leq  CA^2 (1+\tau)^{-2}
\]
and
\[
\mathcal{I}_1 (t,x) =
- \sum_{j=1}^3\int_0^{\f{t}{2} } \int_{\R^3}
\tau (\p_t \p_j G) (t-\theta \tau,x-y)  (B_i u_j)(\tau, y) d \theta  d y d\tau ;
\]
whence
\begin{equation}\nonumber
     \begin{split}
        \|  \mathcal{I}_1(t) \|_{L^p} \leq  & \  \sum_{j=1}^3 \int_0^{\f{t}{2}}
        \int_0^1  \tau \| (\p_t \p_j  G)(t-\theta \tau)  \|_{L^p }
        \|   (B_i u_j )(\tau)  \|_{  L^{1}   } d \theta d \tau
        \\
    \leq  & \  CA^2     \int_0^{\f{t}{2}}
        \int_0^1  \tau   (t-\theta \tau)^{ - \f{3}{2} (1-\f{1}{p})    -\f{3}{2}    } (1+\tau)^{-2  } d \theta d \tau         \\
        \leq  & \ CA^2 t^{    -(1-\f{1}{p})    -\f{3}{2}     }   \int_0^{\f{t}{2}}  (1+\tau)^{-1 } d\tau           \\
        \leq  & \   CA^2 t^{    -\f{3}{2} (1-\f{1}{p})    -\f{3}{2}   } \log(2+t) .   \\
     \end{split}
\end{equation}
Using a change of variable $ x \mapsto t^{\f{1}{2}}  \xi $, it follows that
\[
        \|  \mathcal{I}_2(t) \|_{L^p}
        \leq
t^{-\f{3}{2}  (1-\f{1}{p})-\f{1}{2}}   \sum_{j=1}^3  \int_0^{\f{t}{2}}   \int_{\R^3}
\|\p_j G(1,\cdot - t^{-\f{1}{2}} y ) -\p_j G(1,\cdot) \| _{L^p}
| (B_i u_j) (\tau, y) |       d y  d \tau  .
\]
Observing that $(B_i u_j)  \in L^1(0,\infty; L^1 (\R^3)    )      $ and invoking Lebesgue's dominated convergence theorem, one infers that $\lim_{t\rightarrow \infty }t^{\f{3}{2} (1-\f{1}{p}) +\f{1}{2} } \|  \mathcal{I}_2(t) \|_{L^p}=0  $. Finally,
\begin{equation}\nonumber
     \begin{split}
         \|  \mathcal{I}_3(t) \|_{L^p}+ \|  \mathcal{I}_4(t) \|_{L^p} \leq
         & \      \sum_{j=1}^3  \int_{\f{t}{2}}^{t}   \| \p_j   G(t-\tau) \|_{L^p }        \|  (B_i u_j )(\tau ) \|_{     L^{1}  }    d \tau              \\
   & \  +  \sum_{j=1}^3    \|  \p_j  G(t) \|_{L^p}   \int_{\f{t}{2}}^{\infty}
            \|  (B_i u_j)(\tau ) \|_{     L^{1}    }  d \tau                              \\
     \leq     & \ CA^2  \int_{\f{t}{2}}^{t}   (t-\tau)^{ - \f{3}{2} (1-\f{1}{p})-\f{1}{2} } (1+\tau)^{  -2     } d\tau      \\
     &\ + CA^2     t^{ -\f{3}{2} (1-\f{1}{p}) -\f{1}{2}} \int_{\f{t}{2}}^{\infty}  \tau^{  -2     } d \tau  \\
        \leq  & \   CA^2   t^{    -\f{3}{2} (1-\f{1}{p})    -\f{3}{2}   }.               \\
     \end{split}
\end{equation}
We thus obtain (\ref{hd38}) from the above estimates.      $\Box$

\subsection{Proof of Theorem \ref{thm1}}
This subsection is dedicated to the proof of Theorem \ref{thm1}; it consists of global existence of small solutions in $X^s(\R^3)$, decay estimates of solutions in $L^p$ spaces and the asymptotic profile of solutions.

We begin with the global small solutions in $X^s(\R^3)$. To this end, assume that the integer $s\geq 5$ and $u_0,B_0\in X^s(\R^3)$ obeying $\Grad \cdot u_0=\Grad \cdot B_0=0, |x|B_0(x)\in L^1(\R^3)$, $\| u_0 \|_{H^s}+\| B_0 \|_{H^s} \leq \epsilon$ for $\epsilon>0$ sufficiently small. Let $(u,B)\in C([0,\infty);H^s(\R^3))$ be the global solution ensured by Proposition \ref{pr1}. We conclude from the integral formulations of $(u,B)$, see (\ref{hd4})-(\ref{hd7}), that for any $t>0$
\begin{equation}\label{hd100}
     \begin{split}
         \| (u,B)(t)  &\ \|  _{  L^1_{\xh}(W^{1,1}\cap W^{1,\infty})_{x_3}         }  \\
         \leq  & \
          \| e^{t \Deltah} u_0 \| _{  L^1_{\xh}(W^{1,1}\cap W^{1,\infty})_{x_3}         }
          + \| e^{t \Delta} B_0 \|_{  L^1_{\xh}(W^{1,1}\cap W^{1,\infty})_{x_3}         } \\
     & \  +   C\sum_{m=1}^5  \|   \Nhu_m[u](t)    \| _{ L^1_{\xh} W^{2,1}_{x_3}     }
     + C\sum_{m=1}^3  \|   \Nvu_m[u](t)    \| _{ L^1_{\xh} W^{2,1}_{x_3}     }   \\
         & \   +   C\sum_{m=1}^5  \|   \Nhu_m[B](t)    \| _{ L^1_{\xh} W^{2,1}_{x_3}     }
     + C\sum_{m=1}^3  \|   \Nvu_m[B](t)    \| _{ L^1_{\xh} W^{2,1}_{x_3}     }   \\
         & \  + C\sum_{m=1}^3  \|   \Nhb_m[u,B](t) \|  _{ L^1_{\xh} W^{2,1}_{x_3}     } +
          C\sum_{m=1}^3  \|   \Nvb_m[u,B](t) \|  _{ L^1_{\xh} W^{2,1}_{x_3}     } \\
        \leq   &\ C \| (u_0,B_0) \|_{X^s} +C \int_0^t \|  \p_3(u_3 \uh)(\tau) \|_{ L^1_{\xh} W^{2,1}_{x_3}     }  d\tau  \\
        &\ + C \sum_{k,l=1}^3 \int_0^t (t-\tau)^{ -\f{1}{2}  } \|  (u_k u_l)(\tau) \|_{ L^1_{\xh} W^{2,1}_{x_3}     }  d\tau  \\
        &\ +C   \sum_{k,l=1}^3 \int_0^t (t-\tau)^{ -\f{1}{2}  } \|  (B_k B_l)(\tau) \|_{ L^1_{\xh} W^{2,1}_{x_3}     }  d\tau  \\
         &\ +C   \sum_{k,l=1}^3 \int_0^t (t-\tau)^{ -\f{1}{2}  } \|  (B_k u_l)(\tau) \|_{ L^1_{\xh} W^{2,1}_{x_3}     }  d\tau  \\
         \leq &\   C \| (u_0,B_0) \|_{X^s} + C \int_0^t \| u(\tau) \|^2_{H^3} d\tau
          + C\int_0^t   (t-\tau)^{ -\f{1}{2}  } \| u(\tau) \|^2_{H^2} d\tau \\
          &\ +  C\int_0^t   (t-\tau)^{ -\f{1}{2}  } \| B(\tau) \|^2_{H^2} d\tau
            +  C\int_0^t   (t-\tau)^{ -\f{1}{2}  } \| B(\tau) \|_{H^2} \| u(\tau) \|_{H^2}  d\tau  \\
            \leq &\ C \| (u_0,B_0) \|_{X^s} + C t \| (u_0,B_0) \|^2_{H^3} +C t^{\f{1}{2}} \| (u_0,B_0) \|^2_{H^2} \\
            \leq &\ C(1+t)  \| (u_0,B_0) \|_{X^s},
      \end{split}
\end{equation}
where Lemma \ref{lm1}, Corollary \ref{co1} and the smallness of initial data were invoked. Observe also that the mappings $t \mapsto  \| u(t)   \|  _{  L^1_{\xh}(W^{1,1}\cap W^{1,\infty})_{x_3}         }$, $t \mapsto  \| B(t)   \|  _{  L^1_{\xh}(W^{1,1}\cap W^{1,\infty})_{x_3}         }$ are continuous with respect to time. It follows from (\ref{hd100}) that $u,B\in C([0,\infty); X^s(\R^3))$. This proves the global small smooth solution to the system (\ref{MHD-2}).

We next come to the decay estimates in $L^p$-norms, namely (\ref{ma-1})-(\ref{ma-3}), via the bootstrapping argument. Based on the discussions in subsection \ref{se-l} concerning the linearized equations, we know that there exists a generic positive constant $C_1$ such that
\begin{equation}\label{hd101}
     \begin{split}
      \| \nabla ^{\alpha}   e^{t \Deltah} u_{0,\rm h} (t) \|_{L^p} & \ \leq C_1 t^{ -(1-\f{1}{p}) -\f{|\alphah|}{2}        }  \| (u_0,B_0) \|_{  \overline{ X^{s} } }, \\
    \| \nablah ^{\alphah}   e^{t \Deltah} u_{0,3} (t) \|_{L^p} & \   \leq C_1 t^{ -\f{3}{2}(1-\f{1}{p}) -\f{|\alphah|}{2}        }  \|(u_0,B_0)\|_{  \overline{ X^{s} } }, \\
      \| \nabla ^{\alpha}   e^{t\Delta} B_0 (t) \|_{L^p}    & \  \leq C_1 t^{ -\f{3}{2}(1-\f{1}{p}) -\f{1}{2}-\f{|\alpha|}{2}        }  \| (u_0,B_0) \|_{  \overline{ X^{s} } }, \\
     \end{split}
\end{equation}
for any $t>0$, $1\leq p \leq \infty$ and any $\alpha=(\alphah,\alpha_3)\in (\mathbb{N} \cup \{ 0\})^2 \times (\mathbb{N} \cup \{ 0\})$ with $|\alpha|\leq 1$. Suppose that
\begin{equation}\label{hd102}
     \begin{split}
      \| \nabla ^{\alpha}   \uh (t) \|_{L^p} & \ \leq 2 C_1 t^{ -(1-\f{1}{p}) -\f{|\alphah|}{2}        }  \| (u_0,B_0) \|_{  \overline{ X^{s} } }, \\
    \| \nablah ^{\alphah}   u_3 (t) \|_{L^p} & \   \leq 2 C_1 t^{ -\f{3}{2}(1-\f{1}{p}) -\f{|\alphah|}{2}        }  \|(u_0,B_0)\|_{  \overline{ X^{s} } }, \\
      \| \nabla ^{\alpha}  B (t) \|_{L^p}    & \  \leq 2 C_1 t^{ -\f{3}{2}(1-\f{1}{p}) -\f{1}{2}-\f{|\alpha|}{2}        }  \| (u_0,B_0) \|_{  \overline{ X^{s} } }, \\
     \end{split}
\end{equation}
holds for any $0<t <T$. Upon choosing $A:=2C_1 \| (u_0,B_0) \|_{  \overline{ X^{s} } }$, it follows from Lemma \ref{lem7} that
\begin{equation}\label{hd103}
     \begin{split}
      \| \nabla ^{\alpha}   \uh (t) \|_{L^p} \leq & \   C_1 t^{ -(1-\f{1}{p}) -\f{|\alphah|}{2}        }  \| (u_0,B_0) \|_{  \overline{ X^{s} } } \\
      &\ +C_2 (1+t)^{ -(1-\f{1}{p}) -\f{|\alphah|}{2}        }  \| (u_0,B_0) \|_{  \overline{ X^{s} } }^2, \\
    \| \nablah ^{\alphah}   u_3 (t) \|_{L^p}  \leq & \    C_1 t^{ -\f{3}{2}(1-\f{1}{p}) -\f{|\alphah|}{2}        }  \|(u_0,B_0)\|_{  \overline{ X^{s} } } \\
   &\  + C_2 (1+t)^{ -\f{3}{2}(1-\f{1}{p}) -\f{|\alphah|}{2}        }  \|(u_0,B_0)\|_{  \overline{ X^{s} } }^2,  \\
      \| \nabla ^{\alpha}  B (t) \|_{L^p}   \leq  & \   C_1 t^{ -\f{3}{2}(1-\f{1}{p}) -\f{1}{2}-\f{|\alpha|}{2}        }  \| (u_0,B_0) \|_{  \overline{ X^{s} } }\\
      &\ + C_2 (1+t)^{ -\f{3}{2}(1-\f{1}{p}) -\f{1}{2}-\f{|\alpha|}{2}        }  \| (u_0,B_0) \|_{  \overline{ X^{s} } }^2,
      \\
     \end{split}
\end{equation}
holds for any $0<t <T$ and some $C_2>0$. By setting $\epsilon_1:= \min\{\epsilon, \f{C_1}{2C_2}  \}$ and choosing the initial data sufficiently small such that $\| (u_0,B_0) \|_{  \overline{ X^{s} } } \leq \epsilon_1$, we further deduce from (\ref{hd103}) that
\begin{equation}\label{hd104}
     \begin{split}
      \| \nabla ^{\alpha}   \uh (t) \|_{L^p} & \ \leq \f{3}{2} C_1 t^{ -(1-\f{1}{p}) -\f{|\alphah|}{2}        }  \| (u_0,B_0) \|_{  \overline{ X^{s} } }, \\
    \| \nablah ^{\alphah}   u_3 (t) \|_{L^p} & \   \leq \f{3}{2} C_1 t^{ -\f{3}{2}(1-\f{1}{p}) -\f{|\alphah|}{2}        }  \|(u_0,B_0)\|_{  \overline{ X^{s} } }, \\
      \| \nabla ^{\alpha}  B (t) \|_{L^p}    & \  \leq \f{3}{2} C_1 t^{ -\f{3}{2}(1-\f{1}{p}) -\f{1}{2}-\f{|\alpha|}{2}        }  \| (u_0,B_0) \|_{  \overline{ X^{s} } }, \\
     \end{split}
\end{equation}
for any $0<t <T$; this gives a sharper bound than our initial hypothesis (\ref{hd102}). The bootstrapping argument then assesses that (\ref{hd104}) actually holds for any $0<t<\infty$ and any $1\leq p \leq \infty$. The proof of nonlinear $L^p$ decay estimates is completely finished.

We end the proof by showing the asymptotic profile of solutions. By Lemma \ref{lm1}, Lemma \ref{lem7} and Lemma \ref{lem9}, we obtain for $\uh$ that
\begin{equation}\nonumber
     \begin{split}
     t^{1-\f{1}{p}}  \Big\|    \uh(t,x) &\  -\Gh (t,\xh) \int_{\R^2} u_{0,\rm h}(\yh,x_3) d \yh
       \\
       &\ +  \Gh (t,\xh) \int_0^{\infty} \int_{\R^2} \p_3(u_3 \uh)  (\tau,\yh,x_3) d \yh  d \tau  \\
     &\ -  \Gh (t,\xh) \int_0^{\infty} \int_{\R^2} \p_3(B_3 \Bh)  (\tau,\yh,x_3) d \yh  d \tau \Big \|_{L^p_x}       \\
    \leq & \  t^{1-\f{1}{p}}
  \left\|   e^{t\Deltah} u_{0,\rm h}   -\Gh (t,\xh) \int_{\R^2} u_{0,\rm h}(\yh,x_3) d \yh  \right\|_ {L^p_x}         \\
  &\ + t^{1-\f{1}{p}}
\left\|
\Nhu_1[u] (t,x)+ G_{\rm h}(t,\xh)  \int_0^{\infty} \int_{\R^2}  \p_3 (u_3 \uh) (\tau, y_{\rm h},x_3)
    d y_{\rm h} d \tau
\right\|_{L^p_x}          \\
&\ + t^{1-\f{1}{p}}
\left\|
\Nhu_1[B] (t,x)+ G_{\rm h}(t,\xh)  \int_0^{\infty} \int_{\R^2}  \p_3 (B_3 \Bh) (\tau, y_{\rm h},x_3)
    d y_{\rm h} d \tau
\right\|_{L^p_x}           \\
&\ +\sum_{m=2}^5  t^{1-\f{1}{p}}
\Big( \| \Nhu_m[u] (t)\| _{L^p}     +\| \Nhu_m[B] (t)\| _{L^p}    \Big)        \\
 \leq & \  t^{1-\f{1}{p}}
  \left\|   e^{t\Deltah} u_{0,\rm h}   -\Gh (t,\xh) \int_{\R^2} u_{0,\rm h}(\yh,x_3) d \yh  \right\|_ {L^p_x}         \\
  &\ + t^{1-\f{1}{p}}
\left\|
\Nhu_1[u] (t,x)+ G_{\rm h}(t,\xh)  \int_0^{\infty} \int_{\R^2}  \p_3 (u_3 \uh) (\tau, y_{\rm h},x_3)
    d y_{\rm h} d \tau
\right\|_{L^p_x}          \\
&\ + t^{1-\f{1}{p}}
\left\|
\Nhu_1[B] (t,x)+ G_{\rm h}(t,\xh)  \int_0^{\infty} \int_{\R^2}  \p_3 (B_3 \Bh) (\tau, y_{\rm h},x_3)
    d y_{\rm h} d \tau
\right\|_{L^p_x}           \\
&\ +C \| (u_0,B_0) \|_{   \overline{X^s } }  \Big(
t^{-\f{1}{2} }\log (2+t)+   t^{-\f{1}{2} }
\Big),  \\
     \end{split}
\end{equation}
the right-hand side of which tending to zero as time goes to infinity. This proves (\ref{ma-4}). To proceed, by Lemma \ref{lm2} and Lemma \ref{lem7} we deduce that
\begin{equation}\nonumber
     \begin{split}
     t^{ \f{3}{2} (1-\f{1}{p})  } &\  \Big\|    u_3 (t,x)  -\Gh (t,\xh) \int_{\R^2} u_{0, 3}(\yh,x_3) d \yh
     \Big \|_{L^p_x}
       \\
      \leq  &\
       t^{ \f{3}{2} (1-\f{1}{p})  }  \Big\|    e^{t \Deltah} u_{0,3}(x)   -\Gh (t,\xh) \int_{\R^2} u_{0, 3}(\yh,x_3) d \yh
     \Big \|_{L^p_x}
      \\
      &\ + \sum_{m=1}^3 t^{ \f{3}{2} (1-\f{1}{p})  }  \Big(
      \| \Nvu_m[u] (t)\| _{L^p}     +\| \Nvu_m[B](t) \| _{L^p}
      \Big)
      \\
      \leq &\   C \mathcal{R}_{p,0}(t) ^{\f{1}{p}} \|u_0 \|_{L^1}^{1-\f{1}{p}  }
      +C\| (u_0,B_0) \|_{   \overline{X^s } }  t^{-\f{1}{2p}}     \\
      &\ + C\| (u_0,B_0) \|_{   \overline{X^s } }   t^{-\f{1}{2p} +\f{3}{8} (1-\f{1}{p}) } \log(2+t)    +  C\| (u_0,B_0) \|_{   \overline{X^s } }
       t^{-\f{1}{2}} .  \\
    \end{split}
\end{equation}
Obviously, the right-hand side of the above inequality tends to zero as time approaches to infinity provided that $1\leq p<\infty$. This proves (\ref{ma-5}). Consequently, the proof of Theorem \ref{thm1} is completely finished.

\section{Large time behavior of the MHD system with full dissipation and  horizontal magnetic diffusion}\label{hm}

\subsection{Reformulation of the equations}
In this section we consider the incompressible MHD system with full dissipation and horizontal magnetic diffusion:
\begin{align}\label{hm1}
    \begin{cases}
    \partial_t u -  \Delta u + ( u \cdot \nabla ) u - ( B \cdot \nabla ) B + \nabla p = 0, & t>0 , x \in \mathbb{R}^3,\\
    \partial_t B -  \Deltah B + ( u \cdot \nabla ) B - ( B \cdot \nabla ) u = 0, & t>0 , x \in \mathbb{R}^3,\\
    \nabla \cdot u = \nabla \cdot B = 0 , & t \geqslant 0, x \in \mathbb{R}^3,\\
    u(0,x) = u_0(x), B(0,x) = B_0(x),
    & x \in \mathbb{R}^3.
    \end{cases}
\end{align}
The corresponding integral equation of \eqref{hm1} is then given by
\begin{align}
    &u(t) = e^{ t \Delta} u_0
    - \int_0^t e^{ (t - \tau) \Delta} \mathbb{P} \nabla \cdot ( u \otimes  u - B \otimes B ) (\tau) d\tau,\label{hm2}\\
    &B(t) = e^{ t \Deltah } B_0
    - \int_0^t e^{ (t - \tau) \Deltah } \nabla \cdot ( B \otimes  u - u \otimes B ) (\tau) d\tau,\label{hm3}
\end{align}
where $\mathbb{P}$ is the Helmholtz projection on $\mathbb{R}^3$.
Furthermore, we have the following decomposition of the nonlinear terms of (\ref{hm2})-(\ref{hm3}):
\begin{prop}\label{prop2}
    Let $(u,B)$ be the solution to the integral equation (\ref{hm2})-(\ref{hm3}).
    Then, the integral equation for $u$ is decomposed as
\beq\label{hm4}
        \uh(t) = e^{ t \Delta} u_{0,{\rm h}} + \sum_{m=1}^5 \Big( \Nhe_m[u](t) -  \Nhe_m[B](t) \Big)  ,
\eeq
\beq\label{hm5}
        u_3(t) = e^{ t \Delta} u_{0,3} + \sum_{m=1}^5 \Big( \Nve_m[u](t) -  \Nve_m[B](t) \Big),
\eeq
where
	\begin{align*}
		\Nhe_1[v](t)&:=-\int_0^t    e^{ (t-\tau)\Delta}\partial_{3}(v_3\vh)(\tau)d\tau,\\
		\Nhe_2[v](t)&:=-\int_0^t    e^{ (t-\tau)\Delta}\nablah \cdot (\vh\otimes\vh)(\tau)d\tau,\\
		\Nhe_3[v](t)&:=-\sum_{k,l=1}^2  \int_0^t    \nablah\partial_{k}\partial_{l} N (t-\tau)*(v_k v_l)(\tau)d\tau,\\
			\Nhe_4[v](t)&:=-2 \sum_{k=1}^2  \int_0^t    \nablah\partial_{k}\partial_{3} N (t-\tau)*(v_k v_3)(\tau)d\tau,\\
		\Nhe_5[v](t)&:=- \int_0^t\nablah \p_3 \p_3 N (t-\tau)*(v_3(\tau)^2 ) d\tau , \\
	\end{align*}
	and
	\begin{align*}
		\Nve_1[v](t)&:=  -  \int_0^t  e^{(t-\tau)\Delta}\nablah\cdot(v_3\vh)(\tau)d\tau,\\
			\Nve_2[v](t)&:=  -  \int_0^t  e^{(t-\tau)\Delta}  \p_3 (v_3(\tau)^2 )d\tau,\\
		\Nve_3[v](t)&:=- \sum_{k,l=1}^2\int_0^t
	      \p_3 \p_k \p_l         N(t-\tau)*(v_k v_l)(\tau)d\tau,\\
		\Nve_4[v](t)&:=  -2\sum_{k=1}^2   \int_0^t   \p_3 \p_3 \p_k  N(t-\tau)*(v_k v_3)(\tau)d\tau ,\\
	\Nve_5[v](t)&:=  -   \int_0^t   \p_3 \p_3 \p_3  N(t-\tau)* (v_3(\tau)^2 ) d\tau ,\\
	\end{align*}
	for $v = u, B$.
	Here, $N(t,x)$ is the function with the following representation:
\[
		N(t,x)=\int_0^{\infty}     (4 \pi  (t+s) ) ^{-\f{3}{2}}    e^{-\frac{|x|^2}{4 (t+s)}}    ds
		=\int_0^{\infty}  G( t+s,x)    ds.
\]
    For the integral equation of $B$, it holds
\beq\label{hm6}
        \Bh(t) = e^{ t \Deltah} B_{0,{\rm h}} + \sum_{m=1}^4 \Nhd_m[u,B](t),
\eeq
\beq\label{hm7}
        B_3(t) = e^{ t \Deltah} B_{0,3} + \sum_{m=1}^2 \Nvd_m[u,B](t),
\eeq
where
    \begin{equation}\nonumber
     \begin{split}
         \Nhd_1[u,B](t) & \  :=-\int_0^t e^{(t-\tau)\Deltah} \nablah \cdot ( \Bh \otimes \uh )(\tau)d\tau, \\
     \Nhd_2[u,B](t) & \  :=\int_0^t e^{(t-\tau)\Deltah} \nablah \cdot( \uh \otimes \Bh )(\tau)d\tau, \\
      \Nhd_3[u,B](t)    & \ :=- \int_0^t e^{(t-\tau)\Deltah}  \partial_3 ( u_3 \Bh )(\tau)d\tau,\\
     \Nhd_4[u,B](t)    & \ :=\int_0^t e^{(t-\tau)\Deltah}  \partial_3 ( B_3 \uh )(\tau)d\tau,\\
     \end{split}
\end{equation}
and
    \begin{align*}
        \Nvd_1[u,B](t)&:=-\int_0^t e^{(t-\tau)\Deltah} \nablah \cdot ( B_3  \uh )(\tau)d\tau,\\
        \Nvd_2[u,B](t)&:=\int_0^t e^{(t-\tau)\Deltah} \nablah \cdot (u_3  \Bh )(\tau)d\tau.\\
    \end{align*}
\end{prop}

Observe that the decomposition of magnetic field follows directly from that of integral equation (\ref{hm3}), while the decomposition of velocity field is obtained similar to the anisotropic Navier-Stokes system, cf. \cite{FM}. For the convenience of the reader, we include a short derivation for the latter one in the Appendix, see subsection \ref{de}.

\subsection{Decay estimates for the Duhamel terms}
The goal of this subsection is to show the decay estimates for the Duhamel terms in (\ref{hm4})-(\ref{hm7}). To do this, we assume that $s\in \mathbb{N},0<T\leq \infty, A, A_{\ast} \geq 0$ and make the following ansatzes:
\begin{enumerate}[(i)]
\item {
$u,B\in C([0,\infty);X^s(\R^3)),\Grad \cdot u=\Grad \cdot B=0$ and
\beq\label{hm8}
\| (u,B)(t) \|_{H^s} \leq A,\,\,     \| (u,B)(t) \| _{  L^1_{\xh}(W^{1,1}\cap W^{1,\infty})_{x_3}         }
\leq A(1+t)
\eeq
for any $t>0$.
}
\item {
For any $1\leq p \leq \infty$,
\begin{equation}\label{hm9}
     \begin{split}
      \| \nabla ^{\alpha}   \Bh (t) \|_{L^p} & \ \leq C A t^{ -(1-\f{1}{p}) -\f{|\alphah|}{2}        }  \\
    \| \nablah ^{\alphah}   B_3 (t) \|_{L^p} & \   \leq C A t^{ -\f{3}{2}(1-\f{1}{p}) -\f{|\alphah|}{2}        }  , \\
      \| \nabla ^{\alpha}   u (t) \|_{L^p}    & \  \leq C A  t^{  -\f{9}{8}(1-\f{1}{p})    -\f{ 1+|\alpha|   }{2}             }  \log(2+t) , \\
     \end{split}
\end{equation}
for any $0<t<T$ and any $\alpha=(\alphah,\alpha_3)\in (\mathbb{N} \cup \{ 0\})^2 \times (\mathbb{N} \cup \{ 0\})$ with $|\alpha|\leq 1$.
}
\item{
For any $0<t <T$ and $\alpha=(\alphah,\alpha_3)\in (\mathbb{N} \cup \{ 0\})^2 \times (\mathbb{N} \cup \{ 0\})$ with $|\alpha|\leq 1$, { }
\begin{equation}\label{hm10}
     \begin{split}
        \| \Grad^{\alpha} B (t) \|_{L^{\infty}_{\rm h}  L^{1}_{\rm v}} & \  \leq A t^{-1-\f{|\alphah|}{2}}, \\
     \| \Grad^{\alpha} u (t) \|_{L^{\infty}_{\rm h}  L^{1}_{\rm v}} & \  \leq A t^{-1-\f{1+|\alpha|}{2}}  \log(2+t)  . \\
     \end{split}
\end{equation}
}
\item{
For any $t>0$,  
\begin{equation}\label{hm11}
     \begin{split}
        \| |\xh|\Bh (t,x) \|_{L^{1}_{\rm h}  L^{\infty}_{\rm v}} & \  \leq A_{\ast} (1+t)^{  \f{1}{2}  },  \\
     \| |\xh| B_3 (t,x) \|_{L^{1}_{\rm h}  L^{\infty}_{\rm v}} & \  \leq A_{\ast},  \\
         \| |x| \uh (t,x) \|_{L^{1}_{\rm h}  L^{\infty}_{\rm v}} & \  \leq A_{\ast} (1+t)^{  \f{1}{2}  },  \\
           \| |x| u_3 (t,x) \|_{L^{1}_{\rm h}  L^{\infty}_{\rm v}} & \  \leq A_{\ast} (1+t)^{  \f{1}{2}  }.  \\
     \end{split}
\end{equation}
}
\end{enumerate}

Similar to Section \ref{hd}, one gets the following decay estimates from the ansatzes above. The proof is straightforward and thus omitted.
\begin{lemm}\label{lem12}
Let $(u,B)$ be subject to the ansatzes $(i),(ii)$ with $s\geq 3$. There exists a generic constant $C>0$ such that for any $0<t<T$ and any $1\leq p \leq \infty$ we have the decay estimates for magnetic field:
\begin{equation}\label{hm12}
     \begin{split}
         \| \Bh(t) \|_{L^p},\| \p_3 \Bh (t)\|_{L^p} & \  \leq  CA (1+t)^{-(1-\f{1}{p})}   ,          \\
     \|  \nablah \Bh (t)\|_{L^p}  & \  \leq    \left\{\begin{aligned}
& C A (1+t)^{  -(1-\f{1}{p})-\f{1}{2}           }\,\,\, \text{if  } 2\leq p\leq \infty         \\
& CA t^{  -(1-\f{1}{p})-\f{1}{2}           }\,\,\, \text{if  }   1\leq p <2,         \\
\end{aligned}\right.                           \\
       \| B_3 ( t) \|_{L^p}  & \  \leq      CA (1+t)^{  -\f{3}{2}   (1-\f{1}{p})         }      ,           \\
        \| \p_3 B_3 ( t) \|_{L^p}  & \  \leq      CA (1+t)^{  -   (1-\f{1}{p})  -\f{1}{2}       }      ,           \\
          \|  \nablah B_3 (t)\|_{L^p}  & \   \leq    \left\{\begin{aligned}
& CA (1+t)^{  -\f{3}{2} (1-\f{1}{p})-\f{1}{2}           }\,\,\, \text{if  } 2\leq p\leq \infty         \\
& CA t^{  -  \f{3}{2}  (1-\f{1}{p})-\f{1}{2}           }\,\,\, \text{if  }   1\leq p <2,         \\
\end{aligned}\right.                           \\
     \| \Bh (t) \|_{L^1_{\rm h}  L^{\infty}_{\rm v}}    & \  \leq  CA  ,  \\
     \| B_3 (t) \|_{L^1_{\rm h}  L^{\infty}_{\rm v}}    & \  \leq  C A (1+t)^{-\f{1}{2} },  \\
     \end{split}
\end{equation}
and the decay estimates for velocity field:
\begin{equation}\label{hm13}
     \begin{split}
         \| u(t) \|_{L^p} & \   \leq CA (1+t )^{  -\f{9}{8} (1-\f{1}{p})-\f{1}{2}           } \log(2+t)  ,  \\
        \|  \nablah u (t)\|_{L^p}  & \   \leq    \left\{\begin{aligned}
& CA (1+t)^{  -\f{9}{8} (1-\f{1}{p})-1          } \log(2+t) \,\,\, \text{if  } 2\leq p\leq \infty         \\
& CA t^{  -  \f{9}{8}  (1-\f{1}{p})-1           } \log(2+t) \,\,\, \text{if  }   1\leq p <2,         \\
\end{aligned}\right.                           \\
        \| \p_3 u(t) \|_{L^p} & \   \leq CA (1+t )^{  -\f{9}{8} (1-\f{1}{p})-1           }  \log(2+t),  \\
       \| u(t) \|_{L^1_{\rm h}  L^{\infty}_{\rm v}}   & \ \leq CA (1+t)^{-1} \log(2+t)  .
     \end{split}
\end{equation}
\end{lemm}

In the same spirit of Lemma \ref{lem7}, we next come to the decay estimates of the Duhamel terms. These estimates are crucial in both the nonlinear decay estimates and the asymptotic expansion of solutions. In particular, we shall use extensively Lemma \ref{lem12} without pointing it out each time.
\begin{lemm}\label{lem13}
Let $(u,B)$ be subject to the ansatzes $(i),(ii)$ with $s\geq 5$. There exists a generic constant $C>0$ such that for any $0<t<T$ and any $1\leq p \leq \infty$ we have the decay estimates for the Duhamel terms associated with magnetic equation:
\begin{equation}\label{hm14}
     \begin{split}
         \| \Grad^{\alpha} \Nhd_1[u,B](t)   \|_{L^p}     & \ \leq CA^2 (1+t )^{ -(1-\f{1}{p}) -\f{1+ |\alphah|}{2}   } ,  \\
         \| \Grad^{\alpha} \Nhd_2[u,B](t)   \|_{L^p}     & \ \leq CA^2 (1+t )^{ -(1-\f{1}{p}) -\f{1+|\alphah|}{2}   }  ,        \\
         \| \Grad^{\alpha} \Nhd_3[u,B](t)   \|_{L^p}     & \ \leq CA^2(1+t )^{ -(1-\f{1}{p}) -\f{|\alphah|}{2}   } ,        \\
         \| \Grad^{\alpha} \Nhd_4[u,B](t)   \|_{L^p}     & \ \leq CA^2(1+t )^{ -(1-\f{1}{p}) -\f{|\alphah|}{2}   } ,        \\
     \end{split}
\end{equation}
\begin{equation}\label{hm15}
     \begin{split}
        \| \Grad^{\alpha} \Nvd_1[u,B](t)   \|_{L^p} & \  \leq CA^2 (1+t )^{ -(1-\f{1}{p}) -\f{1+|\alphah|}{2}   }  , \\
     \| \Grad^{\alpha} \Nvd_2[u,B](t)   \|_{L^p} & \  \leq CA^2 (1+t )^{ -(1-\f{1}{p}) -\f{1+|\alphah|}{2}   }   , \\
     \end{split}
\end{equation}
for any $\alpha=(\alphah,\alpha_3)\in (\mathbb{N} \cup \{ 0\})^2 \times (\mathbb{N} \cup \{ 0\})$ with $|\alpha|\leq 1$ and the decay estimates for the Duhamel terms associated with velocity equation:
\begin{equation}\label{hm16}
     \begin{split}
         \| \Grad^{\alpha} \Nhe_1[u](t)   \|_{L^p} & \           \leq CA^2  (1+t )^{  -\f{3}{2}(1-\f{1}{p})    -\f{ 1+|\alpha|   }{2}             }  , \\
    \| \Grad^{\alpha} \Nhe_2[u](t)   \|_{L^p} & \           \leq CA^2 (1+t )^{  -\f{3}{2}(1-\f{1}{p})    -\f{ 1+|\alpha|   }{2}             }  , \\
          \| \Grad^{\alpha} \Nhe_3[u](t)   \|_{L^p} & \           \leq CA^2 (1+t ) ^{  -\f{3}{2}(1-\f{1}{p})    -\f{ 1+|\alpha|   }{2}             } , \\
            \| \Grad^{\alpha} \Nhe_4[u](t)   \|_{L^p} & \           \leq CA^2 (1+t ) ^{  -\f{3}{2}(1-\f{1}{p})    -\f{ 1+|\alpha|   }{2}             } , \\
              \| \Grad^{\alpha} \Nhe_5[u](t)   \|_{L^p} & \           \leq CA^2  (1+t )^{  -\f{3}{2}(1-\f{1}{p})    -\f{ 1+|\alpha|   }{2}             }   , \\
                \| \Grad^{\alpha} \Nhe_1[B](t)   \|_{L^p} & \           \leq CA^2 (1+t ) ^{  -\f{9}{8}(1-\f{1}{p})    -\f{ 1+|\alpha|   }{2}             }  , \\
    \| \Grad^{\alpha} \Nhe_2[B](t)   \|_{L^p} & \           \leq CA^2 (1+t ) ^{  -\f{9}{8}(1-\f{1}{p})    -\f{ 1+|\alpha|   }{2}             }  \log(2+t) , \\
          \| \Grad^{\alpha} \Nhe_3[B](t)   \|_{L^p} & \           \leq CA^2 (1+t ) ^{  -\f{9}{8}(1-\f{1}{p})    -\f{ 1+|\alpha|   }{2}             } \log(2+t)  , \\
            \| \Grad^{\alpha} \Nhe_4[B](t)   \|_{L^p} & \           \leq CA^2 (1+t ) ^{  -\f{3}{2}(1-\f{1}{p})    -\f{ 1+|\alpha|   }{2}             }  , \\
              \| \Grad^{\alpha} \Nhe_5[B](t)   \|_{L^p} & \           \leq CA^2 (1+t ) ^{  -\f{3}{2}(1-\f{1}{p})    -\f{ 1+|\alpha|   }{2}             }  , \\
     \end{split}
\end{equation}
\begin{equation}\label{hm17}
     \begin{split}
         \| \Grad^{\alpha} \Nve_1[u](t)   \|_{L^p} & \           \leq CA^2  (1+t )^{  -\f{3}{2}(1-\f{1}{p})    -\f{ 1+|\alpha|   }{2}             }  , \\
    \| \Grad^{\alpha} \Nve_2[u](t)   \|_{L^p} & \           \leq CA^2 (1+t ) ^{  -\f{3}{2}(1-\f{1}{p})    -\f{ 1+|\alpha|   }{2}             } , \\
          \| \Grad^{\alpha} \Nve_3[u](t)   \|_{L^p} & \           \leq CA^2  (1+t ) ^{  -\f{3}{2}(1-\f{1}{p})    -\f{ 1+|\alpha|   }{2}             }  , \\
           \| \Grad^{\alpha} \Nve_4[u](t)   \|_{L^p} & \           \leq CA^2  (1+t )^{  -\f{3}{2}(1-\f{1}{p})    -\f{ 1+|\alpha|   }{2}             } , \\
            \| \Grad^{\alpha} \Nve_5[u](t)   \|_{L^p} & \           \leq CA^2  (1+t ) ^{  -\f{3}{2}(1-\f{1}{p})    -\f{ 1+|\alpha|   }{2}             }  , \\
              \| \Grad^{\alpha} \Nve_1[B](t)   \|_{L^p} & \           \leq CA^2 (1+t )  ^{  -\f{3}{2}(1-\f{1}{p})    -\f{ 1+|\alpha|   }{2}             } , \\
    \| \Grad^{\alpha} \Nve_2[B](t)   \|_{L^p} & \           \leq CA^2 (1+t ) ^{  -\f{3}{2}(1-\f{1}{p})    -\f{ 1+|\alpha|   }{2}             }  , \\
          \| \Grad^{\alpha} \Nve_3[B](t)   \|_{L^p} & \           \leq CA^2 (1+t )^{  -\f{9}{8}(1-\f{1}{p})    -\f{ 1+|\alpha|   }{2}             }  \log(2+t)  , \\
           \| \Grad^{\alpha} \Nve_4[B](t)   \|_{L^p} & \           \leq CA^2 (1+t )^{  -\f{3}{2}(1-\f{1}{p})    -\f{ 1+|\alpha|   }{2}             }   , \\
            \| \Grad^{\alpha} \Nve_5[B](t)   \|_{L^p} & \           \leq CA^2 (1+t )^{  -\f{3}{2}(1-\f{1}{p})    -\f{ 1+|\alpha|   }{2}             }   , \\
     \end{split}
\end{equation}
for any $\alpha=(\alphah,\alpha_3)\in (\mathbb{N} \cup \{ 0\})^2 \times (\mathbb{N} \cup \{ 0\})$ with $|\alpha|\leq 1$.
\end{lemm}
{\it Proof.} Owing to the interpolation inequality, it suffices to verify all the estimates for $p=1$ and $p=\infty$ respectively. From the expressions of these Duhamel terms we know that the first two estimates in (\ref{hm14}) fall into the following type:
\begin{equation}\label{hm18}
\left\|
\Grad^{\alpha} \int_0^t e^{(t-\tau) \nablah } \nablah (B_k u_l)(\tau)  d \tau
\right\|
\leq CA^2 (1+t )^{ -(1-\f{1}{p}) -\f{1+ |\alphah|}{2}   } ,
\end{equation}
for any $k,l \in \{ 1,2 \}$. If $\alpha_3=0$,
\[
\| B_k u_l (\tau) \|_{L^1} \leq \| B_k (\tau) \|_{L^1} \| u_l (\tau)  \|_{L^{\infty}}
\leq CA^2 (1+\tau) ^{- \f{13}{8}   } \log(2+\tau)        .
\]
If $\alpha_3=1$,
\begin{equation}\nonumber
     \begin{split}
         \| \p_3 (B_k u_l) (\tau) \|_{L^1} \leq & \ \| \p_3 B_k (\tau)  \|_{L^1}  \| u_l(\tau) \|_{L^{\infty}}
+ \|  B_k (\tau)  \|_{L^1}  \| \p_3 u_l(\tau) \|_{L^{\infty}}    \\
  \leq  & \  CA^2 (1+\tau) ^{- \f{13}{8}   } \log(2+\tau)  +CA^2 (1+\tau) ^{- \f{17}{8}   } \log(2+\tau)     \\
       \leq   & \ CA^2 (1+\tau) ^{- \f{13}{8}   } \log(2+\tau)  . \\
     \end{split}
\end{equation}
Thus it holds
\[
\| \p_3^{\alpha_3} (B_k u_l) (\tau) \|_{L^1} \leq  CA^2 (1+\tau) ^{- \f{13}{8}   } \log(2+\tau) .
\]
Similarly,
\begin{equation}\nonumber
     \begin{split}
        \|  \nablah (B_k u_l)(\tau) \|_{L^1} \leq  & \  \|  \nablah B_k (\tau)        \| _{L^{\infty}}     \| u_l (\tau )  \|_{L^1}
 + \|   B_k (\tau)        \|_{L^1} \| \nablah u_l (\tau )  \|_{L^{\infty}}  \\
   \leq  & \  CA^2 (1+\tau)^{-2} \log(2+\tau)  +CA^2 (1+\tau)^{-\f{17}{8}} \log(2+\tau)          \\
      \leq    & \ CA^2 (1+\tau)^{-2} \log(2+\tau) . \\
     \end{split}
\end{equation}
yielding
\[
 \|  \nabla^{\alpha} (B_k u_l)(\tau) \|_{L^1} \leq  CA^2 (1+\tau) ^{- \f{13}{8}   } \log(2+\tau) ;
\]
whence
\begin{equation}\nonumber
     \begin{split}
       \left\|   \Grad^{\alpha} \int_0^t e^{ (t-\tau)\Deltah  }\nablah (B_k u_l ) (\tau  )     d \tau \right\|_{L^1}           \leq  & \  C   \int_0^  {  \f{t}{2} }  \| \nablah ^{\alphah} \nablah \Gh (t-\tau) \|_{L^{1}(\R^2)} \| \p_3^{\alpha_3}  (B_k u_l ) (\tau) \| _{L^1}        d \tau                     \\
         & \ +     C          \int_{  \f{t}{2} } ^t        \| \nablah \Gh (t-\tau) \|_{L^{1}(\R^2)} \| \Grad^{\alpha}  (B_k   u_l) (\tau) \|   _{L^1}          d \tau         \\
     \leq     & \   CA^2 \int_0^  {  \f{t}{2} } (t-\tau)^{ -\f{  1+|\alphah|  }{  2 }  }   (1+\tau) ^{- \f{13}{8}   } \log(2+\tau)  d\tau \\
     &\ + CA^2  \int_{  \f{t}{2} } ^t    (t-\tau)^{ -\f{1}{2}   }  (1+\tau) ^{- \f{13}{8}   } \log(2+\tau)  d\tau \\
\leq      &\  CA^2 t^{ -\f{  1+|\alphah|  }{  2 }  }
 \end{split}
\end{equation}
for $t \geq 1$ and
\begin{equation}\nonumber
     \begin{split}
     \left\|   \Grad^{\alpha}  \int_0^t e^{ (t-\tau)\Deltah  }\nablah (B_k  u_l  ) (\tau  )     d \tau \right\|_{L^{1}}      \leq  & \  C  \int_0^t   \| \nablah \Gh (t-\tau) \|_{L^{1}(\R^2)} \| \Grad^{\alpha}  (B_k   u_l) (\tau) \|   _{L^1}          d \tau                     \\
     \leq     & \   CA^2    \int_0^  {  t }    (t-\tau)^{ -\f{1}{2}   }   (1+\tau) ^{- \f{13}{8}   } \log(2+\tau)  d\tau   \\
     \leq &\ CA^2
     \end{split}
\end{equation}
for $0<t \leq 1$. This gives (\ref{hm18}) with $p=1$. In addition, if $\alpha_3=0$,
\[
\|  (B_k u_l)(\tau) \|_{L^1_{\rm h}  L^{\infty}_{\rm v}}
\leq \| B_k(\tau) \|_{L^{\infty}}  \| u_l (\tau) \|_{L^1_{\rm h}  L^{\infty}_{\rm v}}
\leq CA^2 (1+\tau)^{-2} \log(2+\tau) .
\]
If $\alpha_3=1$,
\begin{equation}\nonumber
     \begin{split}
       \| \p_3 (B_k u_l)(\tau) \|_{L^1_{\rm h}  L^{\infty}_{\rm v}}
  \leq &\    \|\p_3 B_k (\tau)  \|_{L^{\infty}}  \| u_l (\tau) \|_{L^1_{\rm h}  L^{\infty}_{\rm v}} + \|\p_3 u_l (\tau)  \|_{L^{\infty}}  \| B_k (\tau) \|_{L^1_{\rm h}  L^{\infty}_{\rm v}}                \\
     \leq  & \  CA^2 (1+\tau)^{-2} \log(2+\tau) .           \\
     \end{split}
\end{equation}
Thus it holds
\[
 \| \p_3^{\alpha_3} (B_k u_l)(\tau) \|_{L^1_{\rm h}  L^{\infty}_{\rm v}}
  \leq  CA^2 (1+\tau)^{-2} \log(2+\tau) .
\]
Analogously, one sees that
\[
\| (B_k u_l)(\tau)  \|_{L^{\infty}} \leq \| B_k(\tau) \|_{L^{\infty}}   \| u_l (\tau) \|_{L^{\infty}} \leq  CA^2 (1+\tau)^{-\f{21}{8}} \log(2+\tau) ,
\]
\begin{equation}\nonumber
     \begin{split}
      \|\nablah  (B_k u_l)(\tau)  \|_{L^{\infty}} \leq & \ \| \nablah B_k(\tau) \|_{L^{\infty}}   \| u_l (\tau) \|_{L^{\infty}} + \|  B_k(\tau) \|_{L^{\infty}}   \| \nablah u_l (\tau) \|_{L^{\infty}}  \\
   \leq  & \ CA^2 (1+\tau)^{-\f{25}{8}} \log(2+\tau),  \\
     \end{split}
\end{equation}
\begin{equation}\nonumber
     \begin{split}
        \|\p_3  (B_k u_l)(\tau)  \|_{L^{\infty}} \leq  & \ \| \p_3 B_k(\tau) \|_{L^{\infty}}   \| u_l (\tau) \|_{L^{\infty}} + \|  B_k(\tau) \|_{L^{\infty}}   \| \p_3 u_l (\tau) \|_{L^{\infty}}  \\
    \leq  & \  CA^2 (1+\tau)^{-\f{21}{8}} \log(2+\tau)  ,\\
     \end{split}
\end{equation}
yielding
\[
\| \Grad^{\alpha} (B_k u_l)(\tau) \|_{L^{\infty}} \leq
 CA^2 (1+\tau)^{-\f{21}{8}} \log(2+\tau)  .
\]
Consequently,
\begin{equation}\nonumber
     \begin{split}
       \left\|   \Grad^{\alpha} \int_0^t e^{ (t-\tau)\Deltah  }\nablah (B_k u_l ) (\tau  )     d \tau \right\| _{L^{\infty}}              \leq  & \  C   \int_0^  {  \f{t}{2} }  \| \nablah ^{\alphah} \nablah \Gh (t-\tau) \|_{L^{\infty}(\R^2)} \| \p_3^{\alpha_3}  (B_k u_l ) (\tau) \| _{L^1_{\rm h}  L^{\infty}_{\rm v}}        d \tau                     \\
         & \ +     C          \int_{  \f{t}{2} } ^t        \| \nablah \Gh (t-\tau) \|_{L^{1}(\R^2)} \| \Grad^{\alpha}  (B_k   u_l) (\tau) \|   _{L^{\infty}}          d \tau         \\
     \leq     & \   CA^2 \int_0^  {  \f{t}{2} } (t-\tau)^{ -1-\f{  1+|\alphah|  }{  2 }  }    (1+\tau) ^{ -2     } \log(2+\tau) d\tau \\
     &\ + CA^2  \int_{  \f{t}{2} } ^t    (t-\tau)^{ -\f{1}{2}   }   (1+\tau)^{-\f{21}{8}} \log(2+\tau) d\tau \\
\leq      &\  CA^2 t^{-1 -\f{  1+|\alphah|  }{  2 }  }
 \end{split}
\end{equation}
for $t \geq 1$ and
\begin{equation}\nonumber
     \begin{split}
     \left\|   \Grad^{\alpha}  \int_0^t e^{ (t-\tau)\Deltah  }\nablah (B_k  u_l  ) (\tau  )     d \tau \right\|_{L^{\infty}}       \leq  & \  C  \int_0^t   \| \nablah \Gh (t-\tau) \|_{L^{1}(\R^2)} \| \Grad^{\alpha}  (B_k   u_l) (\tau) \|   _{L^{\infty}}          d \tau                     \\
     \leq     & \   CA^2    \int_0^  {  t }    (t-\tau)^{ -\f{1}{2}   }  (1+\tau)^{-\f{21}{8}} \log(2+\tau)  d\tau   \\
     \leq &\ CA^2
     \end{split}
\end{equation}
for $0<t \leq 1$. This gives (\ref{hm18}) with $p=\infty$.

Following the same line as above, we obtain (\ref{hm15}). The details are omitted.

We proceed with (\ref{hm14})$_3$. Observe that
\begin{equation}\nonumber
     \begin{split}
    \| \p_3 (u_3 \Bh) (\tau)\|_{L^1}        \leq  & \    \|\p_3 u_3 (\tau)  \|_{L^{\infty}} \|  \Bh (\tau) \|  _{L^1}   +
    \| u_3 (\tau) \|_{L^{\infty}} \| \p_3 \Bh   (\tau) \|  _{L^1}     \\
    \leq  & \    CA^2 (1+\tau)^{-\f{13}{8}} \log(2+\tau) ,        \\
     \end{split}
\end{equation}
\begin{equation}\nonumber
     \begin{split}
    \| \p_3^2 (u_3 \Bh) (\tau) \|_{L^1}        \leq  & \    \|\p_3^2 u_3 (\tau)  \|_{L^{2}} \|  \Bh (\tau) \|  _{L^2}   +
   2 \| \p_3 u_3 (\tau) \|_{L^{2}} \| \p_3 \Bh (\tau) \|  _{L^2}     \\
   &\ +  \| u_3 (\tau)  \|_{L^{2}} \|   \p_3^2 \Bh (\tau) \|  _{L^2}  \\
   \leq &\  C \Big(  \| \p_3 u_3 (\tau) \|_{L^{2}} ^{\f{2}{3}}   \| \p_3^4  u_3 (\tau) \|_{L^{2}} ^{\f{1}{3}} \|  \Bh (\tau) \|  _{L^2}    +
    \| \p_3 u_3 (\tau) \|_{L^{2}} \| \p_3 \Bh (\tau) \|  _{L^2}             \\
    &\    +  \| u_3 (\tau)  \|_{L^{2}} \|   \p_3^2 \Bh (\tau) \|  _{L^2}   \Big)           \\
 \leq &\   CA^2  (1+\tau) ^{  -\f{17}{16}          }   \log(2+\tau)   ;         \\
     \end{split}
\end{equation}
whence
\begin{equation}\nonumber
     \begin{split}
     \left\|   \Grad^{\alpha} \int_0^t e^{ (t-\tau)\Deltah  }\p_3 (u_3 \Bh ) (\tau  )     d \tau \right\|_{L^1}   \leq  & \  C  \int_0^t  \| \nablah^{\alphah} \Gh (t-\tau) \|_{L^1(\R^2)} \| \p_3^{\alpha_3+1} (u_3 \Bh) (\tau) \|     _{L^1}         d \tau                     \\
  \leq  & \  CA^2          \int_0^  {  \f{t}{2} } (t-\tau)^{  -\f{|\alphah|}{2}         } (1+\tau) ^{  -\f{17}{16}          }   \log(2+\tau)  d \tau                 \\
         & \ +     CA^2          \int_{  \f{t}{2} } ^t (t-\tau)^{  -\f{|\alphah|}{2}         }  (1+\tau) ^{  -\f{17}{16}          }   \log(2+\tau)   d\tau                 \\
     \leq     & \   CA^2 t^{    -\f{|\alphah|}{2}           }
     \end{split}
\end{equation}
for $t \geq 1$ and
\begin{equation}\nonumber
     \begin{split}
     \left\|   \Grad^{\alpha} \int_0^t e^{ (t-\tau)\Deltah  }\p_3 (u_3 \Bh ) (\tau  )     d \tau \right\|_{L^1}   \leq  & \  C  \int_0^t  \| \nablah^{\alphah} \Gh (t-\tau) \|_{L^1(\R^2)} \| \p_3^{\alpha_3+1} (u_3 \Bh) (\tau) \|     _{L^1}         d \tau                     \\
  \leq  & \  CA^2          \int_0^  {  t } (t-\tau)^{  -\f{|\alphah|}{2}         } (1+\tau) ^{  -\f{17}{16}          }   \log(2+\tau)   d \tau                  \\
     \leq     & \   CA^2
     \end{split}
\end{equation}
for $0<t \leq 1$. This proves (\ref{hm14})$_3$ for the case $p=1$. Next we turn to the case $p=\infty$ and divide it into two classes: either $\alpha_3=0$ or $\alpha_3=1$. If $\alpha_3=0$,
\begin{equation}\nonumber
     \begin{split}
        \| \p_3(u_3 \Bh) (\tau) \|_{L^1_{\rm h}  L^{\infty}_{\rm v}}  \leq  & \   \|  \p_3 u_3 (\tau) \| _{L^{\infty}}
        \| \Bh ( \tau ) \|_{L^1_{\rm h}  L^{\infty}_{\rm v}}   +
        \| u_3 ( \tau ) \|_{L^1_{\rm h}  L^{\infty}_{\rm v}}    \|  \p_3  \Bh  (\tau) \| _{L^{\infty}}              \\
    \leq  & \ CA^2 (1+\tau) ^{ -2     }   \log(2+\tau) ,
     \end{split}
\end{equation}
\begin{equation}\nonumber
     \begin{split}
        \| \p_3(u_3 \Bh) (\tau) \|_{L^{\infty}}  \leq  & \   \|  \p_3 u_3 (\tau) \| _{L^{\infty}}
        \| \Bh ( \tau ) \|_{L^{\infty}}  +
        \| u_3 ( \tau ) \|_{L^{\infty}}   \|  \p_3  \Bh  (\tau) \| _{L^{\infty}}              \\
    \leq  & \ CA^2 (1+\tau) ^{ -\f{21}{8}   }  \log(2+\tau ) ;      \\
     \end{split}
\end{equation}
whence
\begin{equation}\nonumber
     \begin{split}
     \left\|   \nablah^{\alphah} \int_0^t e^{ (t-\tau)\Deltah  }\p_3 (u_3\Bh ) (\tau  )     d \tau \right\|_{L^{\infty}}   \leq  & \  C   \int_0^  {  \f{t}{2} }  \| \nablah^{\alphah} \Gh (t-\tau) \|_{L^{\infty}(\R^2)} \| \p_3  (u_3 \Bh) (\tau) \|     _{L^1_{\rm h}  L^{\infty}_{\rm v}}         d \tau                     \\
         & \ +     C          \int_{  \f{t}{2} } ^t        \| \nablah^{\alphah} \Gh (t-\tau) \|_{L^{1}(\R^2)} \| \p_3  (u_3 \Bh) (\tau) \|    _{L^{\infty}}         d \tau         \\
     \leq     & \   CA^2 \int_0^  {  \f{t}{2} } (t-\tau)^{ -1-\f{|\alphah|}{2}   }    (1+\tau) ^{ -2   } \log(2+\tau)      d\tau \\
     &\ + CA^2  \int_{  \f{t}{2} } ^t    (t-\tau)^{ -\f{|\alphah|}{2}   }   (1+\tau) ^{ -\f{21}{8}   }  \log(2+\tau)       d\tau           \\
\leq      &\  CA^2 t ^{ -1-\f{|\alphah|}{2}   }
     \end{split}
\end{equation}
for $t \geq 1$ and
\begin{equation}\nonumber
     \begin{split}
     \left\|  \nablah^{\alphah}  \int_0^t e^{ (t-\tau)\Deltah  }\p_3 (u_3\Bh ) (\tau  )     d \tau \right\|_{L^{\infty}}      \leq  & \  C  \int_0^t  \| \nablah^{\alphah} \Gh (t-\tau) \|_{L^1(\R^2)} \| \p_3 (u_3 \Bh) (\tau) \|     _{L^{\infty}}        d \tau                     \\
     \leq     & \   CA^2    \int_0^  {  t }   (t-\tau)^{  -\f{|\alphah|}{2}         }  (1+\tau) ^{ -\f{21}{8}   }  \log(2+\tau)    d \tau   \\
     \leq &\ CA^2
     \end{split}
\end{equation}
for $0<t \leq 1$. This verifies (\ref{hm14})$_3$ for the case $p=\infty,\alpha_3=0$. If $\alpha_3=1$, we invoke the Gagliardo-Nirenberg inequality to deduce
\begin{equation}\nonumber
     \begin{split}
        \| \p_3^2(u_3 \Bh) (\tau)  \|_{L^1_{\rm h}  L^{\infty}_{\rm v}}  \leq     & \  \| \p_3^2 u_3 (\tau) \|_{L^2_{\rm h}  L^{\infty}_{\rm v}}   \| \Bh (\tau)  \|_{L^2_{\rm h}  L^{\infty}_{\rm v}}
        + 2 \| \p_3 u_3 (\tau) \|_{L^2_{\rm h}  L^{\infty}_{\rm v}}    \| \p_3  \Bh (\tau)  \|_{L^2_{\rm h}  L^{\infty}_{\rm v}}             \\
         & \  +      \|  u_3 (\tau) \|_{L^2_{\rm h}  L^{\infty}_{\rm v}}   \|  \p_3^2 \Bh (\tau)  \|_{L^2_{\rm h}  L^{\infty}_{\rm v}}                             \\
      \leq    & \ C  \Big(   \| \p_3 u_3 (\tau) \|_{L^2} ^{ \f{5}{8}  } \| \p_3^5 u_3 (\tau) \|_{L^2} ^{ \f{3}{8}  }           \|  \Bh(\tau) \| _{L^2}  ^{ \f{1}{2}  }   \| \p_3  \Bh(\tau) \| _{L^2}  ^{ \f{1}{2}  }       \\
         & \  +   \| \p_3 u_3 (\tau) \|_{L^2} ^{ \f{1}{2}  } \| \p_3^2 u_3 (\tau) \|_{L^2} ^{ \f{1}{2}  }                   \| \p_3 \Bh (\tau) \|_{L^2} ^{ \f{1}{2}  } \| \p_3^2 \Bh (\tau) \|_{L^2} ^{ \f{1}{2}  }       \\
         &\ +        \|  u_3 (\tau) \|_{L^2} ^{ \f{1}{2}  } \| \p_3 u_3 (\tau) \|_{L^2} ^{ \f{1}{2}  }                   \| \p_3^2 \Bh (\tau) \|_{L^2} ^{ \f{1}{2}  } \| \p_3^3 \Bh (\tau) \|_{L^2} ^{ \f{1}{2}  }   \Big)       \\
       \leq   &\   CA^2  (1+\tau)^{- \f{33}{32}  } (\log(2+\tau))^{\f{1}{2}}  ,       \\
     \end{split}
\end{equation}
\begin{equation}\nonumber
     \begin{split}
        \| \p_3^2(u_3 \Bh) (\tau)  \|_{L^2_{\rm h}  L^{\infty}_{\rm v}}  \leq     & \  \| \p_3^2 u_3 (\tau) \|_{L^2_{\rm h}  L^{\infty}_{\rm v}}   \| \Bh (\tau)  \|_{L^{\infty}}
        + 2 \| \p_3 u_3 (\tau) \|_{L^{\infty}}     \| \p_3  \Bh (\tau)  \|_{L^2_{\rm h}  L^{\infty}_{\rm v}}             \\
         & \  +      \|  u_3 (\tau) \|_{L^{\infty}}    \|  \p_3^2 \Bh (\tau)  \|_{L^2_{\rm h}  L^{\infty}_{\rm v}}                             \\
      \leq    & \ C  \Big(   \| \p_3 u_3 (\tau) \|_{L^2} ^{ \f{5}{8}  } \| \p_3^5 u_3 (\tau) \|_{L^2} ^{ \f{3}{8}  }              \| \Bh (\tau)  \|_{L^{\infty}}         \\
         & \  +  \| \p_3 u_3 (\tau) \|_{L^{\infty}}                 \| \p_3 \Bh (\tau) \|_{L^2} ^{ \f{1}{2}  } \| \p_3^2 \Bh (\tau) \|_{L^2} ^{ \f{1}{2}  }       \\
         &\ +        \|  u_3 (\tau) \|_{L^{\infty}}                   \| \p_3^2 \Bh (\tau) \|_{L^2} ^{ \f{1}{2}  } \| \p_3^3 \Bh (\tau) \|_{L^2} ^{ \f{1}{2}  }   \Big)                          \\
       \leq   &\   CA^2  (1+\tau)^{- \f{13}{8} }  \log(2+\tau) ;      \\
     \end{split}
\end{equation}
whence
\begin{equation}\nonumber
     \begin{split}
     \left\|   \p_3 \int_0^t e^{ (t-\tau)\Deltah  }\p_3 (u_3\Bh ) (\tau  )     d \tau \right\|_{L^{\infty}}   \leq  & \  C   \int_0^  {  \f{t}{2} }  \|  \Gh (t-\tau) \|_{L^{\infty}(\R^2)} \| \p_3^2  (u_3 \Bh) (\tau) \|     _{L^1_{\rm h}  L^{\infty}_{\rm v}}         d \tau                     \\
         & \ +     C          \int_{  \f{t}{2} } ^t        \|  \Gh (t-\tau) \|_{L^{2}(\R^2)} \| \p_3^2  (u_3 \Bh) (\tau) \|   _{L^2_{\rm h}  L^{\infty}_{\rm v}}         d \tau         \\
     \leq     & \   CA^2 \int_0^  {  \f{t}{2} } (t-\tau)^{ -1  }  (1+\tau)^{- \f{33}{32}  } (\log(2+\tau))^{\f{1}{2}}   d\tau \\
     &\ + CA^2  \int_{  \f{t}{2} } ^t    (t-\tau)^{ -\f{1}{2}   }   (1+\tau)^{- \f{13}{8} }  \log(2+\tau)    d\tau \\
\leq      &\  CA^2 t ^{ -1  }
     \end{split}
\end{equation}
for $t \geq 1$ and
\begin{equation}\nonumber
     \begin{split}
     \left\|  \p_3  \int_0^t e^{ (t-\tau)\Deltah  }\p_3 (u_3\Bh ) (\tau  )     d \tau \right\|_{L^{\infty}}      \leq  & \  C  \int_0^t  \| \Gh (t-\tau) \|_{L^2(\R^2)} \| \p_3^2 (u_3 \Bh) (\tau) \|   _{L^2_{\rm h}  L^{\infty}_{\rm v}}        d \tau                     \\
     \leq     & \   CA^2    \int_0^  {  t }   (t-\tau)^{  -\f{ 1  }{2}         }  (1+\tau)^{- \f{13}{8} }  \log(2+\tau)  d \tau   \\
     \leq &\ CA^2
     \end{split}
\end{equation}
for $0<t \leq 1$. This gives (\ref{hm14})$_3$ with $p=\infty,\alpha_3=1$. Combining the estimates above, we have verified (\ref{hm14})$_3$ completely.

Following the same route, one gets (\ref{hm14})$_4$. Indeed, the proof of (\ref{hm14})$_4$ is easier than (\ref{hm14})$_3$ in the sense that we have better decay estimates of $B_3$. The details are again omitted.

We are now in a position to handle the decay estimates for the Duhamel terms associated with velocity equation. We start with the estimates of $\| \Grad^{\alpha}  \Nhe_i[B](t)  \|_{L^p}, \| \Grad^{\alpha}  \Nve_i[B] (t) \|_{L^p}   $, $i\in \{1,2,3,4,5    \}$.
Notice that
\begin{equation}\nonumber
     \begin{split}
         \| (B_3 \Bh ) (\tau  )  \|_{L^1} \leq  & \  \| B_3(\tau)  \|_{L^{\infty}} \| \Bh(\tau)\|_{L^1}  \leq CA^2  (1+\tau) ^{  -\f{3}{2}          } ,    \\
      \| \p_3 (B_3 \Bh ) (\tau  )  \|_{L^1} \leq & \  \| \p_3 B_3(\tau)  \|_{L^{\infty}} \| \Bh(\tau)\|_{L^1}  +     \|  B_3(\tau)  \|_{L^{\infty}} \| \p_3 \Bh(\tau)\|_{L^1}          \\
        \leq  & \  CA^2  (1+\tau) ^{  -\f{3}{2}          }    ; \\
     \end{split}
\end{equation}
whence
\begin{equation}\nonumber
     \begin{split}
     \left\|   \Grad^{\alpha} \int_0^t e^{ (t-\tau)\Delta  }\p_3 (B_3 \Bh ) (\tau  )     d \tau \right\|_{L^1}  \leq  &\     \int_0^  {  \f{t}{2} }        \|  \Grad^{\alpha}\p_3  G(t-\tau)  \|_{L^1}   \| (B_3 \Bh)(\tau)     \| _{L^1}    d\tau              \\
     &\ +   \int_{  \f{t}{2} } ^t    \|  \Grad^{\alpha}  G(t-\tau)  \|_{L^1}   \| \p_3 (B_3 \Bh)(\tau)     \| _{L^1}       d\tau                    \\
  \leq  & \  CA^2          \int_0^  {  \f{t}{2} } (t-\tau)^{  -\f{1+ |\alpha|}{2}         } (1+\tau) ^{  -\f{3}{2}          }    d \tau                 \\
         & \ +     CA^2          \int_{  \f{t}{2} } ^t (t-\tau)^{  -\f{|\alpha|}{2}         }  (1+\tau) ^{  -\f{3}{2}          }       d\tau                 \\
     \leq     & \   CA^2 t^{    -\f{  1+  |\alpha|}{2}           }
     \end{split}
\end{equation}
for $t \geq 1$ and
\begin{equation}\nonumber
     \begin{split}
    \left\|   \Grad^{\alpha} \int_0^t e^{ (t-\tau)\Delta  }\p_3 (B_3 \Bh ) (\tau  )     d \tau \right\|_{L^1}    \leq  & \  C  \int_0^t   \|  \Grad^{\alpha}  G(t-\tau)  \|_{L^1}   \| \p_3 (B_3 \Bh)(\tau)     \| _{L^1}          d \tau                     \\
  \leq  & \  CA^2          \int_0^  {  t } (t-\tau)^{  -\f{|\alpha|}{2}         } (1+\tau)  ^{  -\f{3}{2}          }    d \tau                  \\
     \leq     & \   CA^2
     \end{split}
\end{equation}
for $0<t \leq 1$. Moreover,
\begin{equation}\nonumber
     \begin{split}
         \| \p_3 (B_3 \Bh ) (\tau  )  \|_{L^4} \leq  & \   \| \p_3 B_3(\tau)  \|_{L^{\infty}} \| \Bh(\tau)\|_{L^4} + \|  B_3(\tau)  \|_{L^{\infty}} \| \p_3 \Bh(\tau)\|_{L^4}       \\
        \leq  & \  CA^2  (1+\tau) ^{  -\f{9}{4}          }    ; \\
     \end{split}
\end{equation}
whence
\begin{equation}\nonumber
     \begin{split}
     \left\|   \Grad^{\alpha} \int_0^t e^{ (t-\tau)\Delta  }\p_3 (B_3 \Bh ) (\tau  )     d \tau \right\|_{L^{\infty}}   \leq  &\     \int_0^  {  \f{t}{2} }        \|  \Grad^{\alpha}\p_3  G(t-\tau)  \|_{L^{\infty}}  \| (B_3 \Bh)(\tau)     \| _{L^1}    d\tau              \\
     &\ +   \int_{  \f{t}{2} } ^t    \|  \Grad^{\alpha}  G(t-\tau)  \|_{L^{ \f{4}{3}} }   \| \p_3 (B_3 \Bh)(\tau)     \| _{L^4}       d\tau                    \\
  \leq  & \  CA^2          \int_0^  {  \f{t}{2} } (t-\tau)^{  -\f{3}{2}-\f{1+ |\alpha|}{2}         } (1+\tau) ^{  -\f{3}{2}          }    d \tau                 \\
         & \ +     CA^2          \int_{  \f{t}{2} } ^t (t-\tau)^{  -\f{3 }{8} -\f{ |\alpha|}{2}            }  (1+\tau) ^{  -\f{9}{4}          }       d\tau                 \\
     \leq     & \   CA^2 t^{  -\f{9}{8}  -\f{  1+  |\alpha|}{2}           }
     \end{split}
\end{equation}
for $t \geq 1$ and
\begin{equation}\nonumber
     \begin{split}
     \left\|  \p_3  \int_0^t e^{ (t-\tau)\Deltah  }\p_3 (u_3\Bh ) (\tau  )     d \tau \right\|_{L^{\infty}}      \leq  & \  C  \int_0^t   \|  \Grad^{\alpha}  G(t-\tau)  \|_{L^{ \f{4}{3}} }   \| \p_3 (B_3 \Bh)(\tau)     \| _{L^4}          d \tau                     \\
     \leq     & \   CA^2    \int_0^  {  t }   (t-\tau)^{  -\f{3 }{8} -\f{ |\alpha|}{2}            }  (1+\tau) ^{  -\f{9}{4}          }    d \tau   \\
     \leq &\ CA^2
     \end{split}
\end{equation}
for $0<t \leq 1$. Combining the above estimates gives the estimate of $\| \Grad^{\alpha}  \Nhe_1 [B](t)  \|_{L^p} $.

To proceed, we see
\begin{equation}\nonumber
     \begin{split}
         \| (\Bh \otimes \Bh ) (\tau  )  \|_{L^1} \leq  & \  \| \Bh (\tau)  \|_{L^{\infty}} \| \Bh(\tau)\|_{L^1}  \leq CA^2  (1+\tau) ^{  -1         } ,    \\
      \| \nablah \cdot (\Bh \otimes \Bh ) (\tau  )  \|_{L^1} \leq & \ C  \| \nablah \Bh(\tau)  \|_{L^{\infty}} \| \Bh(\tau)\|_{L^1}
        \leq    CA^2  (1+\tau) ^{  -\f{3}{2}          }    ; \\
     \end{split}
\end{equation}
whence
\begin{equation}\nonumber
     \begin{split}
     \left\|   \Grad^{\alpha} \int_0^t e^{ (t-\tau)\Delta  }\nablah \cdot (\Bh \otimes \Bh ) (\tau  )     d \tau \right\|_{L^1}  \leq  &\     \int_0^  {  \f{t}{2} }        \|  \Grad^{\alpha}\nablah  G(t-\tau)  \|_{L^1}   \| (\Bh \otimes \Bh)(\tau)     \| _{L^1}    d\tau              \\
     &\ +   \int_{  \f{t}{2} } ^t    \|  \Grad^{\alpha}  G(t-\tau)  \|_{L^1}   \| \nablah \cdot (\Bh \otimes \Bh ) (\tau)     \| _{L^1}       d\tau                    \\
  \leq  & \  CA^2          \int_0^  {  \f{t}{2} } (t-\tau)^{  -\f{1+ |\alpha|}{2}         } (1+\tau) ^{  -1       }    d \tau                 \\
         & \ +     CA^2          \int_{  \f{t}{2} } ^t (t-\tau)^{  -\f{|\alpha|}{2}         }  (1+\tau) ^{  -\f{3}{2}          }       d\tau                 \\
     \leq     & \   CA^2 t^{    -\f{  1+  |\alpha|}{2}           }  \log(2+t)
     \end{split}
\end{equation}
for $t \geq 1$ and
\begin{equation}\nonumber
     \begin{split}
    \left\|   \Grad^{\alpha} \int_0^t e^{ (t-\tau)\Delta  }\nablah \cdot (\Bh \otimes \Bh )  (\tau  )     d \tau \right\|_{L^1}    \leq  & \  C  \int_0^t   \|  \Grad^{\alpha}  G(t-\tau)  \|_{L^1}   \| \nablah \cdot (\Bh \otimes \Bh ) (\tau)     \| _{L^1}          d \tau                     \\
  \leq  & \  CA^2          \int_0^  {  t } (t-\tau)^{  -\f{|\alpha|}{2}         } (1+\tau)  ^{  -\f{3}{2}          }    d \tau                  \\
     \leq     & \   CA^2
     \end{split}
\end{equation}
for $0<t \leq 1$. In addition,
\begin{equation}\nonumber
     \begin{split}
     \left\|   \Grad^{\alpha} \int_0^t e^{ (t-\tau)\Delta  }\nablah \cdot (\Bh \otimes \Bh ) (\tau  )     d \tau \right\|_{L^{\infty}}   \leq  &\     \int_0^  {  \f{t}{2} }        \|  \Grad^{\alpha} \nablah  G(t-\tau)  \|_{L^{\infty}}  \| (\Bh \otimes \Bh)(\tau)     \| _{L^1}    d\tau              \\
     &\ +   \int_{  \f{t}{2} } ^t    \|  \Grad^{\alpha}  G(t-\tau)  \|_{L^{ \f{4}{3}} }   \|  \nablah \cdot (\Bh \otimes \Bh ) (\tau)     \| _{L^4}       d\tau                    \\
  \leq  & \  CA^2          \int_0^  {  \f{t}{2} } (t-\tau)^{  -\f{3}{2}-\f{1+ |\alpha|}{2}         } (1+\tau) ^{  -1         }    d \tau                 \\
         & \ +     CA^2          \int_{  \f{t}{2} } ^t (t-\tau)^{  -\f{3 }{8} -\f{ |\alpha|}{2}            }  (1+\tau) ^{  -\f{9}{4}          }       d\tau                 \\
     \leq     & \   CA^2 t^{  -\f{9}{8}  -\f{  1+  |\alpha|}{2}           }
     \end{split}
\end{equation}
for $t \geq 1$ and
\begin{equation}\nonumber
     \begin{split}
     \left\|  \Grad^{\alpha} \int_0^t e^{ (t-\tau)\Deltah  } \nablah \cdot (\Bh \otimes \Bh ) (\tau  )     d \tau \right\|_{L^{\infty}}      \leq  & \  C  \int_0^t   \|  \Grad^{\alpha}  G(t-\tau)  \|_{L^{ \f{4}{3}} }   \|\nablah \cdot (\Bh \otimes \Bh ) (\tau)     \| _{L^4}          d \tau                     \\
     \leq     & \   CA^2    \int_0^  {  t }   (t-\tau)^{  -\f{3 }{8} -\f{ |\alpha|}{2}            }  (1+\tau) ^{  -\f{9}{4}          }    d \tau   \\
     \leq &\ CA^2
     \end{split}
\end{equation}
for $0<t \leq 1$. These estimates allow us to deduce the estimate of $\| \Grad^{\alpha}  \Nhe_2 [B](t)  \|_{L^p} $.

Next, for any $k,l\in \{  1,2\}$ we infer from Lemma \ref{lm3-1} that
\begin{equation}\nonumber
     \begin{split}
     \Big\|   \Grad^{\alpha}   &\ \int_0^t e^{ (t-\tau)\Delta  } \nablah \p_k \p_l N(t-\tau) \ast (B_k  B_l  ) (\tau  )     d \tau \Big\|_{L^1}  \\
     \leq  &\     \int_0^  {  \f{t}{2} }        \|  \Grad^{\alpha} \nablah \p_k \p_l  N(t-\tau)  \|_{L^1}   \| (B_k  B_l )(\tau)     \| _{L^1}    d\tau         \\
    &\  +   \int_{  \f{t}{2} } ^t    \|  \Grad^{\alpha} \p_k \p_l N(t-\tau)  \|_{L^1}   \| \nablah  (B_k  B_l ) (\tau)     \| _{L^1}       d\tau                    \\
  \leq  & \  CA^2          \int_0^  {  \f{t}{2} } (t-\tau)^{  -\f{3+ |\alpha|}{2} +1        } (1+\tau) ^{  -1       }    d \tau                 \\
         & \ +     CA^2          \int_{  \f{t}{2} } ^t (t-\tau)^{  -\f{2+ |\alpha|}{2} +1        }  (1+\tau) ^{  -\f{3}{2}          }       d\tau                 \\
     \leq     & \   CA^2 t^{    -\f{  1+  |\alpha|}{2}           }  \log(2+\tau)
     \end{split}
\end{equation}
for $t \geq 1$ and
\begin{equation}\nonumber
     \begin{split}
    \Big\|   \Grad^{\alpha} &\   \int_0^t e^{ (t-\tau)\Delta  } \nablah \p_k \p_l N(t-\tau) \ast (B_k  B_l  ) (\tau  )        d \tau \Big\|_{L^1}  \\
    \leq  & \  C  \int_0^t  \|  \Grad^{\alpha} \p_k \p_l N(t-\tau)  \|_{L^1}   \| \nablah  (B_k  B_l ) (\tau)     \| _{L^1}        d \tau                     \\
  \leq  & \  CA^2          \int_0^  {  t }   (t-\tau)^{  -\f{2+ |\alpha|}{2} +1        }  (1+\tau) ^{  -\f{3}{2}          }     d \tau                  \\
     \leq     & \   CA^2
     \end{split}
\end{equation}
for $0<t \leq 1$. Besides,
\begin{equation}\nonumber
     \begin{split}
     \Big\|   \Grad^{\alpha}   &\ \int_0^t e^{ (t-\tau)\Delta  } \nablah \p_k \p_l N(t-\tau) \ast (B_k  B_l  ) (\tau  )     d \tau \Big\|_{L^{\infty}}   \\
     \leq  &\     \int_0^  {  \f{t}{2} }        \|  \Grad^{\alpha} \nablah \p_k \p_l  N(t-\tau)  \|_{L^{\infty}}     \| (B_k  B_l )(\tau)     \| _{L^1}    d\tau         \\
    &\  +   \int_{  \f{t}{2} } ^t    \|  \Grad^{\alpha} \p_k \p_l N(t-\tau)  \|_{L^{ \f{4}{3} }}   \| \nablah  (B_k  B_l ) (\tau)     \| _{L^4}       d\tau                    \\
  \leq  & \  CA^2          \int_0^  {  \f{t}{2} } (t-\tau)^{  -\f{3}{2}-\f{3+ |\alpha|}{2} +1        } (1+\tau) ^{  -1       }    d \tau                 \\
         & \ +     CA^2          \int_{  \f{t}{2} } ^t (t-\tau)^{ -\f{3}{8} -\f{2+ |\alpha|}{2} +1        }  (1+\tau) ^{  -\f{9}{4}          }       d\tau                 \\
     \leq     & \   CA^2 t^{  -\f{9}{8}  -\f{  1+  |\alpha|}{2}           }
     \end{split}
\end{equation}
for $t \geq 1$ and
\begin{equation}\nonumber
     \begin{split}
    \Big\|   \Grad^{\alpha} &\   \int_0^t e^{ (t-\tau)\Delta  } \nablah \p_k \p_l N(t-\tau) \ast (B_k  B_l  ) (\tau  )        d \tau \Big\|_{L^{\infty}}  \\
    \leq  & \  C  \int_0^t  \|  \Grad^{\alpha} \p_k \p_l N(t-\tau)  \|_{L^{ \f{4}{3} }}   \| \nablah  (B_k  B_l ) (\tau)     \| _{L^4}        d \tau                     \\
  \leq  & \  CA^2          \int_0^  {  t }  (t-\tau)^{ -\f{3}{8} -\f{2+ |\alpha|}{2} +1        }  (1+\tau) ^{  -\f{9}{4}          }    d \tau                  \\
     \leq     & \   CA^2
     \end{split}
\end{equation}
for $0<t \leq 1$. The estimate of $\| \Grad^{\alpha}  \Nhe_3 [B](t)  \|_{L^p} $ is thus verified from above.

Next, for any $k\in \{ 1,2 \}$ we have
\begin{equation}\nonumber
     \begin{split}
         \|\nablah (B_k B_3  ) (\tau  )  \|_{L^1} \leq  & \  \| \nablah B_k (\tau)  \|_{L^{\infty}} \| B_3 (\tau)\|_{L^1}        + \| \nablah B_3 (\tau)  \|_{L^{\infty}} \| B_k (\tau)\|_{L^1}              \\
            \leq  & \  CA^2  (1+\tau) ^{  -\f{3}{2}          }    , \\
      \| \nablah (B_k B_3  ) (\tau  )   \|_{L^4} \leq & \  \| \nablah B_k (\tau)  \|_{L^{\infty}} \| B_3 (\tau)\|_{L^4}    +  \| B_k (\tau)\|_{L^4}  \| \nablah B_3 (\tau)  \|_{L^{\infty}}      \\
        \leq  & \  CA^2  (1+\tau) ^{  -\f{3}{2}  -\f{9}{8}         }    ; \\
     \end{split}
\end{equation}
whence
\begin{equation}\nonumber
     \begin{split}
     \Big\|   \Grad^{\alpha}   &\ \int_0^t e^{ (t-\tau)\Delta  } \nablah \p_k \p_3 N(t-\tau) \ast (B_k  B_3  ) (\tau  )     d \tau \Big\|_{L^1}  \\
     \leq  &\     \int_0^  {  \f{t}{2} }        \|  \Grad^{\alpha} \nablah \p_k \p_3  N(t-\tau)  \|_{L^1}   \| (B_k  B_3 )(\tau)     \| _{L^1}    d\tau         \\
    &\  +   \int_{  \f{t}{2} } ^t    \|  \Grad^{\alpha} \p_k \p_3 N(t-\tau)  \|_{L^1}   \| \nablah  (B_k  B_3 ) (\tau)     \| _{L^1}       d\tau                    \\
  \leq  & \  CA^2          \int_0^  {  \f{t}{2} } (t-\tau)^{  -\f{3+ |\alpha|}{2} +1        } (1+\tau) ^{  -\f{3}{2}      }    d \tau                 \\
         & \ +     CA^2          \int_{  \f{t}{2} } ^t (t-\tau)^{  -\f{ |\alpha|}{2}        }  (1+\tau) ^{  -\f{3}{2}          }       d\tau                 \\
     \leq     & \   CA^2 t^{    -\f{  1+  |\alpha|}{2}           }
     \end{split}
\end{equation}
for $t \geq 1$ and
\begin{equation}\nonumber
     \begin{split}
    \Big\|   \Grad^{\alpha} &\   \int_0^t e^{ (t-\tau)\Delta  } \nablah \p_k \p_3 N(t-\tau) \ast (B_k  B_3  ) (\tau  )        d \tau \Big\|_{L^1}  \\
    \leq  & \  C  \int_0^t  \|  \Grad^{\alpha} \p_k \p_3  N(t-\tau)  \|_{L^1}   \| \nablah  (B_k  B_3) (\tau)     \| _{L^1}        d \tau                     \\
  \leq  & \  CA^2          \int_0^  {  t }   (t-\tau)^{  -\f{ |\alpha|}{2}        }  (1+\tau) ^{  -\f{3}{2}          }     d \tau                  \\
     \leq     & \   CA^2
     \end{split}
\end{equation}
for $0<t \leq 1$.
\begin{equation}\nonumber
     \begin{split}
     \Big\|   \Grad^{\alpha}   &\ \int_0^t e^{ (t-\tau)\Delta  } \nablah \p_k \p_3 N(t-\tau) \ast (B_k  B_3  ) (\tau  )     d \tau \Big\|_{L^{\infty}}    \\
     \leq  &\     \int_0^  {  \f{t}{2} }        \|  \Grad^{\alpha} \nablah \p_k \p_3  N(t-\tau)  \|_{L^{\infty}}    \| (B_k  B_3 )(\tau)     \| _{L^1}    d\tau         \\
    &\  +   \int_{  \f{t}{2} } ^t    \|  \Grad^{\alpha} \p_k \p_3 N(t-\tau)  \|_{L^{ \f{4}{3} }}   \| \nablah  (B_k  B_3 ) (\tau)     \| _{L^4}       d\tau                    \\
  \leq  & \  CA^2          \int_0^  {  \f{t}{2} } (t-\tau)^{  -\f{3}{2}-\f{3+ |\alpha|}{2} +1        } (1+\tau) ^{  -\f{3}{2}      }    d \tau                 \\
         & \ +     CA^2          \int_{  \f{t}{2} } ^t (t-\tau)^{ -\f{3}{8} -\f{  |\alpha|}{2}        }  (1+\tau) ^{  -\f{3}{2}    -\f{9}{8}        }       d\tau                 \\
     \leq     & \   CA^2 t^{    -\f{3}{2}  -\f{  1+  |\alpha|}{2}           }
     \end{split}
\end{equation}
for $t \geq 1$ and
\begin{equation}\nonumber
     \begin{split}
    \Big\|   \Grad^{\alpha} &\   \int_0^t e^{ (t-\tau)\Delta  } \nablah \p_k \p_3 N(t-\tau) \ast (B_k  B_3  ) (\tau  )        d \tau \Big\|_{L^1}  \\
    \leq  & \  C  \int_0^t  \|  \Grad^{\alpha} \p_k \p_3  N(t-\tau)  \|_{L^{ \f{4}{3} }}   \| \nablah  (B_k  B_3) (\tau)     \| _{L^4}        d \tau                     \\
  \leq  & \  CA^2          \int_0^  {  t }  (t-\tau)^{ -\f{3}{8} -\f{  |\alpha|}{2}        }  (1+\tau) ^{  -\f{3}{2}    -\f{9}{8}        }     d \tau                  \\
     \leq     & \   CA^2
     \end{split}
\end{equation}
for $0<t \leq 1$. This gives the estimate of $\| \Grad^{\alpha}  \Nhe_4 [B](t)  \|_{L^p} $.

Following the same line, one gets the estimate of $\| \Grad^{\alpha}\Nhe_5 [B](t)  \|_{L^p}$ and that of $\| \Grad^{\alpha}  \Nve_j [B](t)\|_{L^p}$ in (\ref{hm17}) with $j\in \{ 1,2,3,4,5  \}$. The details are omitted.

It remains to show the estimates for $\| \Grad^{\alpha}  \Nhe_j [u](t)\|_{L^p}$ and $\| \Grad^{\alpha}  \Nve_j [u](t)\|_{L^p}$ with $j\in \{ 1,2,3,4,5  \}$. For brevity, we shall present $\| \Grad^{\alpha}  \Nhe_2 [u](t)\|_{L^p}$ for instance.
\begin{equation}\nonumber
     \begin{split}
         \| (\uh \otimes \uh ) (\tau  )  \|_{L^1} \leq  & \  \| \uh (\tau)  \|_{L^{\infty}} \| \uh(\tau)\|_{L^1}  \leq CA^2  (1+\tau) ^{  -\f{17}{8}         } ( \log(2+\tau) )^2,    \\
      \| \nablah \cdot (\uh \otimes \uh ) (\tau  )  \|_{L^1} \leq & \ C  \| \nablah \uh(\tau)  \|_{L^{\infty}} \| \uh(\tau)\|_{L^1}
        \leq    CA^2  (1+\tau) ^{  -\f{21}{8}          }   ( \log(2+\tau) )^2  , \\
        \| \nablah \cdot (\uh \otimes \uh ) (\tau  )  \|_{L^4} \leq & \ C  \| \nablah \uh(\tau)  \|_{L^{\infty}} \| \uh(\tau)\|_{L^4}    \leq CA^2    (1+\tau) ^{  -\f{111}{32}          }   ( \log(2+\tau) )^2     ;
     \end{split}
\end{equation}
whence
\begin{equation}\nonumber
     \begin{split}
     \Big\|   \Grad^{\alpha} &\  \int_0^t e^{ (t-\tau)\Delta  }\nablah \cdot (\uh \otimes \uh ) (\tau  )     d \tau         \Big\|_{L^1}         \\
     \leq  &\     \int_0^  {  \f{t}{2} }        \|  \Grad^{\alpha}\nablah  G(t-\tau)  \|_{L^1}   \| (\uh \otimes \uh)(\tau)     \| _{L^1}    d\tau              \\
     &\ +   \int_{  \f{t}{2} } ^t    \|  \Grad^{\alpha}  G(t-\tau)  \|_{L^1}   \| \nablah \cdot (\uh \otimes \uh ) (\tau)     \| _{L^1}       d\tau                    \\
  \leq  & \  CA^2          \int_0^  {  \f{t}{2} } (t-\tau)^{  -\f{1+ |\alpha|}{2}         }  (1+\tau) ^{  -\f{17}{8}         } ( \log(2+\tau) )^2   d \tau                 \\
         & \ +     CA^2          \int_{  \f{t}{2} } ^t (t-\tau)^{  -\f{|\alpha|}{2}         }  (1+\tau) ^{  -\f{21}{8}          }   ( \log(2+\tau) )^2       d\tau                 \\
     \leq     & \   CA^2 t^{    -\f{  1+  |\alpha|}{2}           }
     \end{split}
\end{equation}
for $t \geq 1$ and
\begin{equation}\nonumber
     \begin{split}
    \left\|   \Grad^{\alpha} \int_0^t e^{ (t-\tau)\Delta  }\nablah \cdot (\uh \otimes \uh )  (\tau  )     d \tau \right\|_{L^1}    \leq  & \  C  \int_0^t   \|  \Grad^{\alpha}  G(t-\tau)  \|_{L^1}   \| \nablah \cdot (\uh \otimes \uh ) (\tau)     \| _{L^1}          d \tau                     \\
  \leq  & \  CA^2          \int_0^  {  t } (t-\tau)^{  -\f{|\alpha|}{2}         } (1+\tau) ^{  -\f{21}{8}          }   ( \log(2+\tau) )^2   d \tau                  \\
     \leq     & \   CA^2
     \end{split}
\end{equation}
for $0<t \leq 1$.
\begin{equation}\nonumber
     \begin{split}
     \left\|   \Grad^{\alpha} \int_0^t e^{ (t-\tau)\Delta  }\nablah \cdot (\uh \otimes \uh ) (\tau  )     d \tau \right\|_{L^{\infty}}   \leq  &\     \int_0^  {  \f{t}{2} }        \|  \Grad^{\alpha} \nablah  G(t-\tau)  \|_{L^{\infty}}  \| (\uh \otimes \uh)(\tau)     \| _{L^1}    d\tau              \\
     &\ +   \int_{  \f{t}{2} } ^t    \|  \Grad^{\alpha}  G(t-\tau)  \|_{L^{ \f{4}{3}} }   \|  \nablah \cdot (\uh \otimes \uh ) (\tau)     \| _{L^4}       d\tau                    \\
  \leq  & \  CA^2          \int_0^  {  \f{t}{2} } (t-\tau)^{  -\f{3}{2}-\f{1+ |\alpha|}{2}         }  (1+\tau) ^{  -\f{17}{8}         } ( \log(2+\tau) )^2    d \tau                 \\
         & \ +     CA^2          \int_{  \f{t}{2} } ^t (t-\tau)^{  -\f{3 }{8} -\f{ |\alpha|}{2}            }   (1+\tau) ^{  -\f{111}{32}          }   ( \log(2+\tau) )^2      d\tau                 \\
     \leq     & \   CA^2 t^{  -\f{3}{2}-\f{1+ |\alpha|}{2}         }
     \end{split}
\end{equation}
for $t \geq 1$ and
\begin{equation}\nonumber
     \begin{split}
     \left\|    \Grad^{\alpha}   \int_0^t e^{ (t-\tau)\Deltah  } \nablah \cdot (\uh \otimes \uh ) (\tau  )     d \tau \right\|_{L^{\infty}}      \leq  & \  C  \int_0^t   \|  \Grad^{\alpha}  G(t-\tau)  \|_{L^{ \f{4}{3}} }   \|\nablah \cdot (\uh \otimes \uh ) (\tau)     \| _{L^4}          d \tau                     \\
     \leq     & \   CA^2    \int_0^  {  t }   (t-\tau)^{  -\f{3 }{8} -\f{ |\alpha|}{2}            }  (1+\tau) ^{  -\f{111}{32}          }   ( \log(2+\tau) )^2   d \tau   \\
     \leq &\ CA^2
     \end{split}
\end{equation}
for $0<t \leq 1$. The estimate of $\| \Grad^{\alpha}  \Nhe_2 [u](t)  \|_{L^p} $ follows immediately.

The proof of Lemma \ref{lem13} is completely finished.                  $\Box$

In the next lemma, we establish the decay estimates for the Duhamel terms in $L^{\infty}_{\rm h}  L^{1}_{\rm v}$-norm.
\begin{lemm}\label{lem14}
Let $(u,B)$ be subject to the ansatzes $(i),(ii),(iii)$ with $s\geq 9$. There exists a generic constant $C>0$ such that for any $0<t<T$ and any $\alpha=(\alphah,\alpha_3)\in (\mathbb{N} \cup \{ 0\})^2 \times (\mathbb{N} \cup \{ 0\})$ with $|\alpha|\leq 1$ we have the decay estimates for the Duhamel terms associated with velocity equation:
\begin{equation}\label{hm19}
     \begin{split}
         \|\Grad^{\alpha}  \Nhe_m[u](t)  \|_{L^{\infty}_{\rm h}  L^{1}_{\rm v}} & \  \leq  CA^2 t^{ -1-\f{1+|\alpha|}{2}  } ,        \\
   \|\Grad^{\alpha}  \Nve_m[u](t)  \|_{L^{\infty}_{\rm h}  L^{1}_{\rm v}} & \  \leq  CA^2 t^{ -1-\f{1+|\alpha|}{2}  }  ,       \\
    \|\Grad^{\alpha}  \Nhe_1[B](t)  \|_{L^{\infty}_{\rm h}  L^{1}_{\rm v}} & \  \leq  CA^2 t^{ -1-\f{1+|\alpha|}{2}  },        \\
    \|\Grad^{\alpha}  \Nhe_2[B](t)  \|_{L^{\infty}_{\rm h}  L^{1}_{\rm v}} & \  \leq  CA^2 t^{ -1-\f{1+|\alpha|}{2}  } \log(2+t) ,        \\
    \|\Grad^{\alpha}  \Nhe_3[B](t)  \|_{L^{\infty}_{\rm h}  L^{1}_{\rm v}} & \  \leq  CA^2 t^{ -1-\f{1+|\alpha|}{2}  } \log(2+t) ,        \\
     \|\Grad^{\alpha}  \Nhe_4[B](t)  \|_{L^{\infty}_{\rm h}  L^{1}_{\rm v}} & \  \leq  CA^2 t^{ -1-\f{1+|\alpha|}{2}  },        \\
      \|\Grad^{\alpha}  \Nhe_5[B](t)  \|_{L^{\infty}_{\rm h}  L^{1}_{\rm v}} & \  \leq  CA^2 t^{ -1-\f{1+|\alpha|}{2}  },        \\
   \|\Grad^{\alpha}  \Nve_1[B](t)  \|_{L^{\infty}_{\rm h}  L^{1}_{\rm v}} & \  \leq  CA^2 t^{ -1-\f{1+|\alpha|}{2}  } ,       \\
      \|\Grad^{\alpha}  \Nve_2[B](t)  \|_{L^{\infty}_{\rm h}  L^{1}_{\rm v}} & \  \leq  CA^2 t^{ -1-\f{1+|\alpha|}{2}  } ,       \\
         \|\Grad^{\alpha}  \Nve_3[B](t)  \|_{L^{\infty}_{\rm h}  L^{1}_{\rm v}} & \  \leq  CA^2 t^{ -1-\f{1+|\alpha|}{2}  } \log(2+t) ,       \\
            \|\Grad^{\alpha}  \Nve_4[B](t)  \|_{L^{\infty}_{\rm h}  L^{1}_{\rm v}} & \  \leq  CA^2 t^{ -1-\f{1+|\alpha|}{2}  } ,       \\
               \|\Grad^{\alpha}  \Nve_5[B](t)  \|_{L^{\infty}_{\rm h}  L^{1}_{\rm v}} & \  \leq  CA^2 t^{ -1-\f{1+|\alpha|}{2}  } ,       \\
     \end{split}
\end{equation}
for any $m\in \{ 1,2,3,4,5 \}$ and the decay estimates for the Duhamel terms associated with magnetic equation:
\begin{equation}\label{hm20}
     \begin{split}
 \|\Grad^{\alpha}  \Nhd_m[u,B](t)  \|_{L^{\infty}_{\rm h}  L^{1}_{\rm v}} & \  \leq  CA^2 t^{ -1-\f{|\alphah|}{2}  } ,       \\
 \|\Grad^{\alpha}  \Nvd_n[u,B](t)  \|_{L^{\infty}_{\rm h}  L^{1}_{\rm v}} & \  \leq  CA^2 t^{ -1-\f{|\alphah|}{2}  } ,       \\
     \end{split}
\end{equation}
for any $m\in \{ 1,2,3,4 \}, n \in \{ 1,2 \}$.
\end{lemm}
{\it Proof.} We begin with the proof of (\ref{hm19}). Observing that
\begin{equation}\nonumber
     \begin{split}
        \| \p_3(B_3 \Bh)(\tau) \|_{L^{\infty}_{\rm h}  L^{1}_{\rm v}}  & \  \leq  \| \p_3 B_3(\tau ) \|_{L^{\infty}}      \| \Bh (\tau )\|_{L^{\infty}_{\rm h}  L^{1}_{\rm v}}      +    \|  B_3(\tau) \|_{L^{\infty}}    \| \p_3  \Bh (\tau) \|     _{L^{\infty}_{\rm h}  L^{1}_{\rm v}}          \\
    & \  \leq  CA^2   (1+\tau)^{-\f{3}{2}} \tau^{-1}     ,        \\
     \end{split}
\end{equation}
yielding
\begin{equation}\nonumber
     \begin{split}
     \Big\|      \nabla^{\alpha}   &\  \int_0^{ t } e^{ (t-\tau)\Delta  }\p_3 (B_3 \Bh ) (\tau  )     d \tau \Big\|_{L^{\infty}_{\rm h}  L^{1}_{\rm v}}   \\
     \leq  & \    \int_0^{\f{t}{2}}  \| \p_3 \nabla^{\alpha} G (t-\tau) \|_{L^{\infty}_{\rm h}  L^{1}_{\rm v}} \|  (B_3 \Bh) (\tau) \|     _{L^1}         d \tau                     \\
     &\ + \int_{  \f{t}{2} } ^t  \| \nabla^{\alpha} G (t-\tau) \|_{L^{1}}  \| \p_3 (B_3 \Bh) (\tau) \|  _{L^{\infty}_{\rm h}  L^{1}_{\rm v}}  d \tau          \\
  \leq  & \  CA^2          \int_0^  {  \f{t}{2} } (t-\tau)^{  -1-\f{1+ |\alpha|}{2}         } (1+\tau) ^{  -\f{3}{2} } d \tau                 \\
         & \ +     CA^2          \int_{  \f{t}{2} } ^t (t-\tau)^{  -\f{|\alpha|}{2}         }    (1+\tau)^{-\f{3}{2}} \tau^{-1}    d\tau                 \\
     \leq     & \   CA^2 t^{  -1  -\f{1+|\alpha|}{2}           }. \\
     \end{split}
\end{equation}
Similarly,
\[
 \| \nablah (\Bh \otimes \Bh ) (\tau  ) \|_{L^{\infty}_{\rm h}  L^{1}_{\rm v}}
 \leq \| \nablah \Bh (\tau)  \| _{L^{\infty}}       \|\Bh(\tau)  \|_{L^{\infty}_{\rm h}  L^{1}_{\rm v}}
 \leq  CA^2  (1+\tau)^{-\f{3}{2}} \tau^{-1}  ;
\]
whence
\begin{equation}\nonumber
     \begin{split}
     \Big\|      \nabla^{\alpha}   &\  \int_0^{ t } e^{ (t-\tau)\Delta  }\nablah \cdot (\Bh \otimes \Bh ) (\tau  )     d \tau \Big\|_{L^{\infty}_{\rm h}  L^{1}_{\rm v}}   \\
     \leq  & \    \int_0^{\f{t}{2}}  \| \nablah \nabla^{\alpha} G (t-\tau) \|_{L^{\infty}_{\rm h}  L^{1}_{\rm v}} \|  (\Bh \otimes \Bh) (\tau) \|     _{L^1}         d \tau                     \\
     &\ + \int_{  \f{t}{2} } ^t  \| \nabla^{\alpha} G (t-\tau) \|_{L^{1}}  \| \nablah \cdot (\Bh \otimes \Bh ) (\tau  )  \|  _{L^{\infty}_{\rm h}  L^{1}_{\rm v}}  d \tau          \\
  \leq  & \  CA^2          \int_0^  {  \f{t}{2} } (t-\tau)^{  -1-\f{1+ |\alpha|}{2}         } (1+\tau) ^{  -1} d \tau                 \\
         & \ +     CA^2          \int_{  \f{t}{2} } ^t (t-\tau)^{  -\f{|\alpha|}{2}         }    (1+\tau)^{-\f{3}{2}} \tau^{-1}    d\tau                 \\
     \leq     & \   CA^2 t^{  -1  -\f{1+|\alpha|}{2}           } \log(2+t). \\
     \end{split}
\end{equation}
Following the same line as above, one gets the estimates of $ \| \Grad^{\alpha} \Nhe_i[B](t)\|  _{L^{\infty}_{\rm h}  L^{1}_{\rm v}} $ for any $i\in \{ 3,4,5 \}$ and $ \| \Grad^{\alpha} \Nve_j[B](t)\|  _{L^{\infty}_{\rm h}  L^{1}_{\rm v}} $ for any $j\in \{1,2,3,4,5  \}$.

Next, we see that
\begin{equation}\nonumber
     \begin{split}
        \| \p_3(u_3 \uh)(\tau) \|_{L^{\infty}_{\rm h}  L^{1}_{\rm v}} \leq  & \    \| \p_3 u_3(\tau ) \|_{L^{\infty}}      \| \uh (\tau )\|_{L^{\infty}_{\rm h}  L^{1}_{\rm v}}      +    \|  u_3(\tau) \|_{L^{\infty}}    \| \p_3  \uh (\tau) \|     _{L^{\infty}_{\rm h}  L^{1}_{\rm v}}          \\
  \leq  & \    CA^2  \Big( (1+\tau)^{-\f{17}{8}} \tau^{-\f{3}{2}}   (\log(2+\tau))^2+
  (1+\tau)^{-\f{13}{8}} \tau^{-2}   (\log(2+\tau))^2
  \Big),        \\
     \end{split}
\end{equation}
implying
\begin{equation}\nonumber
     \begin{split}
     \Big\|      \nabla^{\alpha}   &\  \int_0^{ t } e^{ (t-\tau)\Delta  }\p_3 (u_3 \uh ) (\tau  )     d \tau \Big\|_{L^{\infty}_{\rm h}  L^{1}_{\rm v}}   \\
     \leq  & \    \int_0^{\f{t}{2}}  \| \p_3 \nabla^{\alpha} G (t-\tau) \|_{L^{\infty}_{\rm h}  L^{1}_{\rm v}} \|  (u_3 \uh) (\tau) \|     _{L^1}         d \tau                     \\
     &\ + \int_{  \f{t}{2} } ^t  \| \nabla^{\alpha} G (t-\tau) \|_{L^{1}}  \| \p_3 (u_3 \uh) (\tau) \|  _{L^{\infty}_{\rm h}  L^{1}_{\rm v}}  d \tau          \\
  \leq  & \  CA^2          \int_0^  {  \f{t}{2} } (t-\tau)^{  -1-\f{1+ |\alpha|}{2}         }  (1+\tau)^{-\f{17}{8}}  (\log(2+\tau))^2 d \tau                 \\
         & \ +     CA^2          \int_{  \f{t}{2} } ^t (t-\tau)^{  -\f{|\alpha|}{2}         }   \Big( (1+\tau)^{-\f{17}{8}} \tau^{-\f{3}{2}}   (\log(2+\tau))^2+
  (1+\tau)^{-\f{13}{8}} \tau^{-2}   (\log(2+\tau))^2
  \Big)   d\tau                 \\
     \leq     & \   CA^2 t^{  -1  -\f{1+|\alpha|}{2}           }. \\
     \end{split}
\end{equation}
Analogously,
\begin{equation}\nonumber
     \begin{split}
     \Big\|      \nabla^{\alpha}   &\  \int_0^{ t } e^{ (t-\tau)\Delta  }\nablah \cdot (\uh \otimes \uh ) (\tau  )     d \tau \Big\|_{L^{\infty}_{\rm h}  L^{1}_{\rm v}}   \\
     \leq  & \    \int_0^{\f{t}{2}}  \| \nablah \nabla^{\alpha} G (t-\tau) \|_{L^{\infty}_{\rm h}  L^{1}_{\rm v}} \|  (\uh \otimes \uh) (\tau) \|     _{L^1}         d \tau                     \\
     &\ + \int_{  \f{t}{2} } ^t  \| \nabla^{\alpha} G (t-\tau) \|_{L^{1}}  \| \nablah \cdot (\uh \otimes \uh ) (\tau  )  \|  _{L^{\infty}_{\rm h}  L^{1}_{\rm v}}  d \tau          \\
  \leq  & \  CA^2          \int_0^  {  \f{t}{2} } (t-\tau)^{  -1-\f{1+ |\alpha|}{2}         }   (1+\tau) ^{  -\f{17}{8}         } ( \log(2+\tau) )^2  d \tau                 \\
         & \ +     CA^2          \int_{  \f{t}{2} } ^t (t-\tau)^{  -\f{|\alpha|}{2}         }    (1+\tau)^{-\f{17}{8}} \tau^{-\f{3}{2}}     ( \log(2+\tau) )^2    d\tau                 \\
     \leq     & \   CA^2 t^{  -1  -\f{1+|\alpha|}{2}           } . \\
     \end{split}
\end{equation}
The remaining estimates in (\ref{hm19}) involving $u$ follow as above. The details are omitted.

To finish the proof, it remains to verify (\ref{hm20}). For simplicity, we only provide some typical ones.

If $\alpha_3=0$,
\begin{equation}\nonumber
     \begin{split}
        \| \p_3(u_3 \Bh) \|_{L^{\infty}_{\rm h}  L^{1}_{\rm v}}  & \  \leq  \| \p_3 u_3 \|_{L^{\infty}}    \| \Bh \|_{L^{\infty}_{\rm h}  L^{1}_{\rm v}}   +    \|  u_3 \|_{L^{\infty}}    \| \p_3  \Bh \|_{L^{\infty}_{\rm h}  L^{1}_{\rm v}}        \\
    & \  \leq  CA^2 \left(  (1+\tau)^{-\f{17}{8}}   \tau^{-1} +  (1+\tau)^{-\f{13}{8}} \tau^{-1 }      \right)  \log (2+\tau);        \\
     \end{split}
\end{equation}
whence
\begin{equation}\nonumber
     \begin{split}
     \Big \|   \nablah^{\alphah}  &\  \int_0^{ t } e^{ (t-\tau)\Deltah  }\p_3 (u_3\Bh ) (\tau  )     d \tau \Big\|_{L^{\infty}_{\rm h}  L^{1}_{\rm v}} \\
     \leq  & \    \int_0^{\f{t}{2}}  \| \nablah^{\alphah} \Gh (t-\tau) \|_{L^{\infty}(\R^2)} \| \p_3 (u_3 \Bh) (\tau) \|     _{L^1}         d \tau                     \\
     &\ + \int_{  \f{t}{2} } ^t  \| \nablah^{\alphah} \Gh (t-\tau) \|_{L^{1}(\R^2)}  \| \p_3 (u_3 \Bh) (\tau) \|  _{L^{\infty}_{\rm h}  L^{1}_{\rm v}}  d \tau         \\
  \leq  & \  CA^2          \int_0^  {  \f{t}{2} } (t-\tau)^{  -1-\f{|\alphah|}{2}         } (1+\tau) ^{  -\f{13}{8} }    \log(2+\tau)   d \tau                 \\
         & \ +     CA^2          \int_{  \f{t}{2} } ^t (t-\tau)^{  -\f{|\alphah|}{2}         } \left(  (1+\tau)^{-\f{17}{8}}   \tau^{-1} +  (1+\tau)^{-\f{13}{8}} \tau^{-1 }      \right)  \log (2+\tau)   d\tau                 \\
     \leq     & \   CA^2 t^{  -1  -\f{|\alphah|}{2}           }. \\
     \end{split}
\end{equation}
If $\alpha_3=1$,
\begin{equation}\nonumber
     \begin{split}
        \| \p_3^2(u_3 \Bh) \|_{L^{\infty}_{\rm h}  L^{1}_{\rm v}} \leq  & \    \| \p_3^2 u_3 \|_{L^{\infty}} \|\Bh \| _{L^{\infty}_{\rm h}  L^{1}_{\rm v}}   +2 \| \p_3 u_3 \|_{L^{\infty}}   \| \p_3 \Bh \|_{L^{\infty}_{\rm h}  L^{1}_{\rm v}}          \\
        &\ +  \| \p_3^2 \Bh \|_{L^{\infty}}   \| u_3 \| _{L^{\infty}_{\rm h}  L^{1}_{\rm v}}  \\
       \leq  &\  C \Big( \| u_3 \|_{H^7} ^{\f{2}{7}}    \|\p_3 u_3 \|_{L^{\infty}} ^{\f{5}{7}}
       \|\Bh \| _{L^{\infty}_{\rm h}  L^{1}_{\rm v}}  + \| \p_3 u_3 \|_{L^{\infty}}   \| \p_3 \Bh \|_{L^{\infty}_{\rm h}  L^{1}_{\rm v}}  \\
       &\ +   \| \Bh \|_{H^7} ^{\f{2}{7}}    \|\p_3 \Bh \|_{L^{\infty}} ^{\f{5}{7}}
       \|u_3 \| _{L^{\infty}_{\rm h}  L^{1}_{\rm v}}   \Big)      \\
  \leq  & \    CA^2 \Big(  (1+\tau)^{ - \f{85}{56}   } \tau^{-1}  (\log(2+\tau))^{\f{5}{7}    }   +  (1+\tau)^{-\f{17}{8}} \tau^{-1 } \log(2+\tau) \\
  &\  +   (1+\tau)^{ - \f{5}{7}   }  \tau^{- \f{3}{2} }    \log(2+\tau)   \Big)    ;        \\
     \end{split}
\end{equation}
whence
\begin{equation}\nonumber
     \begin{split}
     \Big\|  \p_3  &\  \int_0^{ t } e^{ (t-\tau)\Deltah  }\p_3 (u_3\Bh ) (\tau  )     d \tau \Big\|_{L^{\infty}_{\rm h}  L^{1}_{\rm v}}  \\
     \leq  & \    \int_0^{\f{t}{2}}  \|  \Gh (t-\tau) \|_{L^{\infty}(\R^2)} \| \p_3^2 (u_3 \Bh) (\tau) \|     _{L^1}         d \tau                     \\
     &\ + \int_{  \f{t}{2} } ^t  \| \Gh (t-\tau) \|_{L^{1}(\R^2)}  \| \p_3^2 (u_3 \Bh) (\tau) \|  _{L^{\infty}_{\rm h}  L^{1}_{\rm v}}  d \tau         \\
  \leq  & \  CA^2          \int_0^  {  \f{t}{2} } (t-\tau)^{  -1     } (1+\tau) ^{  -\f{17}{16} }   \log(2+\tau)    d \tau                 \\
         & \ +     CA^2          \int_{  \f{t}{2} } ^t
       \Big(  (1+\tau)^{ - \f{85}{56}   } \tau^{-1}  (\log(2+\tau))^{\f{5}{7}    }   +  (1+\tau)^{-\f{17}{8}} \tau^{-1 } \log(2+\tau)                  \\
        &\  \qquad  \qquad  \qquad   +   (1+\tau)^{ - \f{5}{7}   }  \tau^{- \f{3}{2} }    \log(2+\tau)   \Big) d\tau    \\
     \leq     & \   CA^2 t^{  -1           }. \\
     \end{split}
\end{equation}
The estimate of $ \|\Grad^{\alpha} \Nhd_3[u,B](t)  \|_{L^{\infty}_{\rm h}  L^{1}_{\rm v}} $ follows immediately.

Seeing that
\begin{equation}\nonumber
     \begin{split}
        \| \Grad^{\alpha} (\uh \otimes \Bh)(\tau) \|_{L^{\infty}_{\rm h}  L^{1}_{\rm v}}  \leq  & \    \|  \nabla ^{\alpha} \uh (\tau)  \| _{L^{\infty}}  \| \Bh \|   _{L^{\infty}_{\rm h}  L^{1}_{\rm v}  }  +  \| \uh (\tau) \|_{L^{\infty}} \|  \Grad^{\alpha} \Bh (\tau)\| _{L^{\infty}_{\rm h}  L^{1}_{\rm v}}     \\
    \leq   & \   CA^2  \Big(     \tau^{-1}   (1+\tau)^{-\f{9}{8}-\f{1+|\alpha|}{2} }   +  \tau^{  -1-\f{|\alphah|}{2}      }  (1+\tau)^{-\f{13}{8} }  \Big)    \log(2+\tau)        ,         \\
     \end{split}
\end{equation}
\begin{equation}\nonumber
     \begin{split}
        \| \p_3^{\alpha_3} ( \uh \otimes \Bh )(\tau) \|_{L^{2}_{\rm h}  L^{1}_{\rm v}}  \leq  & \    \|  \p_3^{\alpha_3} (\uh \otimes \Bh )(\tau) \| _{ L^{1}} ^{ \f{1}{2}  }    \|  \p_3^{\alpha_3} (\uh \otimes \Bh )(\tau) \|_{L^{\infty}_{\rm h}  L^{1}_{\rm v}} ^{ \f{1}{2}  }      \\
    \leq & \   CA^2   \Big(     \tau^{-\f{1}{2}} (1+\tau)^{-\f{30}{16}}    + \tau^{-\f{1}{2}} (1+\tau)^{-\f{26}{16}}          \Big)  \log(2+\tau)       ;        \\
     \end{split}
\end{equation}
whence
\begin{equation}\nonumber
     \begin{split}
     \Big\|   \nabla^{\alpha}  &\   \int_0^{ t } e^{ (t-\tau)\Deltah  } \nablah \cdot (\uh \otimes \Bh  ) (\tau  )     d \tau \Big\|_{L^{\infty}_{\rm h}  L^{1}_{\rm v}}  \\
     \leq  & \    \int_0^{\f{t}{2}}  \| \nablah^{\alphah} \nablah \Gh (t-\tau) \|_{L^{2}(\R^2)} \| \p_3 ^{\alpha_3}(\uh \otimes \Bh ) (\tau) \|   _{L^{2}_{\rm h}  L^{1}_{\rm v}}        d \tau                     \\
     &\ + \int_{  \f{t}{2} } ^t  \| \nablah^{\alphah} \Gh (t-\tau) \|_{L^{1}(\R^2)}  \| \nabla^{\alpha} (\uh \otimes \Bh ) (\tau) \|  _{L^{\infty}_{\rm h}  L^{1}_{\rm v}}  d \tau         \\
  \leq  & \  CA^2          \int_0^  {  \f{t}{2} } (t-\tau)^{  -1-\f{|\alphah|}{2}         } \Big(     \tau^{-\f{1}{2}} (1+\tau)^{-\f{30}{16}}    + \tau^{-\f{1}{2}} (1+\tau)^{-\f{26}{16}}          \Big)  \log(2+\tau)    d \tau                  \\
         & \ +     CA^2          \int_{  \f{t}{2} } ^t (t-\tau)^{  -\f{1}{2}         }  \Big(     \tau^{-1}   (1+\tau)^{-\f{9}{8}-\f{1+|\alpha|}{2} }   +  \tau^{  -1-\f{|\alphah|}{2}      }  (1+\tau)^{-\f{13}{8} }  \Big)    \log(2+\tau)     d\tau                 \\
     \leq     & \   CA^2 t^{  -1  -\f{|\alphah|}{2}           }. \\
     \end{split}
\end{equation}
This gives the estimates of $ \|\Grad^{\alpha} \Nhd_i[u,B](t)  \|_{L^{\infty}_{\rm h}  L^{1}_{\rm v}} $ with $i=1,2$. The remaining estimates in (\ref{hm20}) are verified similarly.                      $\Box$

\subsection{Asymptotic expansions for the Duhamel terms}

The aim of this subsection is to study the asymptotic expansion of the slowly decaying Duhamel terms $\Nhd_3[u,B]$ and $\Nhd_4[u,B]$.

\begin{lemm}\label{lem15}
Let $(u,B)$ be subject to the ansatzes $(i),(ii)$ with $s\geq 3,A>0,T=\infty$. Then there hold the asymptotic limits for any $1\leq p \leq \infty$:
\beq\label{hm21}
\lim_{t\rightarrow \infty} t^{1-\f{1}{p}}
\left\|
\Nhd_3[u,B] (t,x)+ G_{\rm h}(t,\xh)  \int_0^{\infty} \int_{\R^2}  \p_3 (u_3 \Bh) (\tau, y_{\rm h},x_3)
    d y_{\rm h}   d \tau
\right\|_{L^p_x} =0,
\eeq
\beq\label{hm22}
\lim_{t\rightarrow \infty} t^{1-\f{1}{p}}
\left\|
\Nhd_4[u,B] (t,x)-  G_{\rm h}(t,\xh)  \int_0^{\infty} \int_{\R^2}  \p_3 (B_3 \uh) (\tau, y_{\rm h},x_3)
    d y_{\rm h}   d \tau
\right\|_{L^p_x} =0.
\eeq
\end{lemm}
{\it Proof.} We start the proof with (\ref{hm21}). Following \cites{FM,F}, we make the following decomposition
\[
\Nhd_3[u,B] (t,x)+ G_{\rm h}(t,\xh)  \int_0^{\infty} \int_{\R^2}  \p_3 (u_3 \Bh) (\tau, y_{\rm h},x_3)
    d y_{\rm h}   d \tau=-
\sum_{m=1}^4 \mathcal{I}_m (t,x),
\]
with
\begin{equation}\nonumber
     \begin{split}
        \mathcal{I}_1 (t,x) =& \  \int_0^{\f{t}{2}}   \int_{\R^2}
   \Big( \Gh(t-\tau,\xh-\yh)-\Gh(t,\xh-\yh)    \Big) \p_3(u_3 \Bh)(\tau,\yh,x_3)  d y_{\rm h}   d \tau,   \\
   \mathcal{I}_2 (t,x) = & \    \int_0^{\f{t}{2}}   \int_{\R^2}
   \Big( \Gh(t,\xh-\yh)-\Gh(t,\xh)    \Big) \p_3(u_3 \Bh)(\tau,\yh,x_3)  d y_{\rm h}   d \tau,
   \\
       \mathcal{I}_3 (t,x) = & \  \int_{\f{t}{2}}^{t} e^{(t-\tau)\Deltah} \p_3(u_3 \Bh)(\tau,x) d\tau,       \\
       \mathcal{I}_4 (t,x) = & \  -\Gh(t,\xh) \int_{\f{t}{2}}^{\infty} \int_{\R^2} \p_3(u_3 \Bh)(\tau,\yh,x_3)  d y_{\rm h}   d \tau.   \\
     \end{split}
\end{equation}
Thanks to Lemma \ref{lem12}, it holds for any $1\leq p \leq \infty$ that
\begin{equation}\nonumber
     \begin{split}
        \| \p_3(u_3 \Bh)(\tau)  \|_{L^{1}_{\rm h}  L^{p}_{\rm v}} \leq  & \
        \| \p_3 u_3 (\tau) \|_{L^p}   \| \Bh(\tau)   \|_{L^{p'}_{\rm h}  L^{\infty}_{\rm v}}
        +\| \p_3 \Bh (\tau) \|_{L^p}   \| u_3(\tau)   \|_{L^{p'}_{\rm h}  L^{\infty}_{\rm v}}
        \\
   \leq  & \       \| \p_3 u_3 (\tau) \|_{L^p}   \| \Bh(\tau)   \|_{L^{1}_{\rm h}  L^{\infty}_{\rm v}}^{1-\f{1}{p}}  \| \Bh(\tau)   \|_{  L^{\infty} }  ^{\f{1}{p}}                     \\
         & \ +   \| \p_3 \Bh (\tau) \|_{L^p}   \| u_3(\tau)   \|_{L^{1}_{\rm h}  L^{\infty}_{\rm v}}^{1-\f{1}{p}}  \| u_3 (\tau)   \|_{  L^{\infty} }  ^{\f{1}{p}}                     \\
        \leq  & \  CA^2   \Big(   (1+\tau)^{     -\f{17}{8} +\f{1}{8} \f{1}{p}   } \log(2+\tau)  + (1+\tau)^{        -2+ \f{3}{8} \f{1}{p}     }  \log(2+\tau) \Big)  \\
      \leq   &\      CA^2 (1+\tau)^{-\f{3}{2}}.                  \\
     \end{split}
\end{equation}
We are now in a position to estimate $\mathcal{I}_m (m=1,2,3,4)$ suitably. For $\mathcal{I}_1$, we rewrite
\[
 \mathcal{I}_1 (t,x) =  - \int_0^{\f{t}{2}}   \int_{\R^2} \int_0^1
 \tau (\p_t \Gh) (t-\theta\tau,\xh-\yh) \p_3(u_3 \Bh)(\tau,\yh,x_3)   d \theta d y_{\rm h}   d \tau;
\]
whence
\begin{equation}\nonumber
     \begin{split}
        \|  \mathcal{I}_1(t) \|_{L^p} \leq  & \ \int_0^{\f{t}{2}}
        \int_0^1  \tau \| (\p_t \Gh)(t-\theta \tau)  \|_{L^p(\R^2)}
        \|  \p_3 (u_3 \Bh)(\tau)  \|_{L^{1}_{\rm h}  L^{p}_{\rm v}} d \theta d \tau
        \\
    \leq  & \  CA^2     \int_0^{\f{t}{2}}
        \int_0^1  \tau   (t-\theta \tau)^{ -(1-\f{1}{p})    -1    } (1+\tau)^{-\f{3}{2}} d \theta d \tau         \\
        \leq  & \ CA^2 t^{    -(1-\f{1}{p})    -1    }   \int_0^{\f{t}{2}}  (1+\tau)^{-\f{1}{2} } d\tau           \\
        \leq  & \   CA^2 t^{    -(1-\f{1}{p})    -\f{1}{2}    } .   \\
     \end{split}
\end{equation}
Making a change of variable $\xh \mapsto t^{\f{1}{2}}  \xi_{\rm h} $,
\[
        \|  \mathcal{I}_2(t) \|_{L^p}
        \leq
t^{-(1-\f{1}{p})}   \int_0^{\f{t}{2}}   \int_{\R^2}
\| \Gh(1,\cdot - t^{-\f{1}{2}} \yh ) -\Gh(1,\cdot) \| _{L^p(\R^2)}
\| \p_3(u_3 \Bh) (\tau, \yh, \cdot) \|_{L^p(\R)}   d y_{\rm h}   d \tau .
\]
Realizing that $\p_3(u_3 \Bh)  \in L^1(0,\infty;  L^{1}_{\rm h}  L^{p}_{\rm v}(\R^3)    )      $ and invoking Lebesgue's dominated convergence theorem, it follows that $\lim_{t\rightarrow \infty }t^{(1-\f{1}{p})} \|  \mathcal{I}_2(t) \|_{L^p}=0  $. Finally,
\begin{equation}\nonumber
     \begin{split}
         \|  \mathcal{I}_3(t) \|_{L^p}+ \|  \mathcal{I}_4(t) \|_{L^p} \leq
         & \       \int_{\f{t}{2}}^{t}   \|  \Gh(t-\tau) \|_{L^p(\R^2)}        \|  \p_3(u_3 \Bh)(\tau ) \|_{     L^{1}_{\rm h}  L^{p}_{\rm v}   }    d \tau              \\
   & \  +    \|  \Gh(t) \|_{L^p(\R^2)}   \int_{\f{t}{2}}^{\infty}
            \|  \p_3(u_3 \Bh)(\tau ) \|_{     L^{1}_{\rm h}  L^{p}_{\rm v}   }  d \tau                              \\
     \leq     & \ CA^2  \int_{\f{t}{2}}^{t}   (t-\tau)^{ -(1-\f{1}{p}) } (1+\tau)^{  -\f{3}{2}     } d\tau + CA^2     t^{ -(1-\f{1}{p}) } \int_{\f{t}{2}}^{\infty}  (1+\tau)^{  -\f{3}{2}       } d \tau      \\
        \leq  & \   CA^2   t^{    -(1-\f{1}{p})    -\f{1}{2}    }.               \\
     \end{split}
\end{equation}
Combining these estimates, we arrive at (\ref{hm21}).

Exactly in the same manner, one obtains (\ref{hm22}). Indeed, it suffices to observe that
\begin{equation}\nonumber
     \begin{split}
        \| \p_3(B_3 \uh)(\tau)  \|_{L^{1}_{\rm h}  L^{p}_{\rm v}} \leq  & \
        \| \p_3 B_3 (\tau) \|_{L^p}   \| \uh(\tau)   \|_{L^{p'}_{\rm h}  L^{\infty}_{\rm v}}
        +\| \p_3 \uh (\tau) \|_{L^p}   \| B_3(\tau)   \|_{L^{p'}_{\rm h}  L^{\infty}_{\rm v}}
        \\
   \leq  & \       \| \p_3 B_3 (\tau) \|_{L^p}   \| \uh(\tau)   \|_{L^{1}_{\rm h}  L^{\infty}_{\rm v}}^{1-\f{1}{p}}  \| \uh(\tau)   \|_{  L^{\infty} }  ^{\f{1}{p}}                     \\
         & \ +   \| \p_3 \uh (\tau) \|_{L^p}   \| B_3(\tau)   \|_{L^{1}_{\rm h}  L^{\infty}_{\rm v}}^{1-\f{1}{p}}  \| B_3 (\tau)   \|_{  L^{\infty} }  ^{\f{1}{p}}                     \\
        \leq  & \  CA^2   \Big(   (1+\tau)^{     -\f{33}{8} +2 \f{1}{p}   } \log(2+\tau)  + (1+\tau)^{        -\f{33}{8} + \f{13}{8} \f{1}{p}     }  \log(2+\tau) \Big)  \\
      \leq   &\      CA^2 (1+\tau)^{-\f{17}{8}}  \log(2+\tau)     ,                  \\
     \end{split}
\end{equation}
which gives a faster decay rate than $ \| \p_3(u_3 \Bh)(\tau)  \|_{L^{1}_{\rm h}  L^{p}_{\rm v}}$.   $\Box$

In the same spirit of Lemma \ref{lem10}, we are able to show the asymptotic expansions of $\Nvd_1[u,B]$ and $\Nvd_2[u,B]$ by the estimates in Lemma \ref{lem15}. These estimates are crucial for deriving the higher order asymptotic expansion of $B_3(t)$. The details of proof are omitted.
\begin{lemm}\label{lem16}
Let $(u,B)$ be subject to the ansatzes $(i),(ii)$ with $s\geq 3,A>0,T=\infty$. Then there hold the asymptotic limits for any $1\leq p \leq \infty$:
\beq\label{hm25}
\lim_{t\rightarrow \infty} t^{(1-\f{1}{p})+\f{1}{2}}
\left\|
\Nvd_1[u,B] (t,x)+\nablah G_{\rm h}(t,\xh) \cdot  \int_0^{\infty} \int_{\R^2}   (B_3 \uh) (\tau, y_{\rm h},x_3)
    d y_{\rm h}   d \tau
\right\|_{L^p_x} =0,
\eeq
\beq\label{hm26}
\lim_{t\rightarrow \infty} t^{(1-\f{1}{p}) + \f{1}{2}}
\left\|
\Nvd_2[u,B] (t,x)-  \nablah G_{\rm h}(t,\xh) \cdot \int_0^{\infty} \int_{\R^2}   (u_3 \Bh) (\tau, y_{\rm h},x_3)
    d y_{\rm h}   d \tau
\right\|_{L^p_x} =0.
\eeq

\end{lemm}

\subsection{Proof of Theorem \ref{thm2}}

Let $s\geq 5$ be an integer and $u_0,B_0\in X^s(\R^3)$ satisfying $\nabla \cdot u_0=\nabla \cdot B_0=0, |x| u_0(x) \in L^1(\R^3),\nablah u_0 \in L^1(\R^3)$.
We first briefly explain the global small solutions for (\ref{MHD-3}). It follows from Proposition \ref{pr2} that the Cauchy problem (\ref{MHD-3}) admits a unique solution $u,B\in C([0,\infty);H^s(\R^3))$ under the smallness condition that $\| (u_0,B_0) \|_{H^s} \leq \delta$ for some $\delta>0$ sufficiently small. Moreover, similar to the derivation of (\ref{hd100}), one derives
\[
  \| (u,B)(t)   \|  _{  L^1_{\xh}(W^{1,1}\cap W^{1,\infty})_{x_3}         }
  \leq C (1+t) \| (u_0,B_0)  \|_{X^s}
\]
and concludes that the mappings $t \mapsto  \| u(t)   \|  _{  L^1_{\xh}(W^{1,1}\cap W^{1,\infty})_{x_3}         }$, $t \mapsto  \| B(t)   \|  _{  L^1_{\xh}(W^{1,1}\cap W^{1,\infty})_{x_3}         }$ are continuous with respect to time. Consequently, we see that $u,B \in  C([0,\infty);X^s(\R^3))$.

We next come to the decay estimates in $L^p$-norms, namely (\ref{ma-8})-(\ref{ma-10}), via the bootstrapping argument. In accordance with the discussions in subsection \ref{se-l} concerning the linearized equations, we know that there exists generic positive constants $C_0,C_1$ such that
\begin{equation}\label{hm101}
     \begin{split}
      \| \nabla ^{\alpha}   e^{t \Deltah} B_{0,\rm h} (t) \|_{L^p} & \ \leq C_1 t^{ -(1-\f{1}{p}) -\f{|\alphah|}{2}        }  \| (u_0,B_0) \|_{  X^{s}_{\#}  }, \\
    \| \nablah ^{\alphah}   e^{t \Deltah} B_{0,3} (t) \|_{L^p} & \   \leq C_1 t^{ -\f{3}{2}(1-\f{1}{p}) -\f{|\alphah|}{2}        }  \|(u_0,B_0)\|_{  X^{s}_{\#}  }, \\
      \| \nabla ^{\alpha}   e^{t\Delta} u_0 (t) \|_{L^p}    & \  \leq C_0 (1+t)^{ -\f{3}{2}(1-\f{1}{p}) -\f{1}{2}-\f{|\alpha|}{2}        }  \| (u_0,B_0) \|_{  X^{s}_{\#}  } \\
     &\  \leq    C_1  (1+t)^{ -\f{9}{8}(1-\f{1}{p}) -\f{1}{2}-\f{|\alpha|}{2}        } \log(2+t)  \| (u_0,B_0) \|_{  X^{s}_{\#}  } ,        \\
     \end{split}
\end{equation}
for any $t>0$, $1\leq p \leq \infty$ and any $\alpha=(\alphah,\alpha_3)\in (\mathbb{N} \cup \{ 0\})^2 \times (\mathbb{N} \cup \{ 0\})$ with $|\alpha|\leq 1$. Suppose that
\begin{equation}\label{hm102}
     \begin{split}
      \| \nabla ^{\alpha}   \Bh (t) \|_{L^p} & \ \leq 2 C_1 t^{ -(1-\f{1}{p}) -\f{|\alphah|}{2}        }  \| (u_0,B_0) \|_{  X^{s}_{\#}  }, \\
    \| \nablah ^{\alphah}   B_3 (t) \|_{L^p} & \   \leq 2 C_1 t^{ -\f{3}{2}(1-\f{1}{p}) -\f{|\alphah|}{2}        }  \|(u_0,B_0)\|_{  X^{s}_{\#}  }, \\
      \| \nabla ^{\alpha}  u (t) \|_{L^p}    & \  \leq 2 C_1 (1+t)^{ -\f{9}{8}(1-\f{1}{p}) -\f{1}{2}-\f{|\alpha|}{2}        } \log(2+t)  \| (u_0,B_0) \|_{  X^{s}_{\#}  } , \\
     \end{split}
\end{equation}
holds for any $0<t <T$. Upon choosing $A:=2C_1 \| (u_0,B_0) \|_{  X^{s}_{\#}  }$, it follows from Lemma \ref{lem13} that
\begin{equation}\label{hm103}
     \begin{split}
      \| \nabla ^{\alpha}   \Bh (t) \|_{L^p} \leq & \   C_1 t^{ -(1-\f{1}{p}) -\f{|\alphah|}{2}        }  \| (u_0,B_0) \|_{  X^{s}_{\#}  } \\
      &\ +C_2 (1+t)^{ -(1-\f{1}{p}) -\f{|\alphah|}{2}        }  \| (u_0,B_0) \|_{  X^{s}_{\#}  }^2, \\
    \| \nablah ^{\alphah}   B_3 (t) \|_{L^p}  \leq & \    C_1 t^{ -\f{3}{2}(1-\f{1}{p}) -\f{|\alphah|}{2}        }  \|(u_0,B_0)\|_{  X^{s}_{\#}  } \\
   &\  + C_2 (1+t)^{ -\f{3}{2}(1-\f{1}{p}) -\f{1}{2p} -\f{|\alphah|}{2}        }  \|(u_0,B_0)\|_{  X^{s}_{\#}  }^2,  \\
      \| \nabla ^{\alpha}  u (t) \|_{L^p}   \leq  & \   C_1 (1+t)^{ -\f{3}{2}(1-\f{1}{p}) -\f{1}{2}-\f{|\alpha|}{2}        }  \| (u_0,B_0) \|_{  X^{s}_{\#}  }   \\
      &\ + C_2 (1+t)^{ -\f{9}{8}(1-\f{1}{p}) -\f{1}{2}-\f{|\alpha|}{2}        } \log(2+t)  \| (u_0,B_0) \|_{  X^{s}_{\#}  } ^2,
      \\
     \end{split}
\end{equation}
holds for any $0<t <T$ and some $C_2>0$. By setting $\delta_1:= \min\{\delta, \f{C_1}{2C_2}  \}$ and choosing the initial data sufficiently small such that $\| (u_0,B_0) \|_{  X^{s}_{\#}  }  \leq \delta_1$, we further deduce from (\ref{hm103}) that
\begin{equation}\label{hm104}
     \begin{split}
      \| \nabla ^{\alpha}   \Bh (t) \|_{L^p} & \ \leq \f{3}{2} C_1 t^{ -(1-\f{1}{p}) -\f{|\alphah|}{2}        }  \| (u_0,B_0) \|_{  X^{s}_{\#}  }, \\
    \| \nablah ^{\alphah}   B_3 (t) \|_{L^p} & \   \leq \f{3}{2} C_1 t^{ -\f{3}{2}(1-\f{1}{p}) -\f{|\alphah|}{2}        }  \|(u_0,B_0)\|_{  X^{s}_{\#}  }, \\
      \| \nabla ^{\alpha}  u (t) \|_{L^p}    & \  \leq \f{3}{2} C_1 (1+t)^{ -\f{9}{8}(1-\f{1}{p}) -\f{1}{2}-\f{|\alpha|}{2}        } \log(2+t)  \| (u_0,B_0) \|_{  X^{s}_{\#}  }, \\
     \end{split}
\end{equation}
for any $0<t <T$. The estimate (\ref{hm104}) is sharper than our initial hypothesis (\ref{hm102}). The bootstrapping argument then assesses that (\ref{hm104}) actually holds for any $0<t<\infty$ and any $1\leq p \leq \infty$. The proof of nonlinear $L^p$ decay estimates is completely finished.

In the end, we prove the asymptotic profile of solutions, in particular $\Bh(t)$ and $B_3(t)$. By Lemma \ref{lm1}, Lemma \ref{lem13} and Lemma \ref{lem14}, we obtain for any $1\leq p \leq \infty$ that
\begin{equation}\nonumber
     \begin{split}
     t^{1-\f{1}{p}}  \Big\|    \Bh(t,x) &\  -\Gh (t,\xh) \int_{\R^2} B_{0,\rm h}(\yh,x_3) d \yh
       \\
       &\ +  \Gh (t,\xh) \int_0^{\infty} \int_{\R^2} \p_3(u_3 \Bh)  (\tau,\yh,x_3) d \yh  d \tau  \\
     &\ -  \Gh (t,\xh) \int_0^{\infty} \int_{\R^2} \p_3(B_3 \uh)  (\tau,\yh,x_3) d \yh  d \tau \Big \|_{L^p_x}       \\
    \leq & \  t^{1-\f{1}{p}}
  \left\|   e^{t\Deltah} B_{0,\rm h}   -\Gh (t,\xh) \int_{\R^2} B_{0,\rm h}(\yh,x_3) d \yh  \right\|_ {L^p_x}         \\
  &\ + t^{1-\f{1}{p}}
\left\|
\Nhd_3[u,B] (t,x)+ G_{\rm h}(t,\xh)  \int_0^{\infty} \int_{\R^2}  \p_3 (u_3 \Bh) (\tau, y_{\rm h},x_3)
    d y_{\rm h} d \tau
\right\|_{L^p_x}          \\
&\ + t^{1-\f{1}{p}}
\left\|
\Nhd_4[u,B] (t,x)- G_{\rm h}(t,\xh)  \int_0^{\infty} \int_{\R^2}  \p_3 (B_3 \uh) (\tau, y_{\rm h},x_3)
    d y_{\rm h} d \tau
\right\|_{L^p_x}           \\
&\ +\sum_{m=1}^2  t^{1-\f{1}{p}}
 \| \Nhd_m[u,B] (t)\| _{L^p}          \\
 \leq & \  t^{1-\f{1}{p}}
  \left\|   e^{t\Deltah} B_{0,\rm h}   -\Gh (t,\xh) \int_{\R^2} B_{0,\rm h}(\yh,x_3) d \yh  \right\|_ {L^p_x}         \\
  &\ + t^{1-\f{1}{p}}
\left\|
\Nhd_3[u,B] (t,x)+ G_{\rm h}(t,\xh)  \int_0^{\infty} \int_{\R^2}  \p_3 (u_3 \Bh) (\tau, y_{\rm h},x_3)
    d y_{\rm h} d \tau
\right\|_{L^p_x}          \\
&\ + t^{1-\f{1}{p}}
\left\|
\Nhd_4[u,B] (t,x)- G_{\rm h}(t,\xh)  \int_0^{\infty} \int_{\R^2}  \p_3 (B_3 \uh) (\tau, y_{\rm h},x_3)
    d y_{\rm h} d \tau
\right\|_{L^p_x}           \\
&\ +C \| (u_0,B_0) \|_{  X^s_{\# } }
  t^{-\f{1}{2} }.  \\
     \end{split}
\end{equation}
Notice that the right-hand side of the above inequality tends to zero as time goes to infinity. This proves (\ref{ma-10-1}). Next, by Lemma \ref{lm2} and Lemma \ref{lem13} we deduce that
\begin{equation}\nonumber
     \begin{split}
     t^{ \f{3}{2} (1-\f{1}{p})  } &\  \Big\|    B_3 (t,x)  -\Gh (t,\xh) \int_{\R^2} B_{0, 3}(\yh,x_3) d \yh
     \Big \|_{L^p_x}
       \\
      \leq  &\
       t^{ \f{3}{2} (1-\f{1}{p})  }  \Big\|    e^{t \Deltah} B_{0,3}(x)   -\Gh (t,\xh) \int_{\R^2} B_{0, 3}(\yh,x_3) d \yh
     \Big \|_{L^p_x}
      \\
      &\  \qquad + \sum_{m=1}^2 t^{ \f{3}{2} (1-\f{1}{p})  }
      \| \Nvd_m[u,B] (t)\| _{L^p}
      \\
      \leq &\   C \mathcal{R}_{p,0}(t) ^{\f{1}{p}} \|B_0 \|_{L^1}^{1-\f{1}{p}  }
      +C\| (u_0,B_0) \|_{  X^s_{\# }}  t^{-\f{1}{2p}} .    \\
    \end{split}
\end{equation}
It is easily seen that the right-hand side of the above inequality tends to zero as time approaches to infinity provided that $1\leq p<\infty$. This proves (\ref{ma-10-2}). As a consequence, the proof of Theorem \ref{thm2} is completely finished.

\section{Improved estimates and higher order asymptotic expansions} \label{hie}

\subsection{Improved estimates}

Our main goal in this subsection is to confirm, via the bootstrapping argument, that the ansatzes (\ref{hd9-1})-(\ref{hd9-2}) actually hold with $T=\infty,A=C \| (u_0,B_0) \|_{  \overline{ X^{s} } } $ and $A_{\ast}= C\| (u_0,B_0)  \|_{ \widetilde{X^s}  }  $ for some generic $C>0$.
\begin{lemm}\label{lem50}
Let $s \geq 9$ be an integer. Then there exists a positive constant $\epsilon_3=\epsilon_3(s)$ and a generic positive constant $C$ such that for any $u_0,B_0\in X^s(\R^3)$ satisfying $\nabla \cdot u_0=\nabla \cdot B_0=0,|x|B_0(x)\in L^1(\R^3)$, $\| (u_0,B_0) \|_{  \overline{ X^{s} } } \leq \epsilon_3$, the global solution $(u,B)$ to system (\ref{MHD-2}) admits the following decay estimates:
\begin{equation}\label{im70}
     \begin{split}
        \| \Grad^{\alpha} u (t) \|_{L^{\infty}_{\rm h}  L^{1}_{\rm v}} & \  \leq C  \| (u_0,B_0) \|_{  \overline{ X^{s} } } t^{-1-\f{|\alphah|}{2}}, \\
     \| \Grad^{\alpha} B (t) \|_{L^{\infty}_{\rm h}  L^{1}_{\rm v}} & \  \leq C   \| (u_0,B_0) \|_{  \overline{ X^{s} } }    t^{-1-\f{1+|\alpha|}{2}},  \\
     \end{split}
\end{equation}
for any $0<t<\infty$ and $\alpha=(\alphah,\alpha_3)\in (\mathbb{N} \cup \{ 0\})^2 \times (\mathbb{N} \cup \{ 0\})$ with $|\alpha|\leq 1$.
\end{lemm}
{\it Proof.} It follows from Lemma \ref{lm1} and Lemma \ref{lm3-2}\footnote{Indeed, in analogy with the proof of Lemma \ref{lm3-2}, one has
\begin{equation}\nonumber
     \begin{split}
     \| \Grad^{\alpha} e^{t \Delta} B_0 (x)\|_{L^{\infty}_{\rm h}  L^{1}_{\rm v}} \leq & \  C \int_{\R^3}  \int_0^1  \| (\Grad \Grad^{\alpha} G) (t,\cdot-\theta y)  \|_{L^{\infty}_{\rm h}  L^{1}_{\rm v}}  |y| |B_0(y)| d \theta dy       \\
  \leq  & \ C   t^{  -\f{3}{2}-\f{|\alpha|}{2}  }
\| |y| B_0(y) \|_{L^1_y} .              \\
     \end{split}
\end{equation}
} that
\begin{equation}\label{im70-1}
     \begin{split}
        \| \Grad^{\alpha} e^{t \Deltah} u_0 \|_{L^{\infty}_{\rm h}  L^{1}_{\rm v}} & \  \leq C_1 \| u_0 \|_{X^s} t^{-1-\f{|\alphah|}{2}}, \\
     \| \Grad^{\alpha} e^{t \Delta} B_0 \|_{L^{\infty}_{\rm h}  L^{1}_{\rm v}} & \  \leq  C_1 \| B_0 \|_{  \overline{ X^{s} } }   t^{-1-\f{1+|\alpha|}{2}},  \\
     \end{split}
\end{equation}
for any $0<t<\infty$ and $\alpha=(\alphah,\alpha_3)\in (\mathbb{N} \cup \{ 0\})^2 \times (\mathbb{N} \cup \{ 0\})$ with $|\alpha|\leq 1$.
Suppose that
\begin{equation}\label{im70-2}
     \begin{split}
        \| \Grad^{\alpha} u(t) \|_{L^{\infty}_{\rm h}  L^{1}_{\rm v}} & \  \leq 2 C_1 \| (u_0,B_0) \|_{  \overline{ X^{s} } } t^{-1-\f{|\alphah|}{2}}, \\
     \| \Grad^{\alpha} B(t) \|_{L^{\infty}_{\rm h}  L^{1}_{\rm v}} & \  \leq  2 C_1 \| (u_0, B_0) \|_{  \overline{ X^{s} } }    t^{-1-\f{1+|\alpha|}{2}},  \\
     \end{split}
\end{equation}
holds for any $0<t<T$. We further deduce from (\ref{im70-1}) and Lemma \ref{lem8} with $A=C \| (u_0, B_0) \|_{  \overline{ X^{s} } } $ that
\begin{equation}\label{im70-3}
     \begin{split}
      \| \nabla ^{\alpha}   u (t) \|_{L^{\infty}_{\rm h}  L^{1}_{\rm v}} \leq & \   C_1 \| (u_0, B_0) \|_{  \overline{ X^{s} } }   t^{-1-\f{|\alphah|}{2}} \\
      &\ +C_2 \| (u_0,B_0) \|_{  \overline{ X^{s} } }^2   t^{-1-\f{|\alphah|}{2}}          , \\
      \| \nabla ^{\alpha}  B (t) \|_{L^{\infty}_{\rm h}  L^{1}_{\rm v}}  \leq  & \   C_1      \| (u_0, B_0) \|_{  \overline{ X^{s} } }   t^{-1-\f{1+|\alpha|}{2}} \\
      &\ + C_2 \| (u_0,B_0) \|_{X^{s}}^2  t^{-1-\f{1+|\alpha|}{2}} ,
      \\
     \end{split}
\end{equation}
holds for any $0<t <T$ and some $C_2>0$. Thus, by setting $\epsilon_3:= \min\{\epsilon, \f{C_1}{2C_2}  \}$ and choosing the initial data sufficiently small such that $\| (u_0,B_0) \|_{  \overline{ X^{s} } } \leq \epsilon_3$, we further deduce from (\ref{im70-3}) that
\begin{equation}\label{im70-4}
     \begin{split}
        \| \Grad^{\alpha} u(t) \|_{L^{\infty}_{\rm h}  L^{1}_{\rm v}} & \  \leq \f{3}{2} C_1 \| (u_0,B_0) \|_{  \overline{ X^{s} } } t^{-1-\f{|\alphah|}{2}}, \\
     \| \Grad^{\alpha} B(t) \|_{L^{\infty}_{\rm h}  L^{1}_{\rm v}} & \  \leq  \f{2}{2} C_1 \| (u_0, B_0) \|_{  \overline{ X^{s} } }   t^{-1-\f{1+|\alpha|}{2}},  \\
     \end{split}
\end{equation}
holds for any $0<t<T$. In this way we obtain a sharper bound than our initial assumption (\ref{im70-2}). The bootstrapping argument then assesses that (\ref{im70-4}) actually holds for any $0<t<\infty$. The proof of Lemma \ref{lem50} is completely finished.           $\Box$

The next lemma concerns the space-weighted estimates of solutions. The proof is a little tedious and the main idea is a carefully designed integral equation that contains no derivative loss in the vertical direction, together with the contraction mapping principle.
\begin{lemm}\label{lem51}
Let $s \geq 9$ be an integer. Then there exists a positive constant $\epsilon_2=\epsilon_2(s)$ and a generic positive constant $C$ such that for any $u_0,B_0\in X^s(\R^3)$ satisfying $|\xh| u_0(x) \in L^1(\R^2_{\xh};(L^1 \cap L^{\infty})(\R_{x_3})  ), |x| B_0(x)\in L^1(\R^3)$ and $\nabla \cdot u_0=\nabla \cdot B_0=0$, $\| (u_0,B_0) \|_{  \overline{ X^{s} } }  \leq \epsilon_2$, the global solution $(u,B)$ to system (\ref{MHD-2}) admits the following space-weighted estimates:
\begin{equation}\label{im71}
     \begin{split}
        \| |\xh|\uh (t,x) \|_{L^{1}_{\rm h}  L^{\infty}_{\rm v}} & \  \leq C    \| (u_0,B_0)  \|_{ \widetilde{X^s}  }      (1+t)^{  \f{1}{2}  }        ,  \\
     \| |\xh| u_3 (t,x) \|_{L^{1}_{\rm h}  L^{\infty}_{\rm v}} & \  \leq C  \| (u_0,B_0)  \|_{ \widetilde{X^s}  },  \\
         \| |x| \Bh (t,x) \|_{L^{1}_{\rm h}  L^{\infty}_{\rm v}} & \  \leq  C  \| (u_0,B_0)  \|_{ \widetilde{X^s}  }  (1+t)^{  \f{1}{2}  },  \\
           \| |x| B_3 (t,x) \|_{L^{1}_{\rm h}  L^{\infty}_{\rm v}} & \  \leq   C  \| (u_0,B_0)  \|_{ \widetilde{X^s}  }  (1+t)^{  \f{1}{2}  }.  \\
     \end{split}
\end{equation}
\end{lemm}
{\it Proof.} We consider the following integral equation:
\begin{equation}\label{im72}
\left\{\begin{aligned}
&  \vh(t) = e^{ t \Deltah} u_{0,{\rm h}} + \sum_{m=1}^5 \Big( \Nhu_m[u,v](t) -  \Nhu_m[B,W](t) \Big) ,  \\
&  v_3(t) = e^{ t \Deltah} u_{0,3} + \sum_{m=1}^3 \Big( \Nvu_m[u,v](t) -  \Nvu_m[B,W](t) \Big) ,  \\
&    W_{\rm h}(t) = e^{ t \Delta} B_{0,{\rm h}} + \sum_{m=1}^3 \Nhb_m[u,W](t)        ,\\
&       W_3(t) = e^{ t \Delta} B_{0,3} + \sum_{m=1}^3 \Nvb_m[u,W](t)          ,\\
\end{aligned}\right.
\end{equation}
where $v(t,x)=(\vh,v_3)(t,x),W(t,x)=(  W_{\rm h},W_3 )(t,x)$ are the unknown vector fields and $(u,B)$ is the solution to system (\ref{MHD-2}) with initial data $(u_0,B_0)$.

\begin{align*}
		\Nhu_1[u,v](t)&:=-\int_0^t  e^{ (t-\tau)\Deltah}   (   (\p_3 u_3)\vh+ v_3 \p_3 \uh  ) (\tau)    d\tau,\\
		\Nhu_2[u,v](t)&:=-\int_0^t e^{ (t-\tau)\Deltah}\nablah \cdot (\uh \otimes\vh)(\tau)d\tau,\\
		\Nhu_3[u,v](t)&:=\int_0^t\nablah e^{ (t-\tau)\Deltah}( u_3 v_3) (\tau) d\tau,\\
		\Nhu_4[u,v](t)&:=-\sum_{k,l=1}^2\int_0^t\nablah\partial_{k}\partial_{l}K(t-\tau)*(u_k v_l)(\tau)d\tau,\\
		\Nhu_5[u,v](t)&:=2\sum_{k=1}^2\int_0^t\nablah \partial_k(-\Deltah)^{\frac{1}{2}}\widetilde{K}(t-\tau)*(u_3 v_k)(\tau)d\tau\\
		&\qquad +\int_0^t\nablah \Deltah K(t-\tau)*(u_3 v_3) (\tau)  d\tau,
\end{align*}
\begin{align*}
		\Nvu_1[u, v](t)&:=\int_0^te^{(t-\tau)\Deltah}\nablah\cdot(u_3\vh)(\tau)d\tau,\\
		\Nvu_2[u, v](t)&:=\sum_{k,l=1}^2\int_0^t(-\Deltah)^{\frac{1}{2}}\partial_{k}\partial_{l}\widetilde{K}(t-\tau)*(u_k v_l)(\tau)d\tau,\\
		\Nvu_3[u,v](t)&:=2\sum_{k=1}^2\int_0^t\partial_k\Deltah K(t-\tau)*(u_3 v_k)(\tau)d\tau\\
		&\qquad+\int_0^t(-\Deltah)^{\frac{3}{2}}\widetilde{K}(t-\tau)*(u_3 v_3) (\tau)   d\tau ,
\end{align*}
\begin{align*}
		\Nhu_1[B,W](t)&:=-\int_0^t  e^{ (t-\tau)\Deltah}   (   (\p_3 B_3)W_{\rm h}+ W_3 \p_3 \Bh  ) (\tau)    d\tau,\\
		\Nhu_2[B,W](t)&:=-\int_0^t e^{ (t-\tau)\Deltah}\nablah \cdot (\Bh \otimes  W_{\rm h})(\tau)d\tau,\\
		\Nhu_3[B,W](t)&:=\int_0^t\nablah e^{ (t-\tau)\Deltah}( B_3 W_3) (\tau) d\tau,\\
		\Nhu_4[B,W](t)&:=-\sum_{k,l=1}^2\int_0^t\nablah\partial_{k}\partial_{l}K(t-\tau)*(B_k W_l)(\tau)d\tau,\\
		\Nhu_5[B,W](t)&:=2\sum_{k=1}^2\int_0^t\nablah \partial_k(-\Deltah)^{\frac{1}{2}}\widetilde{K}(t-\tau)*(B_3 W_k)(\tau)d\tau\\
		&\qquad +\int_0^t\nablah \Deltah K(t-\tau)*(B_3 W_3) (\tau)  d\tau,
\end{align*}
\begin{align*}
		\Nvu_1[B, W](t)&:=\int_0^te^{(t-\tau)\Deltah}\nablah\cdot(B_3 W_{\rm h} )(\tau)d\tau,\\
		\Nvu_2[B, W](t)&:=\sum_{k,l=1}^2\int_0^t(-\Deltah)^{\frac{1}{2}}\partial_{k}\partial_{l}\widetilde{K}(t-\tau)*(B_k W_l)(\tau)d\tau,\\
		\Nvu_3[B,W](t)&:=2\sum_{k=1}^2\int_0^t\partial_k\Deltah K(t-\tau)*(u_3 v_k)(\tau)d\tau\\
		&\qquad+\int_0^t(-\Deltah)^{\frac{3}{2}}\widetilde{K}(t-\tau)*(B_3 W_3) (\tau)   d\tau ,
\end{align*}
\begin{equation}\nonumber
     \begin{split}
         \Nhb_1[u,W](t) & \  :=-\int_0^t e^{(t-\tau)\Delta} \nabla \cdot ( W_{\rm h} \otimes u )(\tau)d\tau, \\
     \Nhb_2[u,W](t) & \  :=\int_0^t e^{(t-\tau)\Delta} \nablah \cdot( \uh \otimes W_{\rm h} )(\tau)d\tau, \\
     \Nhb_3[u,W](t)    & \ :=\int_0^t e^{(t-\tau)\Delta} \partial_3 ( W_3 \uh )(\tau)d\tau,\\
     \end{split}
\end{equation}
\begin{align*}
        \Nvb_1[u,W](t)&:=-\int_0^t e^{(t-\tau)\Delta} \nabla \cdot ( W_3  u )(\tau)d\tau,\\
        \Nvb_2[u,W](t)&:=\int_0^t e^{(t-\tau)\Delta} \nablah \cdot (u_3  W_{\rm h}  )(\tau)d\tau,\\
        \Nvb_3[u,W](t)&:=\int_0^t e^{(t-\tau)\Delta} \partial_3 ( u_3 W_3)(\tau)d\tau.
\end{align*}

Designing the following space $Y$ with space-weighted norms:
\[
Y= \{  (v,W)\in L^{\infty}( 0,\infty; L^1_{\rm h} L^{\infty}_{\rm v} (\R^3)         ); \|  (v,W) \|_{Y}<\infty       \}
\]
with
\begin{equation}\nonumber
     \begin{split}
      \|  (v,W) \|_{Y}:=  & \ \sup_{t \geq 0} (1+t)^{-\f{1}{2}}   \| |\xh|\vh (t,x) \|_{L^{1}_{\rm h}  L^{\infty}_{\rm v}}  +\sup_{t \geq 0}   \| |\xh| v_3 (t,x) \|_{L^{1}_{\rm h}  L^{\infty}_{\rm v}}            \\
     & \ +   \sup_{t \geq 0} (1+t)^{-\f{1}{2}}   \| |x| W_{\rm h} (t,x) \|_{L^{1}_{\rm h}  L^{\infty}_{\rm v}}    + \sup_{t \geq 0} (1+t)^{-\f{1}{2}}   \| |x| W_{3} (t,x) \|_{L^{1}_{\rm h}  L^{\infty}_{\rm v}}     \\
         & \  +  \sup_{t \geq 0}  \| \vh (t) \|_{  L^1 (\R^2_{\xh} ;W^{1,\infty}(\R_{x_3})   )  }  +  \sup_{t \geq 0} (1+t)^{\f{1}{2}}  \| v_3 (t) \|_{  L^1 ( \R^2_{\xh};W^{1,\infty}(\R_{x_3}) )    }          \\
         & \   +  \sup_{t \geq 0} (1+t)^{\f{1}{2}}  \| W_{\rm h}  (t) \|_{  L^1 ( \R^2_{\xh};W^{1,\infty}(\R_{x_3}) )    }  + \sup_{t \geq 0} (1+t)^{\f{1}{2}}  \| W_{3}  (t) \|_{  L^1 ( \R^2_{\xh};W^{1,\infty}(\R_{x_3}) )    } .   \\
     \end{split}
\end{equation}
We observe that the key issue of the proof lies in the following estimates:
\begin{equation}\label{im73}
     \begin{split}
       \| |\xh|	\Nhu_1[u,v](t,x) \|_{L^{1}_{\rm h}  L^{\infty}_{\rm v}}  \leq C  & \    (1+t)^{\f{1}{2}}  \| (u_0,B_0)     \|_{X^s}  \| (v,B) \|_{Y} ,      \\
   \| 	\Nhu_1[u,v](t,x) \|_{L^{1}_{\rm h}  L^{\infty}_{\rm v}}  \leq C  & \      \| (u_0,B_0)     \|_{X^s}  \| (v,B) \|_{L^{\infty}( 0,\infty; L^1_{\rm h} L^{\infty}_{\rm v}) } ,      \\
        \| \p_3	\Nhu_1[u,v](t,x) \|_{L^{1}_{\rm h}  L^{\infty}_{\rm v}}  \leq C  & \    \| (u_0,B_0)    \|_{X^s}  \| (v,B) \|_{Y} ,      \\
          \| |\xh|	\Nhu_1[B,W](t,x) \|_{L^{1}_{\rm h}  L^{\infty}_{\rm v}}  \leq C  & \    (1+t)^{\f{1}{2}}  \| (u_0,B_0)     \|_{X^s}  \| (v,B) \|_{Y} ,      \\
   \| 	\Nhu_1[B,W](t,x) \|_{L^{1}_{\rm h}  L^{\infty}_{\rm v}}  \leq C  & \      \| (u_0,B_0)     \|_{X^s}  \| (v,B) \|_{L^{\infty}( 0,\infty; L^1_{\rm h} L^{\infty}_{\rm v}) } ,      \\
        \| \p_3	\Nhu_1[B,W](t,x) \|_{L^{1}_{\rm h}  L^{\infty}_{\rm v}}  \leq C  & \    \| (u_0,B_0)    \|_{X^s}  \| (v,B) \|_{Y} ,      \\
     \end{split}
\end{equation}
\begin{equation}\label{im74}
     \begin{split}
       \left\|   |\xh|  \int_0^t  e^{(t-\tau)\Deltah } \nablah (u_k v_l)(\tau,x)d \tau     \right\|_{L^{1}_{\rm h}  L^{\infty}_{\rm v}}   \leq  & \   \left\{\begin{aligned}
&  C  (1+t)^{\f{1}{2}}   \| (u_0,B_0)    \|_{X^s}  \| (v,B) \|_{Y}           \,\,\, \text{if  } k=1,2        \\
& C    \| (u_0,B_0)    \|_{X^s}  \| (v,B) \|_{Y}    \,\,\, \text{if  }   k=3,        \\
\end{aligned}\right.     \\
   \left\|    \int_0^t  e^{(t-\tau)\Deltah } \nablah (u_k v_l)(\tau,x)d \tau     \right\|_{L^{1}_{\rm h}  L^{\infty}_{\rm v}}   \leq  & \   \left\{\begin{aligned}
&  C   \| (u_0,B_0)    \|_{X^s}  \| (v,B) \|_{L^{\infty}( 0,\infty; L^1_{\rm h} L^{\infty}_{\rm v})}           \,\,\, \text{if  } k=1,2        \\
& C   (1+t)^{- \f{1}{2}}   \| (u_0,B_0)    \|_{X^s}  \| (v,B) \|_{ L^{\infty}( 0,\infty; L^1_{\rm h} L^{\infty}_{\rm v})  }        \\
&   \text{if  }  \,\,  k=3,  \\
\end{aligned}\right.     \\
        \left\|  \p_3  \int_0^t  e^{(t-\tau)\Deltah } \nablah (u_k v_l)(\tau,x)d \tau     \right\|_{L^{1}_{\rm h}  L^{\infty}_{\rm v}}   \leq  & \   \left\{\begin{aligned}
&  C    \| (u_0,B_0)    \|_{X^s}  \| (v,B) \|_{Y}           \,\,\, \text{if  } k=1,2        \\
& C   (1+t)^{- \f{1}{2}}   \| (u_0,B_0)    \|_{X^s}  \| (v,B) \|_{Y}    \,\,\, \text{if  }   k=3,        \\
\end{aligned}\right.     \\
   \left\|   |\xh|  \int_0^t  e^{(t-\tau)\Deltah } \nablah (B_k W_l)(\tau,x)d \tau     \right\|_{L^{1}_{\rm h}  L^{\infty}_{\rm v}}   \leq  & \  C    \| (u_0,B_0)    \|_{X^s}  \| (v,B) \|_{Y}  , \,\,  k,l\in \{  1,2,3 \}, \\
   \left\|    \int_0^t  e^{(t-\tau)\Deltah } \nablah (B_k W_l)(\tau,x)d \tau     \right\|_{L^{1}_{\rm h}  L^{\infty}_{\rm v}}   \leq  & \   C   (1+t)^{- \f{1}{2}}   \| (u_0,B_0)    \|_{X^s}  \| (v,B) \|_{ L^{\infty}( 0,\infty; L^1_{\rm h} L^{\infty}_{\rm v})  }      \\
   &\     \,\,  \text{if  }    k,l\in \{  1,2,3 \},        \\
        \left\|  \p_3  \int_0^t  e^{(t-\tau)\Deltah } \nablah (B_k W_l)(\tau,x)d \tau     \right\|_{L^{1}_{\rm h}  L^{\infty}_{\rm v}}   \leq  & \  C   (1+t)^{- \f{1}{2}}   \| (u_0,B_0)    \|_{X^s}  \| (v,B) \|_{Y}  \\
        &\ \,\,  \text{if  }    k,l\in \{  1,2,3 \},          \\
     \end{split}
\end{equation}
\begin{equation}\label{im75}
     \begin{split}
        \left\|   |\xh|  \int_0^t   K_{\beta,\gamma}^{(m) } (t-\tau) \ast (u_k v_l)(\tau,x)            d \tau     \right\|_{L^{1}_{\rm h}  L^{\infty}_{\rm v}}   \leq  & \  C \| (u_0,B_0)    \|_{X^s}  \| (v,B) \|_{Y}       \\
          &\ \,\,  \text{if  }    k,l\in \{  1,2,3 \}, m\in \{  1,2\} ,         \\
    \left\|    \int_0^t   K_{\beta,\gamma}^{(m) } (t-\tau) \ast (u_k v_l)(\tau)            d \tau     \right\|_{L^{1}_{\rm h}  L^{\infty}_{\rm v}}   \leq   & \   C      (1+t)^{- \f{1}{2}}   \| (u_0,B_0)    \|_{X^s}  \| (v,B) \|_{ L^{\infty}( 0,\infty; L^1_{\rm h} L^{\infty}_{\rm v})  }      \\
     &\ \,\,  \text{if  }    k,l\in \{  1,2,3 \},  m\in \{  1,2\} ,         \\
         \left\|   \p_3   \int_0^t   K_{\beta,\gamma}^{(m) } (t-\tau) \ast (u_k v_l)(\tau)            d \tau     \right\|_{L^{1}_{\rm h}  L^{\infty}_{\rm v}}   \leq    & \   C      (1+t)^{- \f{1}{2}}   \| (u_0,B_0)    \|_{X^s}  \| (v,B) \|_{ Y }    \\
           &\ \,\,  \text{if  }    k,l\in \{  1,2,3 \},    m\in \{  1,2\} ,       \\
      \left\|   |x| \int_0^t  e^{ (t-\tau)\Delta } \nablah (u_k W_l)  d \tau     \right\|_{L^{1}_{\rm h}  L^{\infty}_{\rm v}}   \leq   &\       C \| (u_0,B_0)    \|_{X^s}  \| (v,B) \|_{Y}  ,        \\
       &\ \,\,  \text{if  }    k,l\in \{  1,2,3 \},          \\
         \left\|    \int_0^t  e^{ (t-\tau)\Delta } \nablah (u_k W_l)  d \tau     \right\|_{L^{1}_{\rm h}  L^{\infty}_{\rm v}}   \leq   &\       (1+t)^{- \f{1}{2}}   \| (u_0,B_0)    \|_{X^s}  \| (v,B) \|_{ L^{\infty}( 0,\infty; L^1_{\rm h} L^{\infty}_{\rm v})  }           \\
       &\ \,\,  \text{if  }    k,l\in \{  1,2,3 \},          \\
         \left\|   \p_3  \int_0^t  e^{ (t-\tau)\Delta } \nablah (u_k W_l)  d \tau     \right\|_{L^{1}_{\rm h}  L^{\infty}_{\rm v}}   \leq   &\       C      (1+t)^{- \f{1}{2}}   \| (u_0,B_0)    \|_{X^s}  \| (v,B) \|_{ Y }         \\
       &\ \,\,  \text{if  }    k,l\in \{  1,2,3 \},          \\
        \left\|   |x| \int_0^t  e^{ (t-\tau)\Delta } \p_3 (u_k W_l)  d \tau     \right\|_{L^{1}_{\rm h}  L^{\infty}_{\rm v}}   \leq   &\       C \| (u_0,B_0)    \|_{X^s}  \| (v,B) \|_{Y}  ,        \\
       &\ \,\,  \text{if  }    k,l\in \{  1,2,3 \},          \\
         \left\|    \int_0^t  e^{ (t-\tau)\Delta } \p_3 (u_k W_l)  d \tau     \right\|_{L^{1}_{\rm h}  L^{\infty}_{\rm v}}   \leq   &\       (1+t)^{- \f{1}{2}}   \| (u_0,B_0)    \|_{X^s}  \| (v,B) \|_{ L^{\infty}( 0,\infty; L^1_{\rm h} L^{\infty}_{\rm v})  }           \\
       &\ \,\,  \text{if  }    k,l\in \{  1,2,3 \},          \\
         \left\|   \p_3  \int_0^t  e^{ (t-\tau)\Delta } \p_3  (u_k W_l)  d \tau     \right\|_{L^{1}_{\rm h}  L^{\infty}_{\rm v}}   \leq   &\       C      (1+t)^{- \f{1}{2}}   \| (u_0,B_0)    \|_{X^s}  \| (v,B) \|_{ Y }         \\
       &\ \,\,  \text{if  }    k,l\in \{  1,2,3 \}.         \\
     \end{split}
\end{equation}

Taking the estimates (\ref{im73})-(\ref{im75}) for granted, we now proceed to prove Lemma \ref{lem51} via the contraction mapping principle. First of all, we know from Lemmas \ref{lm1}-\ref{lm2} that
\[
\|  ( e^{t \Deltah} u_0 , e^{t \Delta}B_0    )  \|_{Y}
\leq C_1 \|   (u_0,B_0)   \|_{ \widetilde{X^s} }.
\]
For any $(v,W) \in Y$, we define the mapping $\Phi: Y \rightarrow Y$ via
\[
(v,W) \mapsto \Phi (v,W):=( \Phi_{\ast}^{\rm h }[v,W] , \Phi_{\ast}^{\rm v }[v,W]  ,
\Phi_{\ast \ast}^{\rm h }[v,W] , \Phi_{\ast \ast}^{\rm v }[v,W]
)
\]
with
\begin{equation}\nonumber
     \begin{split}
         \Phi_{\ast}^{\rm h }[v,W]  & \  = e^{t \Deltah} u_{0,\rm h} +  \sum_{m=1}^5 \Big( \Nhu_m[u,v](t) -  \Nhu_m[B,W](t) \Big) \\
        &\  =:   e^{t \Deltah} u_{0,\rm h} +   f_{\ast}^{\rm h}(u,v;B,W) ,      \\
   \Phi_{\ast}^{\rm v }[v,W]  & \   = e^{t \Deltah} u_{0,3} + \sum_{m=1}^3 \Big( \Nvu_m[u,v](t) -  \Nvu_m[B,W](t) \Big)  \\
   &\ =: e^{t \Deltah} u_{0,3} +  f_{\ast}^{\rm v}(u,v;B,W) ,              \\
       \Phi_{\ast \ast}^{\rm h }[v,W]    & \ = e^{t \Delta}B_{0,\rm h}+ \sum_{m=1}^3 \Nhb_m[u,W](t)         \\
       &\    =: e^{t \Delta}B_{0,\rm h}+  f_{\ast \ast }^{\rm h}(u,v;B,W) ,        \\
     \Phi_{\ast \ast}^{\rm v }[v,W]     & \ = e^{t \Delta}B_{0,3}+  \sum_{m=1}^3 \Nvb_m[u,W](t)          \\
     &\   =:    e^{t \Delta}B_{0,3}+  f_{\ast \ast }^{\rm v}(u,v;B,W) .   \\
     \end{split}
\end{equation}
It follows from the estimates (\ref{im73})-(\ref{im75}) that
\begin{equation}\nonumber
     \begin{split}
        \|  (  f_{\ast}^{\rm h},   f_{\ast}^{\rm v}, f_{\ast \ast }^{\rm h},  f_{\ast \ast }^{\rm v}    )  (u,v;B,W)       \|_{Y} & \  \leq C_2 \| (u_0,B_0) \|_{X^s} \| (v,W) \|_{Y}   ,            \\
   \|  (  f_{\ast}^{\rm h},   f_{\ast}^{\rm v}, f_{\ast \ast }^{\rm h},  f_{\ast \ast }^{\rm v}    )  (u,q;B,N)       \|_{ L^{\infty}( 0,\infty; L^1_{\rm h} L^{\infty}_{\rm v}          ) } & \  \leq C_2 \| (u_0,B_0) \|_{X^s} \| (q,N) \|_{ L^{\infty}( 0,\infty; L^1_{\rm h} L^{\infty}_{\rm v}          ) }  ,            \\
     \end{split}
\end{equation}
for any $(v,W) \in Y$ and $(q,N)\in L^{\infty}( 0,\infty; L^1_{\rm h} L^{\infty}_{\rm v}    (\R^3)      ) $. Consequently,
\begin{equation}\nonumber
     \begin{split}
         \| \Phi (v,W)  \| _{Y} \leq  & \     \|  ( e^{t \Deltah} u_0 , e^{t \Delta}B_0    )  \|_{Y}  +   \|  (  f_{\ast}^{\rm h},   f_{\ast}^{\rm v}, f_{\ast \ast }^{\rm h},  f_{\ast \ast }^{\rm v}    )  (u,v;B,W)       \|_{Y}                \\
    \leq  & \  C_1 \|   (u_0,B_0)   \|_{ \widetilde{X^s} } +C_2 \| (u_0,B_0) \|_{X^s} \| (v,W) \|_{Y} ,
     \end{split}
\end{equation}
for any $(v,W) \in Y$. This means that $\Phi (v,W)\in Y$ and so the mapping $\Phi$ is well-defined. Next, it is also clear that
\begin{equation}\nonumber
     \begin{split}
        \| \Phi (v^{(1)},W^{(1)})-  \Phi (v^{(2)},W^{(2)})      \| _{Y} \leq & \
        \|       (  f_{\ast}^{\rm h},   f_{\ast}^{\rm v}, f_{\ast \ast }^{\rm h},  f_{\ast \ast }^{\rm v}    )          (u,v^{(1)}-v^{(2)} ; B, W^{(1)}-W^{(2)} )            \|_{Y}
        \\
  \leq  & \ C_2 \| (u_0,B_0) \|_{X^s} \| ( v^{(1)}-v^{(2)},W^{(1)}-W^{(2)} ) \|_{Y}    \\
       \leq  & \  \f{1}{2}   \| (v^{(1)}-v^{(2)},W^{(1)}-W^{(2)}) \|_{Y} ,         \\
     \end{split}
\end{equation}
where we have chosen the initial data sufficiently small, namely $\| (u_0,B_0) \|_{X^s}\leq \min\{   \epsilon_3, \f{1}{2C_2}\}$, in the last step. We may then apply the contraction mapping principle to $\Phi$ so that there exists a unique $( \overline{v},\overline{W}  ) \in Y $ obeying $( \overline{v},\overline{W}  )=\Phi ( \overline{v},\overline{W}  )$. Furthermore, one infers that
\[
\| ( \overline{v},\overline{W}  )  \|_{Y}  \leq 2 C_1 \|   (u_0,B_0)   \|_{ \widetilde{X^s} } .
\]
Next we check that this unique fixed point $( \overline{v},\overline{W}  ) $ coincides with the unique solution $(u,B)$ to the system (\ref{MHD-2}). Upon noticing the fact that $(u,B)$ is also a solution to the integral equation (\ref{im72}),
\begin{equation}\nonumber
     \begin{split}
        \| ( \overline{v},\overline{W}  )- (u,B) \|_{ L^{\infty}( 0,\infty; L^1_{\rm h} L^{\infty}_{\rm v}          ) }  =  & \  \|  \Phi ( \overline{v},\overline{W}  )  -\Phi(u,B)  \| _{ L^{\infty}( 0,\infty; L^1_{\rm h} L^{\infty}_{\rm v}          ) }       \\
    \leq     & \   C_2   \| (u_0,B_0) \|_{X^s}  \| ( \overline{v},\overline{W}  )- (u,B)        \|  _{ L^{\infty}( 0,\infty; L^1_{\rm h} L^{\infty}_{\rm v}          ) }               \\
       \leq   & \  \f{1}{2} \| ( \overline{v},\overline{W}  )- (u,B)        \|  _{ L^{\infty}( 0,\infty; L^1_{\rm h} L^{\infty}_{\rm v}          ) } ,         \\
     \end{split}
\end{equation}
yielding the desired result. Thus $(u,B) \in Y$ and satisfies all estimates of $( \overline{v},\overline{W}  )$ and in particular
\[
\| ( u,B )  \|_{Y}  \leq 2 C_1 \|   (u_0,B_0)   \|_{ \widetilde{X^s} } .
\]
The above estimate directly gives (\ref{im71}).

It thus remains to verify the estimates (\ref{im73})-(\ref{im75}). Notice that the first three estimates in (\ref{im73}), the first three estimates in (\ref{im74}) and the first three estimates in (\ref{im75}) have already been proved in \cite{F}. We only need to prove the remaining estimates. To begin with, we show (\ref{im73})$_4$.
\begin{equation}\nonumber
     \begin{split}
         \| |\xh| &\	\Nhu_1[B,W](t,x) \|_{L^{1}_{\rm h}  L^{\infty}_{\rm v}}  \\
         \leq  & \  \int_0^t \| |\xh| \Gh (t-\tau) \|_{L^1 (\R^2)   }
         \|     ( (\p_3 B_3)W_{\rm h} +W_3 \p_3 \Bh   )(\tau)    \|_{L^{1}_{\rm h}  L^{\infty}_{\rm v}}  d\tau
         \\
     & \  +  \int_0^t  \|  \Gh (t-\tau) \|_{L^1 (\R^2) }
        \|  (  (\p_3 B_3) |\xh|   W_{\rm h} + |\xh|  W_{3}  \p_3 \Bh) (\tau)       \|_{L^{1}_{\rm h}  L^{\infty}_{\rm v}}  d\tau
     \\
      \leq    & \   C \int_0^t   (t-\tau)^{ \f{1}{2}  }  \Big(
      \|  \p_3 B_3(\tau) \|_{L^{\infty}}  \|  W_{\rm h}(\tau)  \|_{L^{1}_{\rm h}  L^{\infty}_{\rm v}}  +  \|  \p_3 \Bh(\tau) \|_{L^{\infty}}  \|  W_{3}(\tau)  \|_{L^{1}_{\rm h}  L^{\infty}_{\rm v}}
      \Big) d \tau  \\
         & \  + C\int_0^t   \Big(     \|  \p_3 B_3(\tau) \|_{L^{\infty}}  \| |\xh| W_{\rm h}(\tau)  \|_{L^{1}_{\rm h}  L^{\infty}_{\rm v}}  +  \|  \p_3 \Bh(\tau) \|_{L^{\infty}}  \| |\xh|  W_{3}(\tau)  \|_{L^{1}_{\rm h}  L^{\infty}_{\rm v}}           \Big) d \tau      \\
      \leq    &\   C\int_0^t     (t-\tau)^{ \f{1}{2}  }  (1+\tau)^{-\f{5}{2} }     d \tau  \|   (u_0,B_0)  \|_{X^s}  \| (v,W) \|_{Y}  \\
      &\   + C\int_0^t (1+\tau)^{-2} d \tau   \|   (u_0,B_0)  \|_{X^s}  \| (v,W) \|_{Y}           \\
     \leq   &\  C (1+t)^{\f{1}{2}  }     \|   (u_0,B_0)  \|_{X^s}  \| (v,W) \|_{Y} .         \\
     \end{split}
\end{equation}
Next, for (\ref{im73})$_5$, it holds that
\begin{equation}\nonumber
     \begin{split}
         \| 	\Nhu_1[B,W] &\ (t,x) \|_{L^{1}_{\rm h}  L^{\infty}_{\rm v}} \\
         \leq  &\
         \int_0^t  \|  \Gh (t-\tau) \|_{L^1 (\R^2) }     \|     ( (\p_3 B_3)W_{\rm h} +W_3 \p_3 \Bh   )(\tau)    \|_{L^{1}_{\rm h}  L^{\infty}_{\rm v}}  d\tau
         \\
   \leq  & \  C    \int_0^t    \Big(
      \|  \p_3 B_3(\tau) \|_{L^{\infty}}  \|  W_{\rm h}(\tau)  \|_{L^{1}_{\rm h}  L^{\infty}_{\rm v}}  +  \|  \p_3 \Bh(\tau) \|_{L^{\infty}}  \|  W_{3}(\tau)  \|_{L^{1}_{\rm h}  L^{\infty}_{\rm v}}
      \Big) d \tau                     \\
        \leq  & \  C    \int_0^t    (1+\tau)^{-\f{5}{2} }   d \tau      \|   (u_0,B_0)  \|_{X^s}  \|   (v,W)   \| _{ L^{\infty}( 0,\infty; L^1_{\rm h} L^{\infty}_{\rm v}          ) }             \\
      \leq   &\ C    \|   (u_0,B_0)  \|_{X^s}  \|   (v,W)   \| _{ L^{\infty}( 0,\infty; L^1_{\rm h} L^{\infty}_{\rm v}          ) }      . \\
     \end{split}
\end{equation}
In a similar manner, we show (\ref{im73})$_6$.
\begin{equation}\nonumber
     \begin{split}
         \| \p_3 &\	\Nhu_1[B,W]  (t,x) \|_{L^{1}_{\rm h}  L^{\infty}_{\rm v}} \\
         \leq  &\
         \int_0^t  \|  \Gh (t-\tau) \|_{L^1 (\R^2) }     \|     ( (\p_3^2 B_3)W_{\rm h} +\p_3 B_3 \p_3 W_{\rm h}+ \p_3 W_{3} \p_3 \Bh  + W_3 \p_3^2 \Bh   )(\tau)    \|_{L^{1}_{\rm h}  L^{\infty}_{\rm v}}  d\tau
         \\
   \leq  & \  C    \int_0^t    \Big(
      \|  \p_3^2 B_3(\tau) \|_{L^{\infty}}  \|  W_{\rm h}(\tau)  \|_{L^{1}_{\rm h}  L^{\infty}_{\rm v}}  +
      \|  \p_3 B_3 \|_{L^{\infty}}   \|  \p_3 W_{\rm h}(\tau)  \|_{L^{1}_{\rm h}  L^{\infty}_{\rm v}}      \\
      &\   +
      \|  \p_3 \Bh(\tau) \|_{L^{\infty}}  \|  \p_3 W_{3}(\tau)  \|_{L^{1}_{\rm h}  L^{\infty}_{\rm v}}
      +   \| W_3 \|_{L^{1}_{\rm h}  L^{\infty}_{\rm v}}   \|  \p_3^2 \Bh (\tau)  \|  _{L^{\infty}}
      \Big) d \tau               \\
     \leq  &\  C    \int_0^t    \Big(          \| B_3(\tau) \|_{H^7}^{ \f{2}{7} }       \|\p_3 B_3(\tau) \|_{L^{\infty}}^{ \f{5}{7} }     \|  W_{\rm h}(\tau)  \|_{L^{1}_{\rm h}  L^{\infty}_{\rm v}}   +
      \|  \p_3 B_3 \|_{L^{\infty}}   \|  \p_3 W_{\rm h}(\tau)  \|_{L^{1}_{\rm h}  L^{\infty}_{\rm v}}            \\
      &\ +    \|  \p_3 \Bh(\tau) \|_{L^{\infty}}  \|  \p_3 W_{3}(\tau)  \|_{L^{1}_{\rm h}  L^{\infty}_{\rm v}}   +\| W_3 \|_{L^{1}_{\rm h}  L^{\infty}_{\rm v}}
       \| \Bh(\tau) \|_{H^7}^{ \f{2}{7} }       \|\p_3 \Bh(\tau) \|_{L^{\infty}}^{ \f{5}{7} }  \Big) d \tau
      \\
        \leq  & \  C    \int_0^t \left(    (1+\tau)^{-\f{5}{2} } +(1+\tau)^{-\f{25}{14} }  \right)    d \tau      \|   (u_0,B_0)  \|_{X^s}  \|   (v,W)   \| _{ Y }             \\
      \leq   &\ C    \|   (u_0,B_0)  \|_{X^s}  \|   (v,W)   \| _{Y}     . \\
     \end{split}
\end{equation}

We proceed with (\ref{im74})$_4$-(\ref{im74})$_6$. For any $k,l\in \{ 1,2,3 \}$, it holds
\begin{equation}\nonumber
     \begin{split}
         \Big\|   |\xh| & \  \int_0^t  e^{(t-\tau)\Deltah } \nablah (B_k W_l)(\tau,x)d \tau     \Big\|_{ L^{1}_{\rm h}  L^{\infty}_{\rm v} }   \\
 \leq  & \  \int_0^t   \| |\xh| \nablah \Gh(t-\tau)     \|_{L^1(\R^2)}
 \| B_k(\tau)  \|_{L^{\infty}}  \|   W_l (\tau)  \|_{ L^{1}_{\rm h}  L^{\infty}_{\rm v} }  d \tau
 \\
         & \  + \int_0^t  \|  \nablah \Gh(t-\tau)     \|_{L^1(\R^2)}   \| B_k(\tau)  \|_{L^{\infty}}     \| |\xh|  W_l (\tau)  \|_{ L^{1}_{\rm h}  L^{\infty}_{\rm v} }  d \tau        \\
     \leq    & \    C  \int_0^t (1+\tau)^{-2}  d \tau  \|   (u_0,B_0)  \|_{X^s}  \|   (v,W)   \| _{Y}            \\
     &\ +    C  \int_0^t (t-\tau)^{ -\f{1}{2} }    (1+\tau) ^{ -\f{3}{2} }    d \tau     \|   (u_0,B_0)  \|_{X^s}  \|   (v,W)   \| _{Y}             \\
      \leq    & \    C  \int_0^t (1+\tau)^{-2}  d \tau  \|   (u_0,B_0)  \|_{X^s}  \|   (v,W)   \| _{Y}            \\
     &\ +    C    (1+t) ^{ -\f{1}{2} }       \|   (u_0,B_0)  \|_{X^s}  \|   (v,W)   \| _{Y}             \\
    \leq  &\   C   \|   (u_0,B_0)  \|_{X^s}  \|   (v,W)   \| _{Y}   ,         \\
     \end{split}
\end{equation}
\begin{equation}\nonumber
     \begin{split}
         \Big\|  \int_0^t  & \   e^{(t-\tau)\Deltah } \nablah (B_k W_l)(\tau,x)d \tau     \Big\|_{ L^{1}_{\rm h}  L^{\infty}_{\rm v} }   \\
 \leq  & \  \int_0^t   \| \nablah \Gh(t-\tau)     \|_{L^1(\R^2)}
 \| B_k(\tau)  \|_{L^{\infty}}  \|   W_l (\tau)  \|_{ L^{1}_{\rm h}  L^{\infty}_{\rm v} }  d \tau
 \\
   \leq  &\     C  \int_0^t (t-\tau)^{ -\f{1}{2} }    (1+\tau) ^{ -2 }    d \tau     \|   (u_0,B_0)  \|_{X^s}  \|   (v,W)   \| _{ L^{\infty}( 0,\infty; L^1_{\rm h} L^{\infty}_{\rm v}          ) }              \\
    \leq  &\   C   (1+t) ^{ -\f{1}{2} }      \|   (u_0,B_0)  \|_{X^s}  \|   (v,W)   \| _{ L^{\infty}( 0,\infty; L^1_{\rm h} L^{\infty}_{\rm v}          ) } ,         \\
     \end{split}
\end{equation}
\begin{equation}\nonumber
     \begin{split}
         \Big\| \p_3 \int_0^t  & \   e^{(t-\tau)\Deltah } \nablah (B_k W_l)(\tau,x)d \tau     \Big\|_{ L^{1}_{\rm h}  L^{\infty}_{\rm v} }   \\
 \leq  & \  \int_0^t   \| \nablah \Gh(t-\tau)     \|_{L^1(\R^2)} \Big(
 \| \p_3 B_k(\tau)  \|_{L^{\infty}}  \|   W_l (\tau)  \|_{ L^{1}_{\rm h}  L^{\infty}_{\rm v} }       +   \| B_k(\tau)  \|_{L^{\infty}}  \|  \p_3 W_l (\tau)  \|_{ L^{1}_{\rm h}  L^{\infty}_{\rm v} }     \Big)       d \tau
 \\
   \leq  &\     C  \int_0^t (t-\tau)^{ -\f{1}{2} }  \left(  (1+\tau) ^{ -2 } +(1+\tau) ^{ -\f{5}{2} } \right)  d \tau     \|   (u_0,B_0)  \|_{X^s}  \|   (v,W)   \| _{Y }              \\
    \leq  &\   C   (1+t) ^{ -\f{1}{2} }      \|   (u_0,B_0)  \|_{X^s}  \|   (v,W)   \| _{Y }             .         \\
     \end{split}
\end{equation}

We end the proof by showing (\ref{im75})$_4$-(\ref{im75})$_9$. We shall only present one of them since the others could be carried out analogously. For any $k,l\in \{ 1,2,3 \}$, it holds
\begin{equation}\nonumber
     \begin{split}
         \Big\|   |x| &\ \int_0^t  e^{ (t-\tau)\Delta } \nablah (u_k W_l)  d \tau     \Big \|_{L^{1}_{\rm h}  L^{\infty}_{\rm v}}    \\
    \leq  & \  \int_0^t   \| |x| \nabla G(t-\tau)     \|_{L^1 }
 \| B_k(\tau)  \|_{L^{\infty}}  \|   W_l (\tau)  \|_{ L^{1}_{\rm h}  L^{\infty}_{\rm v} }  d \tau
 \\
         & \  + \int_0^t  \|  \nabla G(t-\tau)     \|_{L^1}   \| B_k(\tau)  \|_{L^{\infty}}     \| |\xh|  W_l (\tau)  \|_{ L^{1}_{\rm h}  L^{\infty}_{\rm v} }  d \tau        \\
     \leq    & \    C  \int_0^t (1+\tau)^{-2}  d \tau  \|   (u_0,B_0)  \|_{X^s}  \|   (v,W)   \| _{Y}            \\
     &\ +    C  \int_0^t (t-\tau)^{ -\f{1}{2} }    (1+\tau) ^{ -\f{3}{2} }    d \tau     \|   (u_0,B_0)  \|_{X^s}  \|   (v,W)   \| _{Y}             \\
      \leq   &\     C        \|   (u_0,B_0)  \|_{X^s}  \|   (v,W)   \| _{Y}    .          \\
     \end{split}
\end{equation}
The proof of Lemma \ref{lem51} is completely finished.               $\Box$

In the end of this subsection, we point out that the corresponding results to Lemmas \ref{lem50}-\ref{lem51} still hold for the incompressible MHD system with full dissipation and horizontal magnetic diffusion (\ref{MHD-3}). More specifically, we assume that $s\geq 9$ is an integer, $u_0,B_0 \in X^s(\R^3), \Grad \cdot u_0= \Grad \cdot B_0=0, |x|u_0(x)\in L^1(\R^3)$ and $\|  (u_0,B_0) \|_{ X^s_{\#} }$ is sufficiently small. Similar to the proof of Lemma \ref{lem50}, it follows from Lemmas \ref{lm1}, \ref{lm3-2}, \ref{lem14} and the bootstrapping argument that
\begin{equation}\label{im76}
     \begin{split}
        \| \Grad^{\alpha} B (t) \|_{L^{\infty}_{\rm h}  L^{1}_{\rm v}} & \  \leq C  \| (u_0,B_0) \|_{ X^s_{\#}  } t^{-1-\f{|\alphah|}{2}}, \\
     \| \Grad^{\alpha} u (t) \|_{L^{\infty}_{\rm h}  L^{1}_{\rm v}} & \  \leq C   \| (u_0,B_0) \|_{  X^s_{\#}  }    t^{-1-\f{1+|\alpha|}{2}} \log (2+t),  \\
     \end{split}
\end{equation}
for any $0<t<\infty$ and $\alpha=(\alphah,\alpha_3)\in (\mathbb{N} \cup \{ 0\})^2 \times (\mathbb{N} \cup \{ 0\})$ with $|\alpha|\leq 1$. Suppose further that $|\xh| B_0(x) \in L^1(\R^2_{\xh};(L^1 \cap L^{\infty})(\R_{x_3})  )$, one obtains in the same spirit of Lemma \ref{lem51} the following space-weighted estimates:
\begin{equation}\label{im77}
     \begin{split}
        \| |\xh|\Bh (t,x) \|_{L^{1}_{\rm h}  L^{\infty}_{\rm v}} & \  \leq C    \| (u_0,B_0)  \|_{X^s_{\# \# }}      (1+t)^{  \f{1}{2}  }        ,  \\
     \| |\xh| B_3 (t,x) \|_{L^{1}_{\rm h}  L^{\infty}_{\rm v}} & \  \leq C  \| (u_0,B_0)  \|_{X^s_{\# \# }},  \\
         \| |x| \uh (t,x) \|_{L^{1}_{\rm h}  L^{\infty}_{\rm v}} & \  \leq  C  \| (u_0,B_0)  \|_{X^s_{\# \# }}  (1+t)^{  \f{1}{2}  },  \\
           \| |x| u_3 (t,x) \|_{L^{1}_{\rm h}  L^{\infty}_{\rm v}} & \  \leq   C  \| (u_0,B_0)  \|_{X^s_{\# \# }}  (1+t)^{  \f{1}{2}  }.  \\
     \end{split}
\end{equation}
We shall omit the proof of (\ref{im76})-(\ref{im77}) for simplicity.
As a consequence, we conclude that the ansatzes (\ref{hm10})-(\ref{hm11}) actually hold with $A=C \| (u_0,B_0) \|_{X^s_{\#}}, A_{\ast}=C \| (u_0,B_0) \|_{X^s_{\# \# }} $ for some generic $C>0$.

\subsection{Proof of Theorem \ref{thm3}}
Let $s \geq 9$ be an integer and $u_0,B_0\in X^s(\R^3)$ satisfying $\nabla \cdot u_0=\nabla \cdot B_0=0$, $\| (u_0,B_0) \|_{X^s} \leq \epsilon_2$ and $|\xh| u_0(x) \in L^1(\R^2_{\xh};(L^1 \cap L^{\infty})(\R_{x_3})  ), |x| B_0(x)\in L^1(\R^3)$. It follows from Theorem \ref{thm1} and Lemmas \ref{lem50}-\ref{lem51} that there exists a unique global solution $(u,B)$ to (\ref{MHD-2}) with initial data $(u_0,B_0)$ and obey the ansatzes (\ref{hd8})-(\ref{hd9-2}) with $T=\infty,A=C \| (u_0,B_0) \|_{ \overline{X^s } }, A_{\ast}=C \| (u_0,B_0)  \|_{ \widetilde{X^s}  } $ for some suitable generic $C>0$. Based on Lemma \ref{lm1}, Lemma \ref{lem7} and Lemma \ref{lem9}, it holds for any $2<p \leq \infty$ that
\begin{equation}\nonumber
     \begin{split}
          \Big\|  \uh   & \ (t,x) -G_{\rm h}(t,\xh) \int_{\R^2} u_{0,{\rm h}}  (y_{\rm h},x_3)  d y_{\rm h} \\
    & \  + G_{\rm h}(t,\xh)  \int_0^{\infty} \int_{\R^2}  \p_3 (u_3 \uh) (\tau, y_{\rm h},x_3)
    d y_{\rm h} d \tau\\
    & \  - G_{\rm h}(t,\xh)  \int_0^{\infty} \int_{\R^2}  \p_3 (B_3 \Bh) (\tau, y_{\rm h},x_3)
    d y_{\rm h} d \tau  \Big\|_{L^p_x}      \\
   \leq  &\   \Big\|  e^{t \Deltah} u_{0,\rm{h}}(x) -G_{\rm h}(t,\xh) \int_{\R^2} u_{0,{\rm h}}  (y_{\rm h},x_3)  d y_{\rm h} \Big\|_{L^p_x}  \\
   &\ +    \Big\| \Nhu_1[u](t)   +G_{\rm h}(t,\xh)  \int_0^{\infty} \int_{\R^2}  \p_3 (u_3 \uh) (\tau, y_{\rm h},x_3)
    d y_{\rm h} d \tau \Big\|_{L^p_x}         \\
    &\   +    \Big\| \Nhu_1[B](t)   +G_{\rm h}(t,\xh)  \int_0^{\infty} \int_{\R^2}  \p_3 (B_3 \Bh) (\tau, y_{\rm h},x_3)
    d y_{\rm h} d \tau \Big\|_{L^p_x}           \\
    &\ +\sum_{m=2}^5  \|  ( \Nhu_m[u](t),\Nhu_m[B](t)    )          \| _{L^p_x}              \\
   \leq  &\   C  \Big(   \| |\xh|  u_{0,\rm h }\|_{L^1_{\rm h}  L^p_{\rm v}}  t^{ -(1-\f{1}{p})-\f{1}{2}  }  +
     \| (u_0,B_0) \|_{ \overline{X^s } }  \| (u_0,B_0) \|_{\widetilde{  X^s}}   t^{ -(1-\f{1}{p})-\f{1}{2}  } \log (2+t)
   \\
   &\  +      \| (u_0,B_0) \|_{ \overline{X^s } } \| (u_0,B_0) \|_{\widetilde{  X^s}}   t^{ -\f{3}{2}(1-\f{1}{p})-\f{1}{2}  }  +       \| (u_0,B_0) \|_{ \overline{X^s } }  \| (u_0,B_0) \|_{\widetilde{  X^s}}   t^{ -(1-\f{1}{p})-\f{1}{2}  }  (1+t^{-\f{1}{2}})     \Big)     \\
\leq    &\   C     \| (u_0,B_0) \|_{\widetilde{  X^s}}    t^{ -(1-\f{1}{p})-\f{1}{2}  } \log t       \\
     \end{split}
\end{equation}
for any $t \geq 2$. This verifies (\ref{ma11}). To proceed, we deduce from Lemma \ref{lm1} and Lemma \ref{lem7} that for any $ 1\leq p \leq \infty$
\begin{equation}\nonumber
     \begin{split}
       \Big \|
u_3  &\ (t,x)-   G_{\rm h}(t,\xh) \int_{\R^2} u_{0,3} (y_{\rm h},x_3)  d y_{\rm h}
\Big\|_{L^p_x}  \\
\leq  &\       \Big \| e^{t \Deltah} u_{0,3}-   G_{\rm h}(t,\xh) \int_{\R^2} u_{0,3} (y_{\rm h},x_3)  d y_{\rm h}
\Big\|_{L^p_x}  +\sum_{m=1}^3    \|  ( \Nvu_m[u](t),\Nvu_m[B](t)    )          \| _{L^p_x}                   \\
\leq &\  C  \| |\xh|u_{0,3}(x) \|_{L^1_{\rm h}  L^p_{\rm v}}                 t^{ -\f{3}{2} (1-\f{1}{p}) -\f{1}{2p}    }
+C \| (u_0,B_0) \|_{ \overline{X^s } }^2  t^{ -\f{3}{2} (1-\f{1}{p}) -\f{1}{2p}    }
\\
&\ +  C \| (u_0,B_0) \|_{ \overline{X^s } } ^2  t^{ -\f{9}{8} (1-\f{1}{p}) -\f{1}{2}    } \log(2+t)    +  C \| (u_0,B_0) \|_{ \overline{X^s } }^2  t^{ -\f{3}{2} (1-\f{1}{p}) -\f{1}{2}    }            \\
    \leq  & \ \left\{\begin{aligned}
& C  \| (u_0,B_0) \|_{\widetilde{X^s}      }   t^{ -\f{3}{2} (1-\f{1}{p}) -\f{1}{2p}    }      ,\,   1<p \leq \infty       \\
& C  \| (u_0,B_0) \|_{\widetilde{X^s}      } t^{-\f{1}{2}} \log t  ,\,\,\,\, p=1     \\
\end{aligned}\right.
     \end{split}
\end{equation}
for any $t \geq 2$. This proves (\ref{ma12}). Next, in view of Lemma \ref{lm2}, Lemma \ref{lem7} and Lemma \ref{lem10}, we infer for any $1<p\leq \infty$ that
\begin{equation}\nonumber
     \begin{split}
         t^  {  \f{3}{2}(1-\f{1}{p}) +\f{1}{2p} }  \Big\|  & \  u_3 (t,x) -G_{\rm h}(t,\xh) \int_{\R^2} u_{0,3}  (y_{\rm h},x_3)  d y_{\rm h} \\
    & \  + \nablah G_{\rm h}(t,\xh)  \cdot  \int_{\R^2}  \yh u_{0,3} ( y_{\rm h},x_3)
    d y_{\rm h} \\
    & \  -  \nablah  G_{\rm h}(t,\xh)  \cdot \int_0^{\infty} \int_{\R^2}   (u_3 \uh) (\tau, y_{\rm h},x_3)  d y_{\rm h} d \tau          \\
    &\  +  \nablah  G_{\rm h}(t,\xh)  \cdot \int_0^{\infty} \int_{\R^2}   (B_3 \Bh) (\tau, y_{\rm h},x_3)  d y_{\rm h} d \tau     \Big\| _{L^p_x}   \\
   \leq  &\     t^  {  \f{3}{2}(1-\f{1}{p}) +\f{1}{2p} }   \Big\|
   e^{t \Deltah} u_{0,3}(x)-G_{\rm h}(t,\xh) \int_{\R^2} u_{0,3}  (y_{\rm h},x_3)  d y_{\rm h}
   \\
   &\   + \nablah G_{\rm h}(t,\xh)  \cdot  \int_{\R^2}  \yh u_{0,3} ( y_{\rm h},x_3)
    d y_{\rm h}
   \Big\| _{L^p_x} \\
   &\   +     t^  {  \f{3}{2}(1-\f{1}{p}) +\f{1}{2p} } \left\|
\Nvu_1[u] (t,x)-\nablah G_{\rm h}(t,\xh) \cdot  \int_0^{\infty} \int_{\R^2}   (u_3 \uh) (\tau, y_{\rm h},x_3)
    d y_{\rm h}   d \tau
\right\|_{L^p_x}     \\
&\  + t^  {  \f{3}{2}(1-\f{1}{p}) +\f{1}{2p} }     \left\|
\Nvu_1[B] (t,x)-  \nablah G_{\rm h}(t,\xh) \cdot \int_0^{\infty} \int_{\R^2}   (B_3 \Bh) (\tau, y_{\rm h},x_3)
    d y_{\rm h}   d \tau
\right\|_{L^p_x}   \\
&\  +t^  {  \f{3}{2}(1-\f{1}{p}) +\f{1}{2p} }   \sum_{m=1}^3    \|  ( \Nvu_m[u](t),\Nvu_m[B](t)    )          \| _{L^p_x}             \\
\leq &\  t^{-\f{1}{2}+\f{1}{2p}  } \mathcal{R}_{p,1}(t)^{\f{1}{p}}   \|  |\xh| u_0(x)  \|_{L^1_x}^{1-\f{1}{p}}   \\
 &\   +     t^  {  (1-\f{1}{p}) +\f{1}{2} } \left\|
\Nvu_1[u] (t,x)-\nablah G_{\rm h}(t,\xh) \cdot  \int_0^{\infty} \int_{\R^2}   (u_3 \uh) (\tau, y_{\rm h},x_3)
    d y_{\rm h}   d \tau
\right\|_{L^p_x}     \\
&\  + t^  {  (1-\f{1}{p}) +\f{1}{2} }     \left\|
\Nvu_1[B] (t,x)-  \nablah G_{\rm h}(t,\xh) \cdot \int_0^{\infty} \int_{\R^2}   (B_3 \Bh) (\tau, y_{\rm h},x_3)
    d y_{\rm h}   d \tau
\right\|_{L^p_x}   \\
&\ +      C    t^  { - \f{1}{8}(1-\f{1}{p}) }   \log (2+t)
+ C    t^  { - \f{1}{2}(1-\f{1}{p}) } ,   \\
     \end{split}
\end{equation}
where the right-hand side tends to zero as long as $1<p\leq \infty$. This proves (\ref{ma13}). We end the proof by showing the asymptotic expansion of magnetic field. Clearly,
\begin{equation}\nonumber
     \begin{split}
       t ^ {  \f{3}{2}(1-\f{1}{p}) +\f{1}{2} }  \Big\|  B_i (t,x)  & \    +  \nabla G(t,x)  \cdot  \int_{\R^3}   y B_{0,i} ( y)
    d y \\
    & \  +  \sum_{j=1}^3 \p_j  G(t,x) \int_0^{\infty} \int_{\R^3}   (B_i u_j ) (\tau, y )  d y d \tau          \\
     & \  - \sum_{j=1}^3 \p_j  G(t,x) \int_0^{\infty} \int_{\R^3}   (u_i B_j ) (\tau, y )  d y d \tau     \Big\|_{L^p_x}     \\
   \leq    t^  {  \f{3}{2}(1-\f{1}{p}) +\f{1}{2} }  &\        \Big\|    e^{t \Delta}B_{0,i}    +  \nabla G(t,x)  \cdot  \int_{\R^3}   y B_{0,i} ( y)
    d y        \Big\|_{L^p_x}     \\
    &\ +     t^  {  \f{3}{2}(1-\f{1}{p}) +\f{1}{2} }     \Big\|
    \sum_{j=1}^3 \int_0^t  e^{(t-\tau)\Delta} \p_j (B_i u_j) (\tau,x)  d \tau       \\
    &\  -\sum_{j=1}^3 \p_j  G(t,x) \int_0^{\infty} \int_{\R^3}   (B_i u_j ) (\tau, y )  d y d \tau
    \Big\|_{L^p_x}      \\
     &\   +     t^  {  \f{3}{2}(1-\f{1}{p}) +\f{1}{2} }     \Big\|
    \sum_{j=1}^3 \int_0^t  e^{(t-\tau)\Delta} \p_j ( u_i B_j ) (\tau,x)  d \tau      \\
    &\   -
    \sum_{j=1}^3 \p_j  G(t,x) \int_0^{\infty} \int_{\R^3}   (u_i B_j ) (\tau, y )  d y d \tau
    \Big\|_{L^p_x} .      \\
     \end{split}
\end{equation}
From Lemma \ref{lm1} and Lemma \ref{lem11} it follows that the right-hand side tends to zero as time goes to infinity. This verifies (\ref{ma14}). The proof of Theorem \ref{thm3} is completely finished.

\subsection{Proof of Theorem \ref{thm4}}

It follows from Lemma \ref{lm1} and Lemma \ref{lem13} that for any $ 1\leq p \leq \infty$
\begin{equation}\nonumber
     \begin{split}
       \Big \|
B_3  &\ (t,x)-   G_{\rm h}(t,\xh) \int_{\R^2} B_{0,3} (y_{\rm h},x_3)  d y_{\rm h}
\Big\|_{L^p_x}  \\
\leq  &\       \Big \| e^{t \Deltah} B_{0,3}-   G_{\rm h}(t,\xh) \int_{\R^2} B_{0,3} (y_{\rm h},x_3)  d y_{\rm h}
\Big\|_{L^p_x}  +\sum_{m=1}^2    \|   \Nvd_m[u,B](t)     \| _{L^p_x}                   \\
\leq &\  C   \| |\xh| B_{0,3}(x) \|_{L^1_{\rm h}  L^p_{\rm v}}                 t^{ -\f{3}{2} (1-\f{1}{p}) -\f{1}{2p}    }
+C \| (u_0,B_0) \|_{ X^s_{ \#} }^2  t^{ -\f{3}{2} (1-\f{1}{p}) -\f{1}{2p}    }
\\
    \leq  & \ C \| (u_0,B_0) \|_{ X^s_{ \# \# } }  t^{ -\f{3}{2} (1-\f{1}{p}) -\f{1}{2p}    }
     \end{split}
\end{equation}
for any $t \geq 2$. This proves (\ref{ma16}).

Next, we deduce from Lemma \ref{lm2}, Lemma \ref{lem13} and Lemma \ref{lem16} that for any $1\leq p \leq \infty$
\begin{equation}\nonumber
     \begin{split}
         t^  {  \f{3}{2}(1-\f{1}{p}) +\f{1}{2p} }  \Big\|  & \  B_3 (t,x) -G_{\rm h}(t,\xh) \int_{\R^2} B_{0,3}  (y_{\rm h},x_3)  d y_{\rm h} \\
    & \  + \nablah G_{\rm h}(t,\xh)  \cdot  \int_{\R^2}  \yh B_{0,3} ( y_{\rm h},x_3)
    d y_{\rm h} \\
    & \  +  \nablah  G_{\rm h}(t,\xh)  \cdot \int_0^{\infty} \int_{\R^2}   (B_3 \uh) (\tau, y_{\rm h},x_3)  d y_{\rm h} d \tau          \\
    &\  -  \nablah  G_{\rm h}(t,\xh)  \cdot \int_0^{\infty} \int_{\R^2}   (u_3 \Bh) (\tau, y_{\rm h},x_3)  d y_{\rm h} d \tau     \Big\| _{L^p_x}   \\
   \leq  &\     t^  {  \f{3}{2}(1-\f{1}{p}) +\f{1}{2p} }   \Big\|
   e^{t \Deltah} B_{0,3}(x)-G_{\rm h}(t,\xh) \int_{\R^2} B_{0,3}  (y_{\rm h},x_3)  d y_{\rm h}
   \\
   &\   + \nablah G_{\rm h}(t,\xh)  \cdot  \int_{\R^2}  \yh B_{0,3} ( y_{\rm h},x_3)
    d y_{\rm h}
   \Big\| _{L^p_x} \\
   &\   +     t^  {  \f{3}{2}(1-\f{1}{p}) +\f{1}{2p} } \left\|
\Nvd_1[u,B] (t,x)+ \nablah G_{\rm h}(t,\xh) \cdot  \int_0^{\infty} \int_{\R^2}   (B_3 \uh) (\tau, y_{\rm h},x_3)
    d y_{\rm h}   d \tau
\right\|_{L^p_x}     \\
&\  + t^  {  \f{3}{2}(1-\f{1}{p}) +\f{1}{2p} }     \left\|
\Nvd_2[u,B] (t,x)-  \nablah G_{\rm h}(t,\xh) \cdot \int_0^{\infty} \int_{\R^2}   (u_3 \Bh) (\tau, y_{\rm h},x_3)
    d y_{\rm h}   d \tau
\right\|_{L^p_x}   \\
\leq &\  t^{-\f{1}{2}+\f{1}{2p}  } \mathcal{R}_{p,1}(t)^{\f{1}{p}}   \|  |\xh| B_0(x)  \|_{L^1_x}^{1-\f{1}{p}}   \\
 &\   +     t^  {  (1-\f{1}{p}) +\f{1}{2} } \left\|
\Nvd_1[u,B] (t,x)+\nablah G_{\rm h}(t,\xh) \cdot  \int_0^{\infty} \int_{\R^2}   (B_3 \uh) (\tau, y_{\rm h},x_3)
    d y_{\rm h}   d \tau
\right\|_{L^p_x}     \\
&\  + t^  {  (1-\f{1}{p}) +\f{1}{2} }     \left\|
\Nvd_2[u,B] (t,x)-  \nablah G_{\rm h}(t,\xh) \cdot \int_0^{\infty} \int_{\R^2}   (u_3 \Bh) (\tau, y_{\rm h},x_3)
    d y_{\rm h}   d \tau
\right\|_{L^p_x} .  \\
     \end{split}
\end{equation}
It is easily seen that the right-hand side tends to zero as long as $t \rightarrow \infty$ for any $1\leq p \leq \infty$. This proves (\ref{ma17}).

Finally, we come to the higher order asymptotic expansion of velocity field. For any $1\leq p \leq \infty$ and any $0<\sigma < \f{1}{2}$, we deduce from the integral equation (\ref{hm5}) and Lemmas \ref{lm1}, \ref{lem13} that for $t \geq 2$
\begin{equation}\nonumber
     \begin{split}
  t^  {  \f{3}{2}(1-\f{1}{p}) +\f{1}{2} -\sigma }   \Big\| &\  u_3 (t,x)   +  \nabla G(t,x)  \cdot  \int_{\R^3}    y u_{0,3} ( y)
    d y  +  \Nve_3[B] (t)      \Big\|_{L^p_x}             \\
    \leq  & \  t^  {  \f{3}{2}(1-\f{1}{p}) +\f{1}{2} -\sigma }   \Big\|  e^{ t \Delta } u_{0,3} +  \nabla G(t,x)  \cdot  \int_{\R^3}    y u_{0,3} ( y)
    d y    \Big\|_{L^p_x}       \\
    &\  +  t^  {  \f{3}{2}(1-\f{1}{p}) +\f{1}{2} -\sigma }  \Big\|
    \sum_{m=1}^5 \Nve_m[u](t)- \sum_{m=1,2,4,5} \Nve_m[B](t)
    \Big \|_{L^p_x} \\
     \leq  & \  t^  {  \f{3}{2}(1-\f{1}{p}) +\f{1}{2}  }   \Big\|  e^{ t \Delta } u_{0,3} +  \nabla G(t,x)  \cdot  \int_{\R^3}    y u_{0,3} ( y)
    d y    \Big\|_{L^p_x} + t^{-\sigma}  ,    \\
     \end{split}
\end{equation}
where the right-hand side vanishes as $t \rightarrow \infty$. This gives (\ref{ma20}). The asymptotic limits (\ref{ma18})-(\ref{ma19}) are established similarly. The proof of Theorem \ref{thm4} is finished completely.

\section{Conclusion and remarks}\label{con}

In order to reveal the anisotropic mechanism of incompressible MHD equations, we consider two classes of equations in this paper. The first is 3D incompressible MHD system with horizontal dissipation and full magnetic diffusion, the other being 3D incompressible MHD system with full dissipation and horizontal magnetic diffusion. Under suitable assumptions on the initial data, we obtained the decay estimates in the general $L^p$-norms for any $1 \leq p \leq \infty$ and the asymptotic expansions of solutions by making full use of the structure of equations and in particular the anisotropic nature. These results demonstrated the distinct features between the anisotropic MHD equations and the usual MHD equations and further revealed the crucial role played by the magnetic field for 3D anisotropic incompressible MHD equations. More precisely,
\begin{itemize}
\item {
    For the incompressible MHD system with horizontal dissipation and full magnetic diffusion, we mainly focus on the effect of magnetic field in the asymptotic profile of velocity field. It is shown that $\uh(t)$ decays at the rate of $O(t^{-(1-\f{1}{p})})$, like the 2D heat kernel and behaves asymptotically as the combination of linear solution and two integrals resulting from nonlinear effects. We also proved that $u_3(t)$ decays at the rate of $O(t^{-\f{3}{2}(1-\f{1}{p})})$, like the 3D heat kernel and behaves asymptotically as the linear solution. It turns out that the effect of magnetic field only appears in $\uh(t)$, instead of $u_3(t)$, for the leading term of asymptotic expansions. Strikingly, for the higher order asymptotic expansions, we have shown that a nonlinear term arising from magnetic field did affect the second leading term of asymptotic expansion for $u_3(t)$. In addition, under the spatial decay assumptions of initial data we proved the decay rate for the asymptotic expansion of $\uh(t)$, see (\ref{ma11}), for any $2<p\leq \infty$. This is in sharp contrast with the anisotropic incompressible Navier-Stokes system \cite{F}, where the same decay rate was obtained for any $1<p\leq \infty$. Moreover, we have also proved the higher order asymptotic expansion of magnetic field. It turns out that the weak dissipation of velocity field had no essential influence on asymptotic expansion of magnetic field.
    Roughly speaking, although we have the nice decay estimate of magnetic field, it did affect the asymptotic expansions of the velocity field in a very subtle manner; on the contrary, the full magnetic diffusion is robust enough so that the asymptotic expansion of magnetic field remains basically expected.
    }
\item{
    For the incompressible MHD system with full dissipation and horizontal magnetic diffusion, we consider both the decay estimates in $L^p$-norms and the asymptotic expansion of solutions. It is proved that $\Bh(t)$ decays at the rate of $O(t^{-(1-\f{1}{p})})$, like the 2D heat kernel and and behaves asymptotically as the combination of linear solution and two integrals from the nonlinear coupling between magnetic field and velocity field. In addition, for $B_3(t)$ the decay rate was proved to be $O(t^{-\f{3}{2}(1-\f{1}{p})})$, like the 3D heat kernel; the leading term of $B_3(t)$ was given by only the linear solution. Due to the slowly decay estimates of magnetic field and the structure of equations, $u(t)$ decays at the rate of $O(t^{-\f{9}{8}(1-\f{1}{p})  -\f{1}{2}  } \log(2+t)   )$. When it comes to higher order asymptotic expansions, we proved that nonlinear coupling did affect the asymptotic expansion of $B_3(t)$; the nonlinear terms of magnetic field became the second leading term of $u(t)$.
    }

\end{itemize}

Below we give several technical remarks concerning the main theorems.
\begin{itemize}
\item {
    It is possible to relax the regularity index $X^s(\R^3)$ for the initial magnetic field in Theorem \ref{thm1} due to the strong dissipation for magnetic field. Similarly, the regularity index $X^s(\R^3)$ for the initial velocity field may be also relaxed in Theorem \ref{thm2} thanks to the strong dissipation for velocity field. Indeed, it requires a sharper well-posedness theory for the anisotropic incompressible MHD systems than that of Propositions \ref{pr1}-\ref{pr2}. We did not pursue this direction since our main concern in the current paper is the decay estimates in the general $L^p$-norms and the asymptotic expansions of solutions.
    }
    \item{
    For the incompressible MHD system with horizontal dissipation and full magnetic diffusion, the smallness of $\| (u_0,B_0) \|_{ X^s }$ is enough to guarantee the unique global solution $(u,B) \in C ([0,\infty); X^s (\R^3) )$. However, in Theorem \ref{thm1} we require the smallness of $\| (u_0,B_0) \|_{  \overline{X^s}  }$. The smallness of additional part comes from the decay estimates of magnetic field. Moreover, the assumption of $B_0$ could be replaced by $B_0 \in X^s (\R^3), \Grad \cdot B_0 =0$ and
    \[
    \| e^{t \Delta } B_0(x)     \|_{L^1_x} \leq C(B_0)  t^{ -\f{1}{2}}
    \]
    for some $C(B_0) >0$ depending on $B_0$. Notice that our choice of $\int_{\R^3} |x|  | B_0(x )| dx <\infty $ is just a sufficient condition for the above decay estimate to be valid. It simplifies the presentation more or less. The same reason holds for the incompressible MHD system with full dissipation and horizontal magnetic diffusion.
    }
    \item{
    We remark that the decay estimate of $u(t)$ for system (\ref{MHD-3}) is not optimal, namely (\ref{ma-10}). For concreteness, we revisit one typical Duhamel term $\|\Grad^{\alpha} \Nhe_2[B](t) \|_{L^p}$ in Lemma \ref{lem13}. Recall that
    \[
     \left\|   \Grad^{\alpha} \int_0^t e^{ (t-\tau)\Delta  }\nablah \cdot (\Bh \otimes \Bh )  (\tau  )     d \tau \right\|_{L^1}    \leq
     CA^2 (1+t)^{ -\f{1+|\alpha|}{2}} \log(2+t).
    \]
    Observe next that for any $1\leq q \leq \infty$
    \begin{equation}\nonumber
     \begin{split}
     \Big\|   \Grad^{\alpha} &\  \int_0^t e^{ (t-\tau)\Delta  }\nablah \cdot (\Bh \otimes \Bh ) (\tau  )     d \tau \Big\|_{L^{\infty}}     \\
     \leq  &\     \int_0^  {  \f{t}{2} }        \|  \Grad^{\alpha} \nablah  G(t-\tau)  \|_{L^{\infty}}  \| (\Bh \otimes \Bh)(\tau)     \| _{L^1}    d\tau              \\
     &\ +   \int_{  \f{t}{2} } ^t    \|  \Grad^{\alpha}  G(t-\tau)  \|_{L^{ q } }  \|  \nablah \cdot (\Bh \otimes \Bh ) (\tau)     \| _{L^{q'} }       d\tau                    \\
  \leq  & \  CA^2          \int_0^  {  \f{t}{2} } (t-\tau)^{  -\f{3}{2}-\f{1+ |\alpha|}{2}         } (1+\tau) ^{  -1         }    d \tau                 \\
         & \ +     CA^2          \int_{  \f{t}{2} } ^t (t-\tau)^{  -\f{3 }{2}(1-\f{1}{q}) -\f{ |\alpha|}{2}            }  (1+\tau) ^{  -\f{3}{2} -(1-\f{1}{q'})         }       d\tau                 \\
     \leq     & \   CA^2 t^{  -\f{3 }{2}(1-\f{1}{q}) -\f{ |\alpha|}{2}  +1 -\f{3}{2} -\f{1}{q}      }  \\
     = &\    CA^2 t^{  -\f{3}{2}  +\f{1}{2} \f{1}{q}  -\f{1+|\alpha|}{2}   }.
     \end{split}
\end{equation}
}
It should be stressed that in the third inequality we have invoked the following condition for $q$:
\[
0\leq \f{3 }{2}(1-\f{1}{q}) + \f{ |\alpha|}{2}
<1.
\]
Solving the above constraint gives rise to $1\leq q <\f{3}{2}$. Then we realize that the decay estimate in (\ref{ma-10}) corresponds to the case $q=\f{4}{3}$ and the limiting decay estimate of velocity field should be
\[
\| \nabla ^{\alpha}   u (t) \|_{L^p} \leq C (1+t)^{ -\f{7}{6}(1-\f{1}{p}) -\f{1+|\alpha|}{2}        } \log(2+t) \| (u_0,B_0) \|_{X^{s}_{\#}  }, \,\,\,\, 1\leq p \leq \infty.
\]
However, it seems this critical decay rate would never be touched by the current method. Anyway, the key point is that all the Duhamel terms produce strictly slower decay rates than the linear estimates of $u(t)$, thus affecting the decay estimates and asymptotic expansions of velocity field for system (\ref{MHD-3}).

\end{itemize}

\section{Appendix}\label{ap}
\subsection{Global well-posedness theory}
For convenience of the reader and to make the paper self-contained, we present the proof for global well-posedness of \eqref{MHD-2} and \eqref{MHD-3} in this appendix. These results are fundamental in studying the large time behavior of solutions. More precisely, we have
\begin{prop}\label{pr1}
Let $k\geq 2$ be an integer. Assume that
\[
u_0 \in H^k(\R^3),\, B_0 \in H^k(\R^3),\,\, \nabla \cdot u_0=
\nabla \cdot B_0=0.
\]
Then there exists $\epsilon>0$ such that if $\| u_0 \|_{H^k}+\| B_0 \|_{H^k} \leq \epsilon$, the Cauchy problem \eqref{MHD-2} admits a unique global solution satisfying
\[
\| u(t) \|_{H^k}+\| B(t) \|_{H^k} \leq C  \epsilon
\]
for a generic constant $C>0$ and for all $t>0$.
\end{prop}

\begin{prop}\label{pr2}
Let $k\geq 2$ be an integer. Assume that
\[
u_0 \in H^k(\R^3),\, B_0 \in H^k(\R^3),\,\, \nabla \cdot u_0=
\nabla \cdot B_0=0.
\]
Then there exists $\epsilon>0$ such that if $\| u_0 \|_{H^k}+\| B_0 \|_{H^k} \leq \epsilon$, the Cauchy problem \eqref{MHD-3} admits a unique global solution satisfying
\[
\| u(t) \|_{H^k}+\| B(t) \|_{H^k} \leq C  \epsilon
\]
for a generic constant $C>0$ and for all $t>0$.
\end{prop}

Before the proof, we recall a version of anisotropic inequality that will be extensively used below. There exists a generic constant $C>0$ such that
\[
\int_{\R^3}  |fgh|  dx \leq C \| f \|_{L^2}^{ \f{1}{2}       }
\| \p_1f \|_{L^2}^{ \f{1}{2}       }
\| g \|_{L^2}^{ \f{1}{2}       }
\| \p_2 g \|_{L^2}^{ \f{1}{2}       }
\| h \|_{L^2}^{ \f{1}{2}       }
\| \p_3 h \|_{L^2}^{ \f{1}{2}       }
\]
holds if the right-hand side is well-defined. For the proof, one may refer to \cite{WZ}. Moreover, the proof of Proposition \ref{pr2} follows the same line as Proposotion \ref{pr1}.

{\it Proof of Proposition \ref{pr1}.} As is standard in mathematical theory of fluid mechanics, the global solution is essentially constructed by the local solution and the global a priori estimates. To begin with, we explain how to construct the local solutions. Let $P_n$ be the cut-off operator defined by $\mathscr{F}_{\R^3} (P_n f)(\xi):=1_{\{|\xi|\leq n\}}(\xi)  \mathscr{F}_{\R^3} f(\xi)$. Consider the approximated problem
\begin{align}\label{ap-0}
    \begin{cases}
    \partial_t u -  \Deltah u + P_n( u \cdot \nabla  u) -P_n ( B \cdot \nabla  B) + \nabla p = 0, & t>0 , x \in \mathbb{R}^3,\\
    \partial_t B -    \Delta  B + P_n( u \cdot \nabla B) - P_n( B \cdot \nabla  u )= 0, & t>0 , x \in \mathbb{R}^3,\\
    \nabla \cdot u = \nabla \cdot B = 0 , & t \geqslant 0, x \in \mathbb{R}^3,\\
    u(0,x) = P_n u_0(x), B(0,x) = P_n B_0(x),
    & x \in \mathbb{R}^3.
    \end{cases}
\end{align}
Essentially based on an fixed point argument, one obtains the local solution $(u^{(n)},B^{(n)})$ satisfying $(u^{(n)},B^{(n)})=P_n (u^{(n)},B^{(n)})$ on the time interval $[0,T(n))$. In order to pass to the limit $n\rightarrow \infty$ and obtain the global solutions, we are forced to establish sufficient global estimates. Thus, we shall mainly pay attention to global a priori estimates. For brevity, we set
\[
E(t):= \sup_{0\leq s \leq t} ( \| u(s) \|_{H^k}^2+\| B(s) \|_{H^k}^2  )
+\int_0^t     \big( 2 \| \nablah  u(s) \|  _{H^k}^2+
2 \| \nabla B(s) \|  _{H^k}^2
\big)
ds.
\]
It is clear that $\| u \|_{H^k}^2+\| B \|_{H^k}^2  $ is equivalent to
\[
\| (u,B) \|_{L^2}^2 + \sum_{j=1}^3 \| \partial_j ^k u \|_{L^2}^2
+ \sum_{j=1}^3 \| \partial_j ^k B \|_{L^2}^2 .
\]
Dotting \eqref{MHD-2}$_1$ with $u$ and \eqref{MHD-2}$_2$ with $B$ gives rise to
\beq\label{ap-1}
\| (u,B) (t)\|_{L^2}^2+ \int_0^t \left(  2 \|\nablah u (s)  \|_{L^2}^2
+2 \|\nabla B (s)  \|_{L^2}^2
\right) ds
= \| (u_0,B_0) \|_{L^2}^2.
\eeq
To proceed, we apply $\partial_j^k$ to \eqref{MHD-2}$_1$ and \eqref{MHD-2}$_2$, followed by taking $L^2$ inner product, to infer that
\beq\label{ap-2}
\frac{d}{dt}  \sum_{j=1}^3 \|  (\partial_j^k u,\partial_j^k B)    \|_{L^2}^2
+2  \sum_{j=1}^3  \| \nablah \partial_j^k u  \|_{L^2}^2
+2 \sum_{j=1}^3  \| \nabla \partial_j^k B  \|_{L^2}^2 =\sum_{j=1}^4\mathcal{I}_j
\eeq
with
\begin{equation}\label{ap-3}
     \begin{split}
        \mathcal{I}_1=  & \  -2 \sum_{j=1}^3 \int_{\R^3}
         \partial_j^k (u \cdot \nabla u) \cdot \partial_j^k u    dx ,                  \\
   \mathcal{I}_2= & \    2  \sum_{j=1}^3 \int_{\R^3}
      \big( \partial_j^k (B \cdot \nabla B)- B\cdot \nabla \partial_j^k B           \big) \cdot \partial_j^k u   dx  ,                 \\
        \mathcal{I}_3= & \  -2  \sum_{j=1}^3 \int_{\R^3}
        \partial_j^k(u \cdot \nabla B) \cdot \p_j^k B dx ,                             \\
       \mathcal{I}_4=  & \   2  \sum_{j=1}^3 \int_{\R^3}
       \big(  \p_j^k(B\cdot \Grad u )- B \cdot \Grad \p_j^k u          \big) \cdot \p_j^k B  dx
       \\
     \end{split}
\end{equation}
Next we estimate $\mathcal{I}_j\,  (j=1,2,3,4)$ in a suitable manner. For $\mathcal{I}_1$, we infer from $\Grad \cdot u=0$ that
\beq\label{ap-4}
        \mathcal{I}_1=   -2 \sum_{j,i,l=1}^3   \sum_{\alpha=1}^k  \int_{\R^3 }   \p_j^{\alpha} u_{l} \p_{l} \p_j ^{k-\alpha} u_i \p_j^k u_i  dx   .
\eeq
The integrals in \eqref{ap-4} can be classified into two classes: the case that at least one of $j,l$ is equal to $1$ or $2$ and the case that $j=l=3$. For the first case, we may assume, without loss of generality, that $l=1,j=3$ and invoke the anisotropic inequality to see
\begin{equation}\label{ap-5}
     \begin{split}
         \sum_{i=1}^3   \sum_{\alpha=1}^k   & \  \int_{\R^3 } \p_3 ^{\alpha} u_1 \p_1 \p_3^{k-\alpha} u_i \p_3^k u_i dx       \\
    \leq C & \  \| \p_3 ^{\alpha} u_1\|_{L^2}^{\f{1}{2}  } \| \p_1 \p_3 ^{\alpha} u_1\|_{L^2}^{\f{1}{2}  }          \|  \p_1 \p_3^{k-\alpha} u_i \|_{L^2}^{\f{1}{2}  }
      \|  \p_1 \p_3^{k-\alpha+1} u_i \|_{L^2}^{\f{1}{2}  }  \|  \p_3^k u_i \| _{L^2}^{\f{1}{2}  }  \|  \p_2 \p_3^k u_i \| _{L^2}^{\f{1}{2}  }      \\
       \leq C  & \  \| u \|_{H^k}     \|\nablah u \|_{H^k} ^{ 2 } .     \\
     \end{split}
\end{equation}
For the second case, we recall again $\nabla \cdot u=0$ and the anisotropic inequality to obtain
\begin{equation}\label{ap-6}
     \begin{split}
        \sum_{i=1}^3 \sum_{\alpha=1}^k  & \  \int_{\R^3 }  \p_3^{\alpha} u_3 \p_3 \p_3^{k-\alpha} u_i \p_3^k u_i   dx               \\
   = & \  -    \sum_{i=1}^3 \sum_{\alpha=1}^k   \int_{\R^3 }  \p_3^{\alpha-1} (\p_1 u_1+ \p_2 u_2)  \p_3^{k-\alpha+1} u_i \p_3^k u_i dx                     \\
       \leq C   & \   \| \p_3^{\alpha-1}\p_1 u_1 \| _{L^2}^{\f{1}{2}  }      \| \p_3  \p_3^{\alpha-1}\p_1 u_1 \| _{L^2}^{\f{1}{2}  }  \|   \p_3^{k-\alpha+1} u_i \|  _{L^2}^{\f{1}{2}  }
         \| \p_1  \p_3^{k-\alpha+1} u_i \|  _{L^2}^{\f{1}{2}  }   \| \p_3^k u_i \| _{L^2}^{\f{1}{2}  }     \| \p_2  \p_3^k u_i \| _{L^2}^{\f{1}{2}  }                  \\
         & \ +   \| \p_3^{\alpha-1}\p_2 u_2 \| _{L^2}^{\f{1}{2}  }      \| \p_3  \p_3^{\alpha-1}\p_2 u_2 \| _{L^2}^{\f{1}{2}  }  \|   \p_3^{k-\alpha+1} b_i \|  _{L^2}^{\f{1}{2}  }
         \| \p_1  \p_3^{k-\alpha+1} b_i \|  _{L^2}^{\f{1}{2}  }   \| \p_3^k u_i \| _{L^2}^{\f{1}{2}  }     \| \p_2  \p_3^k u_i \| _{L^2}^{\f{1}{2}  }                  \\
          \leq C   & \   \|u  \|_{H^k}   \| \nablah  u \|_{H^k}^{2} .
     \end{split}
\end{equation}
Next we turn to $\mathcal{I}_2$. Obviously,
\beq\label{ap-7}
        \mathcal{I}_2=   -2 \sum_{j,i,l=1}^3   \sum_{\alpha=1}^k  \int_{\R^3 }   \p_j^{\alpha} B_{l} \p_{l} \p_j ^{k-\alpha} B_i \p_j^k u_i  dx   .
\eeq
The integrals in \eqref{ap-7} could be divided into two classes: the case that $j=1$ or $2$ and case that $j=3$. Without loss of generality, it suffices to estimate the following two integrals:
\begin{equation}\label{ap-8}
     \begin{split}
         \sum_{i=1}^3   \sum_{\alpha=1}^k   & \  \int_{\R^3 } \p_1 ^{\alpha} B_3 \p_3 \p_1^{k-\alpha} B_i \p_1^k u_i  dx       \\
    \leq C & \  \| \p_1 ^{\alpha} B_3\|_{L^2}^{\f{1}{2}  } \| \p_3 \p_1 ^{\alpha} B_3\|_{L^2}^{\f{1}{2}  }          \|  \p_3 \p_1^{k-\alpha} B_i \|_{L^2}^{\f{1}{2}  }
      \|  \p_3 \p_1^{k-\alpha+1} B_i \|_{L^2}^{\f{1}{2}  }  \|  \p_1^k u_i \| _{L^2}^{\f{1}{2}  }  \|  \p_2 \p_1^k u_i \| _{L^2}^{\f{1}{2}  }      \\
       \leq C  & \  \| u \|_{H^k} ^{\f{1}{2}  }  \| \nablah u \|_{H^k} ^{\f{1}{2}  }
       \| B \|_{H^k} ^{\f{1}{2}  }    \|\Grad B  \|_{H^k}^{\f{3}{2}  } ;
     \end{split}
\end{equation}
\begin{equation}\label{ap-9}
     \begin{split}
         \sum_{i=1}^3   \sum_{\alpha=1}^k   & \  \int_{\R^3 } \p_3 ^{\alpha} B_2 \p_2 \p_3^{k-\alpha} B_i \p_3^k u_i  dx       \\
    \leq C & \  \| \p_3 ^{\alpha} B_2\|_{L^2}^{\f{1}{2}  } \| \p_2 \p_3 ^{\alpha} B_2\|_{L^2}^{\f{1}{2}  }          \|  \p_2 \p_3^{k-\alpha} B_i \|_{L^2}^{\f{1}{2}  }
      \|  \p_2 \p_3^{k-\alpha+1} B_i \|_{L^2}^{\f{1}{2}  }  \|  \p_3^k u_i \| _{L^2}^{\f{1}{2}  }  \|  \p_1 \p_3^k u_i \| _{L^2}^{\f{1}{2}  }      \\
       \leq C  & \  \| u \|_{H^k} ^{\f{1}{2}  }  \| \nablah u \|_{H^k} ^{\f{1}{2}  }
       \| B \|_{H^k} ^{\f{1}{2}  }    \|\Grad B  \|_{H^k}^{\f{3}{2}  } .
     \end{split}
\end{equation}
To proceed, we easily see
\beq\label{ap-10}
        \mathcal{I}_3=   -2 \sum_{j,i,l=1}^3   \sum_{\alpha=1}^k  \int_{\R^3 }   \p_j^{\alpha} u_{l} \p_{l} \p_j ^{k-\alpha} B_i \p_j^k B_i  dx   .
\eeq
In analogy with before, these integrals could be divided into two classes: the case that $j=1$ or $2$ and case that $j=3$. Likewise,
\begin{equation}\label{ap-11}
     \begin{split}
         \sum_{i=1}^3   \sum_{\alpha=1}^k   & \  \int_{\R^3 } \p_1 ^{\alpha} u_3 \p_3 \p_1^{k-\alpha} B_i \p_1^k B_i  dx       \\
    \leq C & \  \| \p_1 ^{\alpha} u_3\|_{L^2}^{\f{1}{2}  } \| \p_2 \p_1 ^{\alpha} u_3\|_{L^2}^{\f{1}{2}  }          \|  \p_3 \p_1^{k-\alpha} B_i \|_{L^2}^{\f{1}{2}  }
      \|  \p_3 \p_1^{k-\alpha+1} B_i \|_{L^2}^{\f{1}{2}  }  \|  \p_1^k B_i \| _{L^2}^{\f{1}{2}  }  \|  \p_3 \p_1^k B_i \| _{L^2}^{\f{1}{2}  }      \\
       \leq C  & \  \| u \|_{H^k} ^{\f{1}{2}  }  \| \nablah u \|_{H^k} ^{\f{1}{2}  }
       \| B \|_{H^k} ^{\f{1}{2}  }    \|\Grad B  \|_{H^k}^{\f{3}{2}  } ;
     \end{split}
\end{equation}
\begin{equation}\label{ap-12}
     \begin{split}
         \sum_{i=1}^3   \sum_{\alpha=1}^k   & \  \int_{\R^3 } \p_3 ^{\alpha} u_2 \p_2 \p_3^{k-\alpha} B_i \p_3^k B_i  dx       \\
    \leq C & \  \| \p_3 ^{\alpha} u_2\|_{L^2}^{\f{1}{2}  } \| \p_2 \p_3 ^{\alpha} u_2\|_{L^2}^{\f{1}{2}  }          \|  \p_2 \p_3^{k-\alpha} B_i \|_{L^2}^{\f{1}{2}  }
      \|  \p_2 \p_3^{k-\alpha+1} B_i \|_{L^2}^{\f{1}{2}  }  \|  \p_3^k B_i \| _{L^2}^{\f{1}{2}  }  \|  \p_1 \p_3^k B_i \| _{L^2}^{\f{1}{2}  }      \\
       \leq C  & \  \| u \|_{H^k} ^{\f{1}{2}  }  \| \nablah u \|_{H^k} ^{\f{1}{2}  }
       \| B \|_{H^k} ^{\f{1}{2}  }    \|\Grad B  \|_{H^k}^{\f{3}{2}  } .
     \end{split}
\end{equation}
By the same token, we see
\beq\label{ap-13}
        \mathcal{I}_4=   -2 \sum_{j,i,l=1}^3   \sum_{\alpha=1}^k  \int_{\R^3 }   \p_j^{\alpha} B_{l} \p_{l} \p_j ^{k-\alpha} u_i \p_j^k B_i  dx
\eeq
and
\begin{equation}\label{ap-14}
     \begin{split}
         \sum_{i=1}^3   \sum_{\alpha=1}^k   & \  \int_{\R^3 } \p_3 ^{\alpha} B_1 \p_1 \p_3^{k-\alpha} u_i \p_3^k B_i  dx       \\
    \leq C & \  \| \p_3 ^{\alpha} B_1\|_{L^2}^{\f{1}{2}  } \| \p_2 \p_3 ^{\alpha} B_1\|_{L^2}^{\f{1}{2}  }          \|  \p_1 \p_3^{k-\alpha} u_i \|_{L^2}^{\f{1}{2}  }
      \|  \p_1 \p_3^{k-\alpha+1} u_i \|_{L^2}^{\f{1}{2}  }  \|  \p_3^k B_i \| _{L^2}^{\f{1}{2}  }  \|  \p_1 \p_3^k B_i \| _{L^2}^{\f{1}{2}  }      \\
       \leq C  & \  \| u \|_{H^k} ^{\f{1}{2}  }  \| \nablah u \|_{H^k} ^{\f{1}{2}  }
       \| B \|_{H^k} ^{\f{1}{2}  }    \|\Grad B  \|_{H^k}^{\f{3}{2}  } ;
     \end{split}
\end{equation}
\begin{equation}\label{ap-15}
     \begin{split}
         \sum_{i=1}^3   \sum_{\alpha=1}^k   & \  \int_{\R^3 } \p_3 ^{\alpha} B_3 \p_3 \p_3^{k-\alpha} u_i \p_3^k B_i  dx       \\
    \leq C & \  \| \p_3 ^{\alpha} B_3\|_{L^2}^{\f{1}{2}  } \| \p_2 \p_3 ^{\alpha} B_3\|_{L^2}^{\f{1}{2}  }          \|   \p_3^{k-\alpha+1} u_i \|_{L^2}^{\f{1}{2}  }
      \|  \p_1 \p_3^{k-\alpha+1} u_i \|_{L^2}^{\f{1}{2}  }  \|  \p_3^k B_i \| _{L^2}^{\f{1}{2}  }  \|  \p_3 \p_3^k B_i \| _{L^2}^{\f{1}{2}  }      \\
       \leq C  & \  \| u \|_{H^k} ^{\f{1}{2}  }  \| \nablah u \|_{H^k} ^{\f{1}{2}  }
       \| B \|_{H^k} ^{\f{1}{2}  }    \|\Grad B  \|_{H^k}^{\f{3}{2}  } .
     \end{split}
\end{equation}

In view of \eqref{ap-1}-\eqref{ap-2} and recalling the definition of $E(t)$, we find upon using the estimates \eqref{ap-4}-\eqref{ap-15} that
\begin{equation}\label{ap-16}
     \begin{split}
         E(t) \leq  & \ C E(0)+ C \sum_{j=1}^4 \int_0^t \mathcal{I}_j(s)  ds   \\
    \leq  & \ C E(0)+ C  \int_0^t   \|u  \|_{H^k}   \| \nablah  u \|_{H^k}^{2}           ds              \\
         & \  +C \int_0^t \| u \|_{H^k} ^{\f{1}{2}  }  \| \nablah u \|_{H^k} ^{\f{1}{2}  }
       \| B \|_{H^k} ^{\f{1}{2}  }    \|\Grad B  \|_{H^k}^{\f{3}{2}  }   ds   \\
        \leq  & \  C E(0)+ C E(t)^{ \f{3}{2}   }.
     \end{split}
\end{equation}
With the local solution at hand, a direct application of bootstrapping argument to \eqref{ap-16} allows us to pass to the limit $n\rightarrow \infty$ and obtain the global small solutions to the original system \eqref{MHD-2}, see for instance \cite{WZ}. Notice also that uniqueness of solutions is proved via the standard energy method which is omitted here. This completes the proof of Proposition \ref{pr1}. $\Box$

The proof of Proposition \ref{pr1} readily gives the following
\begin{cor}\label{co1}
Let $k\geq 2$ be an integer. Assume that
\[
u_0 \in H^k(\R^3),\, B_0 \in H^k(\R^3),\,\, \nabla \cdot u_0=
\nabla \cdot B_0=0.
\]
Then there exists $\epsilon>0$ such that if $\| u_0 \|_{H^k}+\| B_0 \|_{H^k} \leq \epsilon$, the Cauchy problem \eqref{MHD-2} admits a unique global solution satisfying
\[
\| u(t) \|_{H^j}^2+\| B(t) \|_{H^j}^2 +   \int_0^t   \left(  \|\nablah u (s)  \|_{H^j}^2
+\|\nabla B (s)  \|_{H^j}^2
\right) ds     \leq C \left( \| u_0 \|_{H^j}^2+\| B_0 \|_{H^j}^2  \right)
\]
for a generic constant $C>0$, for all $t>0$ and any $j=0,1,...,k$.
\end{cor}

\subsection{Derivation of the decomposition (\ref{hm4})-(\ref{hm5}) }
For the convenience of the reader, we present a short derivation for the decomposition of velocity field regarding the incompressible MHD system with full dissipation and horizontal magnetic diffusion (\ref{hm1}). First of all, we recall the integral equation of velocity as
\label{de}
\[
u(t) = e^{ t \Delta} u_0
    - \int_0^t e^{ (t - \tau) \Delta} \mathbb{P} \nabla \cdot ( u \otimes  u - B \otimes B ) (\tau) d\tau.
\]
Noticing the identity $\mathbb{P}=(\delta_{kl}+ R_k R_l )_{k,l=1}^3$, where $\{ R_k\}_{k=1}^3$ signifies the Riesz transform in 3D, we see that the $j$-th component of $ e^{(t-\tau)\Delta }   \mathbb{P} \nabla \cdot ( u \otimes  u) (\tau)$ is given by
\begin{equation}\label{ap-17}
     \begin{split}
        \Big( e^{(t-\tau)\Delta } & \  \mathbb{P} \nabla \cdot ( u \otimes  u) (\tau)\Big)_j    \\
    =  & \   \sum_{k=1}^3   e^{(t-\tau)\Delta } \p_k (u_k  u_j)(\tau)
    +
    \sum_{k,l=1}^3  \p_j \p_k \p_l \mathscr{F}^{-1}_{\R^3} \left(
    |\xi|^{-2} e^{  -(t-\tau) |\xi|^2      }
    \right) \ast (u_k u_l )(\tau)
    \\
        =  & \   \sum_{k=1}^3   e^{(t-\tau)\Delta } \p_k (u_k  u_j)(\tau)
    +    \sum_{k,l=1}^3  \p_j \p_k \p_l  N(t-\tau)  \ast (u_k u_l )(\tau),          \\
 =   &\    e^{(t-\tau)\Delta } \p_3 (u_3  u_j)  (\tau) +   \sum_{k=1}^2   e^{(t-\tau)\Delta } \p_k (u_k  u_j)(\tau)   +    \p_j \p_3 \p_3  N(t-\tau)  \ast (u_3 u_3 )(\tau)         \\
 &\   +   \sum_{k,l=1}^2  \p_j \p_k \p_l  N(t-\tau)  \ast (u_k u_l )(\tau)
 +    2    \sum_{k=1}^2      \p_j \p_k \p_3  N(t-\tau)  \ast (u_k u_3 )(\tau)                     \\
     \end{split}
\end{equation}
where we set
\[
N(t,x):= \mathscr{F}^{-1}_{\R^3} \left(
    |\xi|^{-2} e^{  -t |\xi|^2      }
    \right) .
\]
A simple calculation gives another convenient form of $N$. Indeed,
\begin{equation}\label{ap-18}
     \begin{split}
       \mathscr{F}^{-1}_{\R^3} \left(
    |\xi|^{-2} e^{  -t |\xi|^2      }
    \right) (x)=   & \  (2\pi)^{-3} \int_{\R^3}e^{i x \cdot \xi}  |\xi|^{-2} e^{  -t |\xi|^2      } d \xi     \\
  =  & \     (2\pi)^{-3} \int_{\R^3}e^{i x \cdot \xi}  \left(  \int_0^{\infty} e^{ -s |\xi|^2} ds \right)  e^{  -t |\xi|^2      } d \xi                           \\
      =   & \     \int_0^{\infty}    \left(    (2\pi)^{-3} \int_{\R^3}e^{i x \cdot \xi}   e^{  -(t+s) |\xi|^2      }       d \xi \right)                ds                  \\
        =  & \  \int_0^{\infty}  G(t+s,x)             ds   .
     \end{split}
\end{equation}
Here $G(t,x)$ is the 3D Gaussian.

Performing similar calculations yields the decomposition of $ e^{(t-\tau)\Delta }   \mathbb{P} \nabla \cdot ( B \otimes  B) (\tau)$. From these identities we easily obtain the decomposition of velocity field (\ref{hm4})-(\ref{hm5}).

\subsection{On a asymptotic limit of $\Bh(t)$ for system (\ref{MHD-3}) } \label{app01}

In this subsection, we give the explanations why the corresponding estimates for $\Bh(t)$ cannot be achieved in Theorem \ref{thm4}. Indeed, in addition to the assumptions of Lemma \ref{lem15}, 
we suppose further that $(u,B)$ be subject to the ansatzes $(iii),(iv)$, namely (\ref{hm10})-(\ref{hm11}). Then for any $2<p\leq \infty$ there exists $C>0$ depending only on $p$ such that
\begin{equation}\label{hm23}
     \begin{split}
\Big\|
\Nhd_3[u,B] (t,x)+
&\ G_{\rm h}(t,\xh)  \int_0^{\infty} \int_{\R^2}  \p_3 (u_3 \Bh) (\tau, y_{\rm h},x_3)
    d y_{\rm h}  d \tau
\Big\|_{L^p_x} \\
    \leq & \  CA ( A+ A_{\ast} )t^{-(1-\f{1}{p})-\f{1}{2}}  (1+t)^{\f{1}{2}},        \\
     \end{split}
\end{equation}
and for any $2<p \leq \infty$  there exists $C>0$ depending only on $p$ such that
\begin{equation}\label{hm24}
     \begin{split}
\Big\|
\Nhd_4[u,B] (t,x)-
&\ G_{\rm h}(t,\xh)  \int_0^{\infty} \int_{\R^2}  \p_3 (B_3 \uh) (\tau, y_{\rm h},x_3)
    d y_{\rm h}  d \tau
\Big\|_{L^p_x} \\
    \leq & \  \left\{\begin{aligned}
& CA ( A+ A_{\ast} )    t^{-(1-\f{1}{p})-\f{1}{2}}   (1+t)^{ \f{1}{2p}  } \,\,\, \text{if  } 2< p< \infty         \\
&  CA ( A+ A_{\ast} ) t^{-(1-\f{1}{p})-\f{1}{2}}    \log(2+t) \,\,\, \text{if  }   p=\infty.        \\
\end{aligned}\right.
\end{split}
\end{equation}

To show (\ref{hm23}), we notice that
\begin{equation}\nonumber
     \begin{split}
        \|  |\yh| \p_3(u_3 \Bh)(\tau ) \|_{L^1}  \leq & \   \| \p_3 u_3(\tau) \| _{     L^{\infty}_{\rm h}  L^{1}_{\rm v}   }
        \| |\yh|  \Bh (\tau)  \| _{     L^{1}_{\rm h}  L^{\infty}_{\rm v}   }       \\
    & \ +   \| \p_3 \Bh (\tau) \| _{     L^{\infty}_{\rm h}  L^{1}_{\rm v}   }
        \| |\yh|  u_3 (\tau)  \| _{     L^{1}_{\rm h}  L^{\infty}_{\rm v}   }      \\
       \leq   & \  CA A_{\ast}  \Big( \tau^{-2} (1+\tau)^{     \f{1}{2}  } \log(2+\tau)   +\tau^{-1} (1+\tau)^{     \f{1}{2}  }     \Big), \\
     \end{split}
\end{equation}
\begin{equation}\nonumber
     \begin{split}
        \|  |\yh| \p_3(u_3 \Bh)(\tau ) \|_{     L^{1}_{\rm h}  L^{\infty}_{\rm v}   }  \leq & \   \| \p_3 u_3(\tau) \| _{   L^{\infty}  }
        \| |\yh|  \Bh (\tau)  \| _{     L^{1}_{\rm h}  L^{\infty}_{\rm v}   }       \\
    & \ +   \| \p_3 \Bh (\tau) \| _{     L^{\infty}   }
        \| |\yh|  u_3 (\tau)  \| _{     L^{1}_{\rm h}  L^{\infty}_{\rm v}   }      \\
       \leq   & \  CA A_{\ast}   (1+\tau)^{-\f{1}{2} } , \\
     \end{split}
\end{equation}
which immediately yields
\begin{equation}\nonumber
     \begin{split}
         \|  |\yh| \p_3(u_3 \Bh)(\tau ) \|_{     L^{1}_{\rm h}  L^{p}_{\rm v}   }  \leq       & \      \|  |\yh| \p_3(u_3 \Bh)(\tau ) \|_{L^1} ^{ \f{1}{p} }     \|  |\yh| \p_3(u_3 \Bh)(\tau ) \|_{     L^{1}_{\rm h}  L^{\infty}_{\rm v}   }  ^{ 1-\f{1}{p} }             \\
    \leq  & \  CA A_{\ast}  \Big(  \tau^{-\f{2}{p}   }  (1+\tau)^{  -\f{1}{2}+  \f{1}{p}       } (\log(2+\tau))^{   \f{1}{p}   } +
    \tau^{-\f{1}{p}   }   (1+\tau)^{  -\f{1}{2}+  \f{1}{p}       } \Big)
    .    \\
     \end{split}
\end{equation}
Observe that $\mathcal{I}_2$ is equivalent to
\[
\mathcal{I}_2(t,x)=
- \int_0^{\f{t}{2}  } \int_{\R^2} \int_0^1
(\nablah \Gh) (t,\xh-\theta \yh) \cdot \yh \p_3(u_3 \Bh)(\tau,\yh,x_3)
d \theta d \yh d \tau;
\]
whence
\begin{equation}\nonumber
     \begin{split}
         \|  \mathcal{I}_2(t) \|_{L^p} \leq  & \
         \int_0^{\f{t}{2}  } \int_{\R^2} \int_0^1
         \|   (\nablah \Gh) (t,\cdot-\theta \yh)      \|_{L^p(\R^2)}
         \|    |\yh| \p_3(u_3 \Bh)(\tau,\yh,\cdot)   \|_{L^p(\R)}d \theta d \yh d \tau
         \\
  \leq  & \  C t^{  -(1-\f{1}{p})-\f{1}{2}     }  \int_0^{\f{t}{2}  }
   \|    |\yh| \p_3(u_3 \Bh)(\tau)   \|_{     L^{1}_{\rm h}  L^{p}_{\rm v}   }  d \tau
  \\
        \leq  & \  CA A_{\ast}t^{  -(1-\f{1}{p})-\f{1}{2}     }
        \int_0^{\f{t}{2}  }  \Big(  \tau^{-\f{2}{p}   }  (1+\tau)^{  -\f{1}{2}+  \f{1}{p}       } (\log(2+\tau))^{   \f{1}{p}   } +
    \tau^{-\f{1}{p}   }   (1+\tau)^{  -\f{1}{2}+  \f{1}{p}       } \Big)   d \tau
        \\
     \leq     & \   CA A_{\ast}t^{  -(1-\f{1}{p})-\f{1}{2}     }  (1+t)^{ \f{1}{2} } .    \\
     \end{split}
\end{equation}
It should be pointed out that the assumption $2<p\leq \infty$ was invoked in the last integral. Recalling $\mathcal{I}_m$ with $m=1,3,4$ in Lemma \ref{lem15}, we arrive at (\ref{hm23}).

The proof of (\ref{hm24}) is carried out analogously. Indeed, it suffices to notice that
\begin{equation}\nonumber
     \begin{split}
        \|  |\yh| \p_3(B_3 \uh)(\tau ) \|_{L^1}  \leq & \   \| \p_3 B_3(\tau) \| _{     L^{\infty}_{\rm h}  L^{1}_{\rm v}   }
        \| |\yh|  \uh (\tau)  \| _{     L^{1}_{\rm h}  L^{\infty}_{\rm v}   }       \\
    & \ +   \| \p_3 \uh (\tau) \| _{     L^{\infty}_{\rm h}  L^{1}_{\rm v}   }
        \| |\yh|  B_3 (\tau)  \| _{     L^{1}_{\rm h}  L^{\infty}_{\rm v}   }      \\
       \leq   & \  CA A_{\ast}  \Big( \tau^{-1} (1+\tau)^{     \f{1}{2}  }   +\tau^{-2}  \log(2+\tau)    \Big), \\
     \end{split}
\end{equation}
\begin{equation}\nonumber
     \begin{split}
        \|  |\yh| \p_3(B_3 \uh)(\tau ) \|_{     L^{1}_{\rm h}  L^{\infty}_{\rm v}   }  \leq & \   \| \p_3 B_3(\tau) \| _{   L^{\infty}  }
        \| |\yh|  \uh (\tau)  \| _{     L^{1}_{\rm h}  L^{\infty}_{\rm v}   }       \\
    & \ +   \| \p_3 \uh (\tau) \| _{     L^{\infty}   }
        \| |\yh|  B_3 (\tau)  \| _{     L^{1}_{\rm h}  L^{\infty}_{\rm v}   }      \\
       \leq   & \  CA A_{\ast}   (1+\tau)^{- 1} ; \\
     \end{split}
\end{equation}
whence
\begin{equation}\nonumber
     \begin{split}
         \|  |\yh| \p_3(B_3 \uh)(\tau ) \|_{     L^{1}_{\rm h}  L^{p}_{\rm v}   }  \leq       & \      \|  |\yh| \p_3(B_3 \uh)(\tau ) \|_{L^1} ^{ \f{1}{p} }     \|  |\yh| \p_3(B_3 \uh)(\tau ) \|_{     L^{1}_{\rm h}  L^{\infty}_{\rm v}   }  ^{ 1-\f{1}{p} }             \\
    \leq  & \  CA A_{\ast}  \Big(  \tau^{-\f{1}{p}   }  (1+\tau)^{  -1+  \f{3}{2} \f{1}{p}       }  +
    \tau^{-\f{2}{p}   }   (1+\tau)^{  -1+  \f{1}{p}       }  (\log(2+\tau)) ^{ \f{1}{p}}         \Big)
    .    \\
     \end{split}
\end{equation}
Thus,
\[
\int_0^{\f{t}{2}}  \tau^{-\f{2}{p}   }   (1+\tau)^{  -1+  \f{1}{p}       }  (\log(2+\tau)) ^{ \f{1}{p}}        d \tau
\leq \log(2+t)
\]
for any $2<p \leq \infty$ and
\[
\int_0^{\f{t}{2}} \tau^{-\f{1}{p}   }  (1+\tau)^{  -1+  \f{3}{2} \f{1}{p}       }d \tau \leq   \left\{\begin{aligned}
& (1+t)^{ \f{1}{2p}  } \,\,\, \text{if  } 1< p< \infty         \\
& \log(2+t) \,\,\, \text{if  }   p=\infty.         \\
\end{aligned}\right.
\]

\clearpage\newpage

{\centerline{\bf{Acknowledgements}}}
The research of Yang Li is supported by National Natural Science Foundation of China under grant No. 12001003. The author sincerely thanks Dr. Mikihiro Fujii for many fruitful discussions. He also thanks Professor Yongzhong Sun and Professor Boqing Dong for constant encouragement and guidance.

{\centerline{\bf{Data Availability Statement}}}

All data generated or analysed during this study are included in this published article.

{\centerline{\bf{Conflict of Interests}}}
The author hereby confirms that there is no conflict of
interest.


\end{document}